\documentclass[11pt]{article}
\usepackage{amsfonts,amssymb,graphicx,bm}
\usepackage{psfrag}

\newcommand{\Ga}{\Gamma}
\newcommand{\ga}{\gamma}
\newcommand{\si}{\sigma}
\newcommand{\Si}{\Sigma}
\newcommand{\Th}{\Theta}
\newcommand{\al}{\alpha}
\newcommand{\rh}{\rho}
\newcommand{\de}{\delta}
\newcommand{\be}{\beta}
\newcommand{\ta}{\tau}
\newcommand{\io}{\iota}
\newcommand{\ep}{\epsilon}
\newcommand{\Om}{\Omega}
\newcommand{\om}{\omega}
\newcommand{\N}{\mathbb{N}}
\newcommand{\C}{\mathbb{C}}
\newcommand{\R}{\mathbb{R}}
\newcommand{\Proj}{\mathbb{P}}
\newcommand{\cO}{\mathcal{O}}
\newcommand{\cD}{\mathcal{D}}
\newcommand{\cS}{\mathcal{S}}
\newcommand{\cH}{\mathcal{H}}
\newcommand{\cC}{\mathcal{C}}
\def\<{\left\langle} \def\>{\right\rangle}
\newcommand{\ux}{\mathbf{x}}
\newcommand{\uy}{\mathbf{y}}
\newcommand{\uz}{\mathbf{z}}
\newcommand{\gS}{\mathfrak{S}}
\newcommand{\qed}{\hfill\rule{3mm}{3mm}}
\newtheorem{Lemma}{Lemma}
\newtheorem{Remark}{Remark}
\newtheorem{Theorem}{Theorem}
\newtheorem{Proposition}{Proposition}
\newtheorem{Corollary}{Corollary}
\newtheorem{Conjecture}{Conjecture}
\newtheorem{Definition}{Definition}

\begin{document}
\title{On the volume conjecture for classical spin networks}
%\begin{center}
%{\it Dedicated to the memory of Pierre Leroux}
%\end{center}
\author{Abdelmalek Abdesselam}
\maketitle

\begin{center}
{\it Department of Mathematics,
P. O. Box 400137,
University of Virginia,
Charlottesville, VA 22904-4137, USA }\\
email: \texttt{malek@virginia.edu}
\end{center}

\bigskip
\begin{center}
{\it In memoriam Pierre Leroux}
\end{center}

\bigskip

\parbox{11.8cm}{ \small
{\bf Abstract.}
We prove an upper bound for the evaluation of all classical $SU_2$ spin
networks
conjectured by Garoufalidis and van der Veen. This implies one half of the
analogue of the volume conjecture which they proposed for classical spin
networks. We are also able to obtain the other half, namely, an exact
determination of
the spectral radius, for the special class of generalized
drum graphs. Our proof uses
a version of Feynman diagram calculus which we developed as a tool
for the interpretation of the symbolic method
of classical invariant theory, in a manner which is rigorous yet true to the
spirit of the classical literature.
}

\bigskip

\parbox{12cm}{\small
}

\tableofcontents

\section{Introduction}

\subsection{Motivation}
The volume conjecture~\cite{Kashaev,MurakamiM} is one of the most
important open problems in low-dimensional topology, with vast ramifications
in many active areas of mathematics and theoretical physics (see,
e.g.,~\cite{Roberts2,DijkgraafF,DimofteGLZ} and references therein). 
The conjecture relates the exponential growth rate of the colored Jones
polynomial of a hyperbolic knot evaluated at a root of unity to the
hyperbolic volume of the knot complement.
There is a close connection between this problem and that of analysing
asymptotics of $q$-deformed or quantum spin
networks~\cite{Roberts2,MurakamiTet,TaylorW,Costantino1,Costantino2,vdVeen}.
Perhaps as a simpler setting for studying such asymptotics,
a vigorous program for the systematic study of large angular momentum
asymptotics of classical spin networks (CSN), with precisely formulated
conjectures, was presented in~\cite{GvdV}.
It is the main source of inspiration for the present article.
As a further piece of motivation for the study of CSN's, one can
note that spin networks feature in a great variety of topics, e.g.,
quantum computation~\cite{MarzuoliR}, Tyurin's approach to nonabelian
theta functions~\cite{Tyurin}, shell models of turbulence~\cite{Eyink},
to only cite some of the perhaps less well known.

\subsection{History}
The introduction of CSN's is usually attributed to
Roger Penrose~\cite{Penrose1,Penrose2} who used them in an attempt to 
combinatorially quantize gravity.
This is in similar spirit to Regge's
calculus~\cite{Regge}. Indeed, by analysing the asymtotics of CSN's
when all angular momenta are large, Ponzano and Regge~\cite{PonzanoR}
made a strong case for the emergence of the Regge action of 3d gravity
from this semiclassical regime. This is based on a precise
asymptotic formula for the 6-j symbol conjectured in~\cite{PonzanoR}
but proved much later by Roberts in~\cite{Roberts1}.
Recently spin networks and their generalizations such as the
Barrett-Crane model~\cite{BarrettC}
and spin foams~\cite{Perez}
have become a staple food at the table
of loop quantum gravity (see~\cite{Rovelli} and references therein).
A central issue in this approach to quantum gravity is the understanding
of this semiclassical regime
(see, e.g., \cite{BaezCE,BarrettS,FriedelL,Gurau}).

Earlier, CSN's essentially appeared in the works of 
the Lituanian School of quantum angular momentum theory
(QAMT)~\cite{YutsisLV}. Indeed some sort of graphical notation
becomes indispensable when calculating complicated 3n-j symbols.
However, the beginnings of QAMT and the theory of CSN's go much further back
to the classical invariant theory (CIT) of binary forms, although there is some
effort in mathematical translation needed in order to see the connection.
The ideal tool for this translation is Feynman diagram calculus (FDC),
namely a diagrammatic representation of contractions of tensors as in
the remarkable book~\cite{Cvitanovic}. Enough elements of such a translation
were presented in~\cite{AC1,AC2,AC3}
to suit the needs of these articles. More is required here in order
to translate the
modern CSN formalism into the framework of CIT, and this is the object
of \S\ref{CGsection}.
A piece of data which appears naturally in the CIT/FDC picture is
the notion of smooth orientations of a cubic graph, according to
the terminology of~\cite{KiralyS}. This is key to our solution of some
of the questions raised in~\cite{GvdV}.

There is a huge physical literature on QAMT where such objects
as Clebsch-Gordan, Clebsch-Gordon (sic), Clebsh-Gordon (sic),\ldots
coefficients are ubiquitous.
Yet, with only a few exceptions such as~\cite{Dowker},
this literature shows almost no sign of awareness or aknowledgement
of the work of Alfred Clebsch and Paul Gordan.
Aside from introducing some foundational material needed in the
subsequent proofs, we attemp to correct this injustice in \S\ref{CGsection}.
Our hope is that by allowing the users of QAMT to tap into the vast
and most often perfectly
rigorous 19th century literature on CIT, and through
cross-fertilization with modern theories in mathematics and physics,
new ideas and unexpected connections will emerge.

\subsection{Definitions and statement of results}\label{introdefsec}
We will follow the definitions of~\cite{GvdV} where a CSN is defined
as a pair
$(\Ga,\ga)$ consisting of an abstract cubic ribbon graph $\Ga$ together with
a decoration $\ga$ of the edges by natural integers.
By cubic ribbon graph we mean a trivalent regular graph
equipped with a cyclic ordering of the edges
at each vertex.
These are also called rotation systems, or fat graphs, although one has to be
careful as definitions may vary in the literature~\cite[\S2.1]{BollobasR}.
The notion we use here is that of pure rotation systems
as in~\cite[Chap. 3]{GrossT}. It is well known that such a ribbon graph
defines (up to orientation preserving equivalence of imbeddings)
a unique imbedding of the underlying abstract graph into a compact
orientable Riemann surface~\cite[Thm. 3.2.3]{GrossT}.
Note that we allow multiple edges and loops, i.e., edges with both ends
attached to the same vertex. We also allow $\Ga$ to be disconnected.
Finally, although this means a slight arm twisting on the usual
definition of a graph, we allow trivial components without vertices
which are made of a single edge closing upon itself.
\[
\parbox{2.5cm}{
\begin{center}
\includegraphics[width=0.9cm]{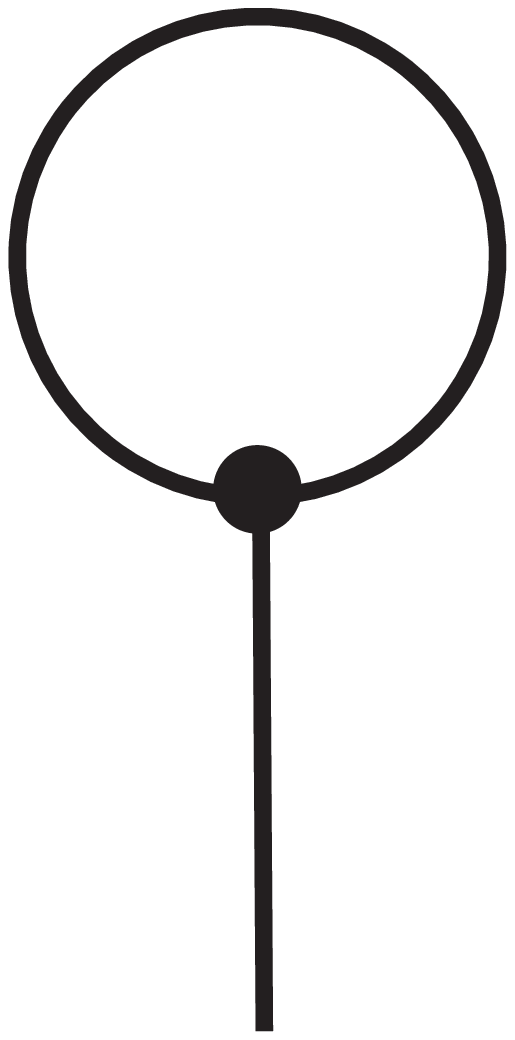}\\
{\rm a\ loop}
\end{center}}
\qquad,\qquad
\parbox{2.5cm}{
\begin{center}
\includegraphics[width=0.9cm]{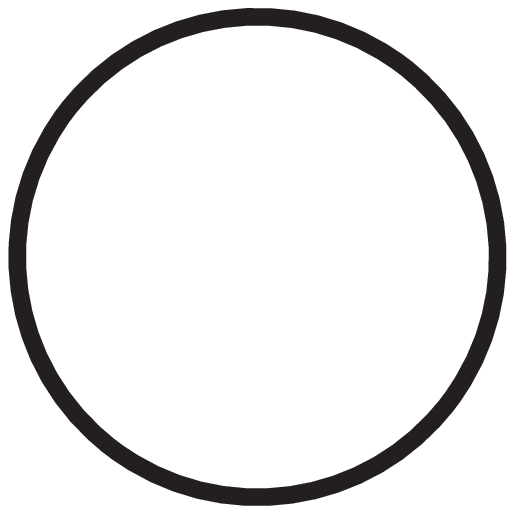}\\
{\ }\\
{\rm a\ trivial}\\
{\rm component}
\end{center}}
\]
Note that the decorations $\ga=(\ga(e))_{\ga\in E(\Ga)}$
in $\N^{E(\Ga)}$ where $E(\Ga)$ is the edge set of $\Ga$ must satisfy
the following admissibility conditions.
For every vertex, the three integers $a,b,c$ associated to the
incident edges are such that $a+b+c$ is even and the triangle
inquality $|a-b|\le c\le a+b$ holds.
Note that for a loop vertex
\[
\parbox{1.5cm}{
\psfrag{a}{${\scriptstyle a}$}
\psfrag{b}{${\scriptstyle b}$}
\includegraphics[width=1.5cm]{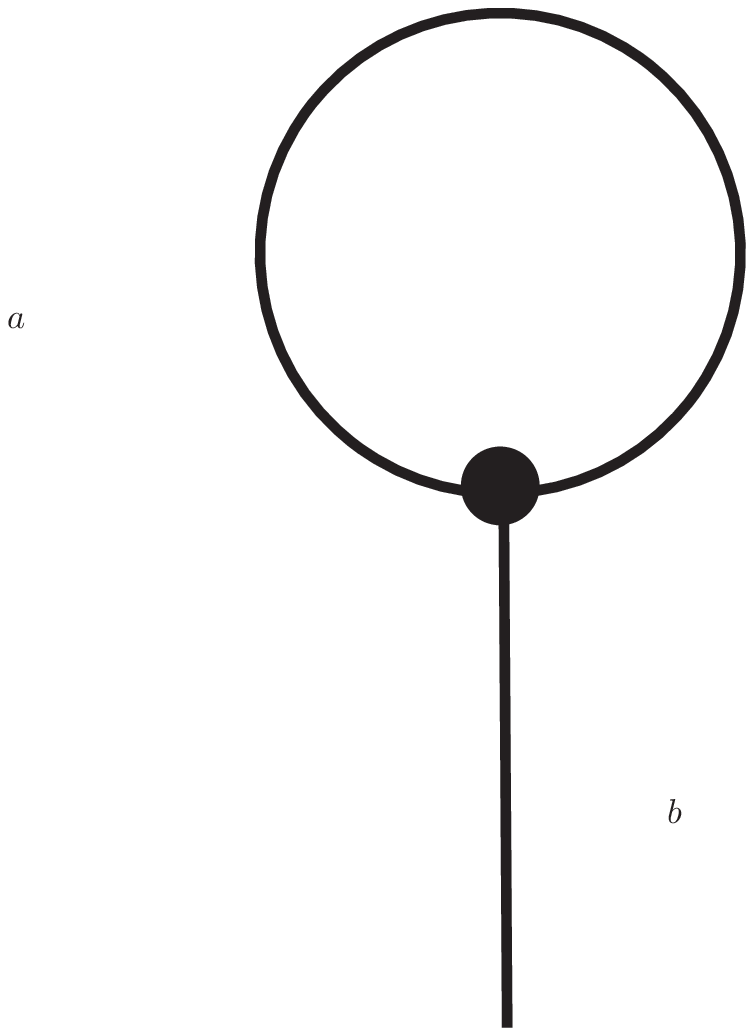}
}
\]
the decoration $a$ is counted twice so the constraint reduces to:
$b$ is even and $0\le b\le 2a$.
For a trivial component
\[
\parbox{1.2cm}{
\psfrag{a}{${\scriptstyle a}$}
\includegraphics[width=1.5cm]{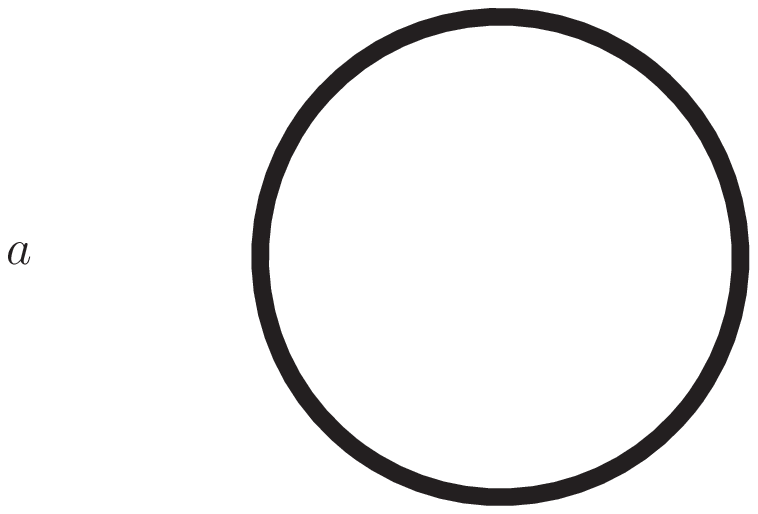}
}
\]
the decoration $a$ can be any nonnegative integer.

We can now define as in~\cite{GvdV} the Penrose evaluation $\<\Ga,\ga\>^P$
of such a spin network:
\begin{enumerate}
\item
Use the imbedding into the surface $\Si$ to thicken the vertices
into discs and the edges into bands.
\item
On an edge carrying the decoration $a$, draw $a$ parallel strands
and a perpendicular bar $\parbox{0.7cm}{
\psfrag{P}{${\scriptstyle P}$}
\includegraphics[width=0.7cm]{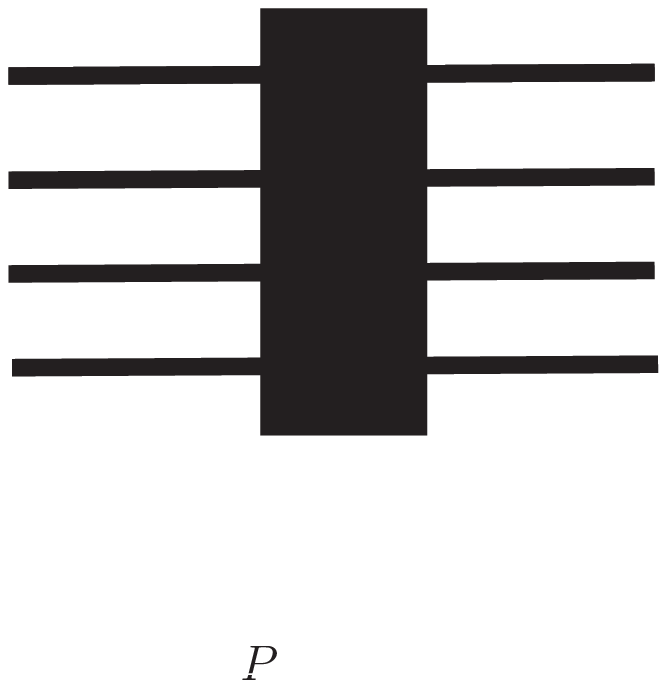}}$ which we call a Penrose bar:
\[
\parbox{1.5cm}{
\psfrag{a}{${\scriptstyle a=3}$}
\includegraphics[width=1.5cm]{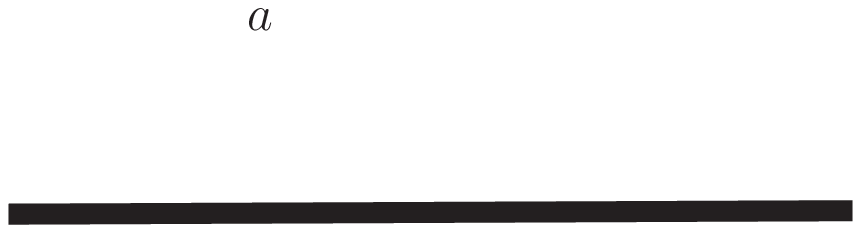}}
\qquad\longrightarrow\qquad
\parbox{1.7cm}{
\psfrag{P}{${\scriptstyle P}$}
\includegraphics[width=1.7cm]{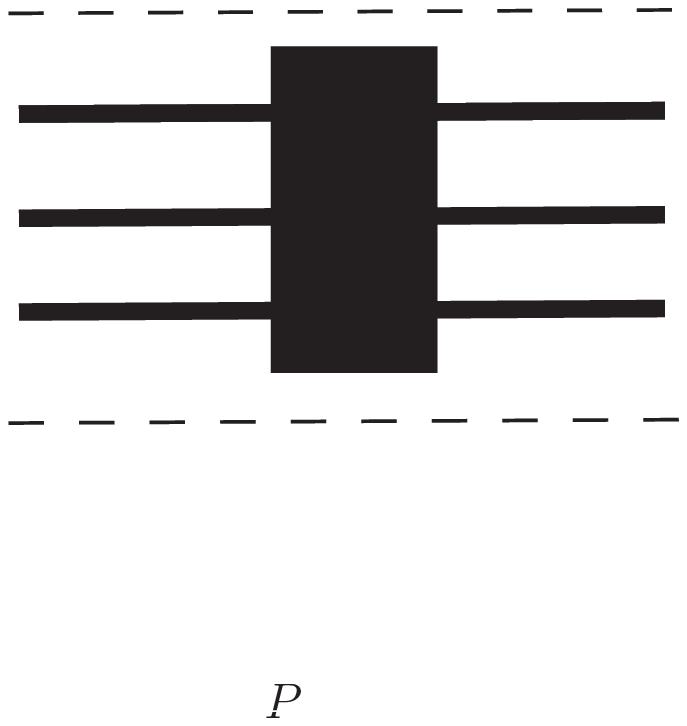}}\qquad .
\]
\item
Connect the strands at each vertex as in the picture
\[
\parbox{2cm}{
\psfrag{a}{${\scriptstyle a=6}$}
\psfrag{b}{${\scriptstyle b=7}$}\psfrag{c}{${\scriptstyle c=5}$}
\includegraphics[width=2cm]{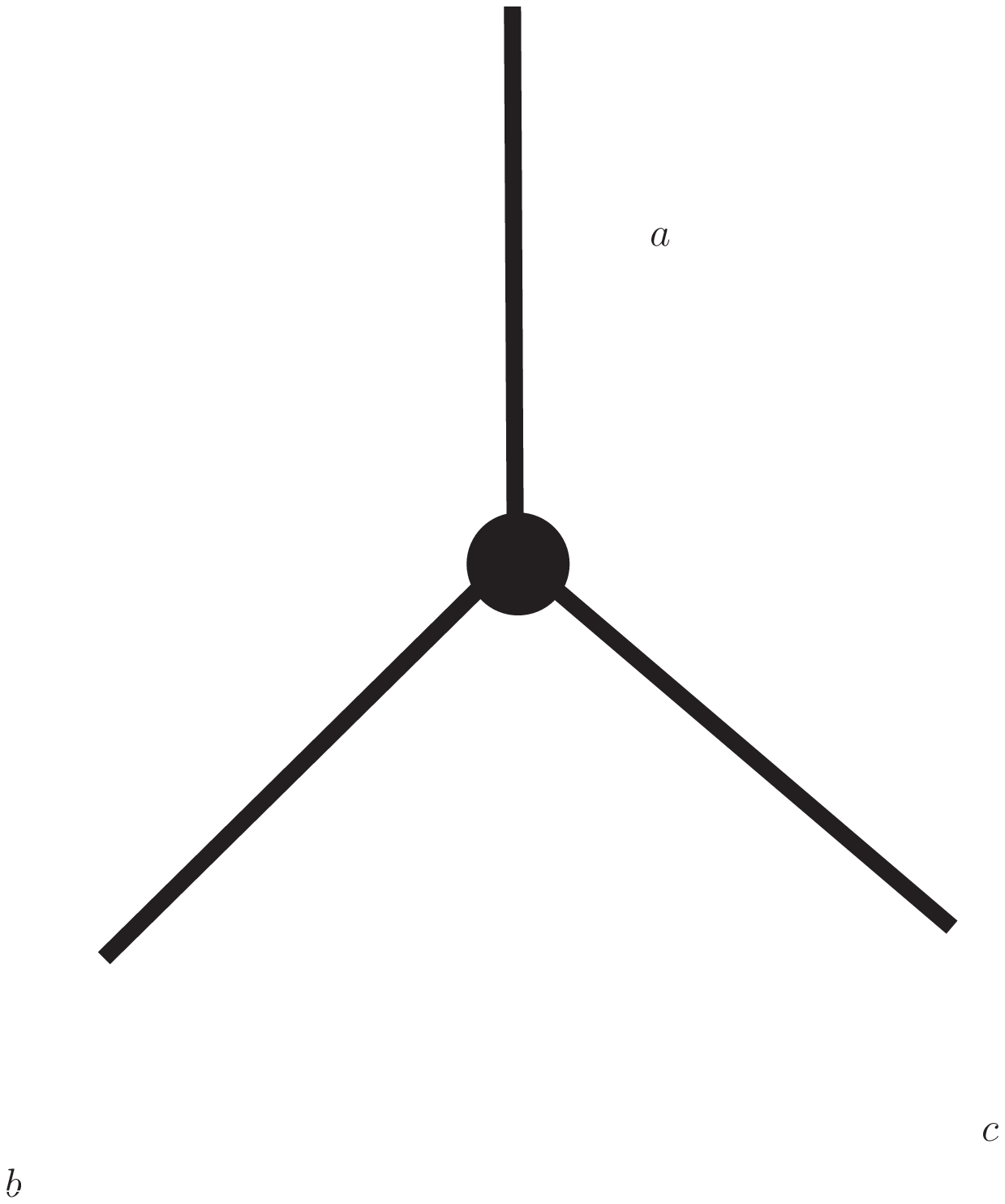}}
\qquad\longrightarrow\qquad
\parbox{3cm}{
\includegraphics[width=3cm]{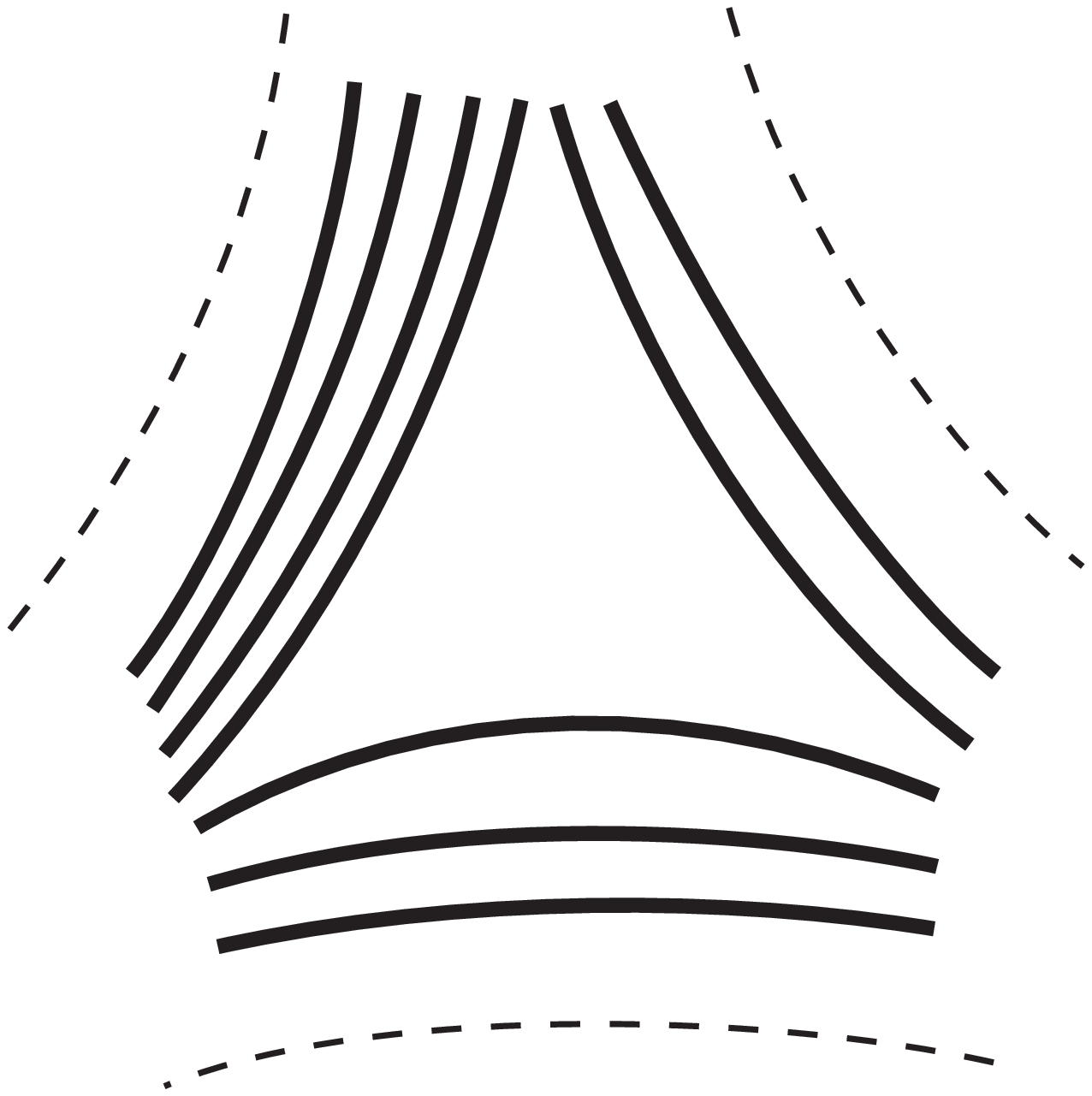}}
\]
seen from outside $\Si$. This connection is made possible by the admissibility
constraints. It is also essentially unique since the number
of strands connecting the $a$ edge to the $b$ edge is $\frac{a+b-c}{2}$, etc.
\item
Replace each Penrose bar $\parbox{0.7cm}{
\psfrag{P}{${\scriptstyle P}$}\psfrag{a}{${\scriptstyle a}$}
\includegraphics[width=0.7cm]{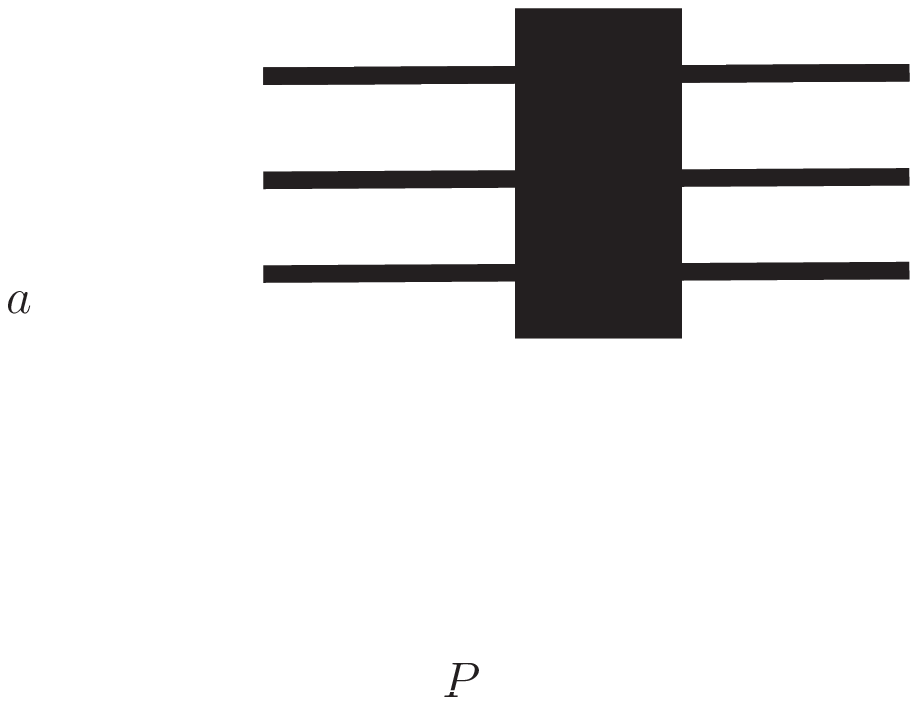}}$ by an alternating sum over
permutations $\si\in\gS_a$ as in
\[
\parbox{1.7cm}{
\psfrag{P}{${\scriptstyle P}$}\psfrag{a}{${\scriptstyle a=3}$}
\includegraphics[width=1.7cm]{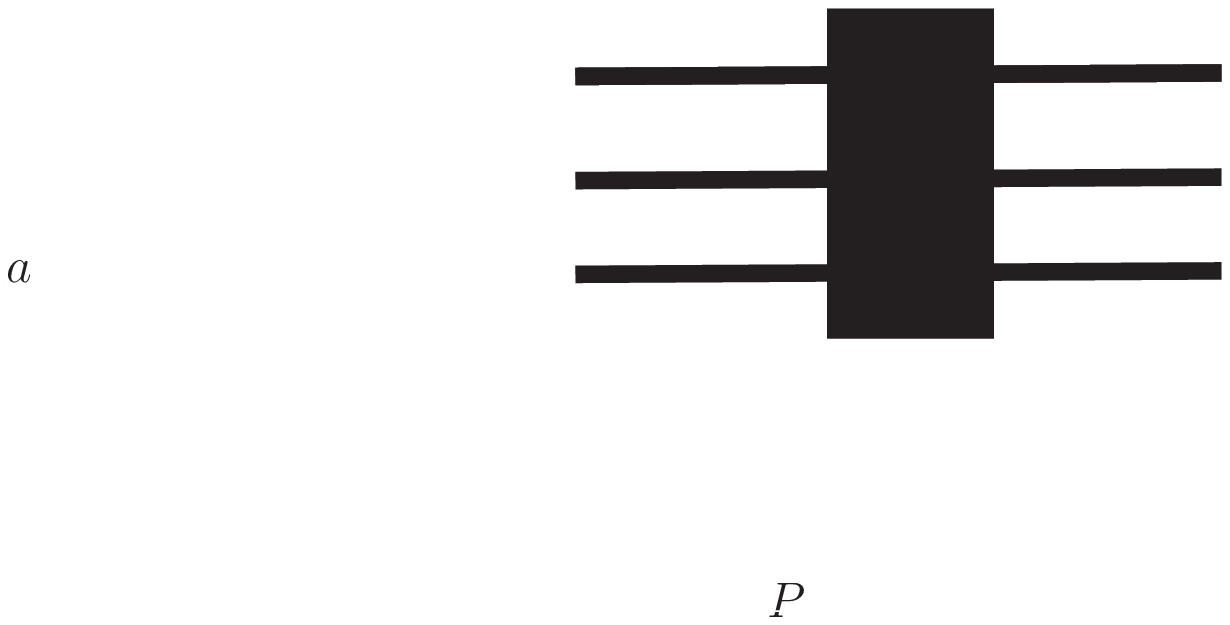}}
\longrightarrow
\parbox{0.8cm}{
\includegraphics[width=0.8cm]{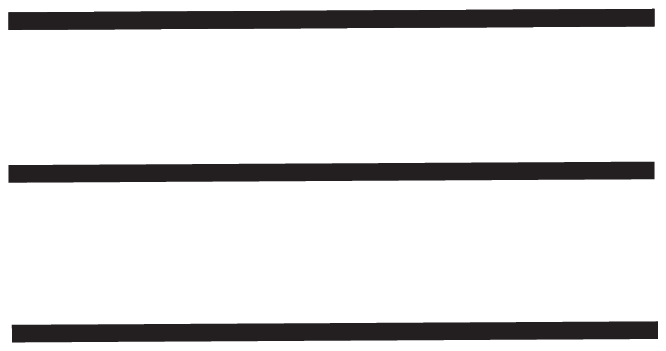}}
-\parbox{0.8cm}{
\includegraphics[width=0.8cm]{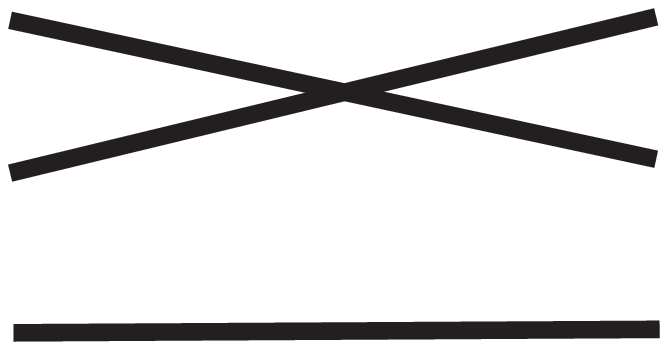}}
-\parbox{0.8cm}{
\includegraphics[width=0.8cm]{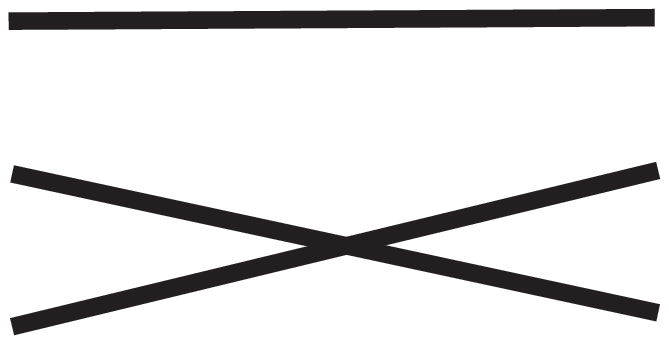}}
-\parbox{0.8cm}{
\includegraphics[width=0.8cm]{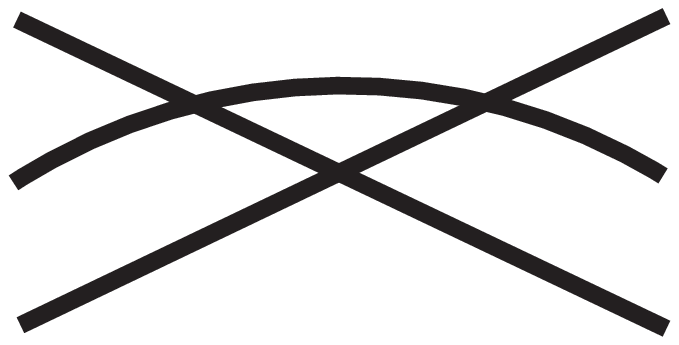}}
+\parbox{0.8cm}{
\includegraphics[width=0.8cm]{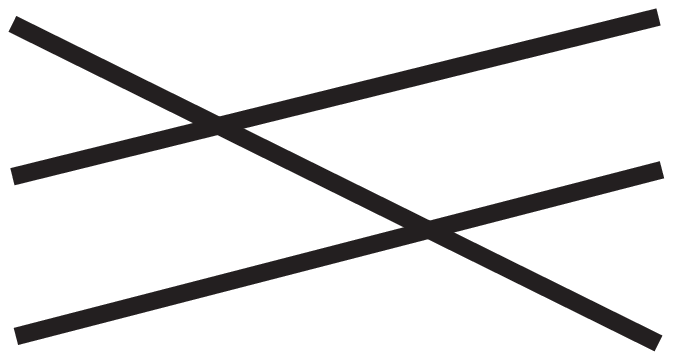}}
+\parbox{0.8cm}{
\includegraphics[width=0.8cm]{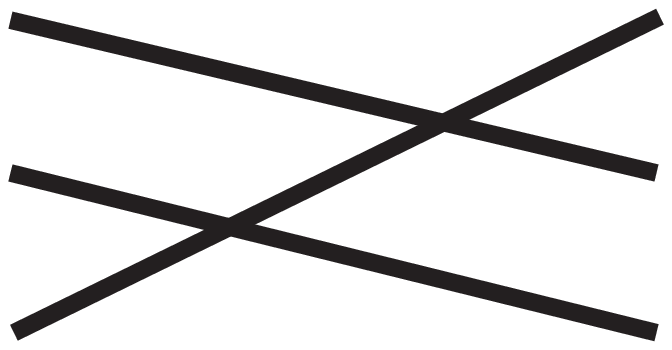}}
\qquad .
\]
\item
Finally each term in the sum over permutations will produce
a collection of closed curves drawn on the surface $\Si$. One associates
a factor (-2) to each such curve.
\item
The Penrose evaluation $\<\Ga,\ga\>^P$ is the result of summing the
corresponding (-2) to the power of the number of curves times
(-1) to the power of the number of crossings, over
all possible states or connection schemes specified by the permutations
$\si$ at each edge.
\end{enumerate}

\noindent{\bf Example:}
\[
\<\parbox{1.3cm}{\psfrag{2}{$\scriptstyle{2}$}
\includegraphics[width=1.3cm]{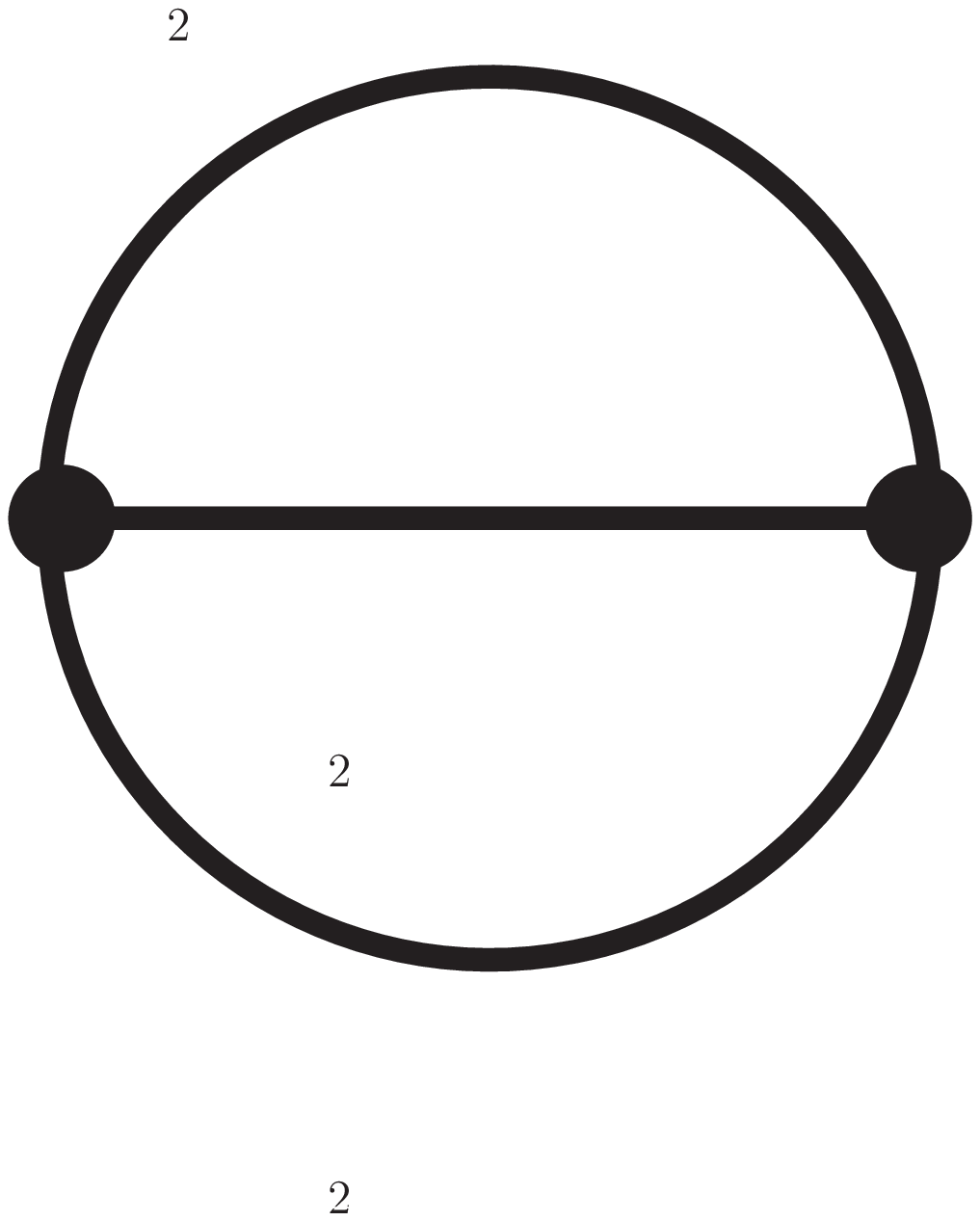}}
\>^P=\parbox{1.7cm}{\psfrag{P}{$\scriptstyle{P}$}\psfrag{p}{$\scriptscriptstyle{P}$}
\includegraphics[width=1.7cm]{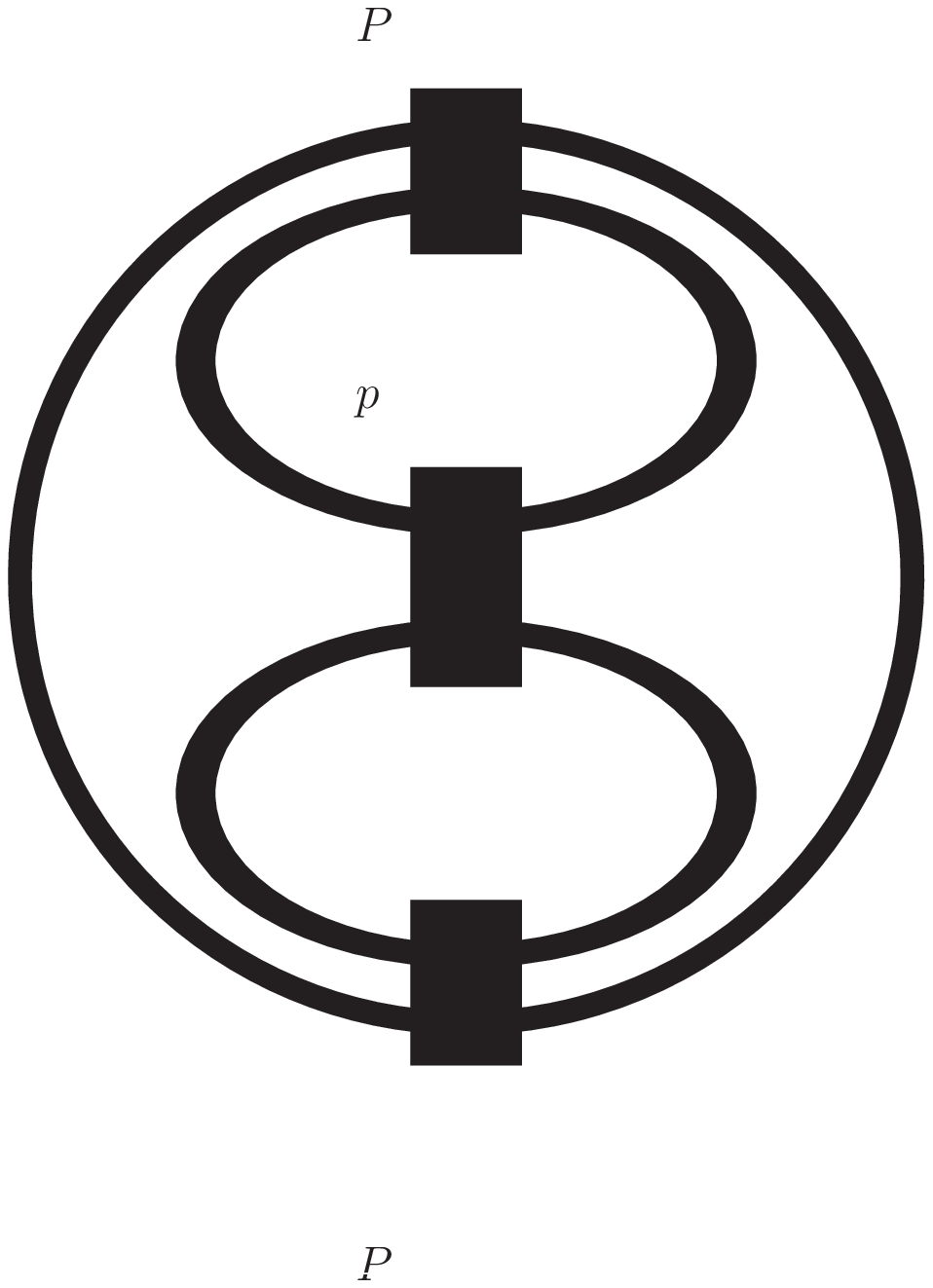}}
\]
\[
=\parbox{1.7cm}{
\includegraphics[width=1.7cm]{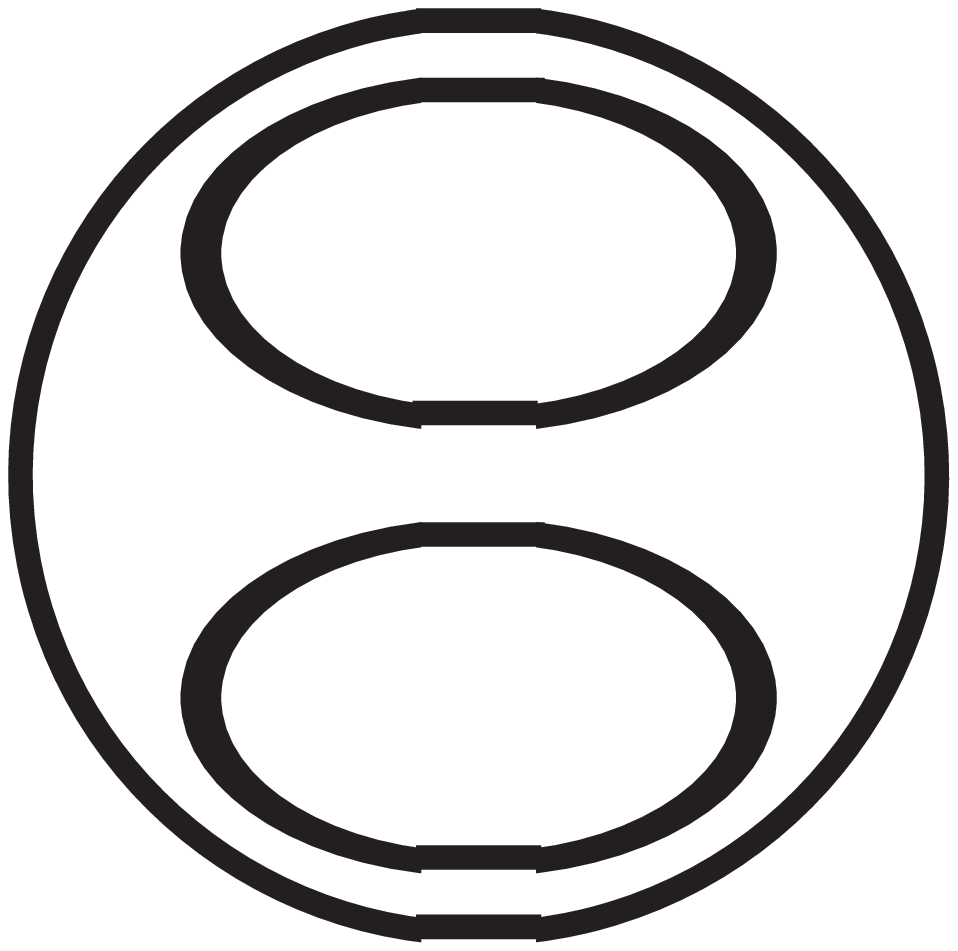}}
-\parbox{1.7cm}{
\includegraphics[width=1.7cm]{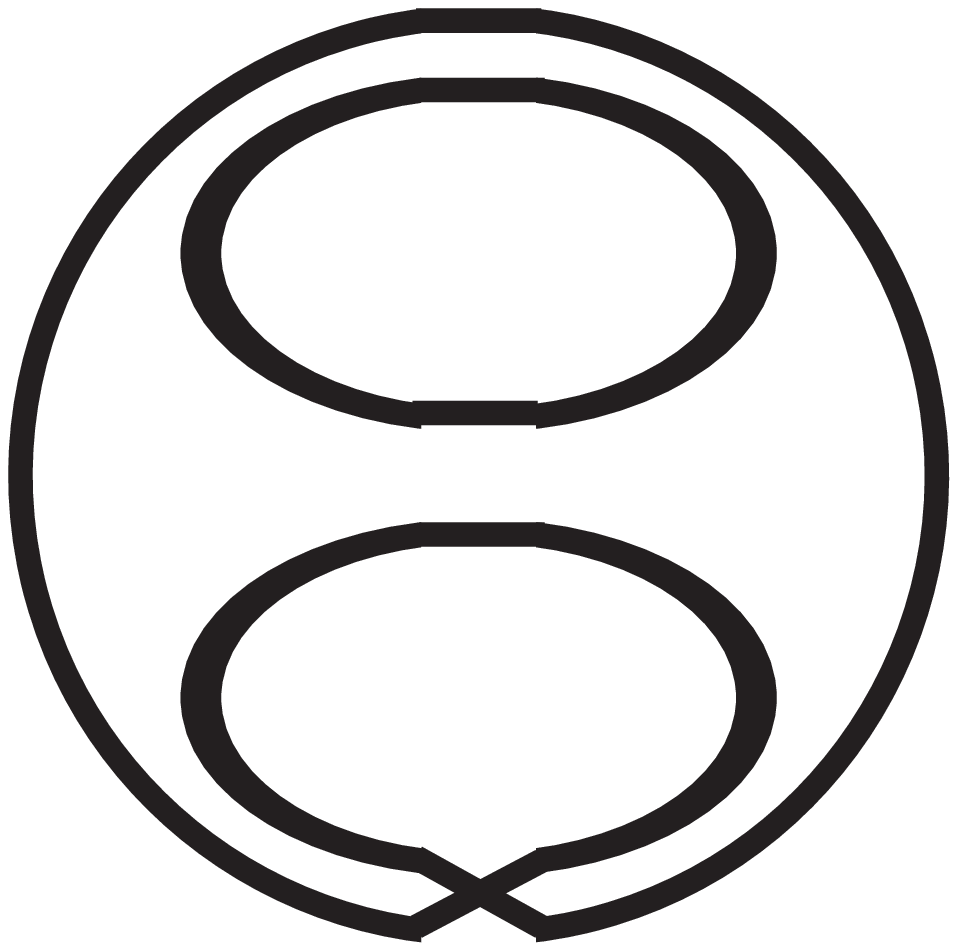}}
-\parbox{1.7cm}{
\includegraphics[width=1.7cm]{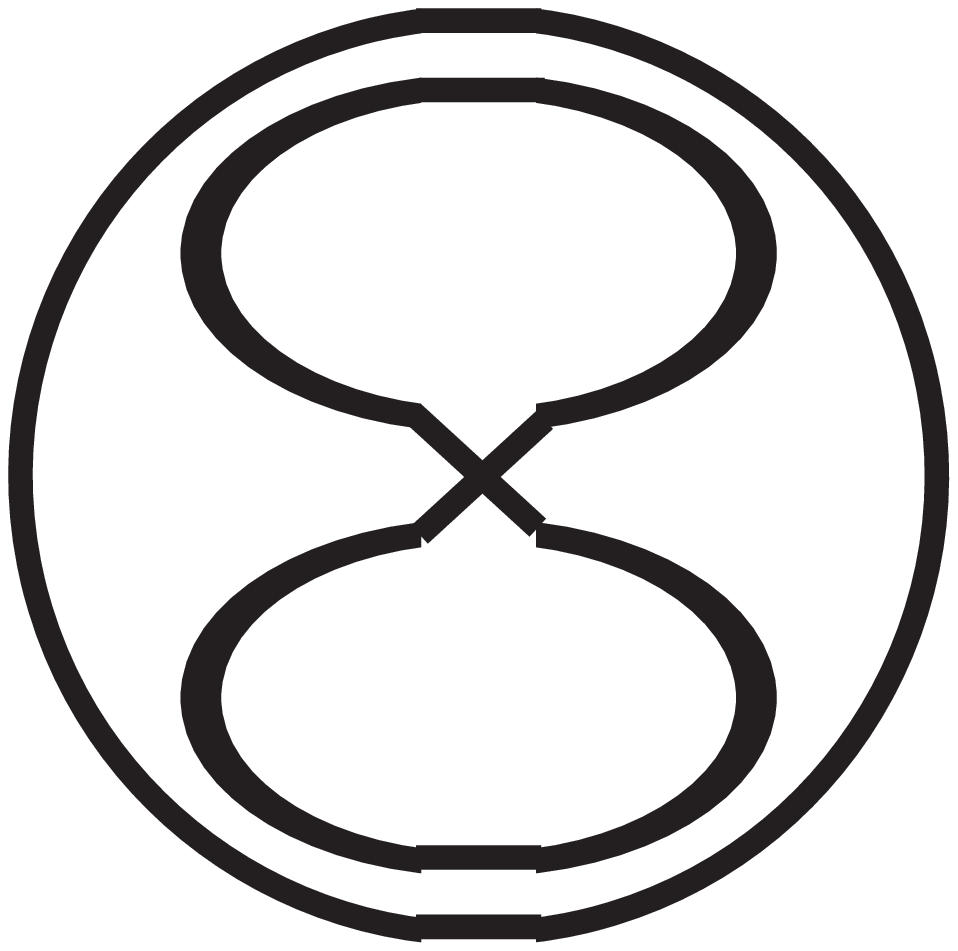}}
+\parbox{1.7cm}{
\includegraphics[width=1.7cm]{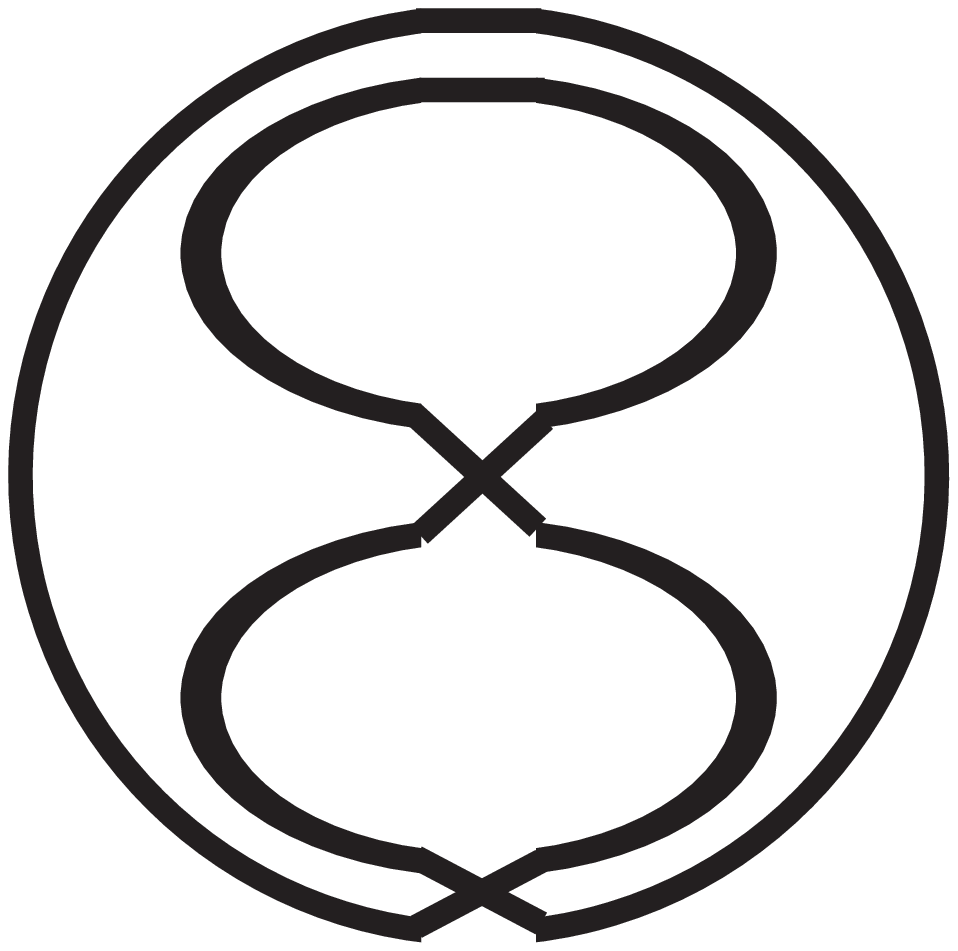}}
\]
\[
-\parbox{1.7cm}{
\includegraphics[width=1.7cm]{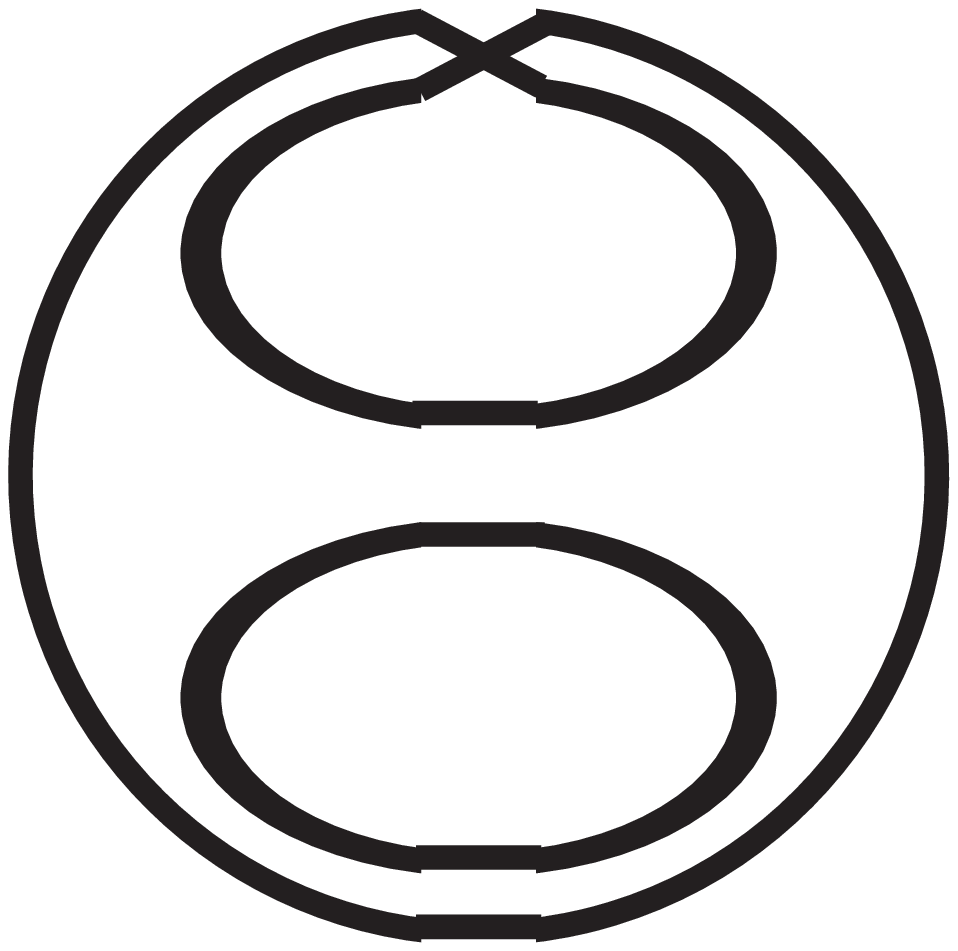}}
+\parbox{1.7cm}{
\includegraphics[width=1.7cm]{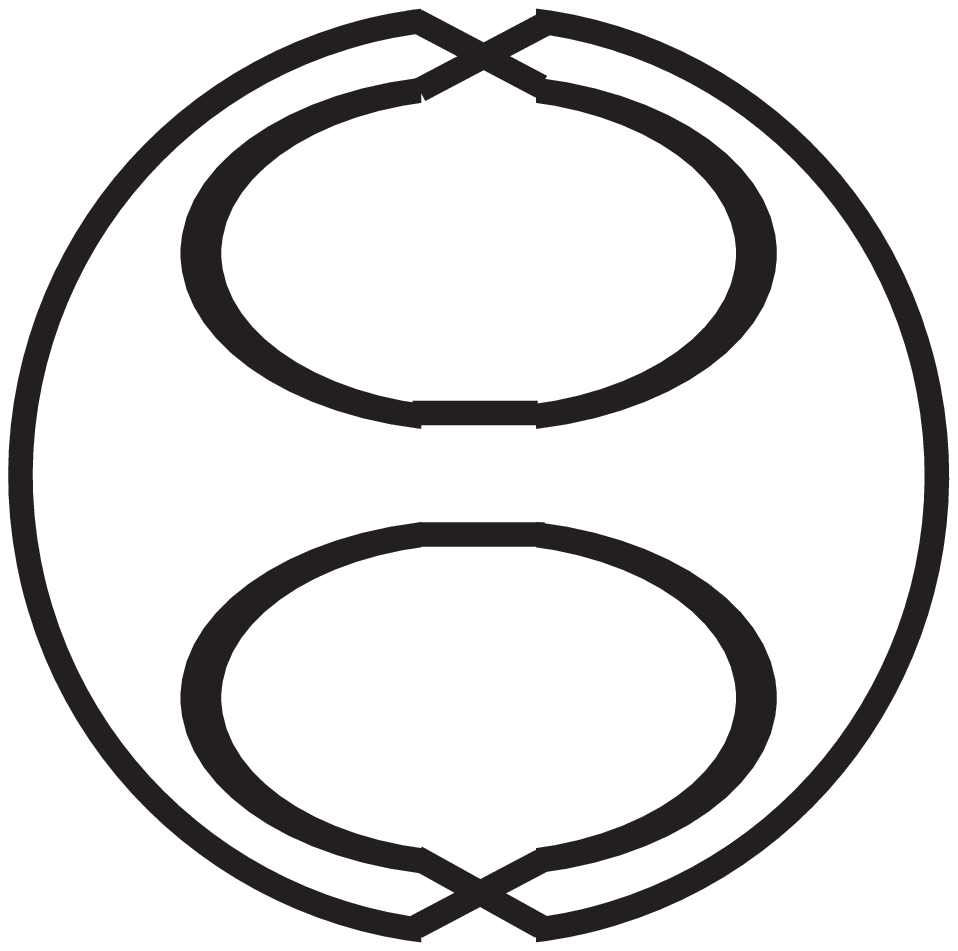}}
+\parbox{1.7cm}{
\includegraphics[width=1.7cm]{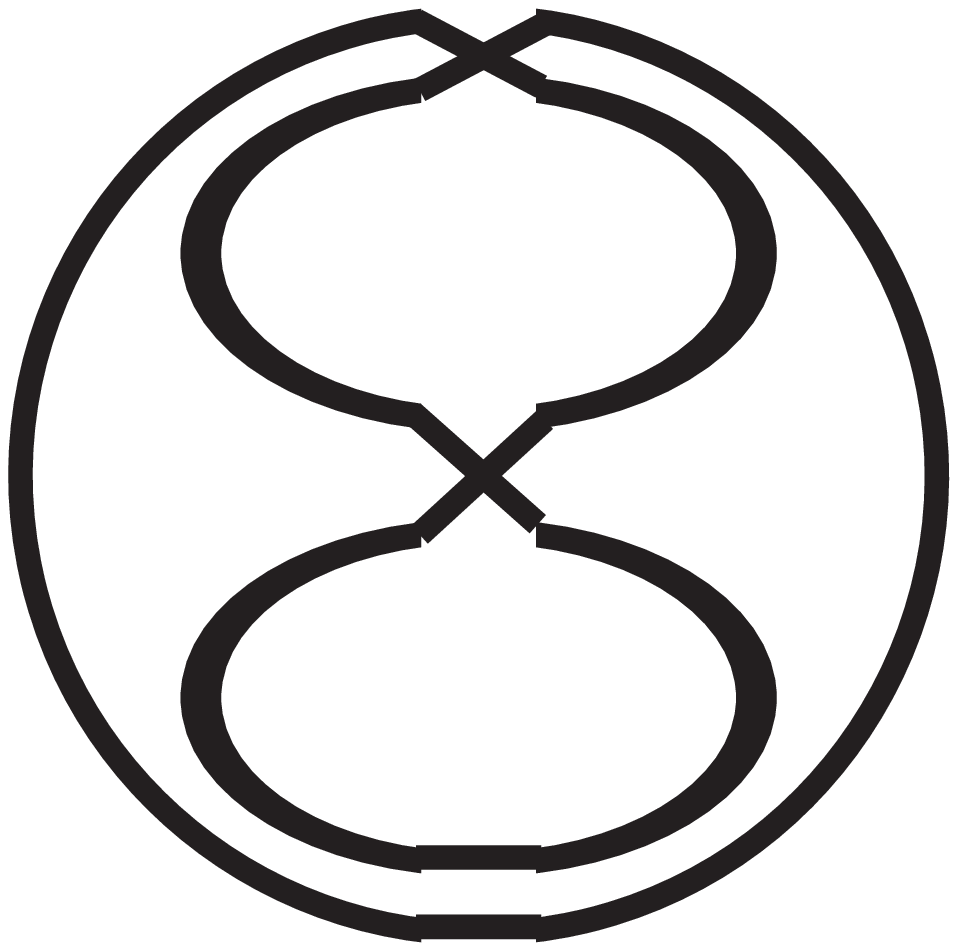}}
-\parbox{1.7cm}{
\includegraphics[width=1.7cm]{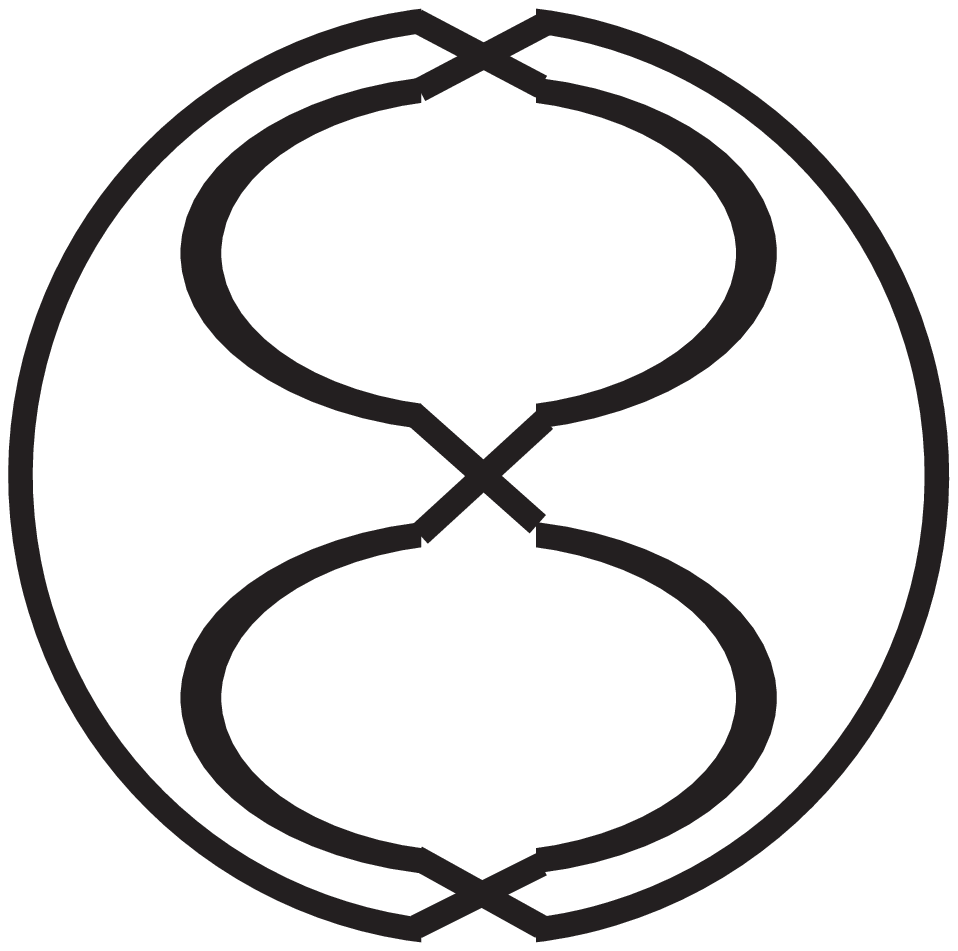}}
\]
\[
= (-2)^3-(-2)^2-(-2)^2+(-2)
\]
\[
-(-2)^2+(-2)+(-2)-(-2)
\]
\[
=-24\ .
\]
Note that we used the counterclockwise cyclic ordering at the two
vertices and the drawing on the plane of the written page, or equivalently
on the sphere.

The article~\cite{GvdV} also considers
the standard evaluation $\<\Ga,\ga\>^S$ of a spin network
\begin{equation}
\<\Ga,\ga\>^S=\<\Ga,\ga\>^P\times
\prod_{v\in V(\Ga)}
\left\{
\left(\frac{a_v+b_v-c_v}{2}\right)!
\left(\frac{a_v+c_v-b_v}{2}\right)!
\left(\frac{b_v+c_v-a_v}{2}\right)!
\right\}^{-1}
\label{Sdef}
\end{equation}
where $V(\Ga)$ is the vertex set of $\Ga$ and the
$a_v,b_v,c_v$ denote the decorations of the edges incident to vertex $v$.
One again counts a decoration twice in the case of a loop vertex.
One also defines the unitary evaluation $\<\Ga,\ga\>^U$
by
\begin{equation}
\<\Ga,\ga\>^U=\<\Ga,\ga\>^S\times
\prod_{v\in V(\Ga)}
\Th(a_v,b_v,c_v)^{-\frac{1}{2}}
\label{Udef}
\end{equation}
where
\[
\Th(a,b,c)=\frac{\left(\frac{a+b+c}{2}+1\right)!}
{\left(\frac{a+b-c}{2}\right)!
\left(\frac{a+c-b}{2}\right)!
\left(\frac{b+c-a}{2}\right)!}\ \ .
\]
\begin{Remark}
In~\cite[Def. 9.4]{GvdV}, due to a typo, (\ref{Udef}) is incorrectly
stated with $\<\Ga,\ga\>^P$ instead of $\<\Ga,\ga\>^S$.
\end{Remark}

In~\cite[Lem. 6.1]{GvdV} the following property is proved.

\begin{Lemma}
Changing the cyclic orientations at the vertices modifies any of the
previous evaluations by a sign.
More precisely, if one changes the cyclic ordering at a vertex with decorations
$a,b,c$, the resulting sign factor is
\[
(-1)^{\frac{a(a-1)+b(b-1)+c(c-1)}{2}}\ \ .
\]
\end{Lemma}

The case of trivial components, which is a good warm-up exercise on
Penrose evaluations, is dispensed with in the next easy lemma.

\begin{Lemma}\label{trivcomplemma}
\[
\<\parbox{1cm}{\psfrag{a}{$\scriptstyle{a}$}
\includegraphics[width=1cm]{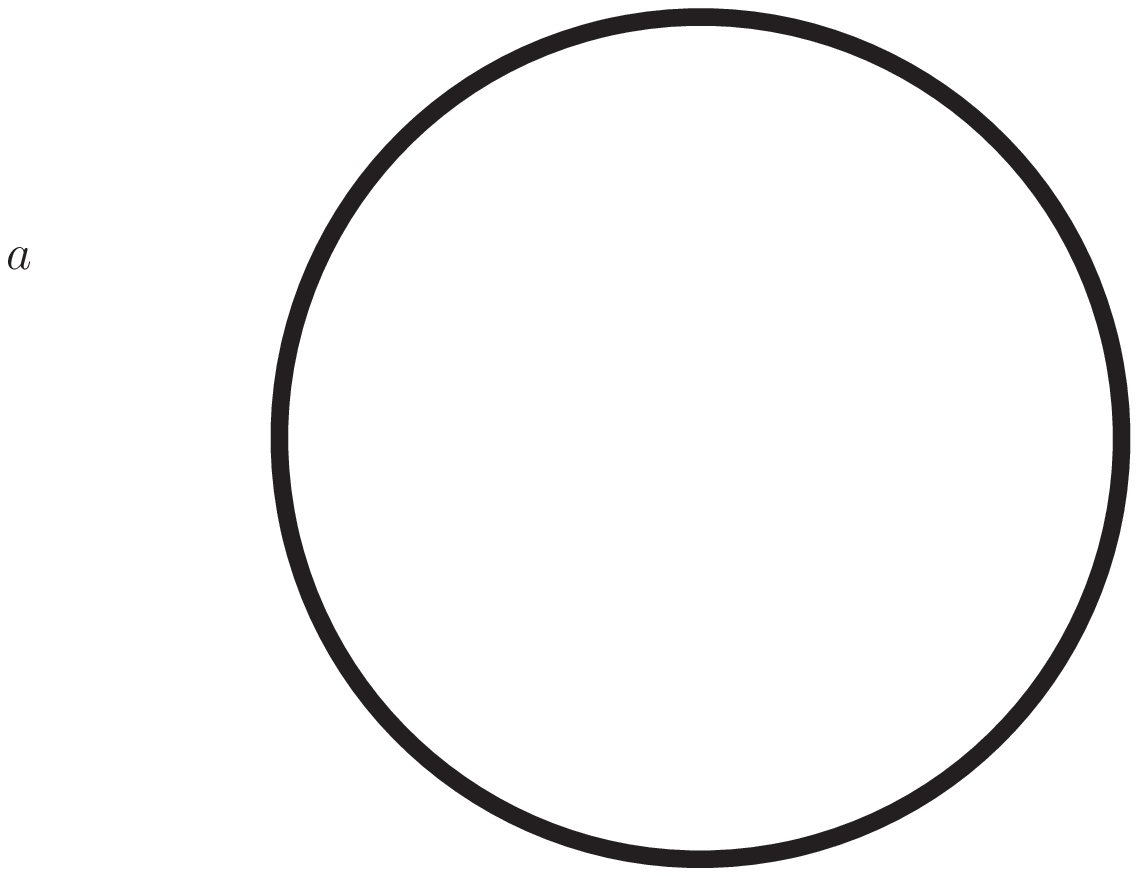}}
\>^P=\<\parbox{1cm}{\psfrag{a}{$\scriptstyle{a}$}
\includegraphics[width=1cm]{Fig28.eps}}
\>^S=\<\parbox{1cm}{\psfrag{a}{$\scriptstyle{a}$}
\includegraphics[width=1cm]{Fig28.eps}}
\>^U=(-1)^a (a+1)!\ \ .
\]
\end{Lemma}
\noindent{\bf Proof:}
One can use the chromatic method of Penrose and Moussouris and
a recursion such as~\cite[Eq. 6.20]{Cvitanovic},
but we prefer to use the definition.
With self-explanatory notations:
\[
\<\parbox{1cm}{\psfrag{a}{$\scriptstyle{a}$}
\includegraphics[width=1cm]{Fig28.eps}}\>^P=
\sum_{\si\in\gS_a} {\rm sign}(\si)\ 
\parbox{1.7cm}{\psfrag{a}{$\scriptstyle{\sigma}$}
\includegraphics[width=1.7cm]{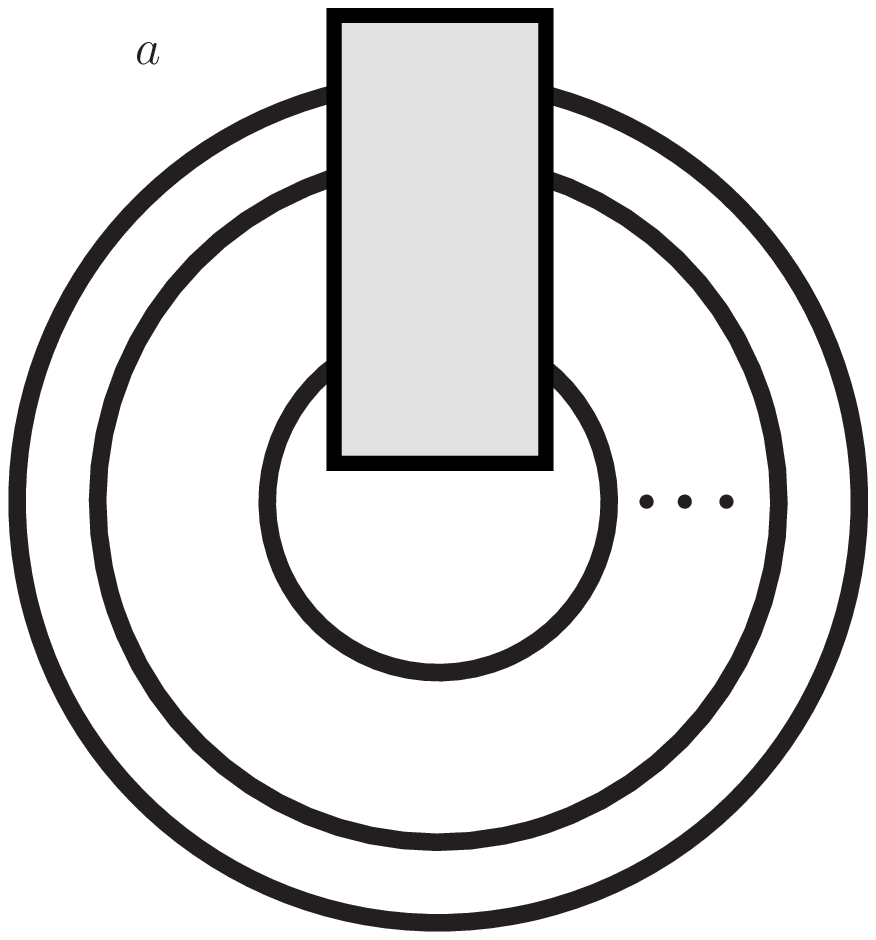}}
\]
\[
=\sum_{\si\in\gS_a} (-1)^{a-c(\si)} (-2)^{c(\si)}
\]
where $c(\si)$ is the number of cycles of the permutation $\si$.
Therefore,
\[
\<\parbox{1cm}{\psfrag{a}{$\scriptstyle{a}$}
\includegraphics[width=1cm]{Fig28.eps}}
\>^P=(-1)^a
\sum\limits_{k=0}^{a} 2^k c(a,k) 
\]
where $c(a,k)$ is the number of permutations in $\gS_a$ with exactly
$k$ cycles.
It is related to the Stirling number of the first kind $s(a,k)$ by
$c(a,k)=(-1)^{a-k}s(a,k)$.
By~\cite[Prop. 1.3.4]{Stanley}
one has
\[
\sum\limits_{k=0}^{a} 2^k c(a,k) =2(2+1)\cdots(2+a-1)=(a+1)!\ \ .
\]
\qed

Another trivial consequence of the previous definitions is the following.
\begin{Lemma}\label{factolemma}
The $\<\cdots\>^P$, $\<\cdots\>^S$ and $\<\cdots\>^U$ evaluations factorize over the connected
components of $\Ga$. 
\end{Lemma}

Given an admissible spin network $(\Ga,\ga)$, for any $n\in\N$,
the dilation $(\Ga,n\ga)$ where each decoration gets multiplied
by $n$ is also admissible.
The main problem addressed in~\cite{GvdV} and which goes back to~\cite{PonzanoR}
is the study of the asymptotics of evaluations of $(\Ga,n\ga)$ as
$n$ goes to infinity.
To this end, Garoufalidis and van der Veen introduced the power series
\[
F_{\Ga,\ga}(z)=\sum\limits_{n=0}^{\infty}
\<\Ga,n\ga\>^S z^n
\]
and defined the spectral radius $\rh_{\Ga,\ga}\in[0,\infty]$
of the spin network as the inverse of the radius of convergence of the series
$F_{\Ga,\ga}$.
Among the challenge problems mentioned in~\cite{GvdV},
Problem 2 therein is the statement that $\<\Ga,\ga\>^U$
is bounded in absolute value by one.
Problem 3 therein is the statement that for uniform decorations
$\ga(e)=2$ for every edge $e$, the spectral radius is exactly
$3^{\frac{3|V(\Ga)|}{2}}$ where we used
notation $|\cdot|$ for the cardinality of finite sets.
This is a hopefully easier analogue for CSN's of the volume conjecture
in knot theory~\cite{Kashaev,MurakamiM}. Indeed the decorations play
a role similar to that of the colors of the colored Jones polynomial.
We prefer the ``decoration'' terminology rather than the ``coloring'' one
used in~\cite{GvdV}.

Clearly, trivial components violate both statements completely
and must be excluded.
Another interesting practice example is the dumbell which evaluates to
\begin{equation}
\<\parbox{2.2cm}{\psfrag{a}{$\scriptstyle{a}$}
\psfrag{b}{$\scriptstyle{b}$}\psfrag{c}{$\scriptstyle{c}$}
\includegraphics[width=2.2cm]{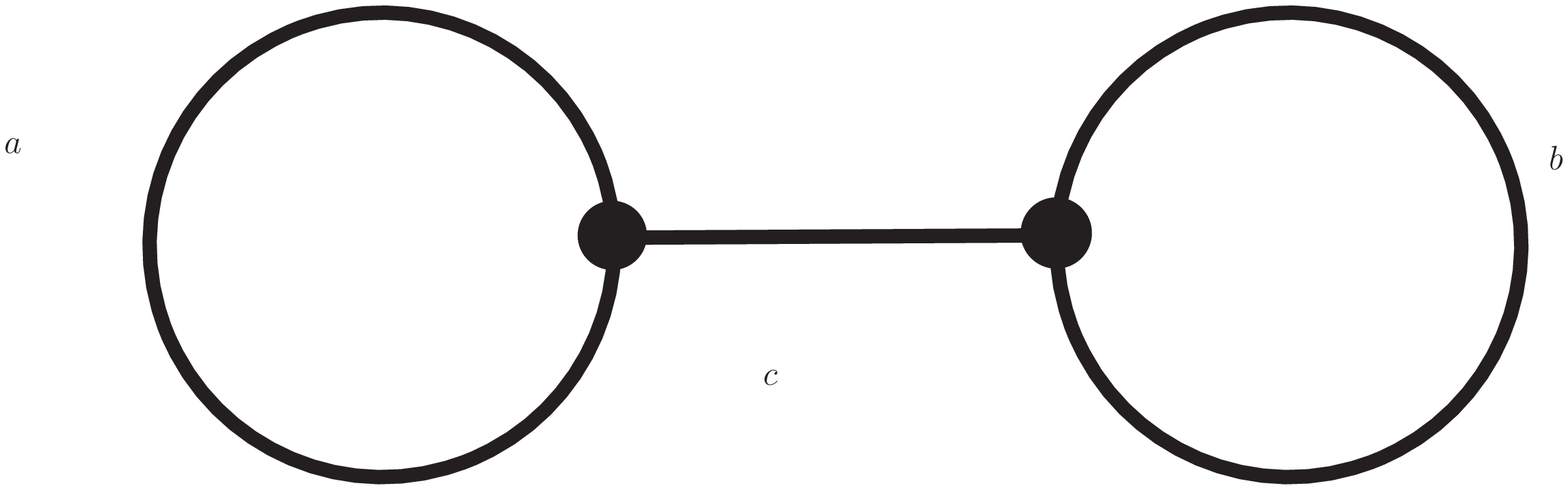}}
\>^U=\de_{c,0} (-1)^{a+b}\sqrt{(a+1)(b+1)}\ \ ,
\label{dumbell}
\end{equation}
see Remark \ref{dumbellrk}.
More generally, loops which carry sufficiently large decorations
can cause the failure of the bound $|\<\Ga,\ga\>^U|\le 1$.
Bridges also spoil the conjecture for the value
of the spectral radius because of the Kronecker delta factor as
in (\ref{dumbell}).

Our main result is the solution of Problem 2 in~\cite{GvdV}.

\begin{Theorem}\label{mainthm}
For any spin network $(\Ga,\ga)$ without trivial components
and without loops, one has
\[
|\<\Ga,\ga\>^U|\le 1\ .
\]
\end{Theorem}

For the case of graphs with loops one has the following weaker
statement.

\begin{Theorem}\label{thmwithloops}
For any spin network $(\Ga,\ga)$ without trivial components,
$|\<\Ga,n\ga\>^U|$ grows at most polynomially with $n$.
\end{Theorem}

As to Problem 3 in~\cite{GvdV}, we can state the following
corollary which via Stirling's formula is an easy consequence of
Theorem \ref{thmwithloops}.

\begin{Corollary}
For any spin network $(\Ga,\ga)$ without trivial components,
the spectral radius satisfies the bound
\[
\rh_{\Ga,\ga}\le \prod_{v\in V(\Ga)}\sqrt{\be(a_v,b_v,c_v)}
\]
where
\[
\be(a,b,c)=
\frac{\left(\frac{a+b+c}{2}\right)^{\left(\frac{a+b+c}{2}\right)}}
{
\left(\frac{a+b-c}{2}\right)^{\left(\frac{a+b-c}{2}\right)}
\left(\frac{a+c-b}{2}\right)^{\left(\frac{a+c-b}{2}\right)}
\left(\frac{b+c-a}{2}\right)^{\left(\frac{b+c-a}{2}\right)}
}
\]
and using the convention $0^0=1$ so
that degenerate cases are covered as well.
\end{Corollary}

Note that $\be(2,2,2)=9$ so in the case $\ga \equiv 2$
we have `half' of
the volume conjecture for CSN's.

\begin{Corollary}
For any spin network $(\Ga,\ga)$ without trivial components,
\[
\rh_{\Ga,\ga\equiv 2}\le 3^{\frac{3|V(\Ga)|}{2}}\ \ .
\]
\end{Corollary}

The full conjecture should be as follows.

\begin{Conjecture}\label{mainconj}
For any spin network $(\Ga,\ga)$ without trivial components,
and without bridges, i.e., which is 2-edge-connected
\[
\rh_{\Ga,\ga\equiv 2}= 3^{\frac{3|V(\Ga)|}{2}}\ \ .
\]
\end{Conjecture}

We were able to prove the equality for graphs of the form
\[
\parbox{3cm}{\includegraphics[width=3cm]{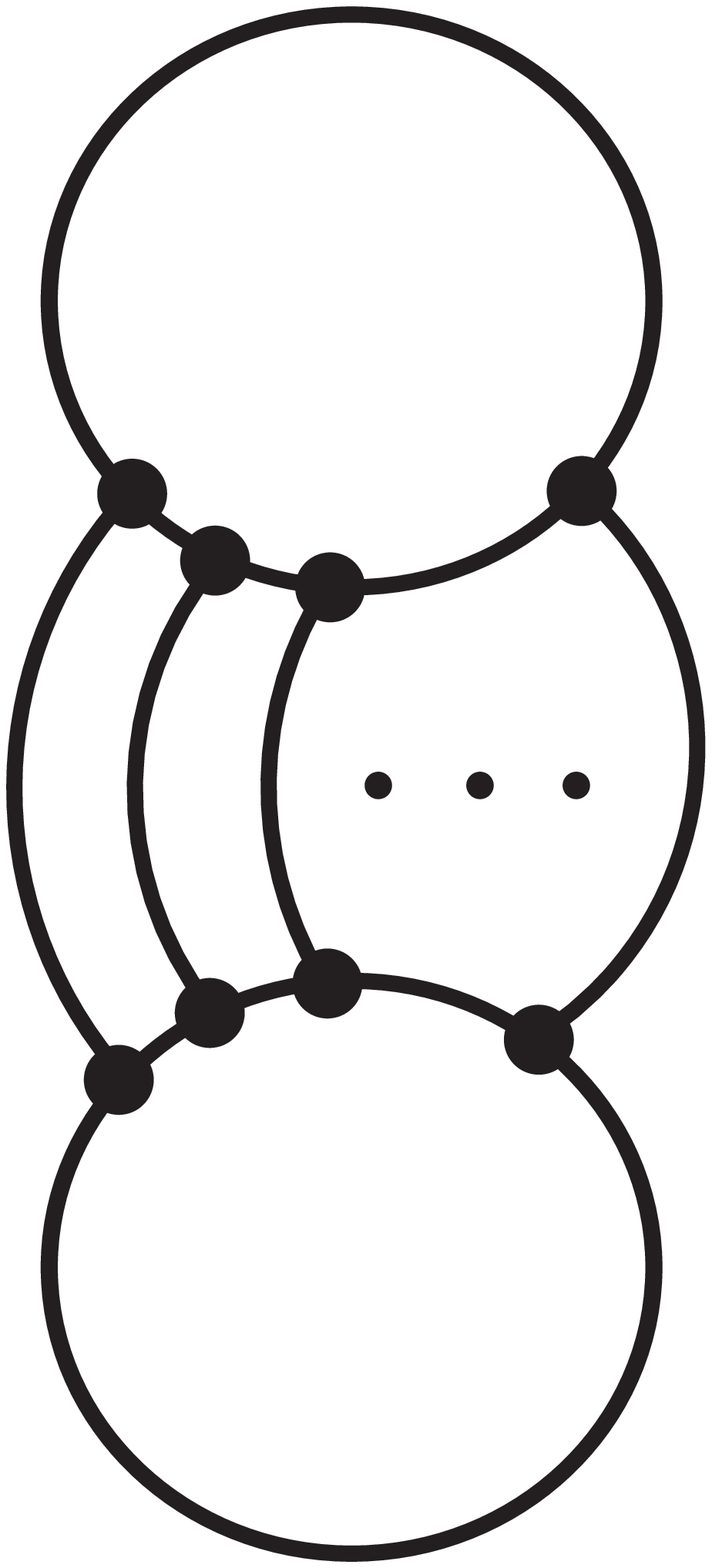}}
\]
with $s$ edges between the two circles. These are called generalized drum
graphs and are denoted by ${\rm Drum}_s$ in~\cite{GvdV}.

\begin{Theorem}\label{drumthm}
For any $s\ge 2$, 
\[
\rh_{\rm Drum_s,\ga\equiv 2}=3^{3s}\ \ .
\]
\end{Theorem}

\noindent
This of course covers the case of the cube ($s=4$)
which is new.

\section{Clebsch-Gordan networks and a brief tour of classical invariant theory}
\label{CGsection}
Let $G$ be a cubic graph, possibly disconnected, which may contain multiple
edges and loops.
However, we exclude trivial vertex-less components in this section.
A smooth orientation $\cO$ of $G$
is an orientation of the edges of $G$ such that the resulting digraph
only has two types of vertices:
\[
\parbox{1.6cm}{\includegraphics[width=1.6cm]{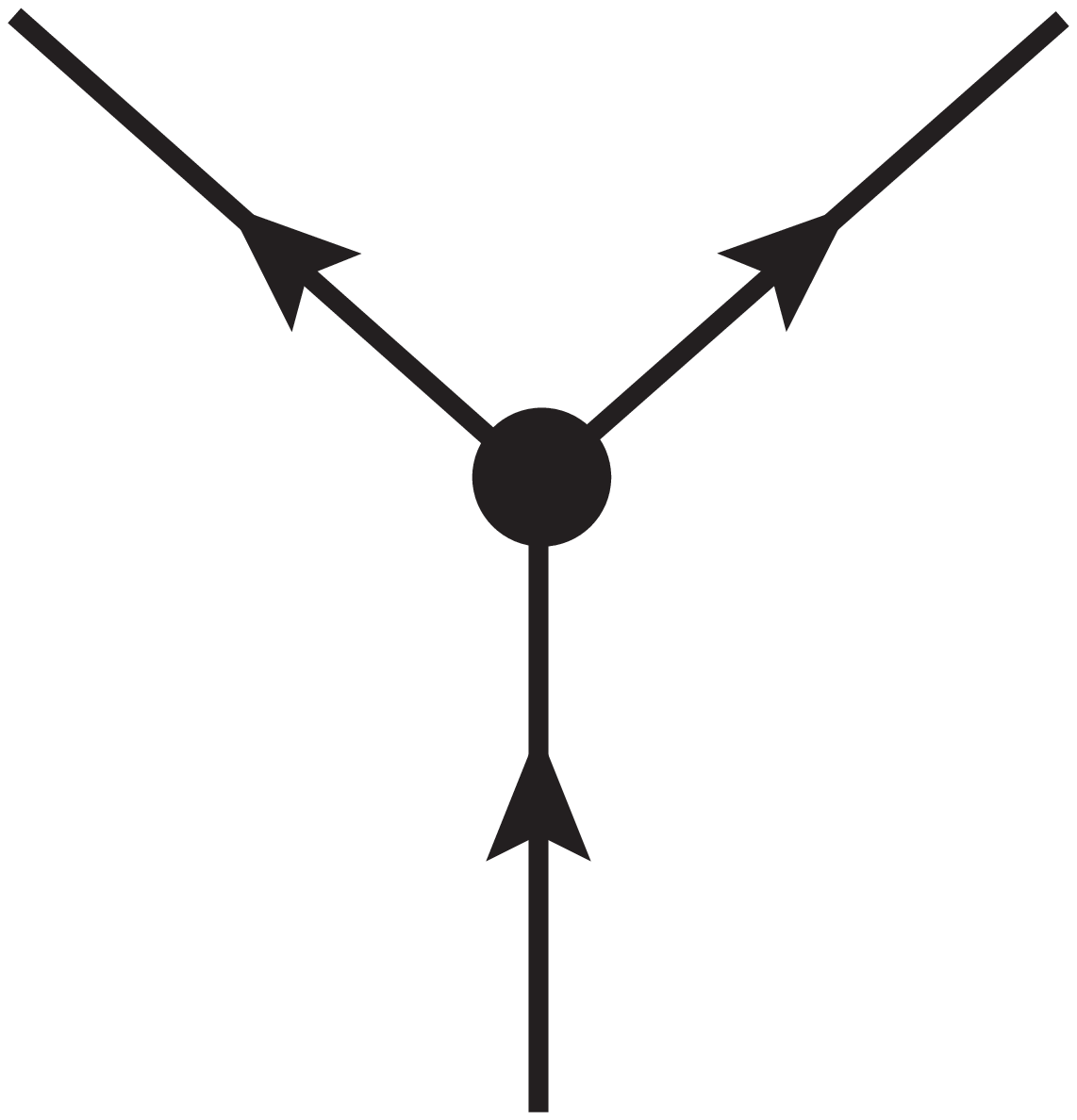}}
\qquad {\rm and}\qquad
\parbox{1.6cm}{\includegraphics[width=1.6cm]{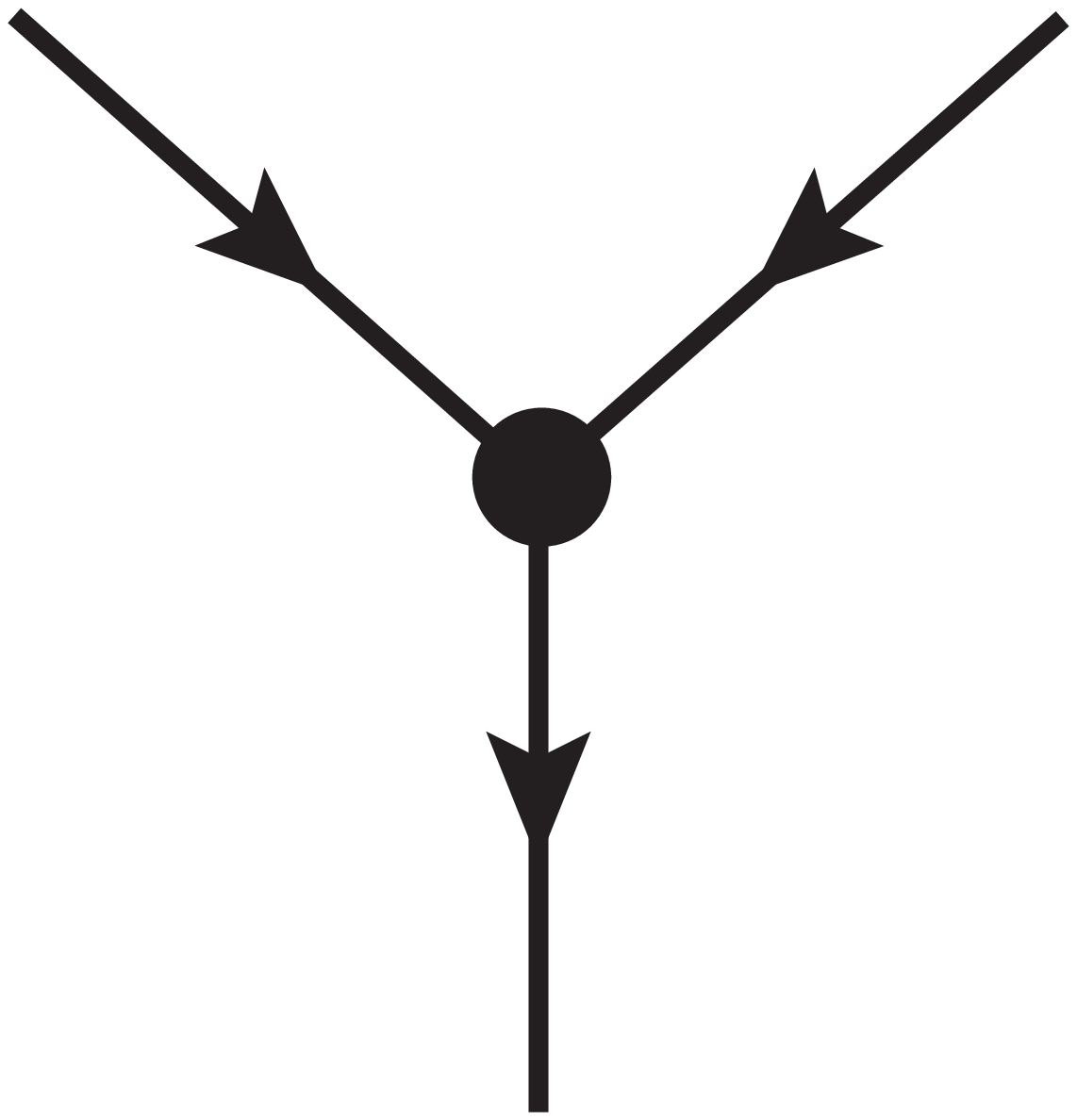}}\qquad .
\]
A gate signage $\ta$ corresponds to an ordering for every vertex of the two
edges (or rather half-edges in order to cover the loop case as well)
which share the same direction.
Such pair of half-edges is called a gate.
The ordering is indicated by a small curved arrow as in:
\[
\parbox{1.6cm}{\includegraphics[width=1.6cm]{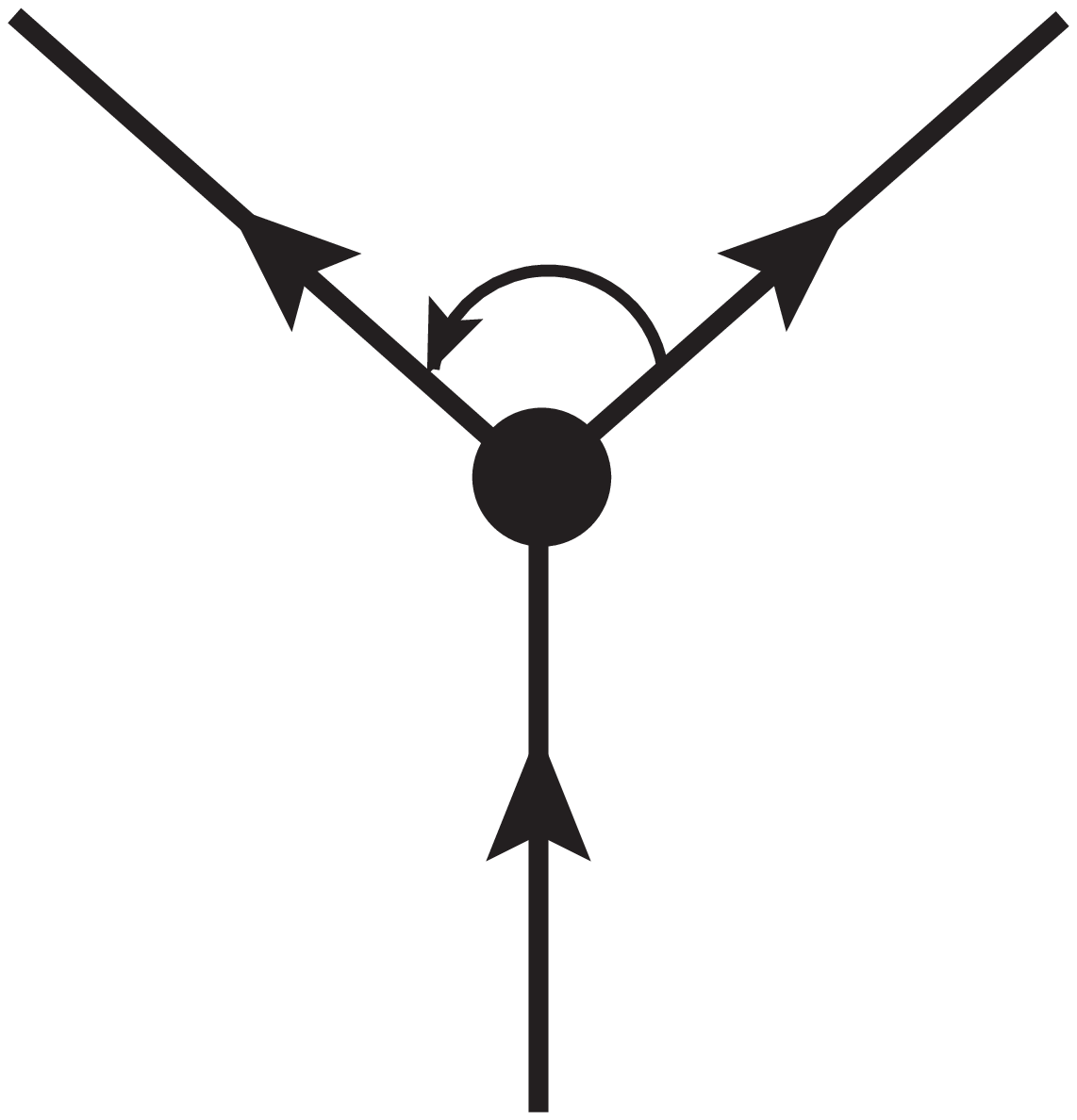}}
\qquad{\rm or}\qquad
\parbox{1.6cm}{\includegraphics[width=1.6cm]{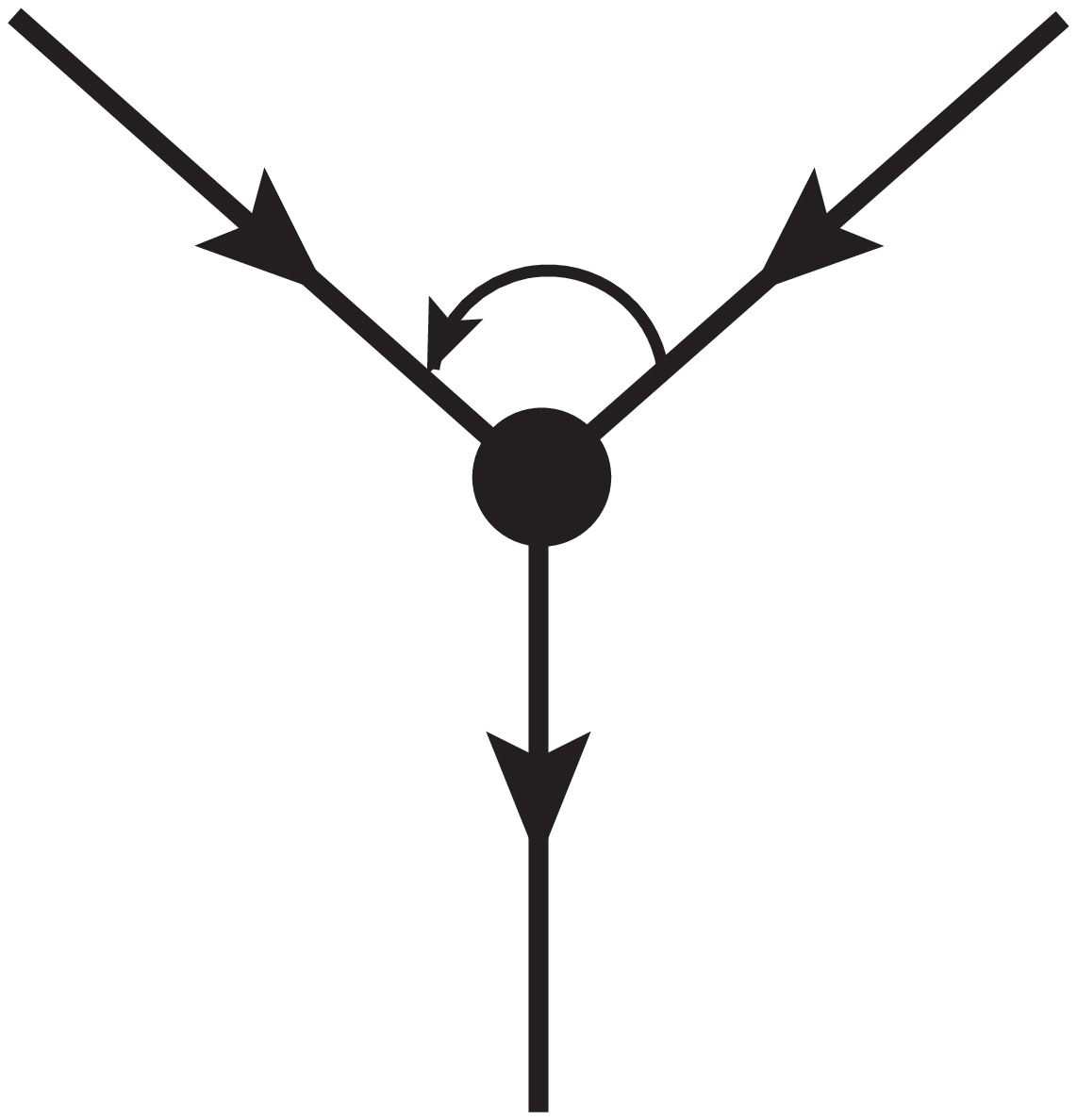}}\qquad .
\]
We also consider a decoration $\ga$ of the edges by nonnegative integers
satisfying the same admissibility conditions as in \S\ref{introdefsec}.
To the data $(G,\cO,\ta,\ga)$ which we call a Clebsch-Gordan or CG network
we will associate a number $\<G,\cO,\ta,\ga\>^{CG}$
using FDC, i.e., a graphical encoding of tensor contractions
(see, e.g., \cite{AbdesselamSLC,AC1}).
These tensors typically belong to spaces of the the form
$V\otimes V\otimes V^\ast\otimes V^\ast\otimes V\cdots$
where the fundamental vector space is $V=\C^2$.
The FDC has no need for a `coordinate free' approach, since it is
to modern tensor algebra what matrix algebra is to abstract linear algebra.
We think of tensors simply as arrays or `matrices' of numbers $T_{i_1 i_2\cdots}$
with any number of indices taking their values in $\{1,2\}$.
One can of course extend the scalars and allow entries which are polynomials
in formal letters, but here we will mostly work with complex numbers.

\begin{Remark}
Note that the FDC used in quantum field theory is concerned with the
situation where $V$ is infinite-dimensional say $V=L^2(\R^d)$,
tensors become integral kernels and sums over indices become
integration over $\R^d$.
Since the kernels involved are singular, this poses a nontrivial
problem of mathematical analysis which is the object of
renormalization theory. It is a far-reaching extension of Schwartz's
kernel theory which would be enough if the `sums over indices' or rather
$L^2(\R^d)$ pairings were always between $S(\R^d)$ and $S'(\R^d)$
tensor factors.
\end{Remark}

A vector $\ux=(x_1,x_2)$ in $V$ is denoted graphically by
\[
x_i=
\parbox{0.5cm}{\psfrag{x}{$\scriptstyle{x}$}\psfrag{i}{$\scriptstyle{i}$}
\includegraphics[width=0.5cm]{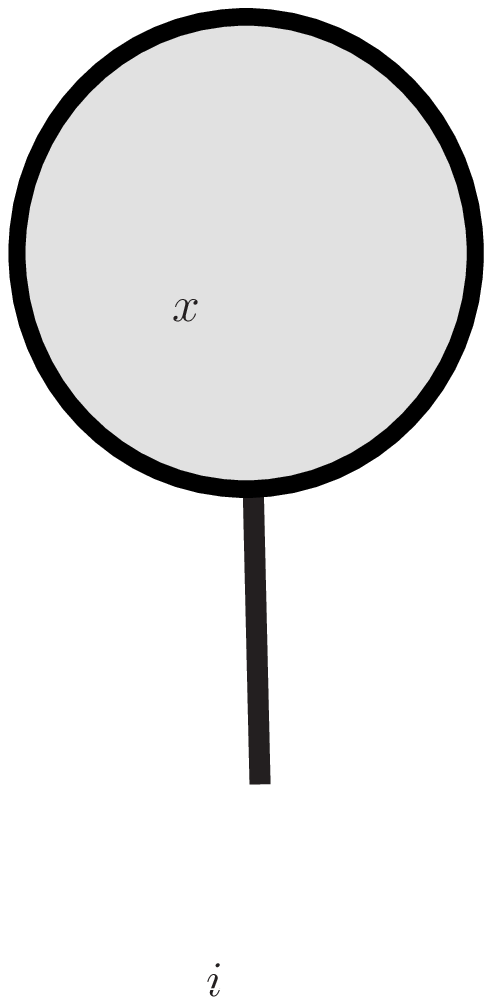}}
\qquad,\qquad i=1,2\ .
\]
A binary form $F$ of order $r$ is a polynomial in $\ux$,
homogeneous of degree $r$, which following Cayley we write using binomial
coefficients as:
\[
F(\ux)=\sum\limits_{p=0}^{r} \left(
\begin{array}{c}
r\\
p
\end{array}
\right)
f_p\ x_1^{r-p} x_2^p\ .
\]
Equivalently,
\[
F(\ux)=\sum\limits_{i_1,\ldots,i_r=1}^{2}
F_{i_1\ldots i_r}\ x_{i_1}\ldots x_{i_r}
\]
where the tensor $(F_{i_1\ldots i_r})_{1\le i_1,\ldots,i_r\le 2}$
is completely symmetric in its
$r$ indices and is related to the previous description by
$F_{i_1\ldots i_r}=f_p$
where $p$ is the number of indices which happen to be equal to 2.
Introduce the graphical notation
\[
\parbox{1.2cm}{\psfrag{F}{$\scriptstyle{F}$}\psfrag{1}{$\scriptstyle{i_1}$}
\psfrag{2}{$\scriptstyle{i_2}$}\psfrag{r}{$\scriptstyle{i_r}$}
\includegraphics[width=1.2cm]{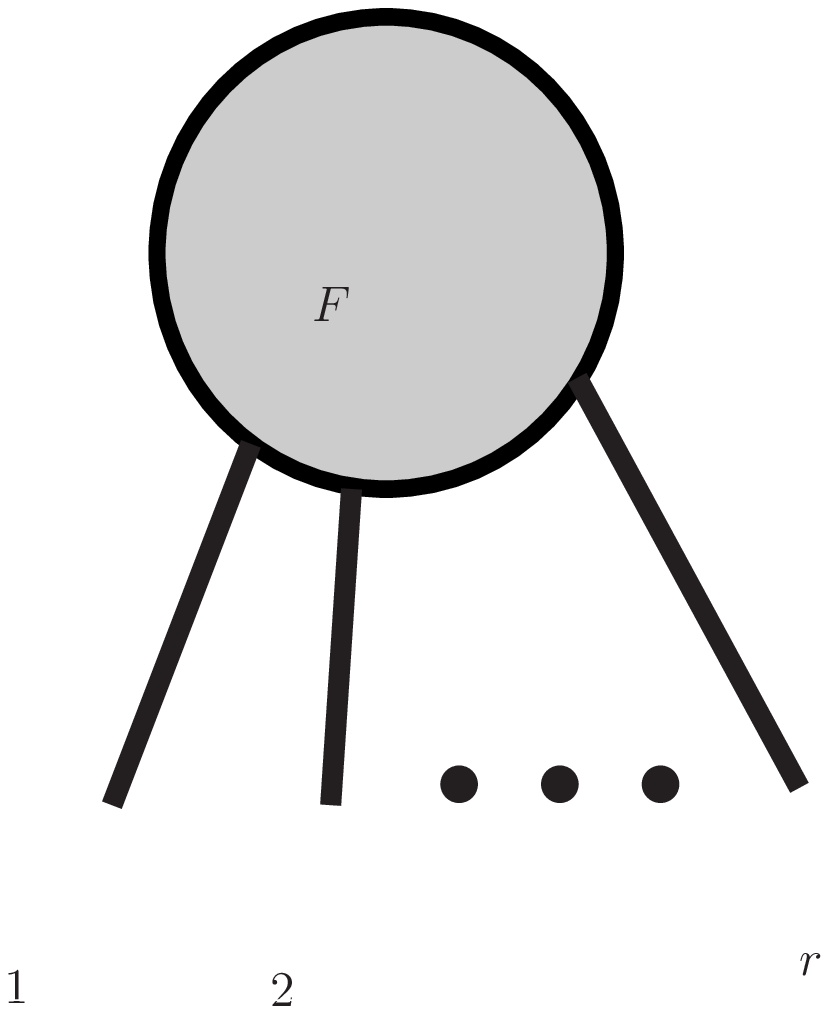}}
= F_{i_1\ldots i_r}\ \ .
\]
Then
\[
F(\ux)=\underbrace{
\parbox{1.5cm}{\psfrag{F}{$\scriptstyle{F}$}\psfrag{x}{$\scriptstyle{x}$}
\includegraphics[width=1.5cm]{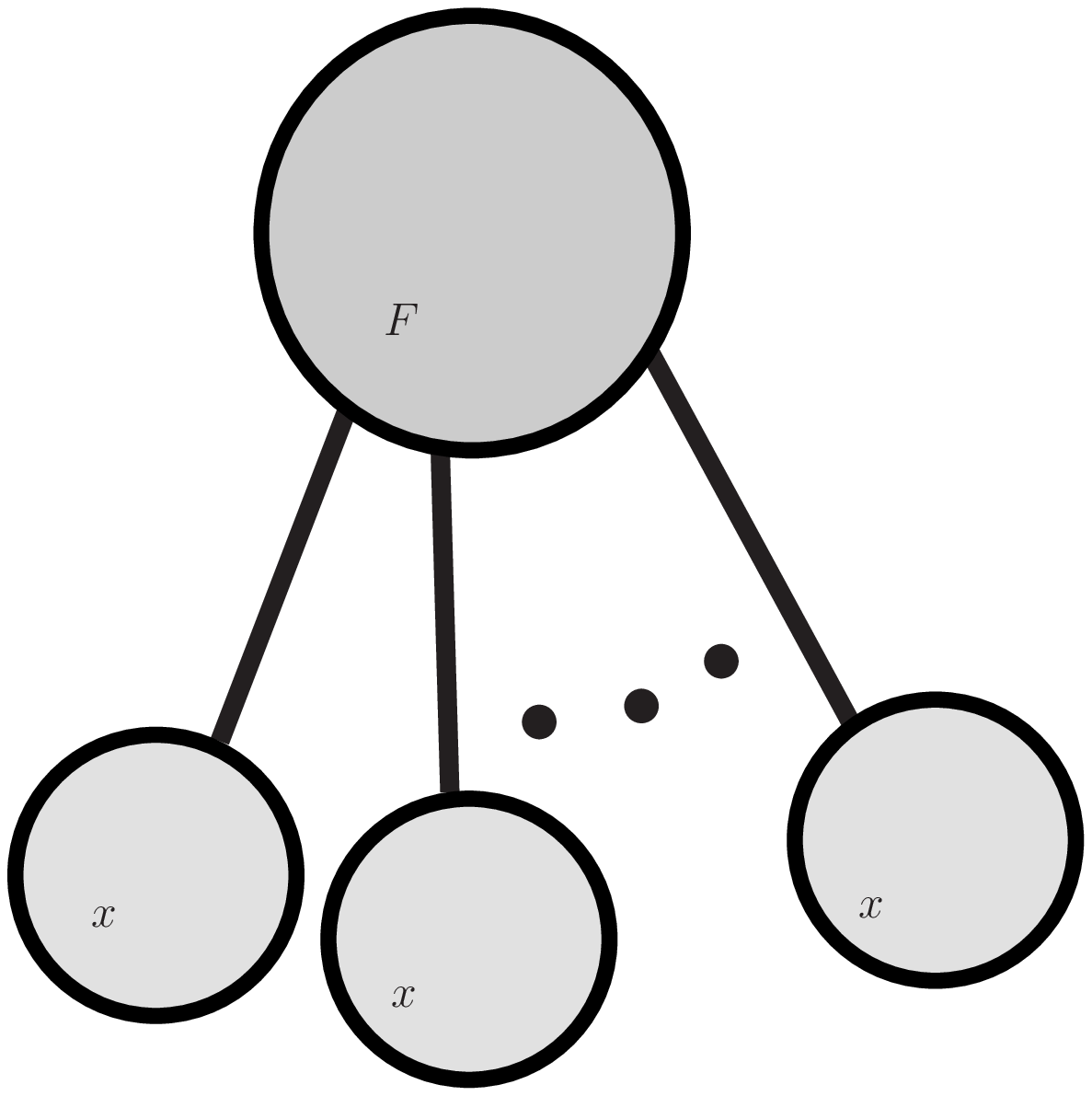}}
}_{r}
\]
using the main recipe of FDC:
\begin{quote}
The evaluation of a diagram obtained
from basic pieces say $\parbox{0.4cm}{\psfrag{x}{$\scriptstyle{x}$}
\includegraphics[width=0.4cm]{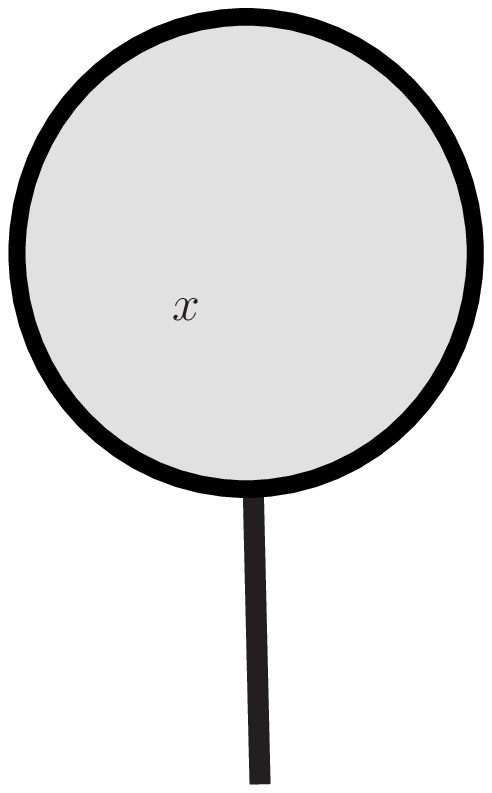}}$ and $\parbox{0.7cm}{\psfrag{F}{$\scriptstyle{F}$}
\includegraphics[width=0.7cm]{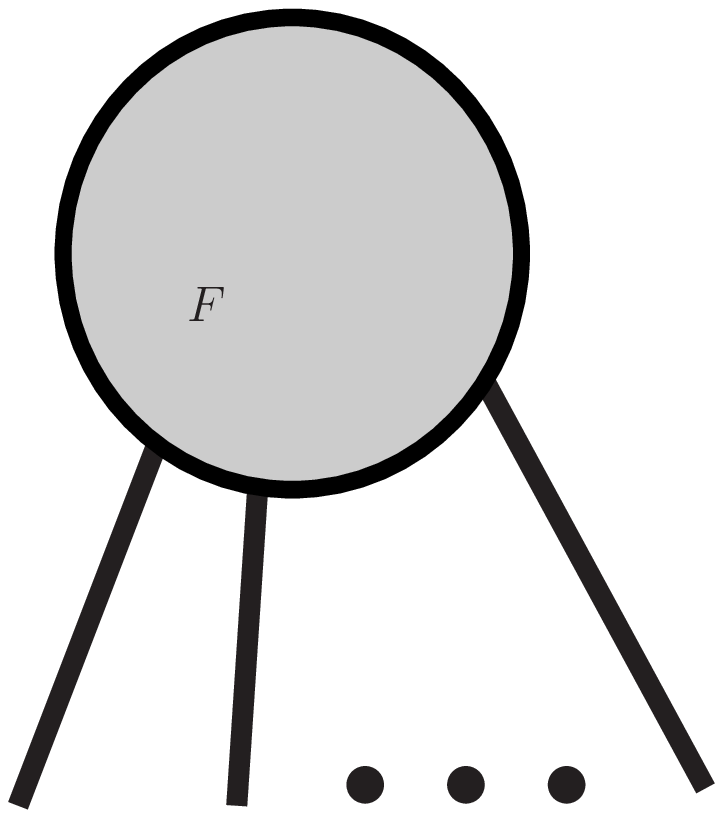}}$
by a gluing of the half-lines (or legs) is the result of
assigning an index to each pair of glued half-lines, taking the product
of corresponding tensor elements and summing over indices.
\end{quote}
We need more pieces to continue playing this game.
The Kronecker delta is
\[
\parbox{1.5cm}{\psfrag{j}{$\scriptstyle{j}$}\psfrag{i}{$\scriptstyle{i}$}
\includegraphics[width=1.5cm]{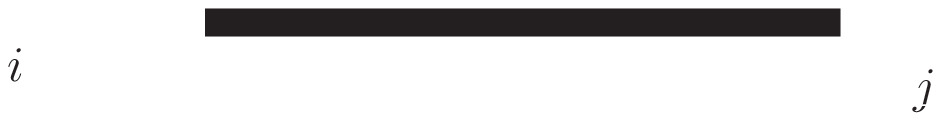}}
=\de_{ij}\ \ .
\]
The epsilon or Levi-Civita tensor is
\[
\parbox{1.5cm}{\psfrag{j}{$\scriptstyle{j}$}\psfrag{i}{$\scriptstyle{i}$}
\includegraphics[width=1.5cm]{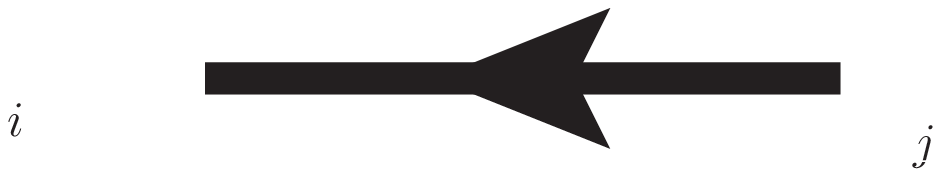}}
=\ep_{ij}
\]
where $\ep=(\ep_{ij})_{1\le i,j\le 2}$ is the antisymetric matrix
$\left(
\begin{array}{cc}
0 & 1\\
-1 & 0
\end{array}
\right)$.
The symmetrizer is
\[
\parbox{1.2cm}{\psfrag{a}{$\scriptstyle{j_1}$}\psfrag{b}{$\scriptstyle{j_2}$}
\psfrag{c}{$\scriptstyle{j_a}$}\psfrag{d}{$\scriptstyle{i_1}$}
\psfrag{e}{$\scriptstyle{i_2}$}\psfrag{f}{$\scriptstyle{i_a}$}
\includegraphics[width=1.2cm]{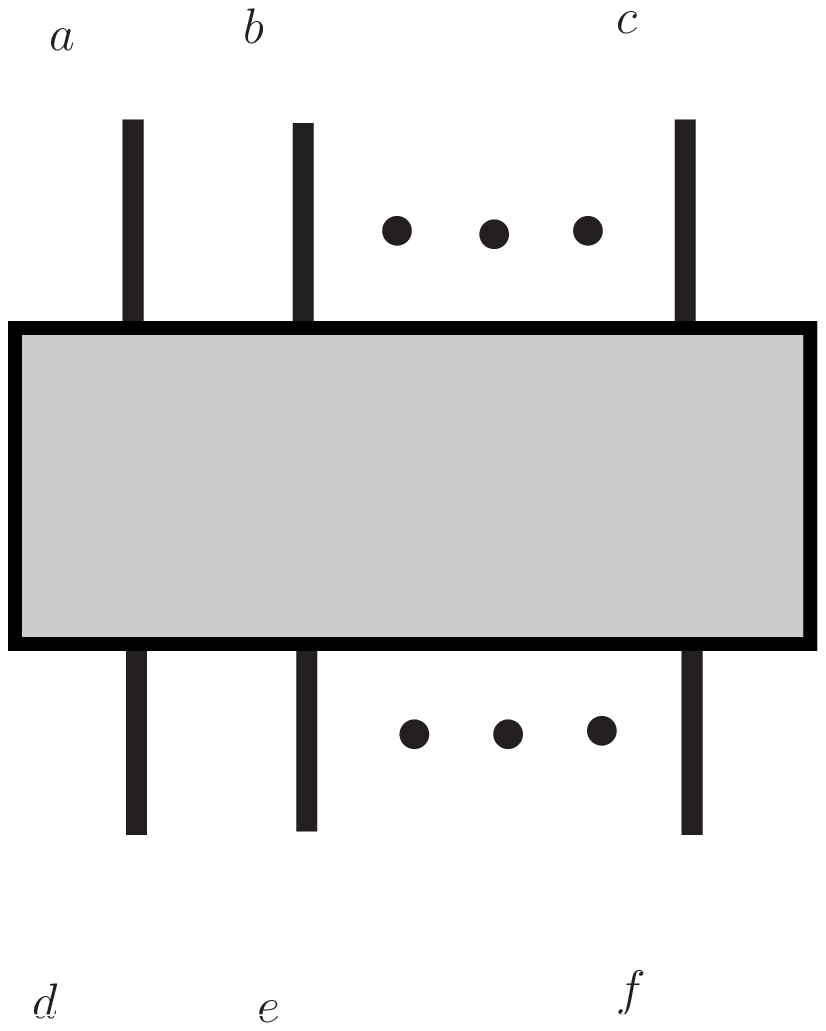}}
=\frac{1}{a!}\sum\limits_{\si\in\gS_a}
\de_{i_1 j_{\si(1)}}\ldots \de_{i_a j_{\si(a)}}\ \ .
\]
For example
\[
\parbox{1.2cm}{\psfrag{a}{$\scriptstyle{l}$}\psfrag{b}{$\scriptstyle{m}$}
\psfrag{c}{$\scriptstyle{n}$}\psfrag{d}{$\scriptstyle{i}$}
\psfrag{e}{$\scriptstyle{j}$}\psfrag{f}{$\scriptstyle{k}$}
\includegraphics[width=1.2cm]{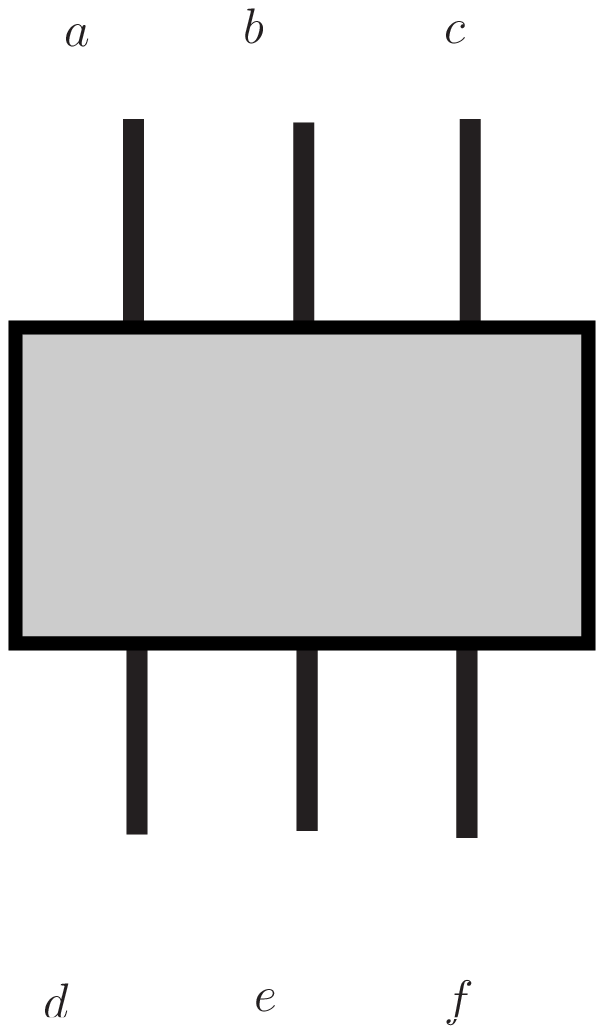}}
=\frac{1}{6}
\left[
\parbox{1cm}{\psfrag{a}{$\scriptstyle{l}$}\psfrag{b}{$\scriptstyle{m}$}
\psfrag{c}{$\scriptstyle{n}$}\psfrag{d}{$\scriptstyle{i}$}
\psfrag{e}{$\scriptstyle{j}$}\psfrag{f}{$\scriptstyle{k}$}
\includegraphics[width=1cm]{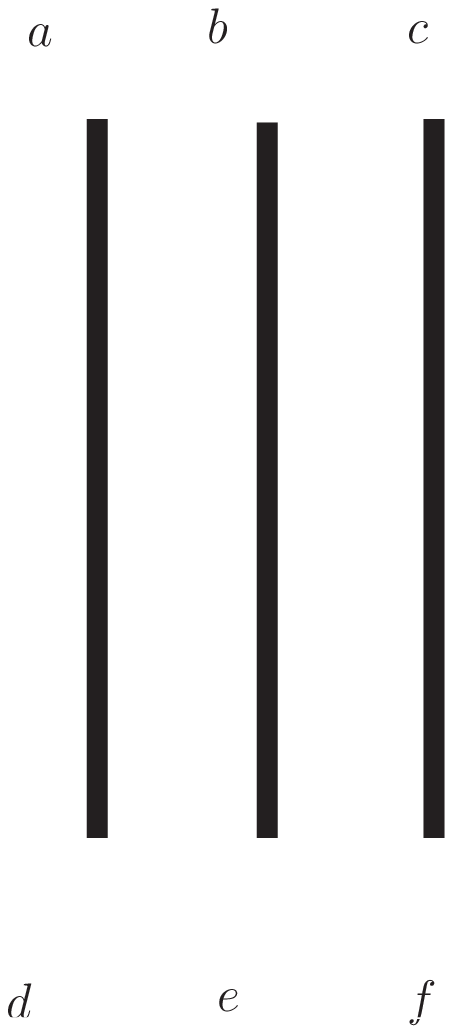}}
+
\parbox{1cm}{\psfrag{a}{$\scriptstyle{l}$}\psfrag{b}{$\scriptstyle{m}$}
\psfrag{c}{$\scriptstyle{n}$}\psfrag{d}{$\scriptstyle{i}$}
\psfrag{e}{$\scriptstyle{j}$}\psfrag{f}{$\scriptstyle{k}$}
\includegraphics[width=1cm]{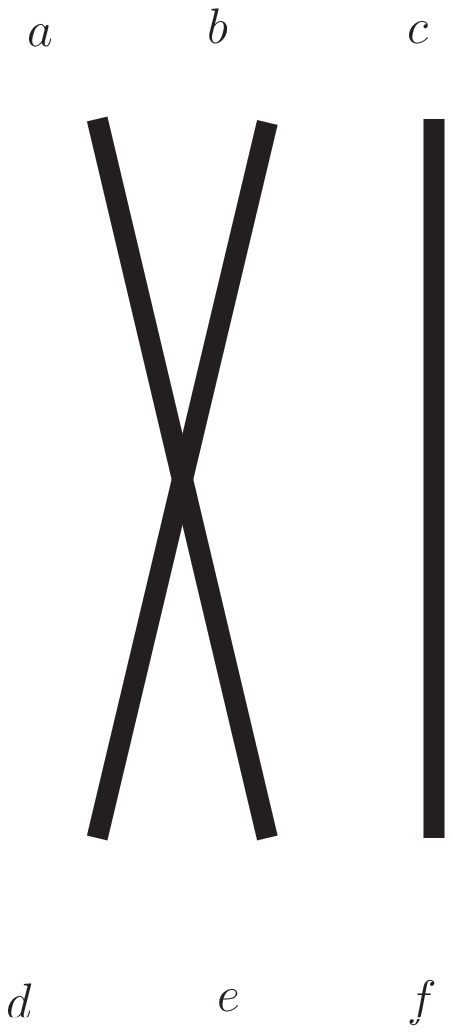}}
+
\parbox{1cm}{\psfrag{a}{$\scriptstyle{l}$}\psfrag{b}{$\scriptstyle{m}$}
\psfrag{c}{$\scriptstyle{n}$}\psfrag{d}{$\scriptstyle{i}$}
\psfrag{e}{$\scriptstyle{j}$}\psfrag{f}{$\scriptstyle{k}$}
\includegraphics[width=1cm]{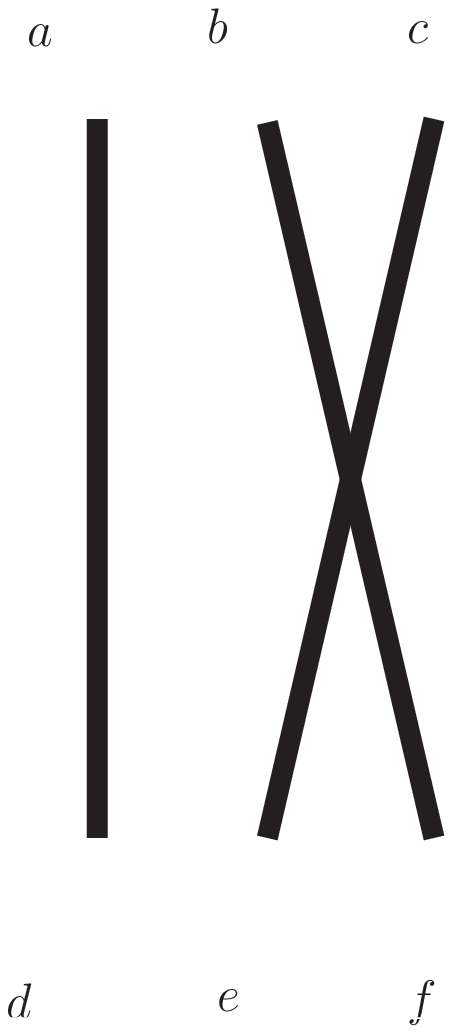}}
+
\parbox{1cm}{\psfrag{a}{$\scriptstyle{l}$}\psfrag{b}{$\scriptstyle{m}$}
\psfrag{c}{$\scriptstyle{n}$}\psfrag{d}{$\scriptstyle{i}$}
\psfrag{e}{$\scriptstyle{j}$}\psfrag{f}{$\scriptstyle{k}$}
\includegraphics[width=1cm]{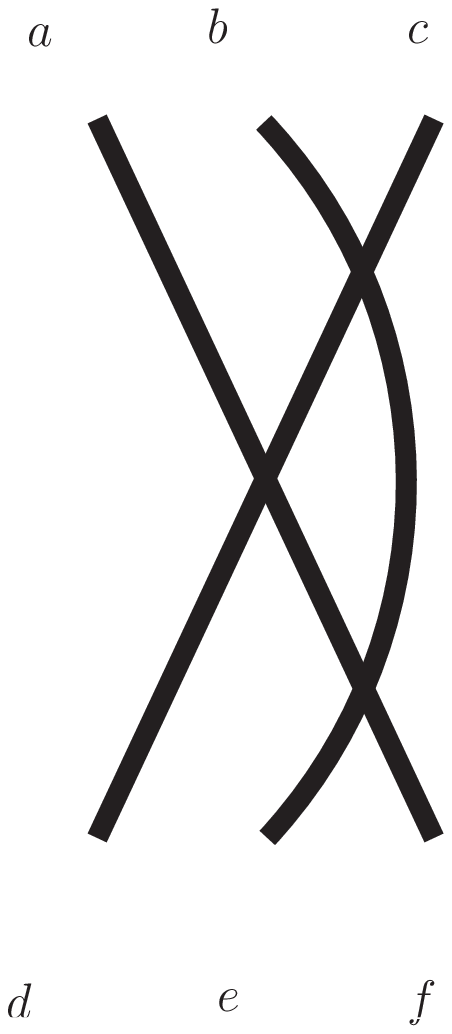}}
+
\parbox{1cm}{\psfrag{a}{$\scriptstyle{l}$}\psfrag{b}{$\scriptstyle{m}$}
\psfrag{c}{$\scriptstyle{n}$}\psfrag{d}{$\scriptstyle{i}$}
\psfrag{e}{$\scriptstyle{j}$}\psfrag{f}{$\scriptstyle{k}$}
\includegraphics[width=1cm]{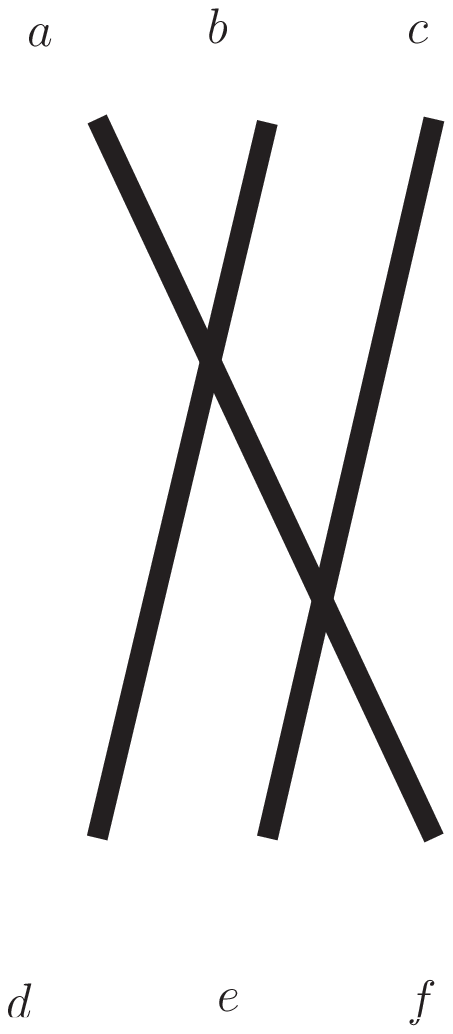}}
+
\parbox{1cm}{\psfrag{a}{$\scriptstyle{l}$}\psfrag{b}{$\scriptstyle{m}$}
\psfrag{c}{$\scriptstyle{n}$}\psfrag{d}{$\scriptstyle{i}$}
\psfrag{e}{$\scriptstyle{j}$}\psfrag{f}{$\scriptstyle{k}$}
\includegraphics[width=1cm]{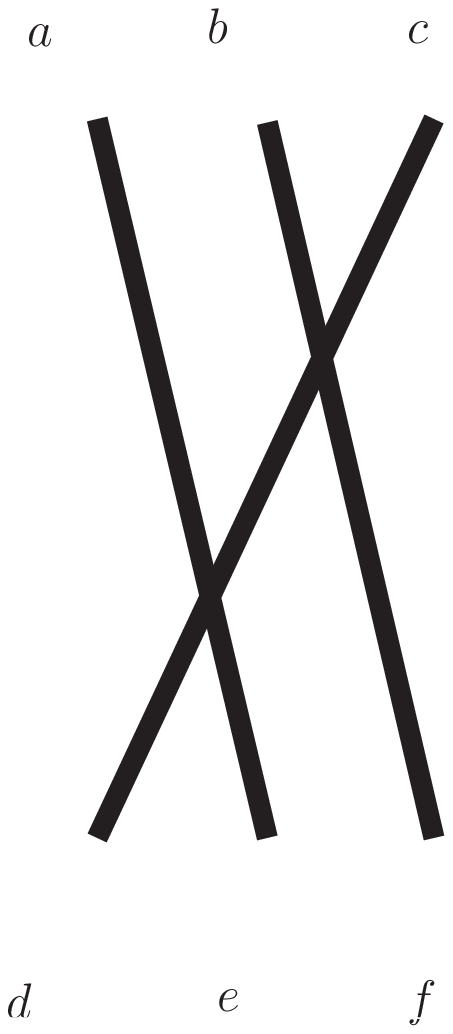}}
\right]\ \ .
\]
A $2\times 2$ matrix $g=(g_{ij})_{1\le i,j\le 2}$ is represented
by
\[
\parbox{1.6cm}{\psfrag{i}{$\scriptstyle{i}$}\psfrag{j}{$\scriptstyle{j}$}
\psfrag{g}{$\scriptstyle{g}$}
\includegraphics[width=1.6cm]{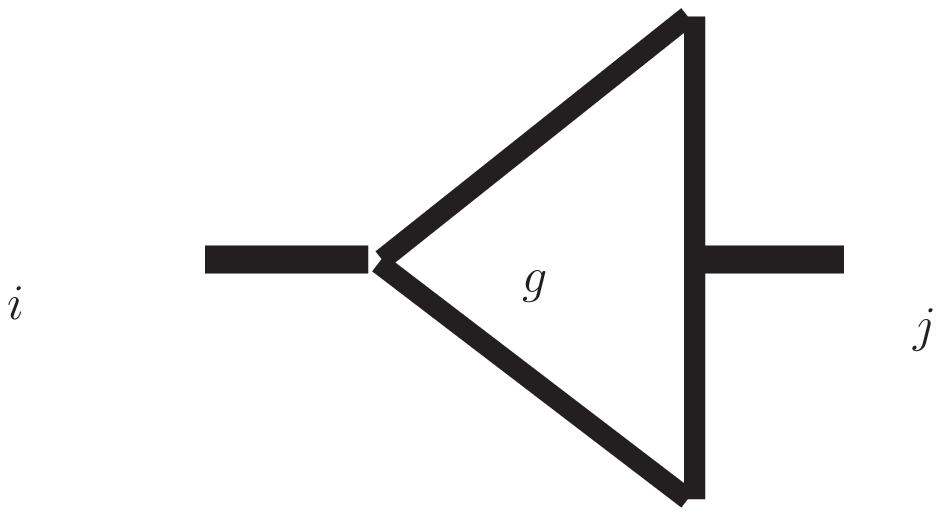}}
=g_{ij}\ \ .
\]
One can also use this graphical representation
for differential operators
\[
\parbox{0.6cm}{\psfrag{x}{$\scriptstyle{\partial x}$}\psfrag{i}{$\scriptstyle{i}$}
\includegraphics[width=0.6cm]{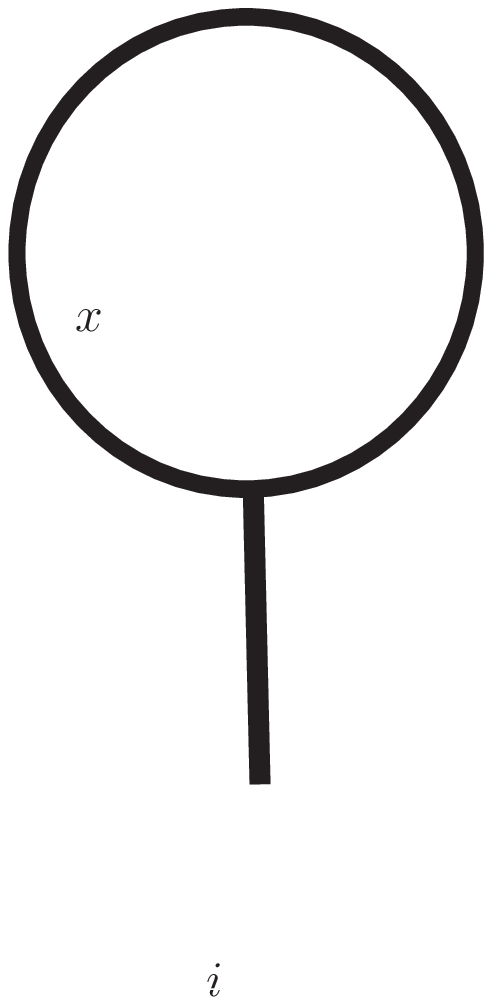}}
=\frac{\partial}{\partial x_i}
\]
the only difference is that one needs to indicate a direction of reading, since
the $\frac{\partial}{\partial x}$'s do not commute with the $x$'s.
For instance one has the identity
\begin{equation}
\frac{1}{r!}
\underbrace{
\parbox{1.8cm}{\psfrag{F}{$\scriptstyle{F}$}\psfrag{x}{$\scriptstyle{\partial a}$}
\includegraphics[width=1.8cm]{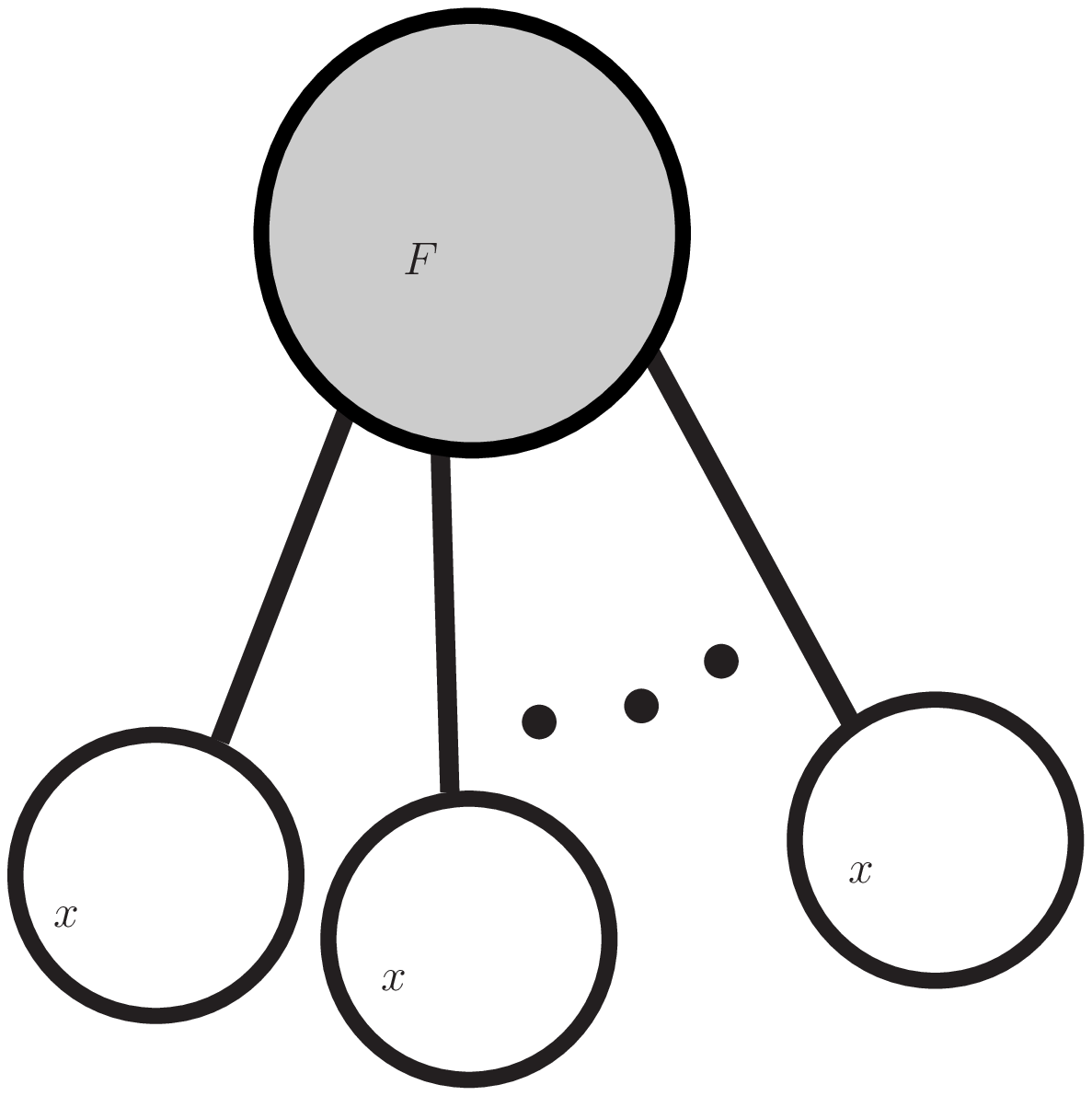}}}_{r}
\qquad
\underbrace{
\parbox{0.4cm}{\psfrag{a}{$\scriptstyle{a}$}\psfrag{x}{$\scriptstyle{x}$}
\includegraphics[width=0.4cm]{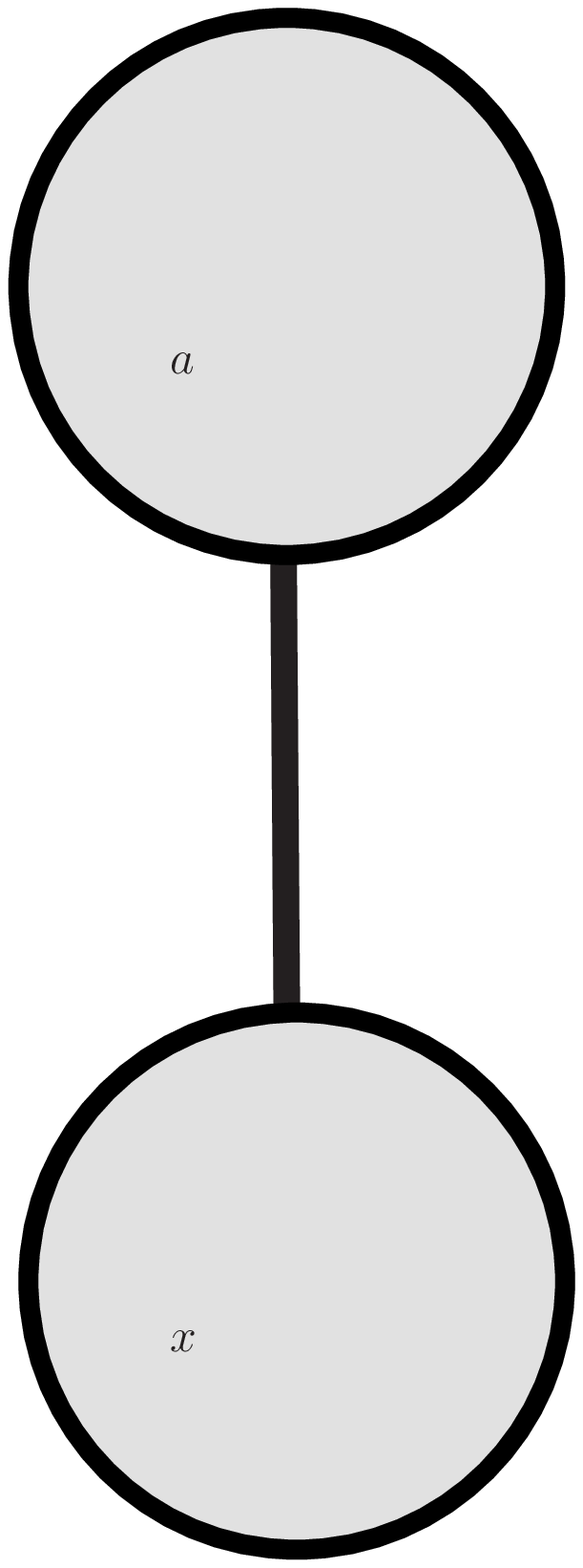}}\ 
\parbox{0.4cm}{\psfrag{a}{$\scriptstyle{a}$}\psfrag{x}{$\scriptstyle{x}$}
\includegraphics[width=0.4cm]{Fig54.eps}}
\cdots
\parbox{0.4cm}{\psfrag{a}{$\scriptstyle{a}$}\psfrag{x}{$\scriptstyle{x}$}
\includegraphics[width=0.4cm]{Fig54.eps}}
}_{r}
=
\underbrace{
\parbox{1.8cm}{\psfrag{F}{$\scriptstyle{F}$}\psfrag{x}{$\scriptstyle{x}$}
\includegraphics[width=1.8cm]{Fig38.eps}}}_{r}
\label{graphdiff}
\end{equation}
where the formula on the left-hand side is read from left to right.
The familiar Liebnitz's rule can be interpreted by saying
that each $\parbox{0.8cm}{\psfrag{x}{$\scriptstyle{\partial a}$}
\includegraphics[width=0.8cm]{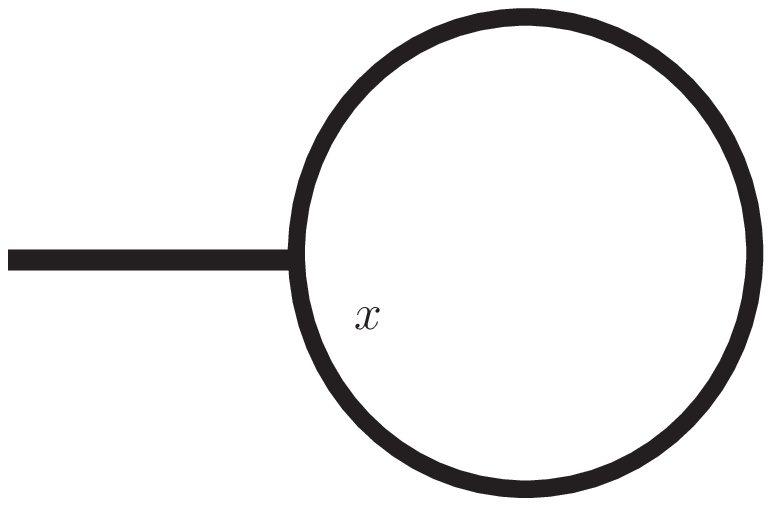}}$ selects an
$\parbox{0.8cm}{\psfrag{x}{$\scriptstyle{a}$}
\includegraphics[width=0.8cm]{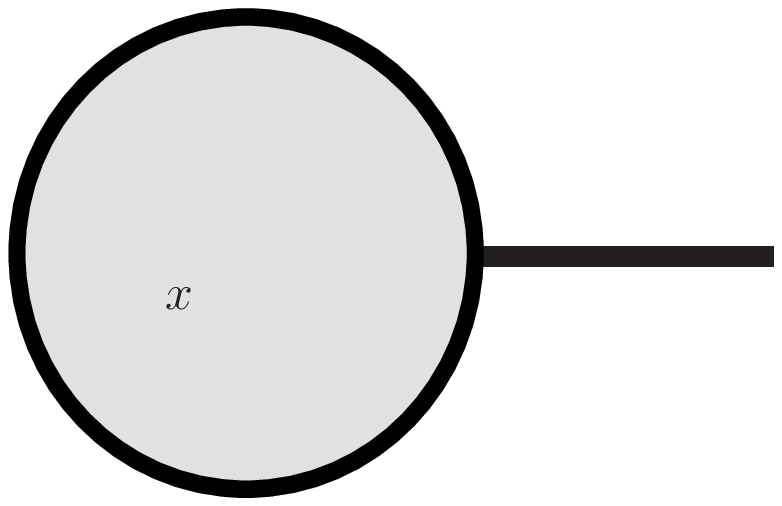}}$ to annihilate
and glues whatever it is attached to together with
what the selected
$\parbox{0.8cm}{\psfrag{x}{$\scriptstyle{a}$}
\includegraphics[width=0.8cm]{Fig56.eps}}$
is attached to. One also needs to sum over
all possible ways of making these selections.
See~\cite{BarNatanGRT} for a nice
presentation along this line of FDC and Feynman diagram expansions
intuitively seen as the result of a `chemical reaction' between
$\parbox{0.8cm}{\psfrag{x}{$\scriptstyle{\partial a}$}
\includegraphics[width=0.8cm]{Fig55.eps}}$
graphs and 
$\parbox{0.8cm}{\psfrag{x}{$\scriptstyle{a}$}
\includegraphics[width=0.8cm]{Fig56.eps}}$
graphs.
The previous identity
says nothing more than the trivial fact
\[
\frac{1}{r!}
\sum\limits_{i_1,\ldots,i_r=1}^{2}
F_{i_1\ldots i_r}\frac{\partial}{\partial a_{i_1}}
\ldots\frac{\partial}{\partial a_{i_r}}
\sum\limits_{j_1,\ldots,j_r=1}^{2}
a_{j_1} x_{j_1}\ldots a_{j_r} x_{j_r}=F(\ux)\ .
\]
However, it is the basis of a correct interpretation of the very powerful
classical symbolic method, as we will shortly see.

\begin{Remark}
In our opinion, the most pedagogically useful way of introducing
the so-called `Wick Theorem' which generates Feynman diagram expansions
from integration with respect to a Gaussian measure with covariance $C$,
is by following the same philosophy as above. Namely, integration
amounts in this case to applying a differential operator
\[
\exp\left(\frac{1}{2}
\parbox{1.8cm}{\psfrag{x}{$\scriptstyle{\partial x}$}\psfrag{C}{$\scriptstyle{C}$}
\includegraphics[width=1.8cm]{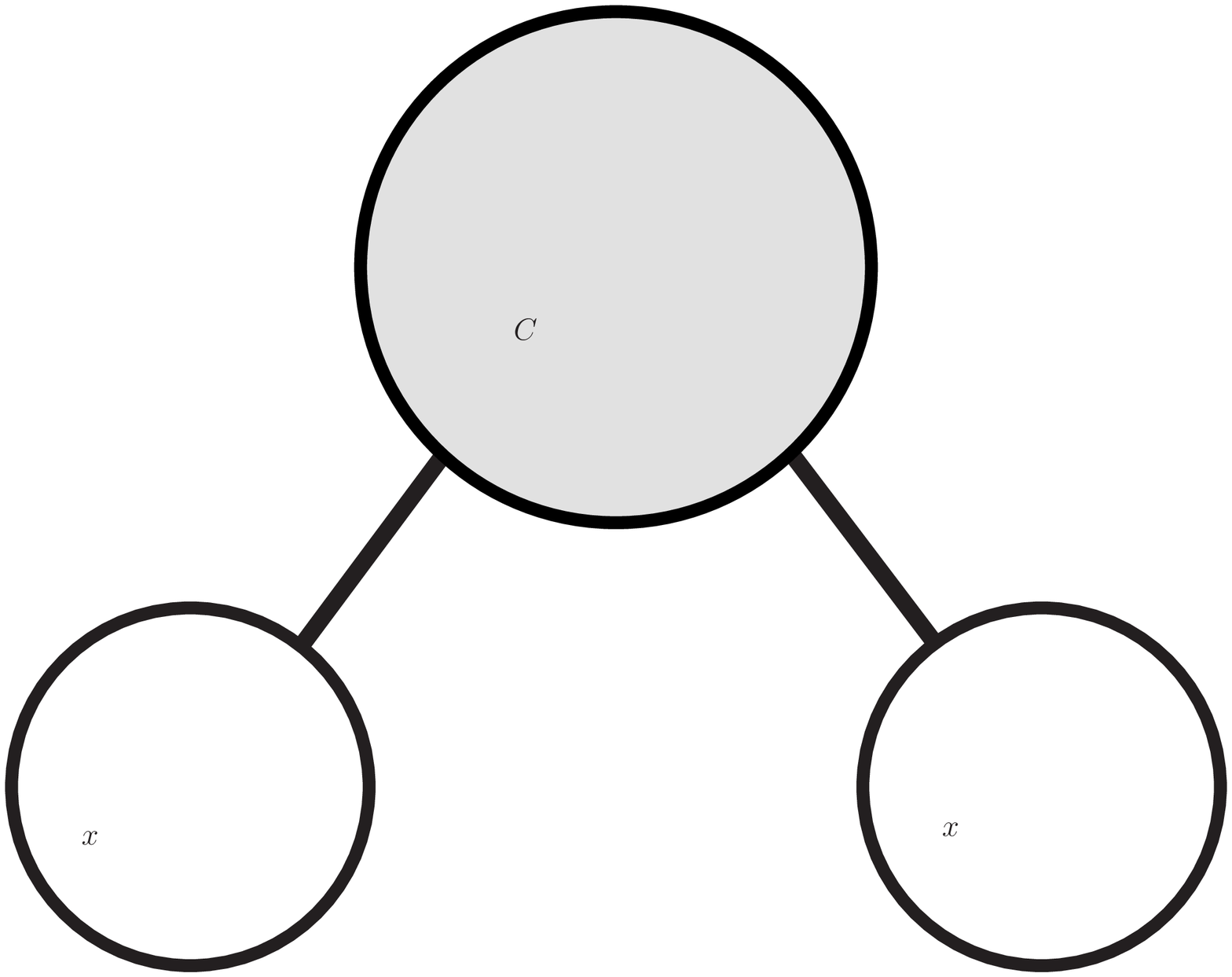}}
\right)
\]
albeit of infinite order.
The immediate follow up is that one does not need a Gaussian measure
in order to generate such a diagram expansion, and the differential operator
in the exponential needs not be of second order. One does not
need to generate graphs but could also generate higher-dimensional
objects as in random matrix theory.
As is well-known, traditional Feynman diagram expansions provide the explicit
form of the stationary phase asymptotic series in the presence of
nondegenerate critical points. Trying to generalize this
to degenerate critical points~\cite{ArnoldGV} is largely open
and is related to the problem addressed in~\cite{EtingofKP}.
We believe this investigation would most likely benefit
from the classical invariant theoretic point of view.
Also note that generating graphs by applying differential operators
in this way, was already known to Arthur Cayley, one of the first
developers of CIT~\cite{Cayley2}. 
Also, Wick's Theorem which expresses the moments of a Gaussian
measure was discovered by Leon Isserlis~\cite{Isserlis1} who in addition
to his expertise in statistics was well versed in the techniques
of CIT~\cite{Isserlis2,Irwin}.
\end{Remark}

The fundamental propertry of $SL_2(\C)$ invariance of diagrams built
with the previous pieces is the identity
\begin{equation}
\parbox{1.8cm}{\psfrag{i}{$\scriptstyle{i}$}\psfrag{j}{$\scriptstyle{j}$}
\psfrag{g}{$\scriptstyle{g}$}
\includegraphics[width=1.8cm]{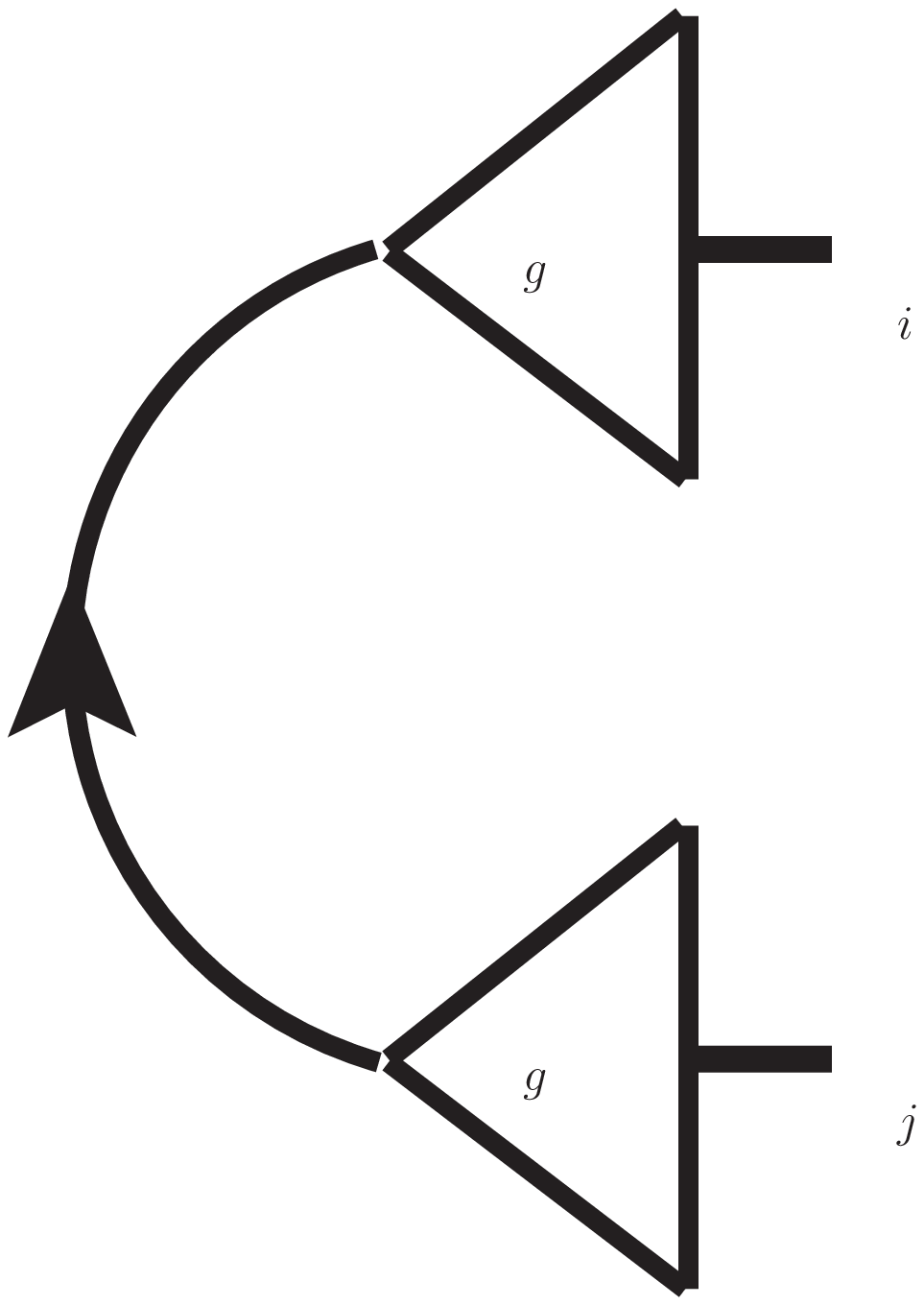}}
=({\rm det}\ g)\ 
\parbox{1cm}{\psfrag{i}{$\scriptstyle{i}$}\psfrag{j}{$\scriptstyle{j}$}
\includegraphics[width=1cm]{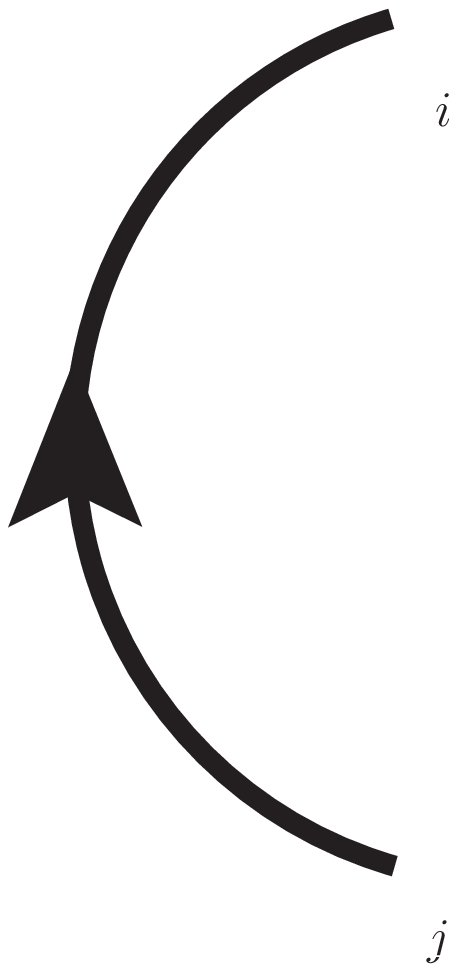}}\ \ .
\label{epsilonident}
\end{equation}
Indeed, if we decide that $GL_2(\C)$ acts on vectors by
\begin{equation}
\parbox{0.4cm}{\psfrag{x}{$\scriptstyle{x}$}
\includegraphics[width=0.4cm]{Fig39.eps}}
\longmapsto
\parbox{0.6cm}{\psfrag{x}{$\scriptstyle{gx}$}
\includegraphics[width=0.6cm]{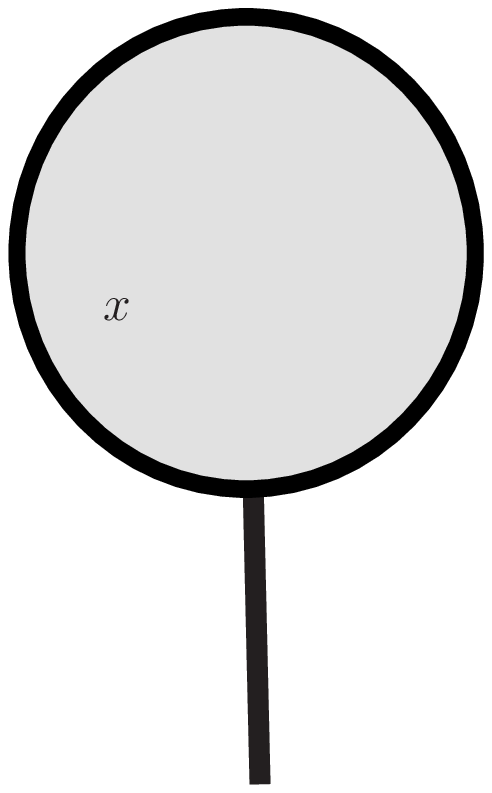}}
\ =\ \parbox{0.5cm}{\psfrag{x}{$\scriptstyle{x}$}\psfrag{g}{$\scriptstyle{g}$}
\includegraphics[width=0.5cm]{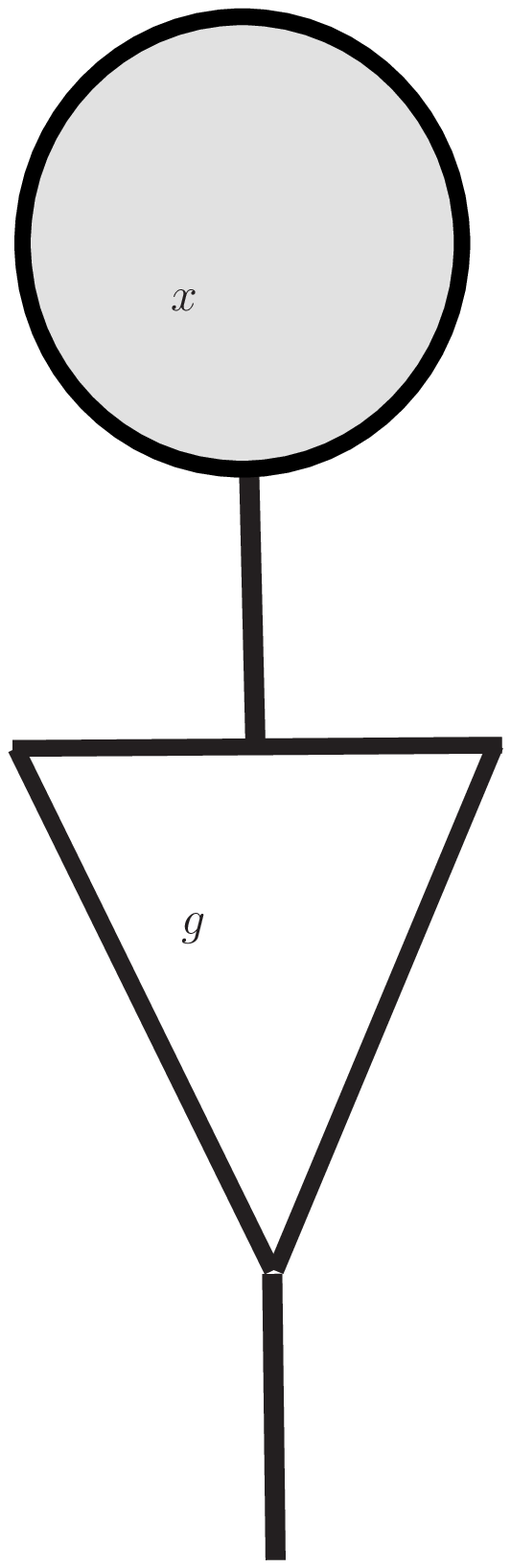}}
\label{xtransf}
\end{equation}
then it will also act on binary forms by
\begin{equation}
\parbox{0.8cm}{\psfrag{F}{$\scriptstyle{F}$}
\includegraphics[width=0.8cm]{Fig40.eps}}
\longmapsto
\parbox{1cm}{\psfrag{F}{$\scriptstyle{gF}$}
\includegraphics[width=1cm]{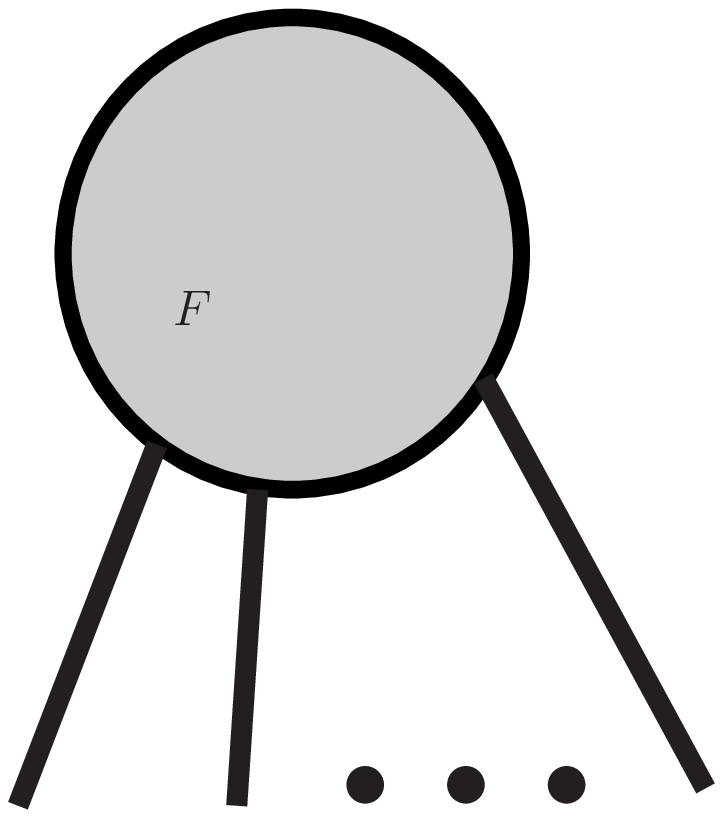}}
\ =\ 
\parbox{4cm}{\psfrag{F}{$F$}\psfrag{g}{$\scriptscriptstyle{g^{\scriptscriptstyle{-1}}}$}
\includegraphics[width=4cm]{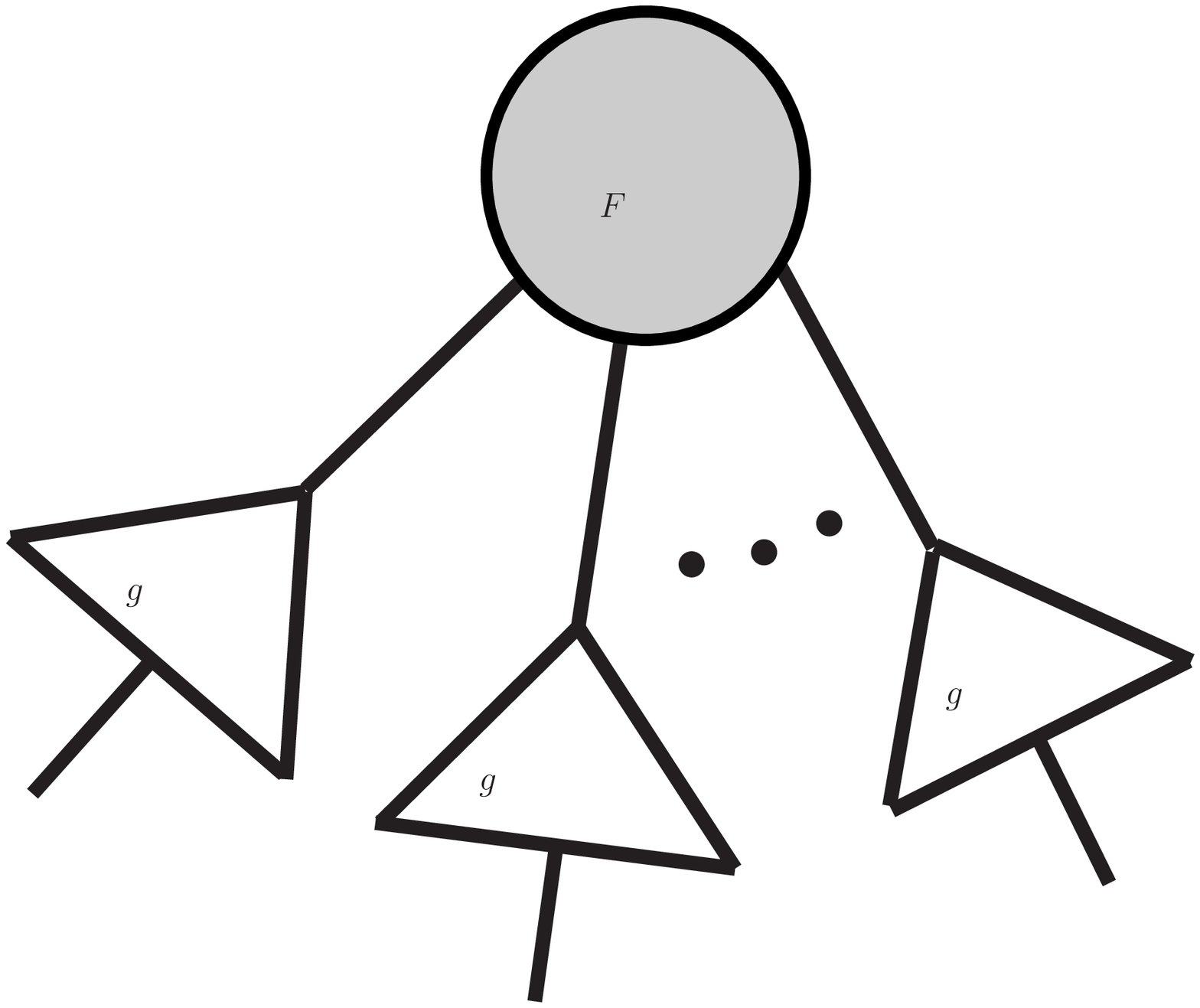}}\ \ .
\label{Ftransf}
\end{equation}
A (homogeneous) classical invariant of degree $d$ of a binary form $F$
of order $r$ is a polynomial
$I(F)=I(f_0,f_1,\ldots,f_r)$ which is homogeneous of degree $d$ in the
coefficients of $F$ and satisfies
\[
I(gF)=I(F)\times ({\rm det}\ g)^{-w}
\]
for any $g\in GL_2(\C)$.
The power $w$ is called the weight of the invariant and it is given by
$w=\frac{dr}{2}$.
More generally, a classical covariant is a polynomial
\[
C(F,\ux)=C(f_0,f_1,\ldots,f_r; x_1,x_2)
\]
separately homogeneous in $F$ and $\ux$ which satisfies
\[
C(gF,g\ux)=C(F,\ux)\times ({\rm det}\ g)^{-w}
\]
where the weight is $w=\frac{dr-n}{2}$.
Here $d$ is the degree of the covariant, i.e., the degree in the coefficients
of $F$. Whereas $n$ is the order of the covariant, i.e., the degree in $\ux$.
Invariants are covariants of order $0$.

The so-called First Fundamental Theorem of CIT for $SL_2$ says:

\begin{enumerate}
\item
Consider Feynman diagrams obtained by assembling pieces
taken among 
$\parbox{0.4cm}{\psfrag{x}{$\scriptstyle{x}$}
\includegraphics[width=0.4cm]{Fig39.eps}}$ , $\parbox{0.7cm}{\psfrag{F}{$\scriptstyle{F}$}
\includegraphics[width=0.7cm]{Fig40.eps}}$
and
$\parbox{0.7cm}{\includegraphics[width=0.7cm]{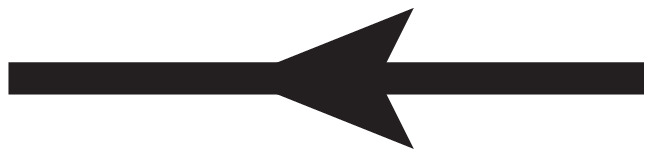}}$, with the condition
that direct connections
between two $F$'s or two $x$'s are forbidden, and
that $\ep$ arrows must join two $F$ vertices or `blobs'.
The evaluation of such a diagram is a covariant.
The degree is the number of $F$ blobs, the order
is the number of $x$ blobs, and the weight
is the number of $\ep$ arrows.
\item
Any covariant is a linear combination of such diagrams (or evaluations
thereof).
\end{enumerate}

Part 1) as well as the conceptual framework in which expressing
this statement was made possible are due to Cayley~\cite{Cayley1}.
This property is a consequence of (\ref{epsilonident}). Note that
the meaning of the word ``hyperdeterminant'' used by Cayley in~\cite{Cayley1}
is somewhat different from that in some of his earlier work
and the modern understanding~\cite{GelfandKZ}.
It refers to {\em any} polynomial one can obtain by application of the rules
in Part 1).
The easy proof of Part 1) is best seen on an example.
If $F$ is a binary quintic, its canonisant is the polynomial
\begin{equation}
C(F,\ux)=
\parbox{2.7cm}{\psfrag{F}{$F$}\psfrag{x}{$x$}
\includegraphics[width=2.7cm]{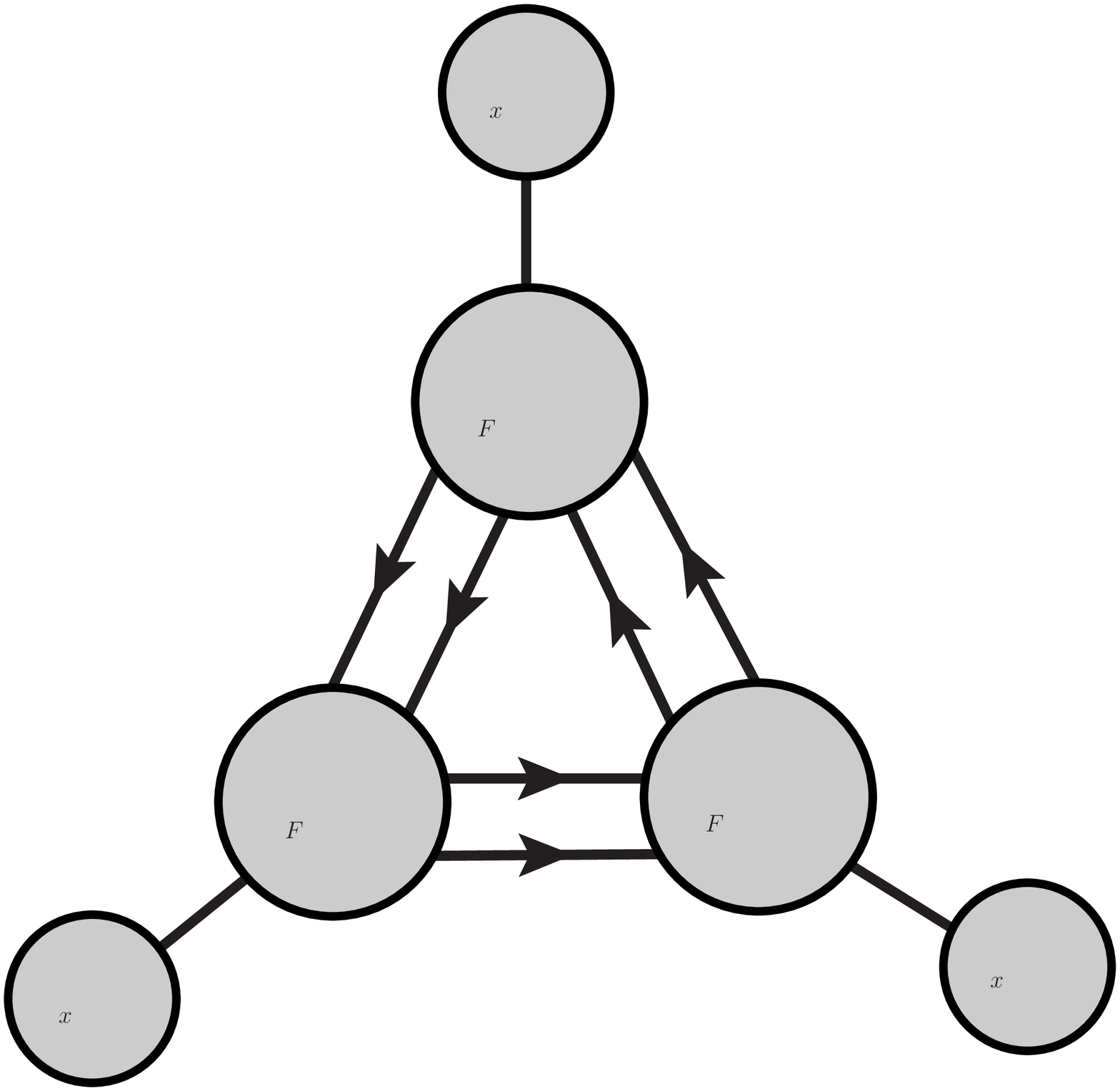}}\ \ .
\label{canonisant}
\end{equation}
Now, applying the definitions (\ref{xtransf}) and (\ref{Ftransf})
\[
C(gF,g\ux)=
\parbox{8cm}{\psfrag{F}{$F$}\psfrag{x}{$x$}
\psfrag{1}{$\scriptscriptstyle{+}$}\psfrag{2}{$\scriptscriptstyle{-}$}
\includegraphics[width=3.7cm]{Fig66.eps}}
\]
\[
\begin{array}{c}
\ \\
\ \\
\ 
\end{array}
\]
where ``$+$'' refers to the matrix $g$ and ``$-$'' to the matrix $g^{-1}$. Therefore
\[
C(gF,g\ux)
= [{\rm det}\ (g^{-1})]^6\times C(F,\ux)
\]
because of (\ref{epsilonident})
and $\parbox{2.2cm}{\psfrag{1}{$\scriptstyle{+}$}
\psfrag{2}{$\scriptstyle{-}$}
\includegraphics[width=2.2cm]{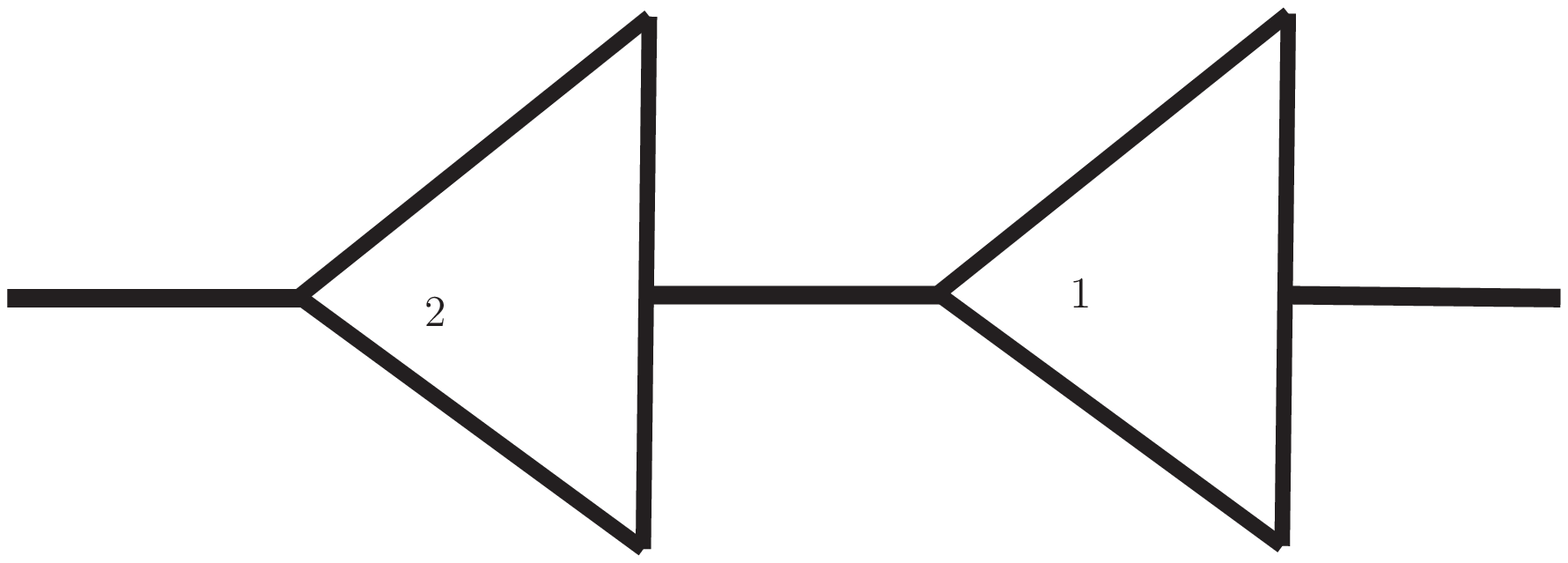}}
=\parbox{1cm}{\includegraphics[width=1cm]{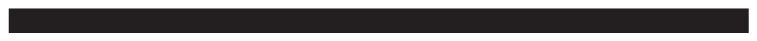}}$, i.e., $g^{-1} g={\rm Id}$.
Note in passing the use of `substitutions of diagrams into blobs'.
This is related to plethysm and was exploited extensively by the classics.

The proof of Part 2), for $SL_n$, is due to Clebsch~\cite[\S3]{Clebsch1}.
With the hindsight provided by the FDC devised by Michael Creutz for
$SU_n$ integration~\cite{Creutz1,Creutz2}, i.e.,
the misnamed Reynolds operator,
one can paraphrase {\em Clebsch's original proof}, in
the simpler case of invariants, as follows.
An invariant $I$ also has its blob:
\begin{equation}
I(F)=
\underbrace{
\parbox{2.2cm}{\psfrag{F}{$\scriptstyle{F}$}\psfrag{I}{$I$}
\includegraphics[width=2.2cm]{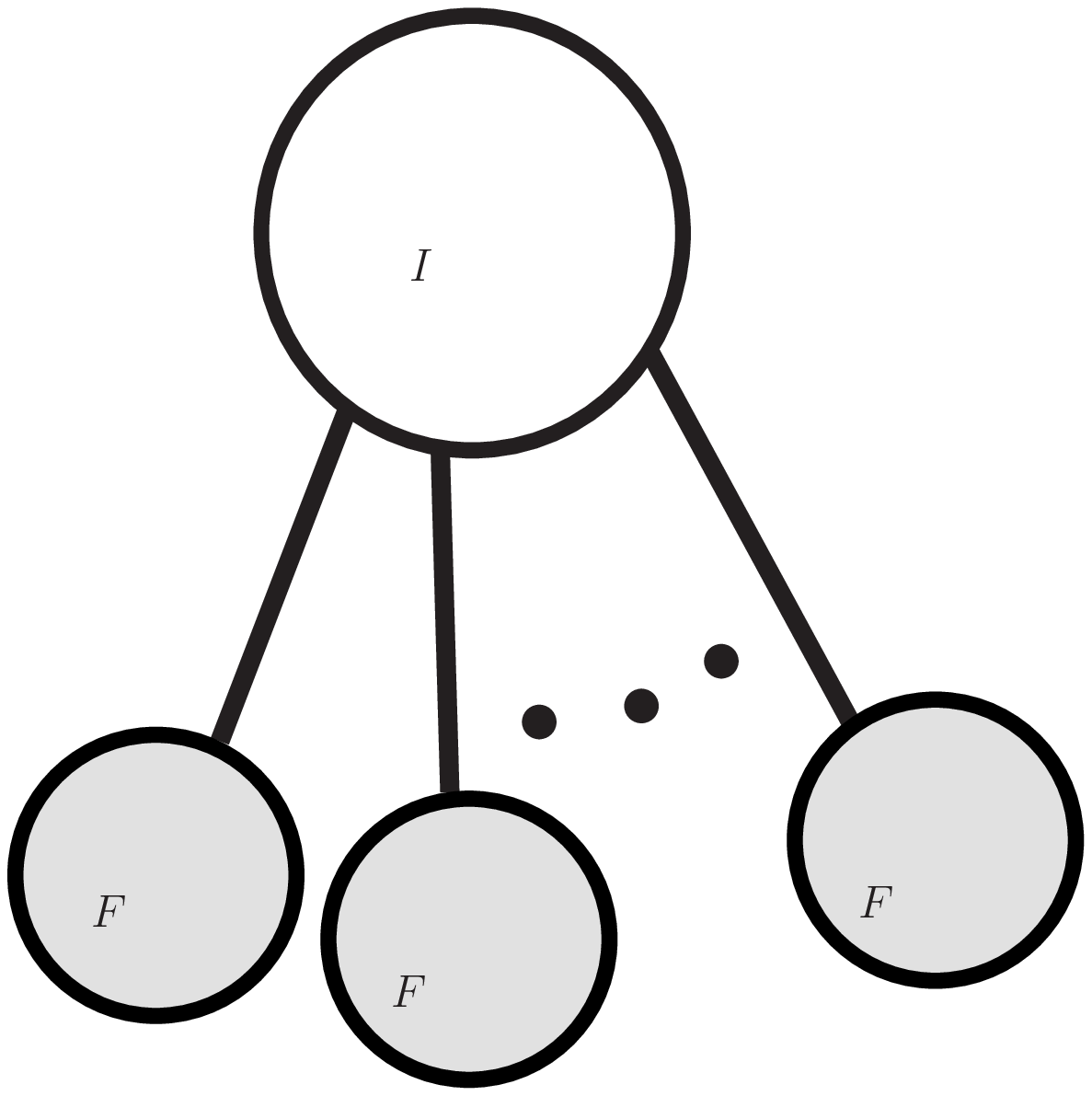}}
}_{d}
\label{Iblobdef}
\end{equation}
except that indices run over the set with $r+1$ elements
which labels a basis of ${\rm Sym}^d(V^\ast)$.
These are the equivalence classes of tuples $(i_1,\ldots,i_r)\in
\{1,2\}^r$ under permutations.
One can revert back to the description in terms of binary indices
and rewrite (\ref{Iblobdef}) `microscopically' as
\[
I(F)=
\parbox{3.5cm}{\psfrag{F}{$\scriptstyle{F}$}\psfrag{I}{$I$}
\psfrag{U}{$\overbrace{\parbox{2.5cm}{\ }}^{d}$}\psfrag{D}{$\scriptstyle{r}\underbrace{\ \ }$}
\includegraphics[width=3.5cm]{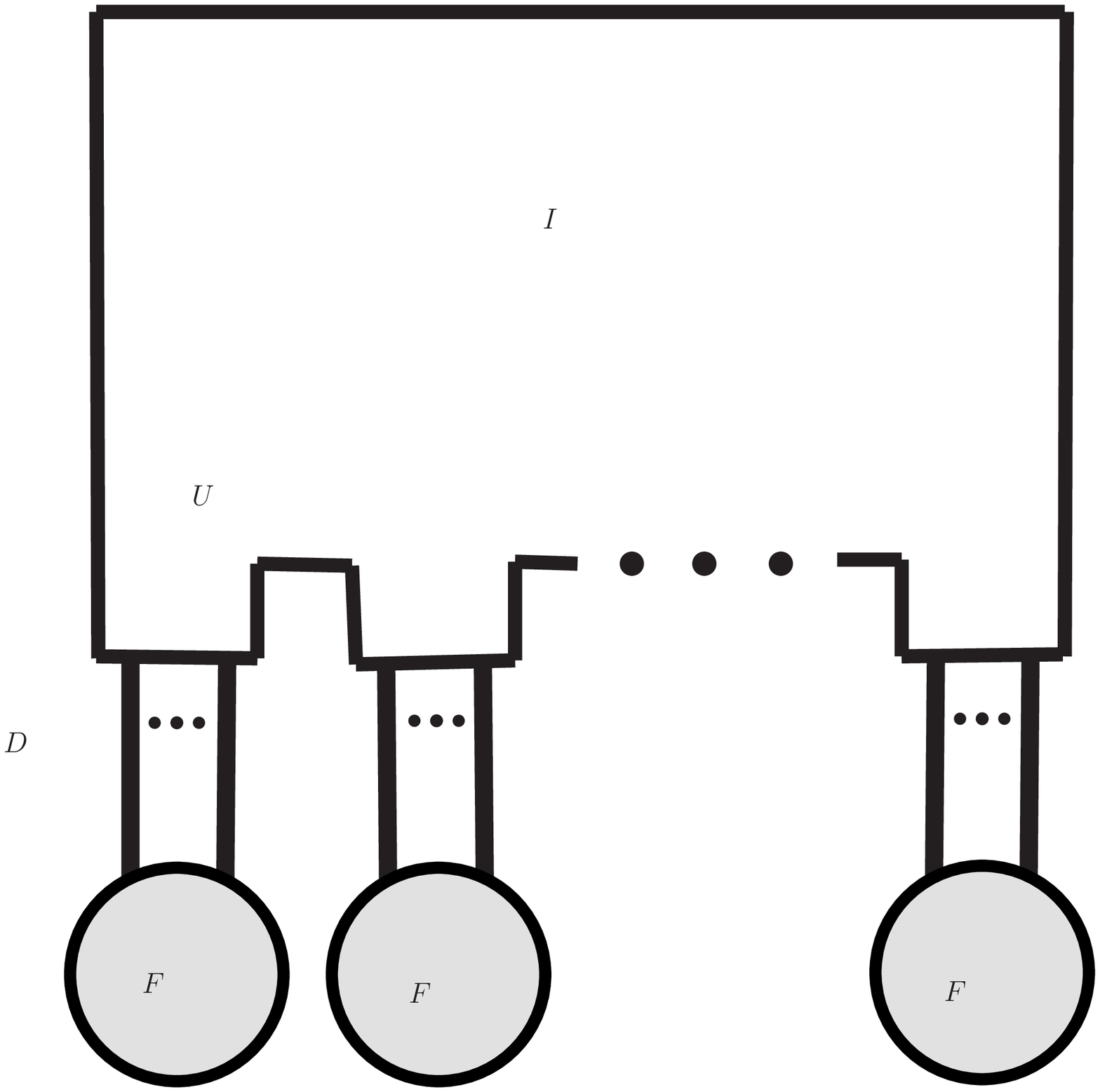}}
\]
with indices in $\{1,2\}$.
This is an example of ``categorification'' in the sense of~\cite{BaezD},
since we sum over tuples instead of equivalence classes.
The tensor $I$ has $dr$ indices, arranged in $d$ groups of $r$, taking values
in $\{1,2\}$.
The tensor, or corresponding blob, must be symmetric with respect
to permutation of the groups as well as permutations of indices within groups.
In other words, $I$ lives in ${\rm Sym}^d({\rm Sym}^r(V))$.

By $SL_2(\C)$ and therefore $SU_2(\C)$ invariance, $I(F)=I(g^{-1}F)$
for any $g\in SU_2(\C)$.
Thus, denoting the normalized Haar measure on $SU_2(\C)$ by ${\rm d}\mu$,
\[
I(F)=\int_{SU_2(\C)} {\rm d}\mu(g)\ I(g^{-1}F)
= {\rm Cst}\times
{\rm det}\left(\frac{\partial}{\partial g}\right)^{\frac{dr}{2}}\ I(g^{-1}F)
\]
\[
= {\rm Cst}\times
\underbrace{
\parbox{2cm}{\psfrag{d}{$\scriptscriptstyle{\partial g}$}
\includegraphics[width=2cm]{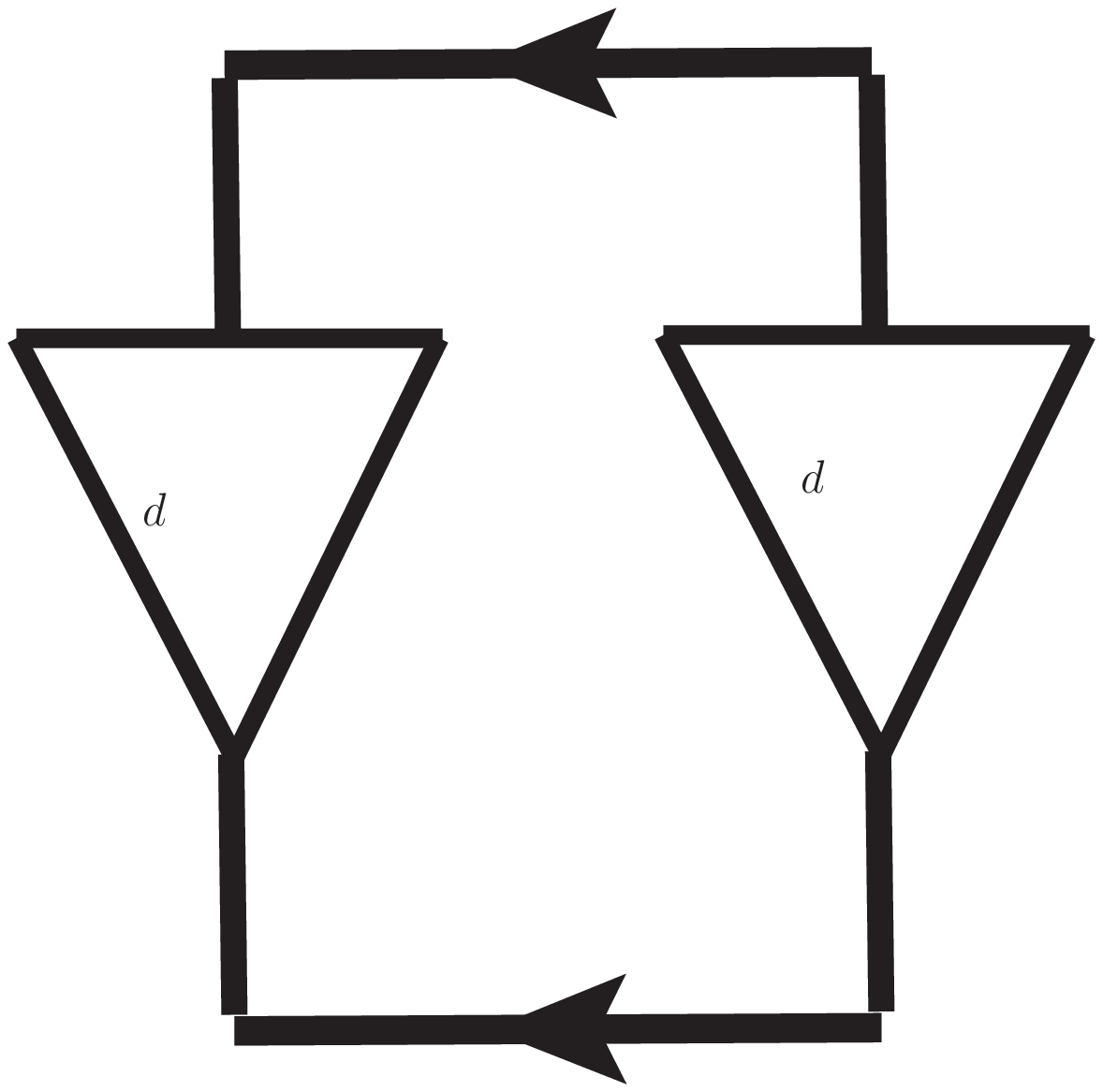}}
\cdots
\parbox{2cm}{\psfrag{d}{$\scriptscriptstyle{\partial g}$}
\includegraphics[width=2cm]{Fig71.eps}}
}_{\frac{dr}{2}\ {\rm times}}
\qquad
\parbox{3.5cm}{\psfrag{F}{$\scriptstyle{F}$}\psfrag{I}{$I$}
\psfrag{g}{$\scriptscriptstyle{g}$}
\includegraphics[width=3.5cm]{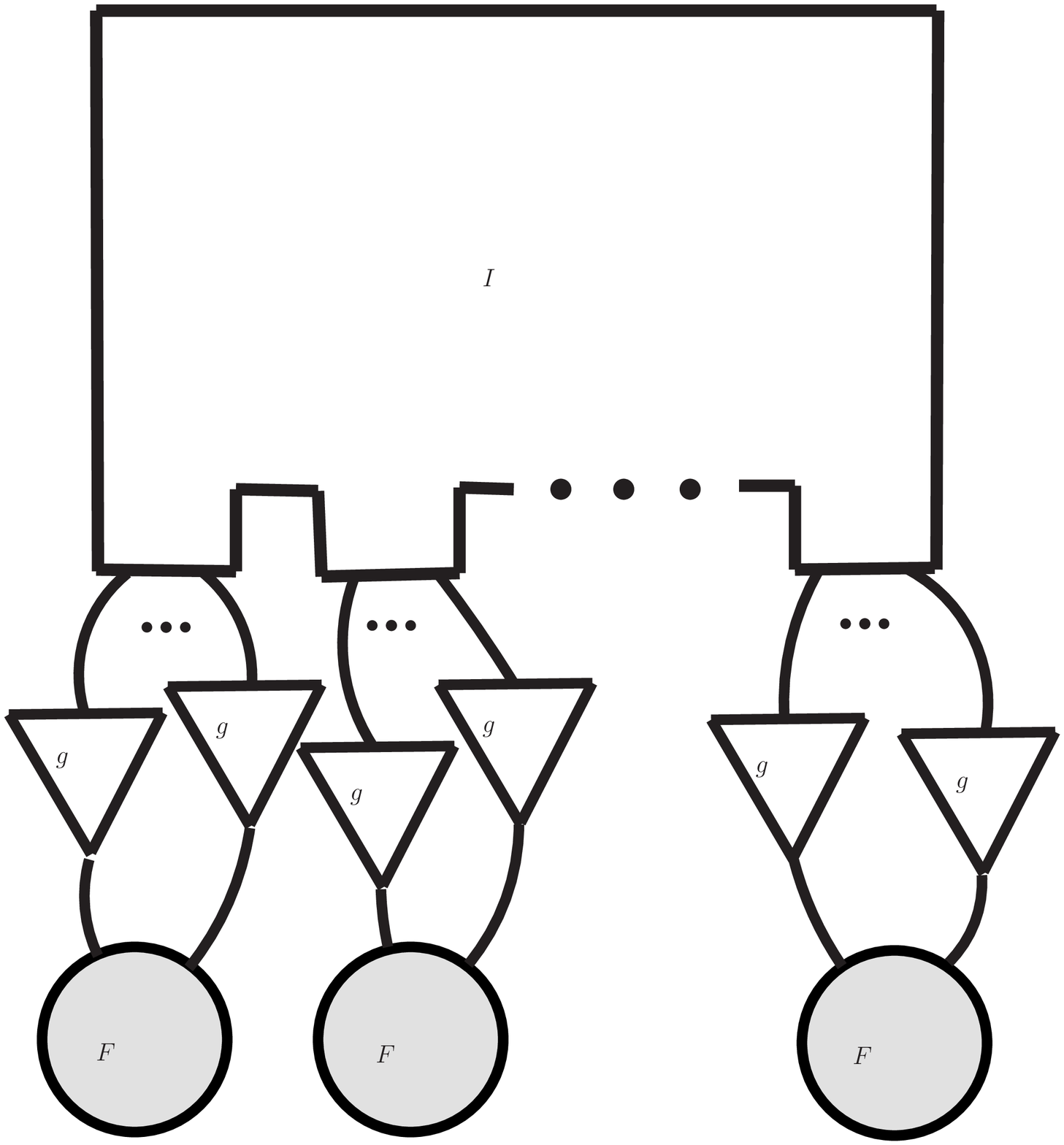}}\ \ .
\]
Now the reader who understood (\ref{graphdiff})
sees what is happening.
Each $\partial g$ square selects, in all possible ways,
a pair of $g$'s to contract to as in
\[
\parbox{2cm}{\psfrag{d}{$\scriptscriptstyle{\partial g}$}
\includegraphics[width=2cm]{Fig71.eps}}
\qquad
\parbox{0.8cm}{\psfrag{g}{$\scriptscriptstyle{g}$}
\includegraphics[width=0.8cm]{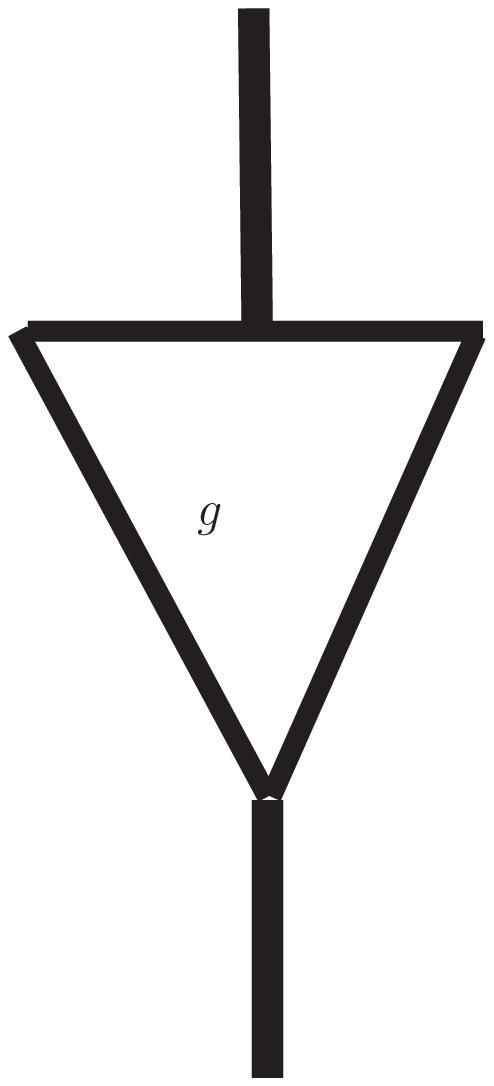}}
\ \ 
\parbox{0.8cm}{\psfrag{g}{$\scriptscriptstyle{g}$}
\includegraphics[width=0.8cm]{Fig73.eps}}
=2\ 
\parbox{1.6cm}{\includegraphics[width=1.6cm]{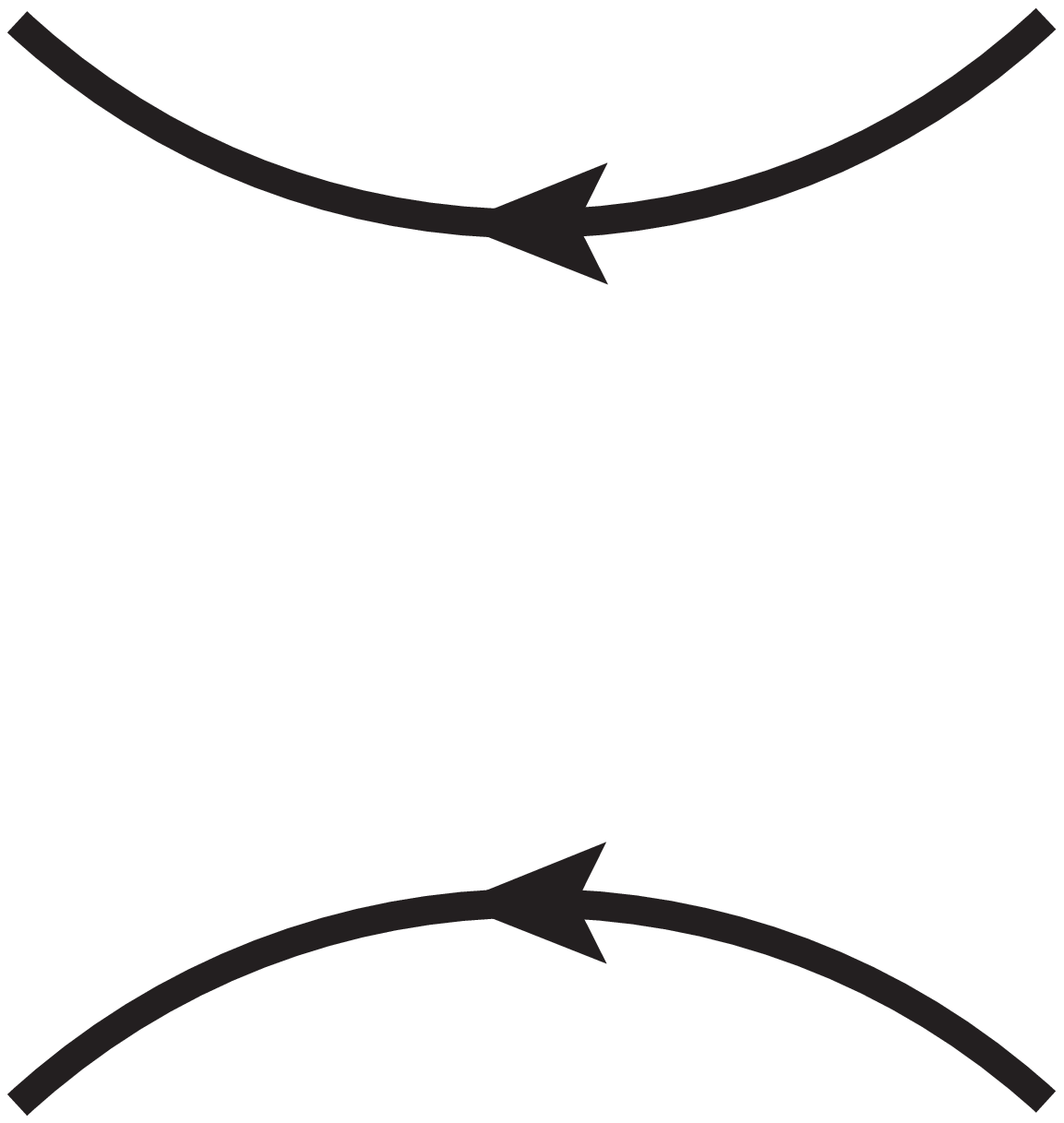}}\ \ .
\]
One has severed the line of communication between
the $I$ blob and the $F$'s and one obtains a sum of diagrams as described
in Part 1) times pure scalars (or `reduced tensor elements')
of the form
\[
\parbox{3.5cm}{\psfrag{I}{$I$}
\includegraphics[width=3.5cm]{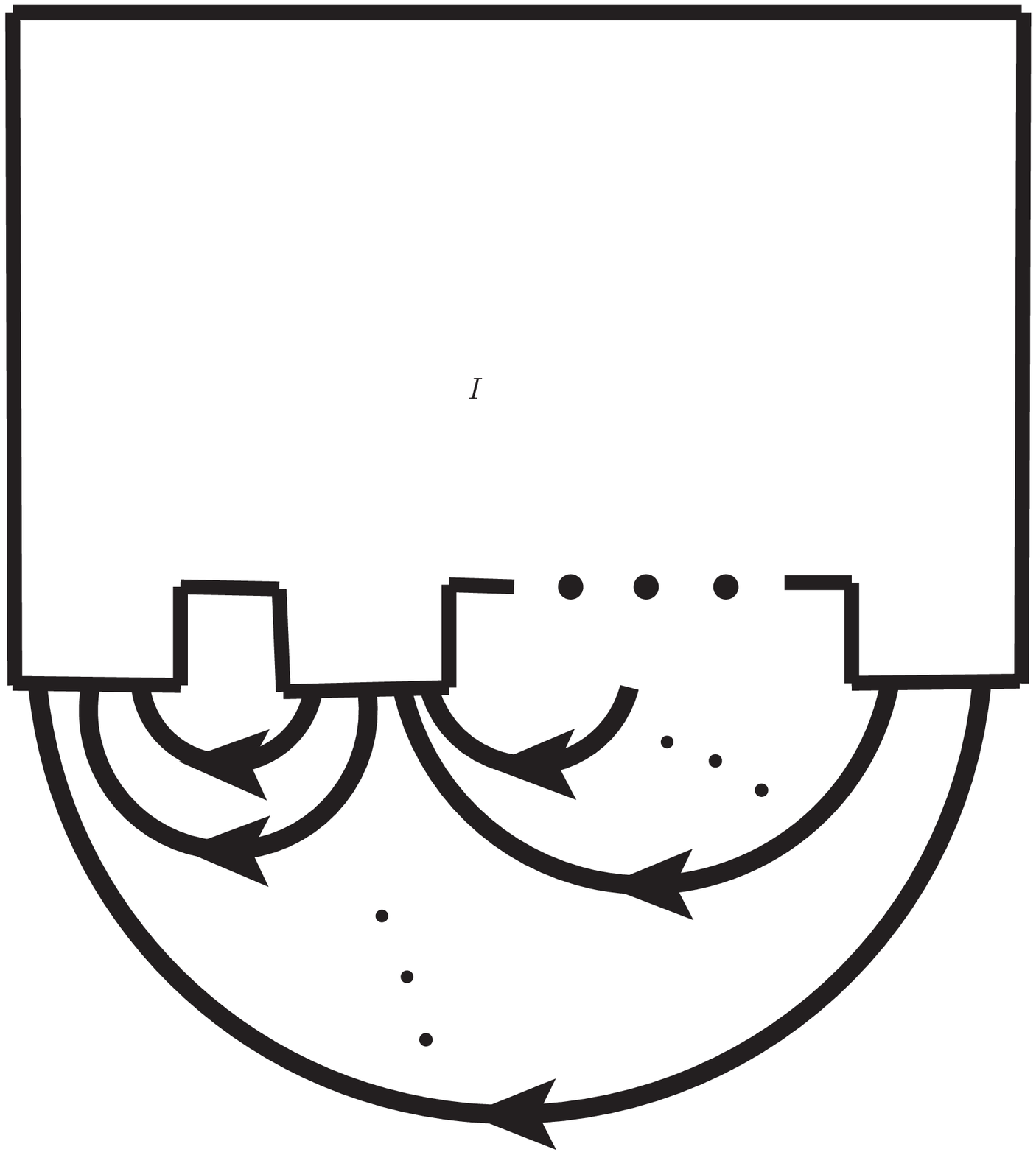}}
\]
where all legs are contracted by $\ep$ arrows.
This is essentially what physicists call the Wigner-Eckart Theorem,
or rather a microscopic version of it (see~\cite[\S5.3]{Cvitanovic}).

At this point one could object to this rewriting of the history of 19th century
invariant theory using post-WW II
Feynman diagrams, and ask where are these graphs to be found in the
classical literature?
One can answer: in~\cite{Cayley2,Clifford,Sylvester,Buchheim,Kempe1,Kempe2},
see also~\cite{Sabidussi}.
A more important point is that classics went beyond graphs and devised
a formalism in order to encode graphs by compact algebraic expressions:
the symbolic method.
The availability of this algebraic formalism is another
possible explanation to be added to those listed in~\cite[\S4.9]{Cvitanovic}
as per
the relative rarity of graphs in the printed CIT literature.

\begin{Conjecture} (For historians of mathematics)
If one could get a hold of handwritten notes and papers by Paul
Gordan, especially from the time when he wrote~\cite{Gordan1},
one should find many `birdtracks' in the sense of~\cite{Cvitanovic}.
\end{Conjecture}

Going back to the example of the canonisant of a binary quintic
(\ref{canonisant}), one can rewrite it using (\ref{graphdiff}) as
\[
\parbox{2.7cm}{\psfrag{F}{$F$}\psfrag{x}{$x$}
\includegraphics[width=2.7cm]{Fig65.eps}}
=\frac{1}{5!^3}\times
\]
\[
\parbox{2.4cm}{\psfrag{F}{$F$}\psfrag{d}{$\scriptscriptstyle{\partial a}$}
\includegraphics[width=2.4cm]{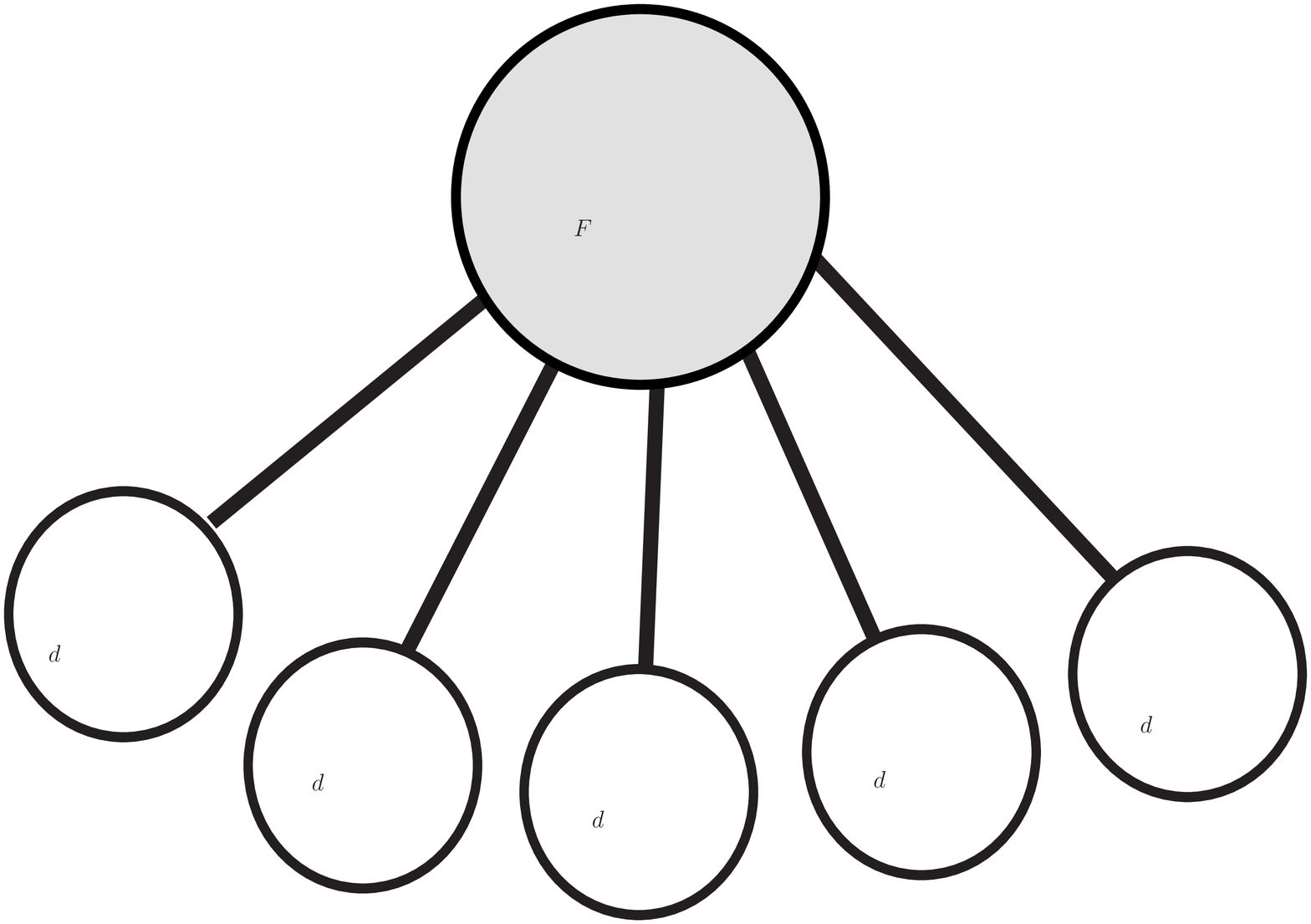}}\ 
\parbox{2.4cm}{\psfrag{F}{$F$}\psfrag{d}{$\scriptscriptstyle{\partial b}$}
\includegraphics[width=2.4cm]{Fig76.eps}}\ 
\parbox{2.4cm}{\psfrag{F}{$F$}\psfrag{d}{$\scriptscriptstyle{\partial c}$}
\includegraphics[width=2.4cm]{Fig76.eps}}
\parbox{3.7cm}{\psfrag{x}{$\scriptstyle{x}$}
\psfrag{a}{$\scriptstyle{a}$}\psfrag{b}{$\scriptstyle{b}$}
\psfrag{c}{$\scriptstyle{c}$}
\includegraphics[width=2.7cm]{Fig77.eps}}
\]
\[
\begin{array}{c}
\ \\
\ 
\end{array}
\]
\[
=\cD\ \cS
\]
where $\cD$ is the differential operator
\[
\frac{1}{5!^3}
\ F\left(\frac{\partial}{\partial a}\right)
\ F\left(\frac{\partial}{\partial b}\right)
\ F\left(\frac{\partial}{\partial c}\right)
\]
acting on the symbolic expression
\[
\cS=(ab)^2(ac)^2(bc)^2 a_{\ux} b_{\uy} c_{\ux}\ .
\]
Here a bracket factor $(ab)$ is shorthand for 
$\left|
\begin{array}{cc}
a_1 & b_1\\
a_2 & b_2
\end{array}
\right|$,
and $a_{\ux}=a_1 x_1+a_2 x_2$, etc. 
The symbolic expression $\cS$ is a bonafide polynomial involving
auxiliary variables $a,b,\ldots$ which play the same role as dummy
variables of integration~\cite[\S4.3]{AC2}.
To make things harder for the modern reader, the classics did not bother
writing the operator $\cD$ and wrote instead an {\em equality}
between a covariant
and its symbolic expression
\[
C(F,\ux)=(ab)^2(ac)^2(bc)^2 a_{\ux} b_{\uy} c_{\ux}\ .
\]
with the provision that the right-hand side must `interpreted' according to
the recipe:
\begin{quote}
Expand the right-hand side as a polynomial in $a,b,c,\ux$.
Keep the $x$'s as they are. Turn $a_1^{5-p} a_2^p$
into the coefficient $f_p$ of $F$, and likewise for $b_1^{5-p} b_2^p$
and $c_1^{5-p} c_2^p$.
\end{quote}
For instance, the expansion of the symbolic expression for the canonisant
produces terms such as
\[
-a_1 b_2 a_2 b_1 a_1 c_2 a_1 c_2 b_1 c_2 b_1 c_2 a_1 x_1 b_2 x_2 c_2 x_2
\]
\[
= -a_1^4 a_2\ b_1^3 b_2^2\ c_2^5\ x_1 x_2^2
\]
which is to be interpreted as
$-f_1f_2f_5\ x_1 x_2^2$.
Note that the weight is $6=(1+2+5)-2$, i.e., the sum of the
$f$ subscripts minus the number of $x_2$'s.
In particular, invariants are isobaric.

This symbolic recipe can be made precise using the umbral calculus
of~\cite{KungR}. We believe however that the interpretation using
differential operators $\cD$ (or, even better, using the integral
notation of~\cite[\S4.3]{AC2})
is conceptually much simpler and also more powerful: one can treat
some variables as symbolic, others as `actual', one can `iterate' the symbolic
representation as when substituting diagrams into blobs, etc.

Let $\cH_n={\rm Sym}^n(V^\ast)$
be the space of binary forms of order $n$.
It has a natural finite dimensional Hilbert space structure given by the
inner product
\begin{equation}
\langle F|G\rangle=\parbox{1.8cm}{\psfrag{F}{$\overline{F}$}
\psfrag{G}{$G$}
\includegraphics[width=1.8cm]{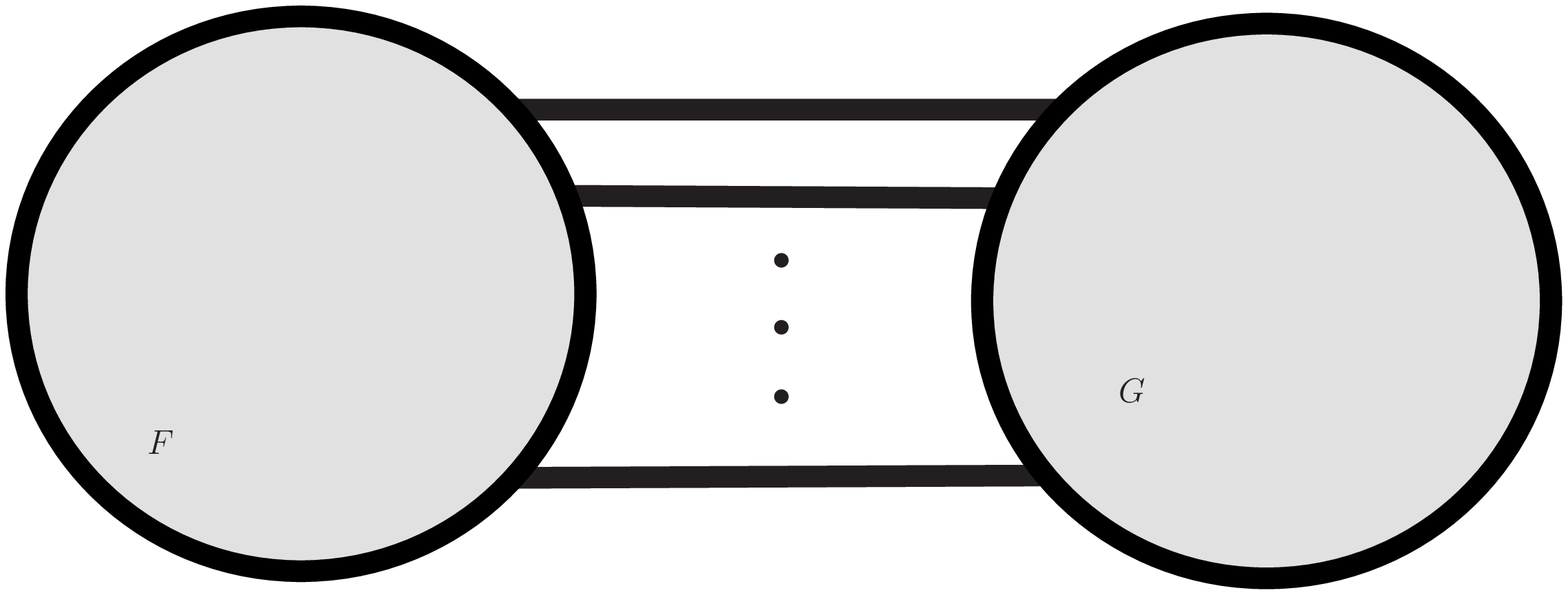}}
\label{innerprod}
\end{equation}
where $\parbox{1.2cm}{\psfrag{F}{$\scriptstyle{\overline{F}}$}\psfrag{1}{$\scriptstyle{i_1}$}
\psfrag{2}{$\scriptstyle{i_2}$}\psfrag{r}{$\scriptstyle{i_n}$}
\includegraphics[width=1.2cm]{Fig37.eps}}$
is the complex conjugate of
$\parbox{1.2cm}{\psfrag{F}{$\scriptstyle{F}$}\psfrag{1}{$\scriptstyle{i_1}$}
\psfrag{2}{$\scriptstyle{i_2}$}\psfrag{r}{$\scriptstyle{i_n}$}
\includegraphics[width=1.2cm]{Fig37.eps}}$.
This is the spin $\frac{n}{2}$ irreducible representation of $SU_2(\C)$
where the group action is that of (\ref{Ftransf}).
The tensor product $\cH_m \otimes \cH_n$
can be seen as the space of bihomogeneous forms
$B(\ux,\uy)$ of degree $m$ in $\ux=(x_1,x_2)$ and degree $n$
in $\uy=(y_1,y_2)$.
Graphically
\[
B(\ux,\uy)=
\parbox{2.2cm}{\psfrag{B}{$B$}\psfrag{1}{$\scriptstyle{i_1}$}
\psfrag{x}{$\scriptstyle{x}$}\psfrag{y}{$\scriptstyle{y}$}
\psfrag{m}{$\underbrace{\qquad\ }_{\scriptstyle{m}}$}
\psfrag{n}{$\underbrace{\qquad\ \ \ }_{\scriptstyle{n}}$}
\includegraphics[width=2.2cm]{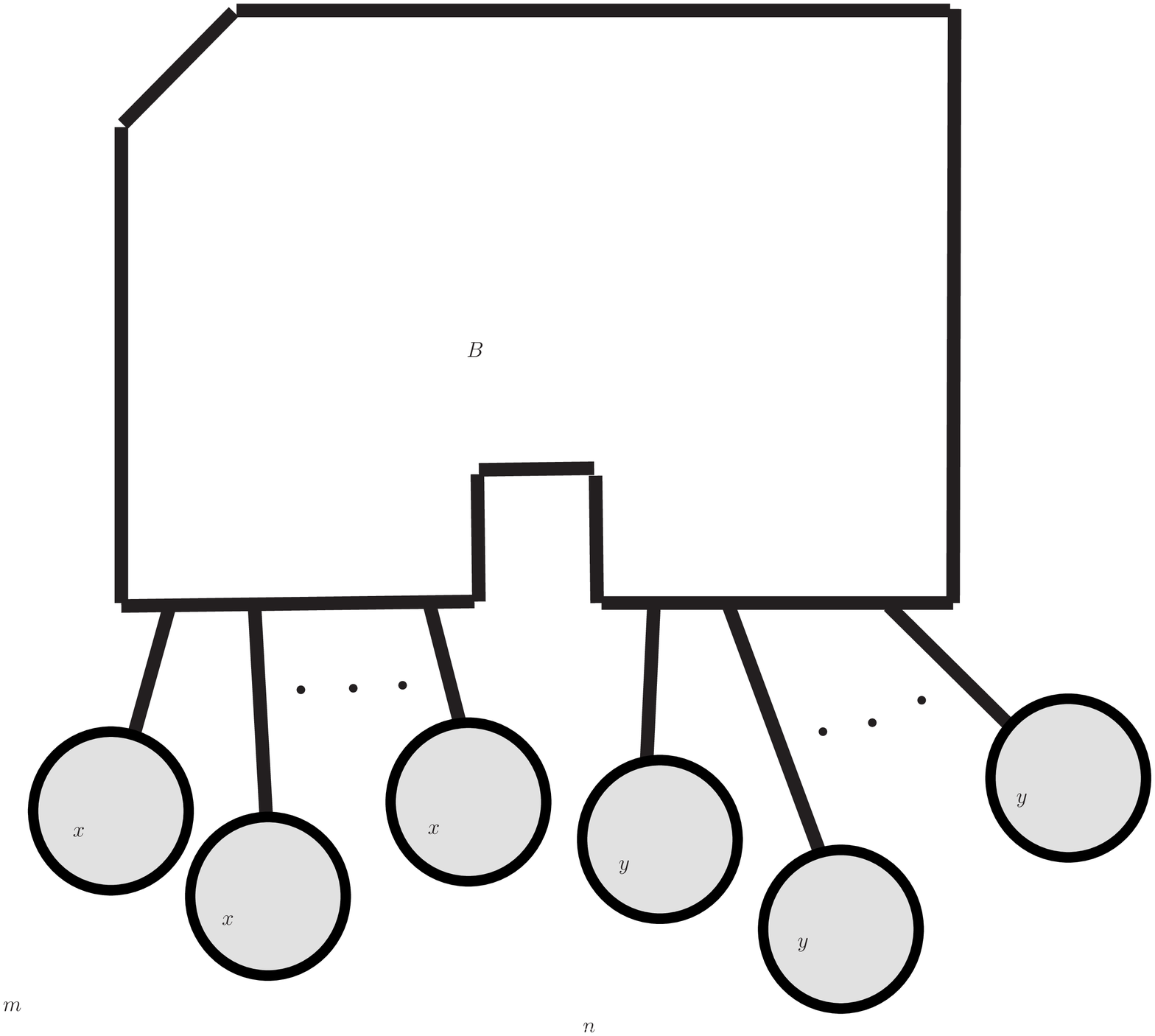}}\ \ .
\]
\[
\ 
\]
The blob of such a form $B$ is symmetric in its first $m$ indices and in its 
last $n$ indices.
One also has a natural inner product on $\cH_m \otimes \cH_n$
defined as in (\ref{innerprod}).

For any integer $k$, $0\le k\le \min(m,n)$, one has a natural equivariant map
$\cH_m \otimes \cH_n\rightarrow \cH_{m+n-2k}$ given by
\[
\parbox{2.2cm}{\psfrag{B}{$B$}
\includegraphics[width=2.2cm]{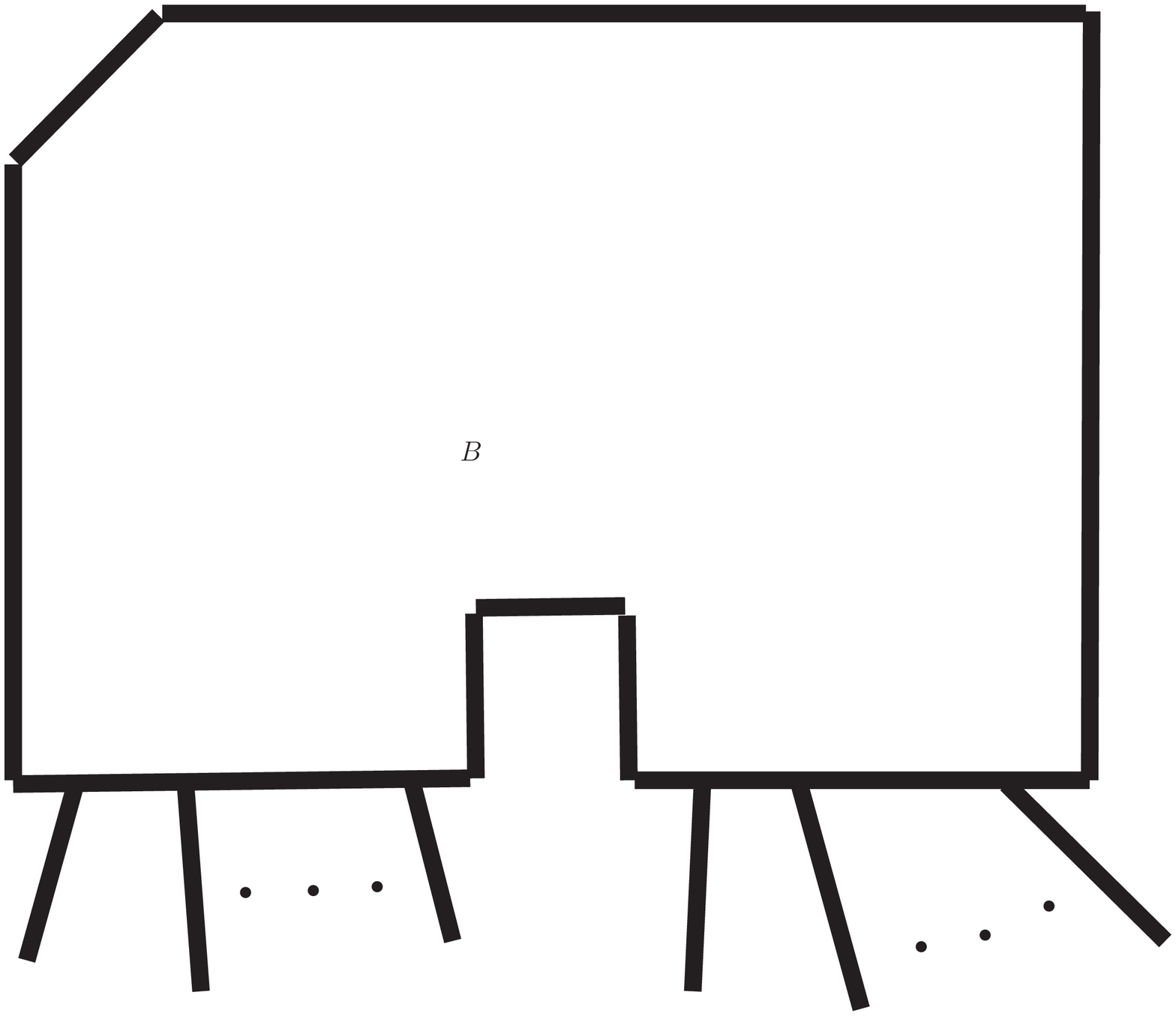}}
\longmapsto
\parbox{2.2cm}{\psfrag{B}{$B$}\psfrag{k}{$\scriptscriptstyle{k}$}
\includegraphics[width=2.2cm]{Fig83.eps}}
\]
where $k$ epsilon arrows as well as one symmetrizer are used.
This applied to a decomposable element of the form
$F\otimes G$
is called the $k$-th transvectant (or Ueberschibung)
of the binary forms $F(\ux)$ and $G(\ux)$,
denoted by
\begin{equation}
(F,G)_k=
\parbox{2.4cm}{\psfrag{F}{$F$}\psfrag{G}{$G$}
\psfrag{x}{$\scriptstyle{x}$}\psfrag{k}{$\scriptscriptstyle{k}$}
\includegraphics[width=2.4cm]{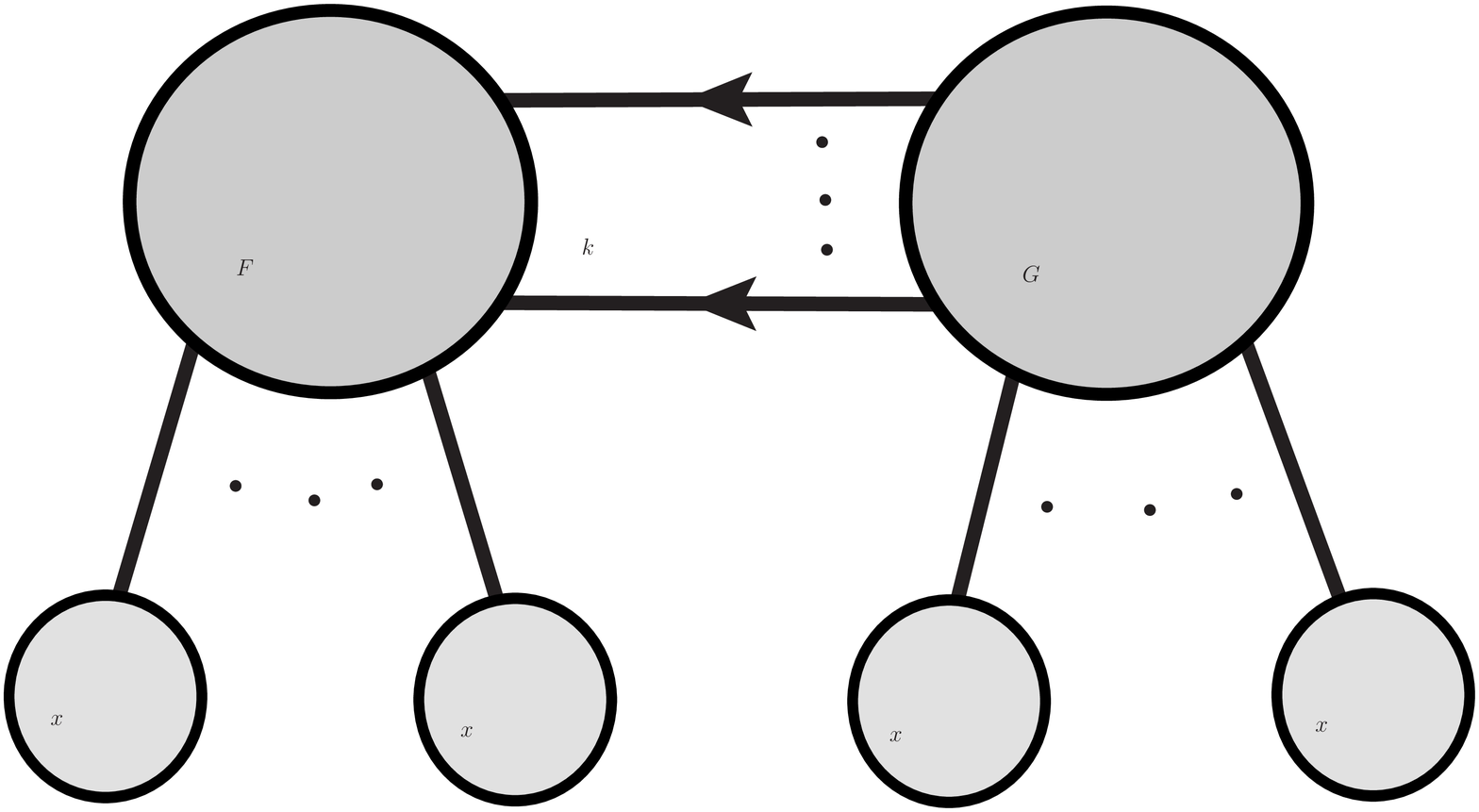}}
\label{transvectdef}
\end{equation}
or $(ab)^k a_{\ux}^{m-k}\ b_{\ux}^{n-k}$ in classical symbolic notation, where
$a,b$ are symbolic letters for $F,G$ respectively.
The beginning of $SU_2$ recoupling theory for quantum
angular momentum is the discovery by Paul Gordan and Alfred Clebsch
of the fundamental identity
\begin{equation}
\parbox{1.6cm}{\psfrag{u}{$\overbrace{\ \ }^{\scriptstyle{m}}$}
\includegraphics[width=1.6cm]{Fig85.eps}}
\ 
\parbox{1.6cm}{\psfrag{u}{$\overbrace{\ \ }^{\scriptstyle{n}}$}
\includegraphics[width=1.6cm]{Fig85.eps}}
=
\sum\limits_{k=0}^{\min(m,n)}
\frac{
\left(
\begin{array}{c}
m\\
k
\end{array}
\right)
\left(
\begin{array}{c}
n\\
k
\end{array}
\right)
}{
\left(
\begin{array}{c}
m+n-k+1\\
k
\end{array}
\right)
}
\parbox{3.7cm}{\psfrag{m}{$\scriptstyle{m}$}
\psfrag{n}{$\scriptstyle{n}$}\psfrag{k}{$\scriptstyle{k}$}
\psfrag{l}{$\scriptstyle{m+n-2k}$}
\includegraphics[width=3cm]{Fig86.eps}}
\label{Gordanseries}
\end{equation}
\[
\begin{array}{c}
\ \\
\ 
\end{array}
\]
which holds for any assignment of the $2m+2n$ indices
attached to the legs of the diagrams, consistently on both sides
of the equation.
The existence of such an identity with undetermined numerical coefficients
was proved by Gordan in~\cite[\S2]{Gordan1}.
Gordan's argument is very elegant and based on a sort of Taylor
expansion with respect to the diagonal $\ux=\uy$
of $\Proj^1\times\Proj^1$ (see also~\cite{Wahl}).
The explicit formula with binomial coefficients seems to be due
to Clebsch~\cite{Clebsch2}.
The key to the determination of these coefficients is the identity
\begin{equation}
\parbox{3.7cm}{\psfrag{m}{$\scriptstyle{m}$}
\psfrag{n}{$\scriptstyle{n}$}\psfrag{k}{$\scriptscriptstyle{k}$}
\psfrag{l}{$\scriptscriptstyle{l}$}
\psfrag{p}{$\scriptstyle{p}$}\psfrag{q}{$\scriptstyle{q}$}
\includegraphics[width=3cm]{Fig87.eps}}
=
\de_{pq}
\frac{k!\ (m+n-k+1)!\ (m-k)!\ (n-k)!}
{m!\ n!\ (m+n-2k+1)!}
\parbox{0.8cm}{\psfrag{p}{$\scriptstyle{p}$}
\psfrag{q}{$\scriptstyle{q}$}
\includegraphics[width=0.8cm]{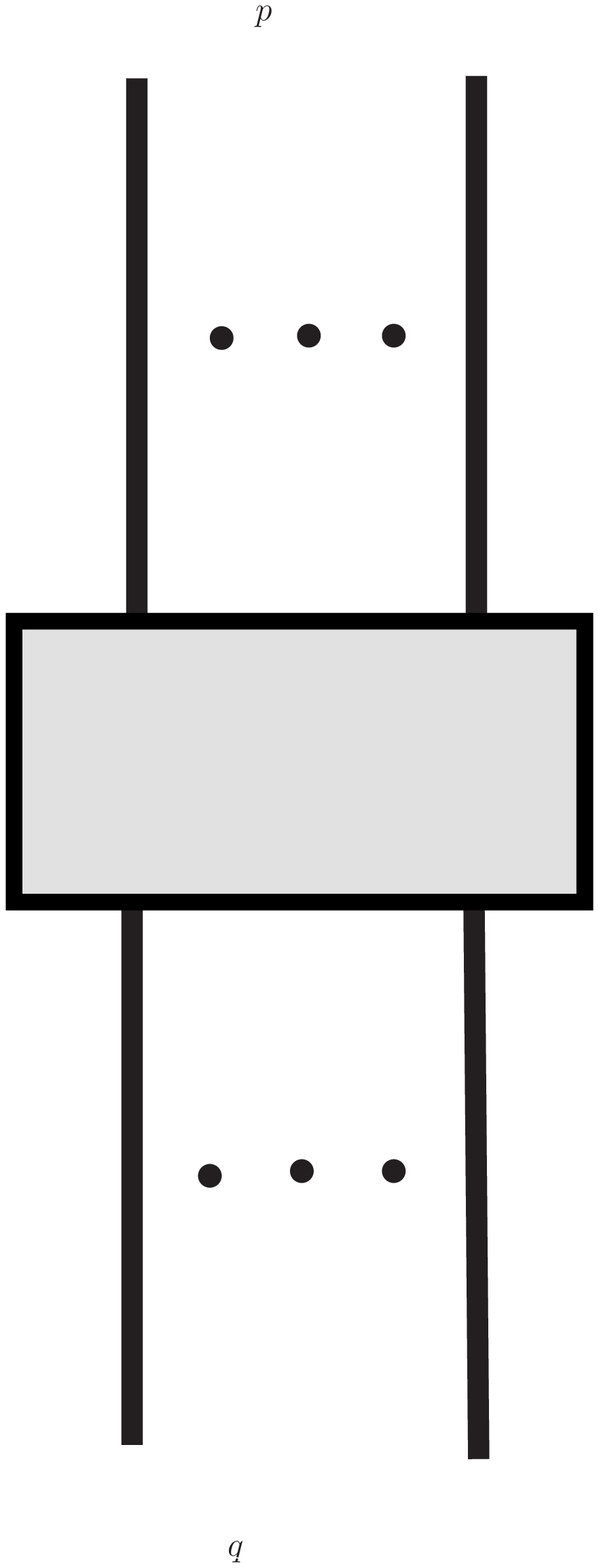}}
\label{Clebschident}
\end{equation}
\[
\begin{array}{c}
\ \\
\ 
\end{array}
\]
with $p=m+n-2k$ and $q=m+n-2l$. Indeed,
\[
\parbox{3.7cm}{\psfrag{m}{$\scriptstyle{m}$}
\psfrag{n}{$\scriptstyle{n}$}\psfrag{k}{$\scriptscriptstyle{k}$}
\psfrag{l}{$\scriptscriptstyle{l}$}
\psfrag{p}{$\scriptstyle{p}$}\psfrag{q}{$\scriptstyle{q}$}
\psfrag{a}{$\scriptstyle{a}$}\psfrag{x}{$\scriptstyle{x}$}
\includegraphics[width=3cm]{Fig89.eps}}
\ \ =\ \ 
\parbox{3cm}{\psfrag{m}{$\scriptstyle{m}$}
\psfrag{n}{$\scriptstyle{n}$}\psfrag{k}{$\scriptscriptstyle{k}$}
\psfrag{l}{$\scriptscriptstyle{l}$}
\psfrag{1}{$\overbrace{\qquad}^{\scriptstyle{m-k}}$}
\psfrag{2}{$\overbrace{\qquad}^{\scriptstyle{n-k}}$}
\psfrag{3}{$\underbrace{\qquad}_{\scriptstyle{m-l}}$}
\psfrag{4}{$\underbrace{\qquad}_{\scriptstyle{n-l}}$}
\psfrag{a}{$\scriptstyle{a}$}\psfrag{x}{$\scriptstyle{x}$}
\includegraphics[width=3cm]{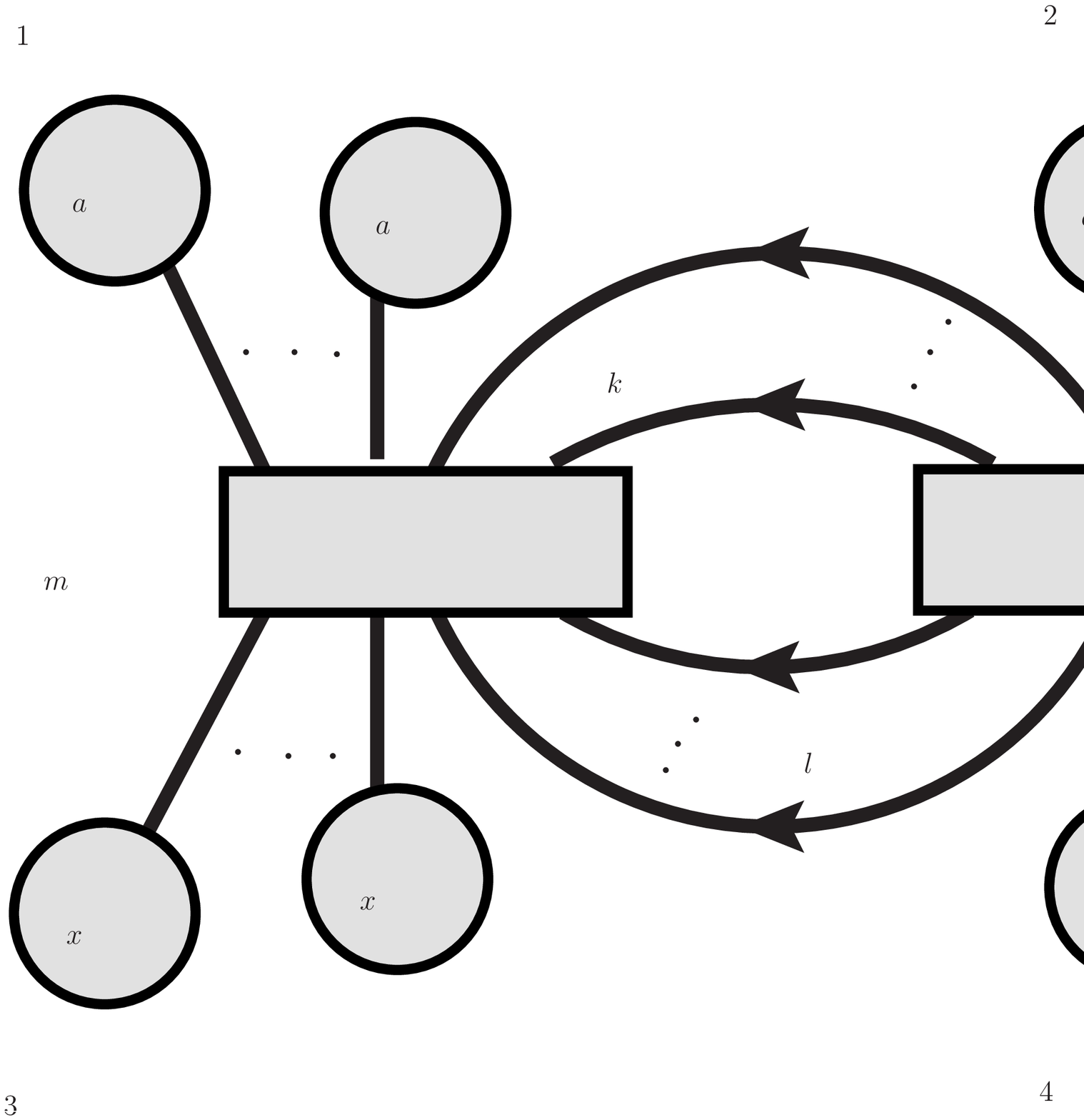}}
\ \ =\ \ 
\left.
\parbox{3cm}{\psfrag{m}{$\scriptstyle{m}$}
\psfrag{n}{$\scriptstyle{n}$}\psfrag{k}{$\scriptscriptstyle{k}$}
\psfrag{l}{$\scriptscriptstyle{l}$}
\psfrag{1}{$\overbrace{\qquad}^{\scriptstyle{m-k}}$}
\psfrag{2}{$\overbrace{\qquad}^{\scriptstyle{n-k}}$}
\psfrag{3}{$\underbrace{\qquad}_{\scriptstyle{m-l}}$}
\psfrag{4}{$\underbrace{\qquad}_{\scriptstyle{n-l}}$}
\psfrag{a}{$\scriptstyle{a}$}\psfrag{x}{$\scriptstyle{x}$}
\psfrag{y}{$\scriptstyle{y}$}
\includegraphics[width=3cm]{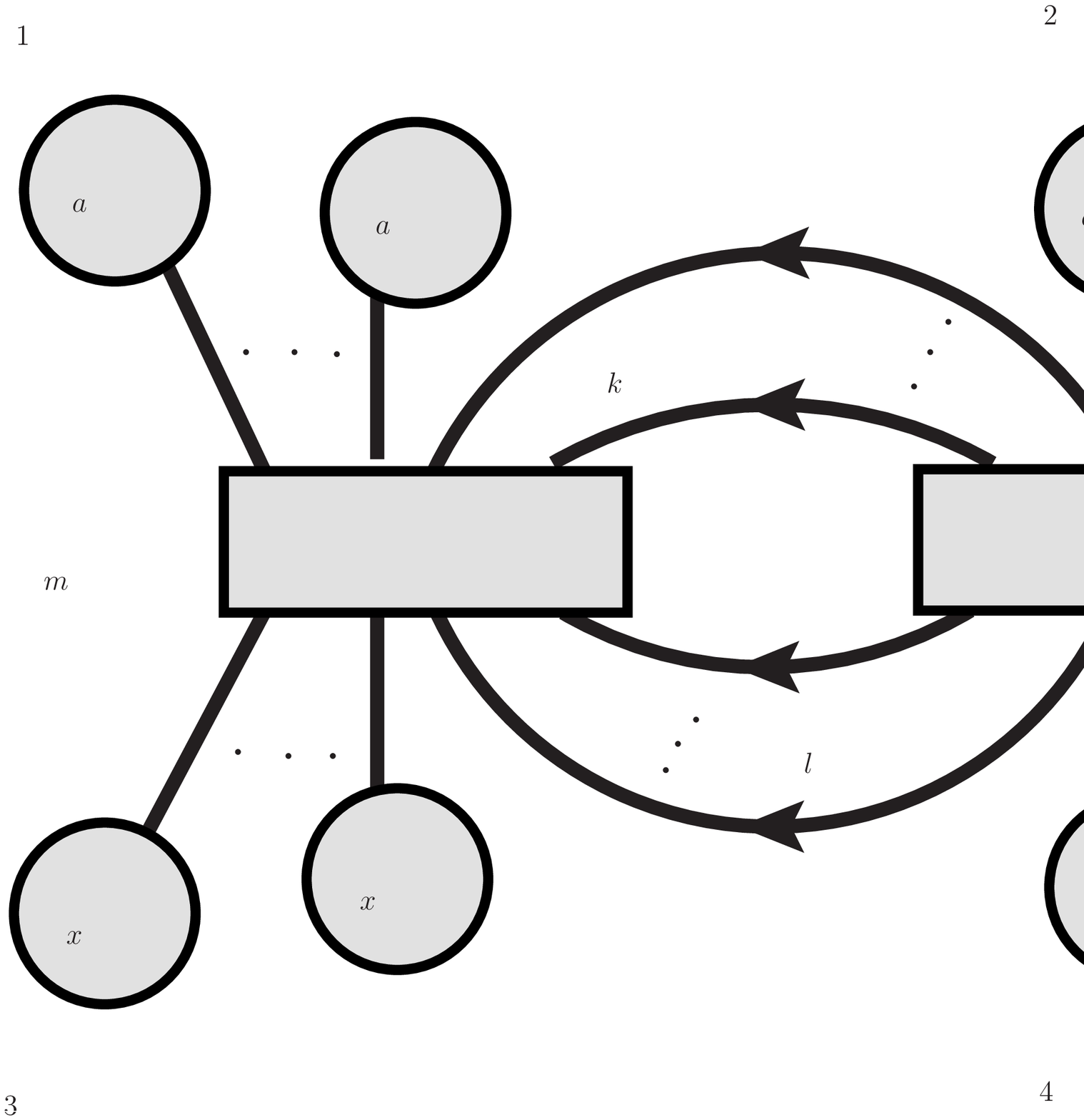}}
\ \ \right|_{\uy=\ux}
\]
\[
\begin{array}{c}
\ \\
\ 
\end{array}
\]
\[
=\frac{(m-l)!\ (n-l)!}{m!\ n!}
\left.\left\{\quad
\parbox{2cm}{\psfrag{x}{$\scriptstyle{\partial x}$}
\psfrag{y}{$\scriptstyle{\partial y}$}
\psfrag{l}{$l$}
\includegraphics[width=2cm]{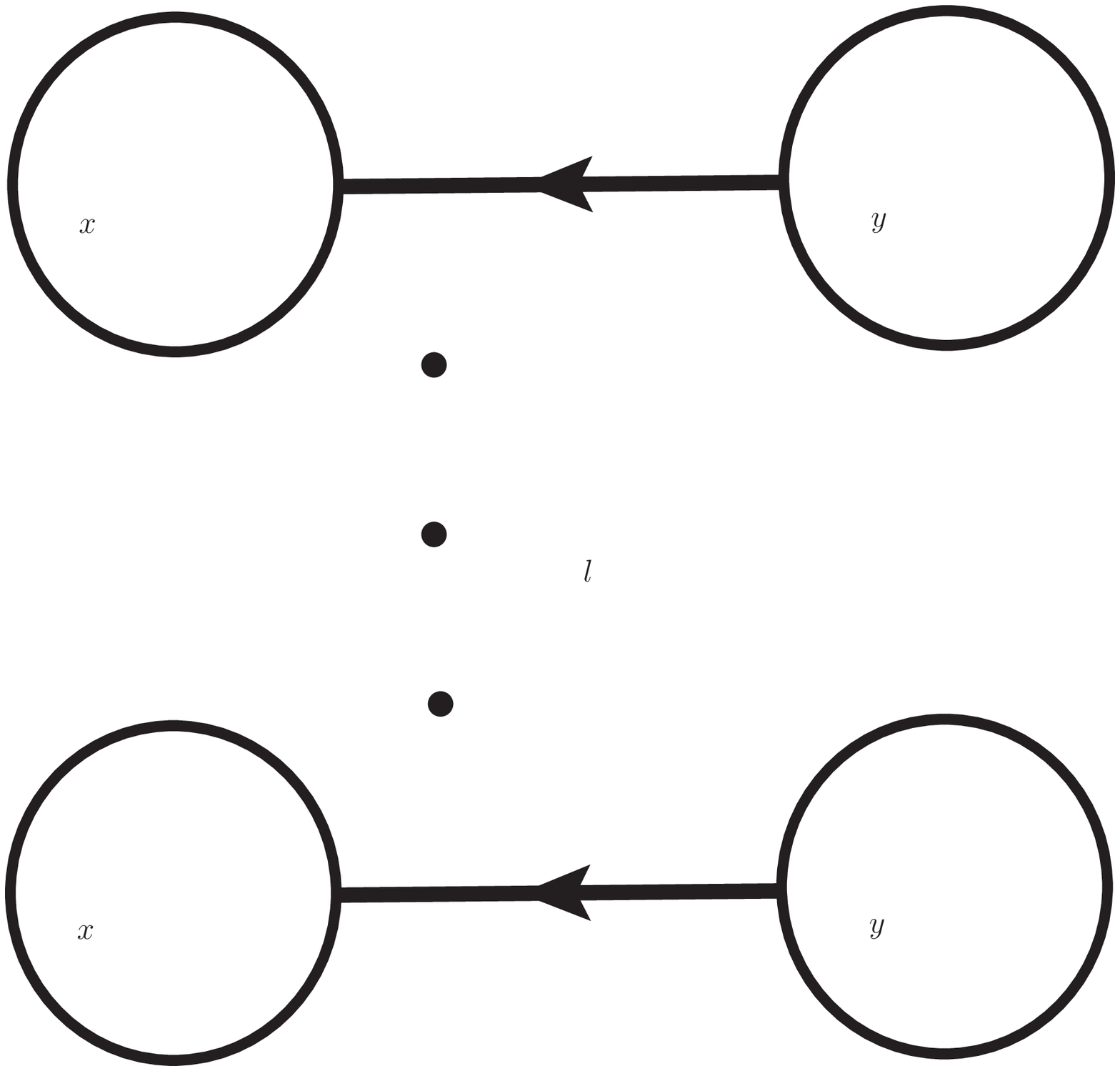}}
\qquad
\parbox{3cm}{\psfrag{m}{$\scriptstyle{m}$}
\psfrag{n}{$\scriptstyle{n}$}\psfrag{k}{$\scriptscriptstyle{k}$}
\psfrag{1}{$\overbrace{\qquad}^{\scriptstyle{m-k}}$}
\psfrag{2}{$\overbrace{\qquad}^{\scriptstyle{n-k}}$}
\psfrag{a}{$\scriptstyle{a}$}\psfrag{x}{$\scriptstyle{x}$}
\psfrag{y}{$\scriptstyle{y}$}
\includegraphics[width=3cm]{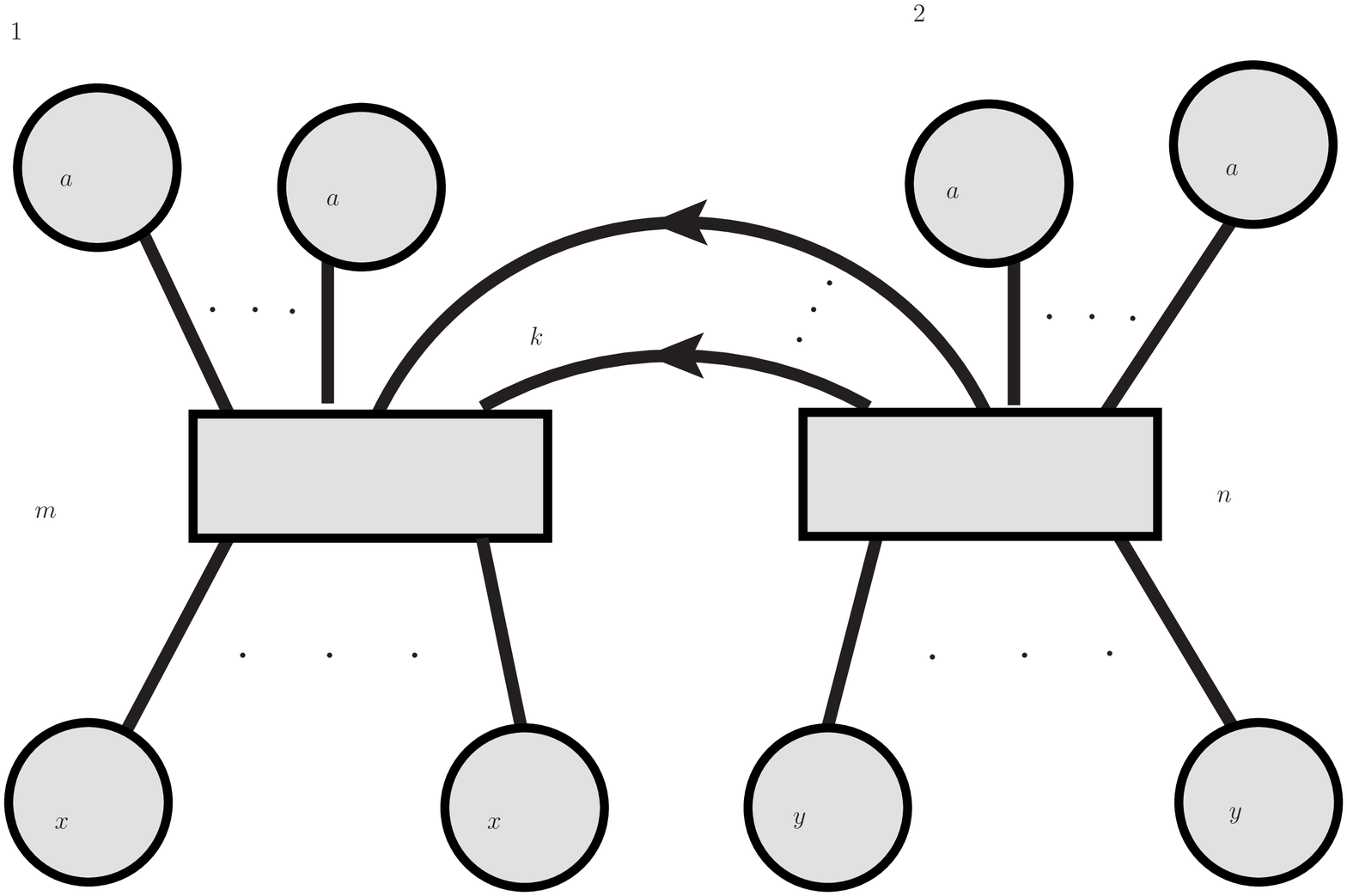}}
\quad\right\}\right|_{\uy = \ux}
\]
\[
=\frac{(m-l)!\ (n-l)!}{m!\ n!}
\left.\left\{
\Om_{\ux \uy}^l (\ux \uy)^k a_{\ux}^{m-k}\ a_{\uy}^{n-k}
\right\}\ \right|_{\uy =\ux}
\]
in classical notation where
\[
\Om_{\ux \uy}=
\parbox{2cm}{\psfrag{x}{$\scriptstyle{\partial x}$}
\psfrag{y}{$\scriptstyle{\partial y}$}
\includegraphics[width=2cm]{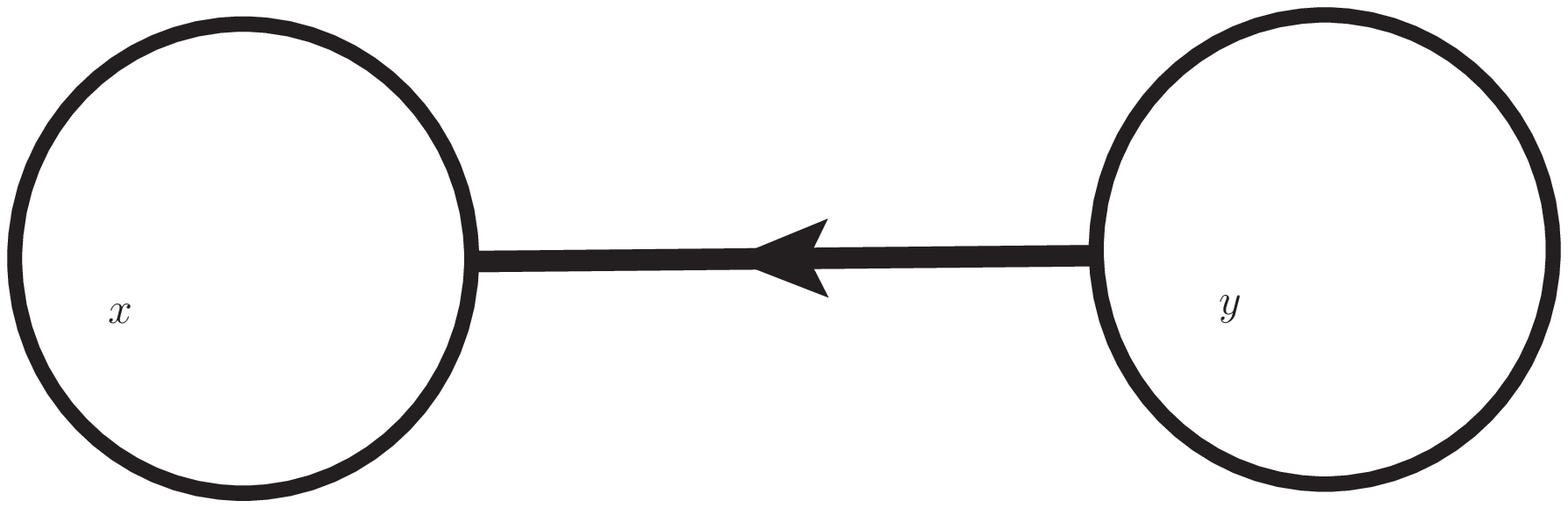}}
=\left|
\begin{array}{cc}
\frac{\partial}{\partial x_1} & \frac{\partial}{\partial y_1}\\
\frac{\partial}{\partial x_2} & \frac{\partial}{\partial y_2}
\end{array}
\right|
\]
is Cayley's Omega Operator.
Formula (\ref{Clebschident}) is thus an immediate consequence
of the easy identity
\[
\Om_{\ux \uy} \left[
(\ux \uy)^k a_{\ux}^{m-k} a_{\uy}^{n-k}\right]
=\left\{
\begin{array}{ll}
0 & {\rm if}\ k=0 ,\\
k(m+n-k+1)\ (\ux \uy)^{k-1}\ a_{\ux}^{m-k}\ a_{\uy}^{n-k} & {\rm if}\ k\ge 1 .
\end{array}
\right.
\]
Without the $a$'s this is the simplest nontrivial
case of the so-called
Cayley Identity which is nowhere to be found in Cayley's
work~\cite{AbdesselamCS}.

\begin{Remark}
We are not aware of a higher-dimensional generalization
of the identity with the $a$'s. This is related to the study
of 3-j symbols (the Theta graph rather than the 3-jm Wigner symbols)
for the group $U_n$~\cite{ElvangCK}. The classical references
related to this issue can be found in~\cite{Littlewood}.
Note that even the admissibility condition becomes nontrivial
since it is given by the Littlewood-Richardson rule.
\end{Remark}

The existence of identity (\ref{Gordanseries}) was proved by Gordan in
order to show that every covariant of a binary form $F$ is a
linear combination of compounded transvectants
of the form
$(\cdots((F,F)_{k_1},F)_{k_2},\ldots,F)_{k_p}$
These can be depicted graphically using the FDC introduced above.
For instance if $F$ is of order $d$
\[
(((F,F)_{k_1},F)_{k_2},F)_{k_3}=
\parbox{4cm}{\psfrag{x}{$x$}
\psfrag{F}{$F$}
\psfrag{d}{$\scriptstyle{d}$}
\psfrag{1}{$\scriptstyle{k_1}$}
\psfrag{2}{$\scriptstyle{k_2}$}
\psfrag{3}{$\scriptstyle{k_3}$}
\psfrag{4}{$\scriptstyle{2d-2k_1}$}
\psfrag{5}{$\scriptstyle{3d-2k_1-2k_2}$}
\psfrag{6}{$\scriptstyle{4d-2(k_1+k_2+k_3)}$}
\includegraphics[width=4cm]{Fig95.eps}}\qquad\qquad .
\]
\[
\begin{array}{c}
\ \\
\ \\
\ 
\end{array}
\]
This can be represented `macroscopically' or in shorthand
notation as
\[
\parbox{5cm}{\psfrag{x}{$x$}
\psfrag{F}{$F$}
\psfrag{d}{$\scriptstyle{d}$}
\psfrag{4}{$\scriptstyle{2d-2k_1}$}
\psfrag{5}{$\scriptstyle{3d-2k_1-2k_2}$}
\psfrag{6}{$\scriptstyle{4d-2(k_1+k_2+k_3)}$}
\includegraphics[width=5cm]{Fig96.eps}}\qquad\qquad .
\]
\[
\begin{array}{c}
\ \\
\ 
\end{array}
\]
In other words, classical covariants of binary forms expressed
using iterated transvectants are examples of
open-ended spin networks as considered in~\cite{Penrose1}
where external legs here correspond to the $F$ blobs
or a collection of $x$ blobs.
The orientation of edges keeps track of the distinction between
$V$ and $V^\ast$, even though they
are the same as $SL_2$ representations, which is special to the
2d situation.
Here we will mainly be concerned with vacuum diagrams without external legs
or what we earlier called a CG network.

An example of such a structure $(G,\cO,\ta,\ga)$ is the Theta network:
\[
\parbox{3cm}{
\psfrag{a}{$\scriptstyle{a}$}
\psfrag{b}{$\scriptstyle{b}$}
\psfrag{c}{$\scriptstyle{c}$}
\includegraphics[width=3cm]{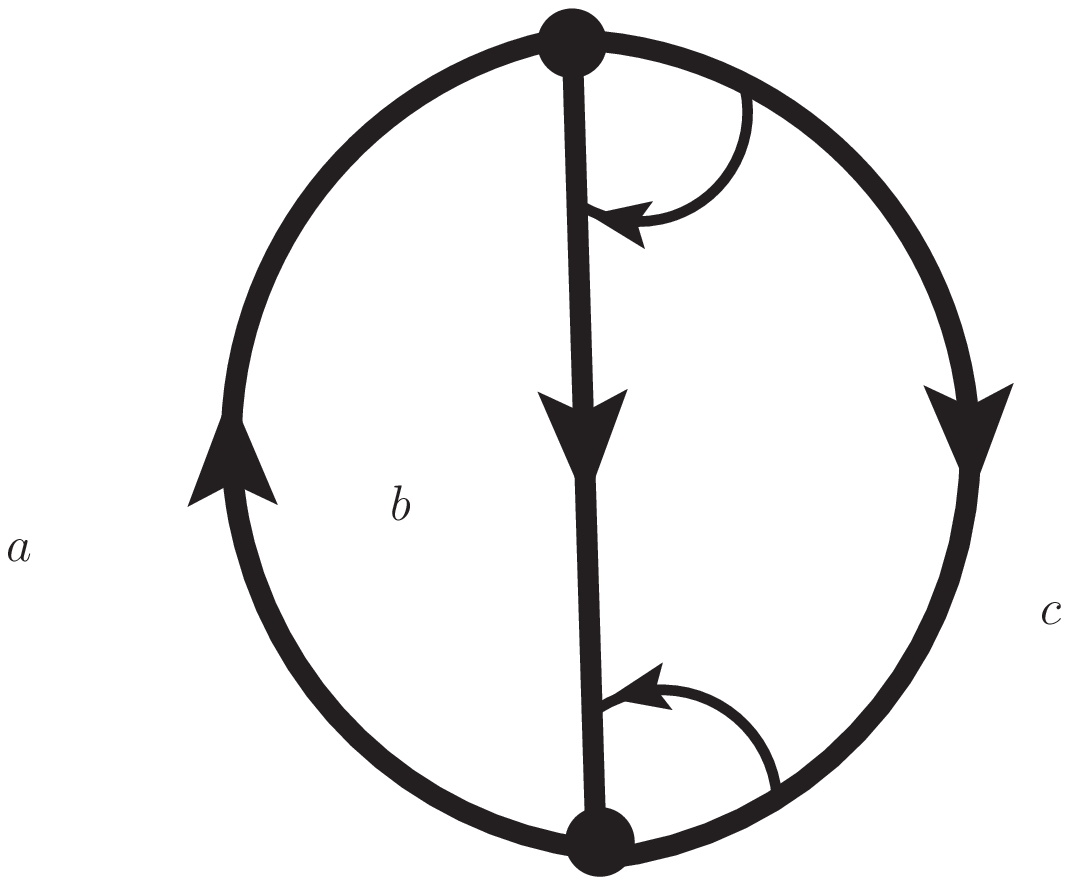}}
\]
where $G=\parbox{1cm}{
\includegraphics[width=1cm]{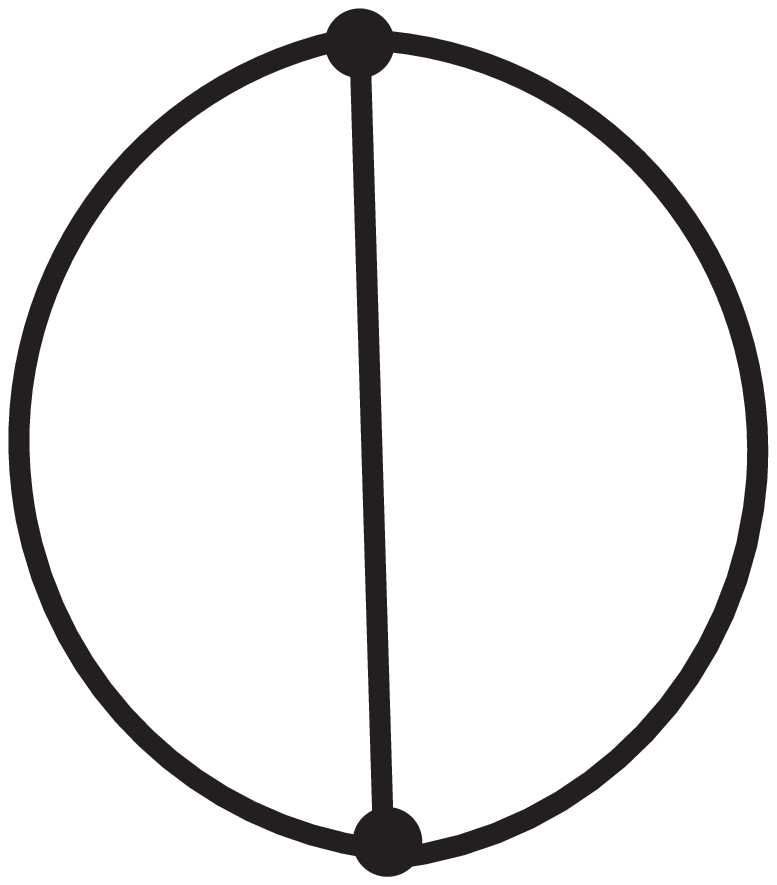}}$ is the underlying cubic graph,
$\cO$ is the given orientation of the edges of $G$.
The gate signage $\ta$ has been indicated by the two small curved arrows.
Finally the decoration $\ga$ is indicated by the nonnegative
integers $a,b,c$ which satisfy $a+b+c\in 2\N$ and
$|a-b|\le c\le a+b$.
Now, by definition,
what we call the Clebsch-Gordan evaluation of this CG network
is
\[
\<
\parbox{2.5cm}{
\psfrag{a}{$\scriptstyle{a}$}
\psfrag{b}{$\scriptstyle{b}$}
\psfrag{c}{$\scriptstyle{c}$}
\includegraphics[width=2.5cm]{Fig98.eps}}
\>^{CG}=
\parbox{3.5cm}{
\psfrag{a}{$\scriptstyle{a}$}
\psfrag{b}{$\scriptstyle{b}$}
\psfrag{c}{$\scriptstyle{c}$}
\psfrag{k}{$\scriptscriptstyle{k}$}
\includegraphics[width=3.5cm]{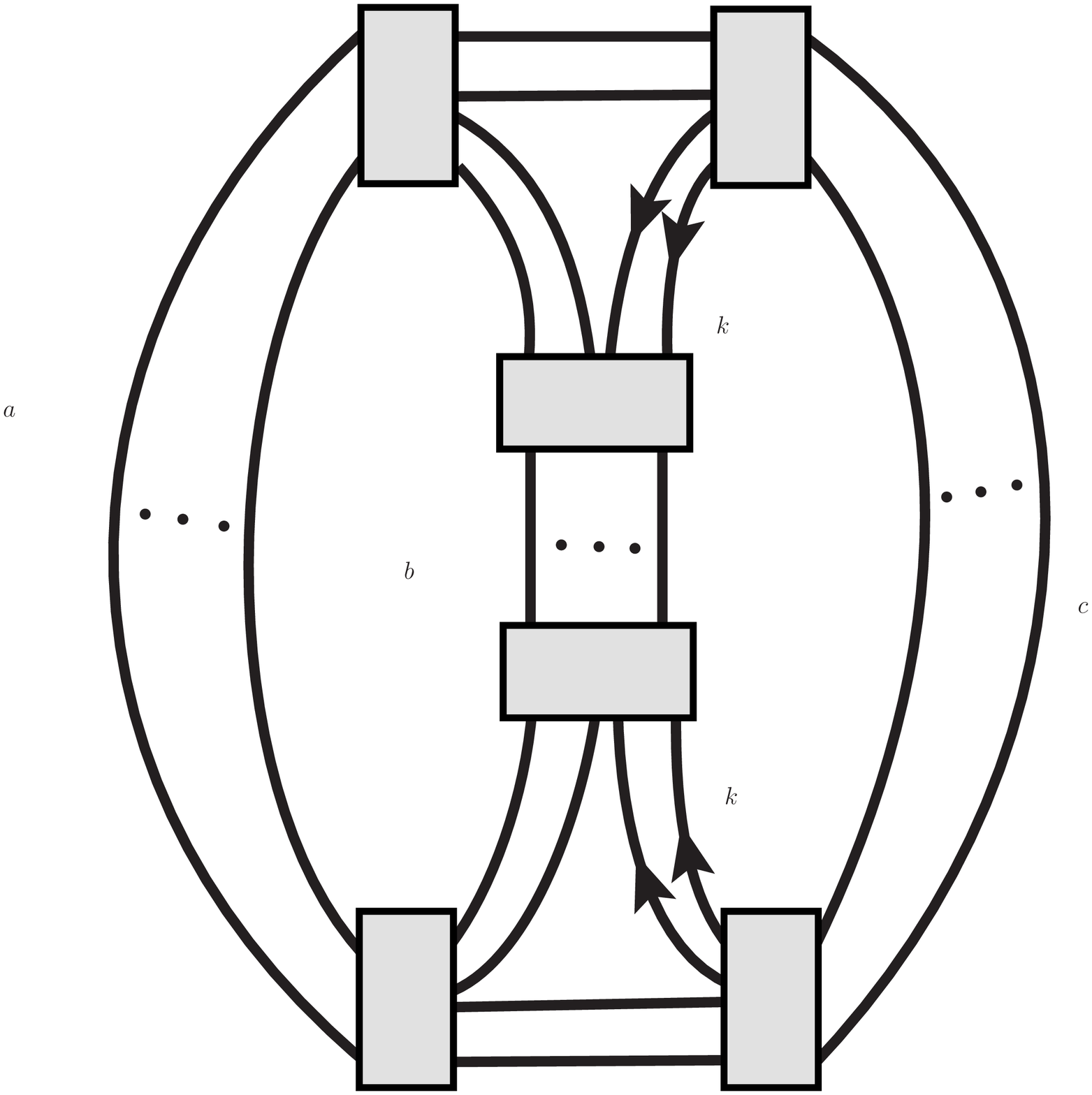}}\ \ .
\]
Here $a,b,c$ indicate the number of Kronecker delta strands
to be used.
The number of $\ep$ arrows at each trivalent vertex is
$k=\frac{b+c-a}{2}$.
The right-hand side is meant as the FDC evaluation of the given diagram
made of $\de$'s, $\ep$'s and symmetrizers as explained above.
Using the identity (\ref{Clebschident})
and the idempotence of symetrizers one immediately obtains
\[
\<
\parbox{2.5cm}{
\psfrag{a}{$\scriptstyle{a}$}
\psfrag{b}{$\scriptstyle{b}$}
\psfrag{c}{$\scriptstyle{c}$}
\includegraphics[width=2.5cm]{Fig98.eps}}
\>^{CG}=
\frac{
\left(\frac{b+c-a}{2}\right)!
\ \left(\frac{a+b+c}{2}+1\right)!
\ \left(\frac{a+b-c}{2}\right)!
\ \left(\frac{a+c-b}{2}\right)!
}{b!\ c!\ (a+1)!}
\times
\parbox{1.3cm}{
\psfrag{a}{$\scriptstyle{a}$}
\includegraphics[width=1.3cm]{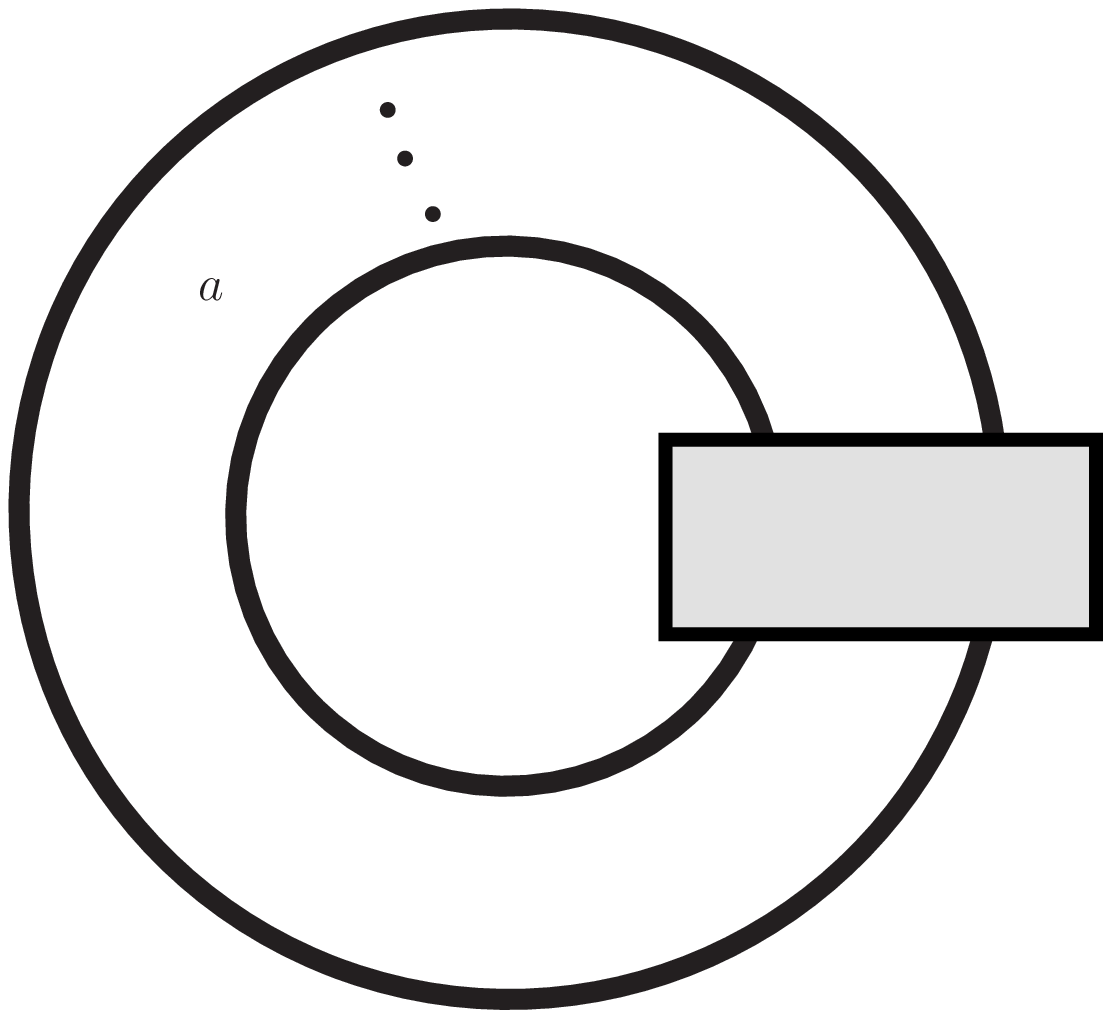}}\ .
\]
Using the same calculation and notation as in Lemma \ref{trivcomplemma}
one can write
\begin{equation}
\parbox{1.3cm}{
\psfrag{a}{$\scriptstyle{a}$}
\includegraphics[width=1.3cm]{Fig101.eps}}
=
\frac{1}{a!}\sum\limits_{\si\in\gS_a} 2^{c(\si)}
=\frac{1}{a!}\sum\limits_{k=0}^{a} 2^k c(a,k)=a+1
\label{loopevalCG}
\end{equation}
as it should be. Indeed, this is the trace of the identity
operator on ${\rm Sym}^a(V^\ast)$, i.e.,
the dimension of this representation.
Finally, the evaluation of the Theta graph in the CG formalism
is
\[
\<
\parbox{2.5cm}{
\psfrag{a}{$\scriptstyle{a}$}
\psfrag{b}{$\scriptstyle{b}$}
\psfrag{c}{$\scriptstyle{c}$}
\includegraphics[width=2.5cm]{Fig98.eps}}
\>^{CG}=
\frac{
\left(\frac{a+b+c}{2}+1\right)!
\ \left(\frac{a+b-c}{2}\right)!
\ \left(\frac{a+c-b}{2}\right)!
\ \left(\frac{b+c-a}{2}\right)!
}{a!\ b!\ c!}\ \ .
\]

More generally, we have the following definition.
\begin{Definition}\label{CGevaldef}
The CG evaluation $\<G,\cO,\ta,\ga\>^{CG}$ of a CG network
$(G,\cO,\ta,\ga)$ is obtained by:
\begin{enumerate}
\item
replacing each edge $e$ carrying a decoration $\ga(e)$ by
a number $\ga(e)$ of Kronecker delta strands,
\item
replacing a vertex
\[
\parbox{1.6cm}{\psfrag{a}{$\scriptstyle{a}$}
\psfrag{b}{$\scriptstyle{b}$}\psfrag{c}{$\scriptstyle{c}$}
\includegraphics[width=1.6cm]{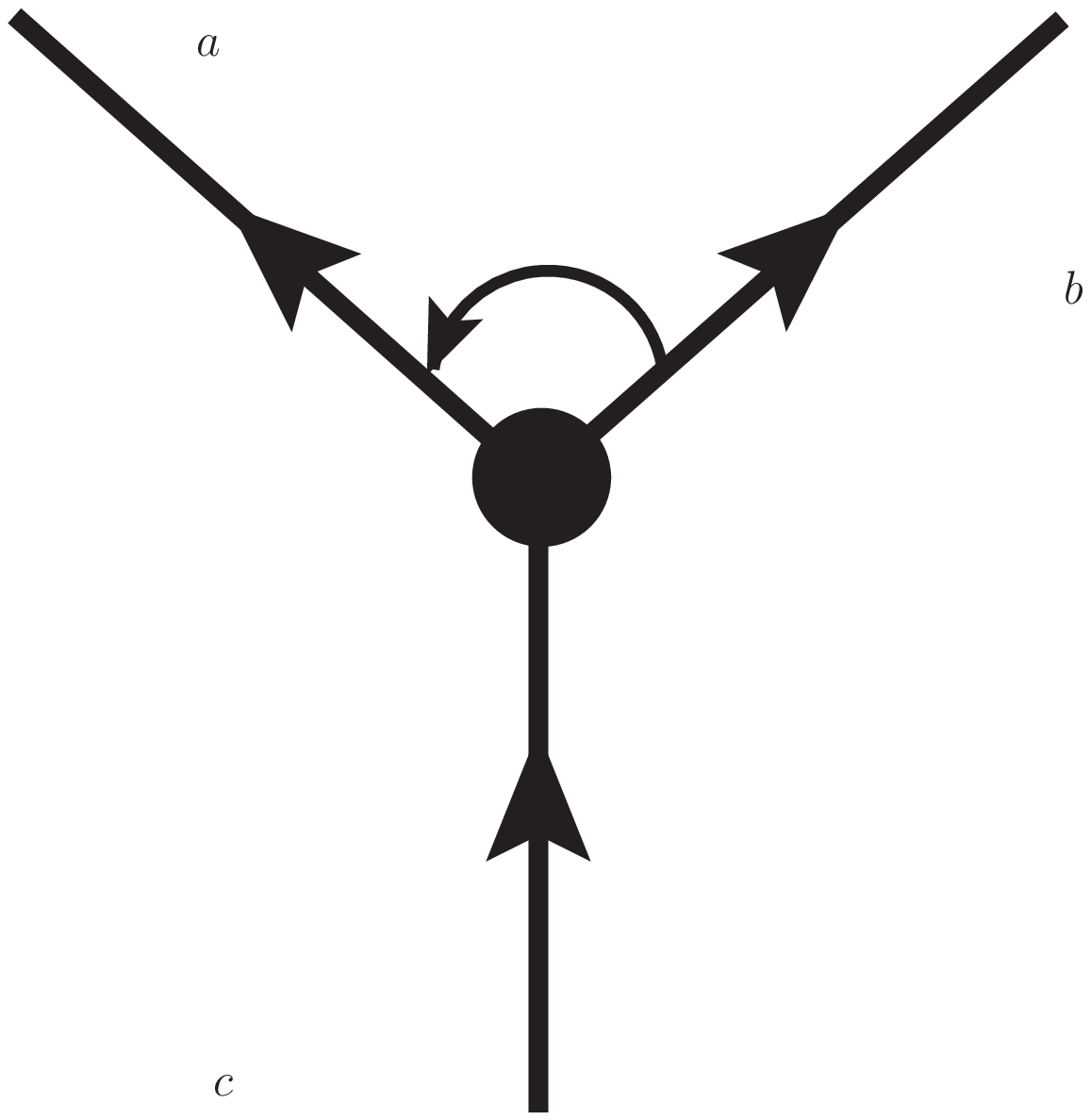}}
\qquad {\rm or}\qquad
\parbox{1.6cm}{\psfrag{a}{$\scriptstyle{a}$}
\psfrag{b}{$\scriptstyle{b}$}\psfrag{c}{$\scriptstyle{c}$}
\includegraphics[width=1.6cm]{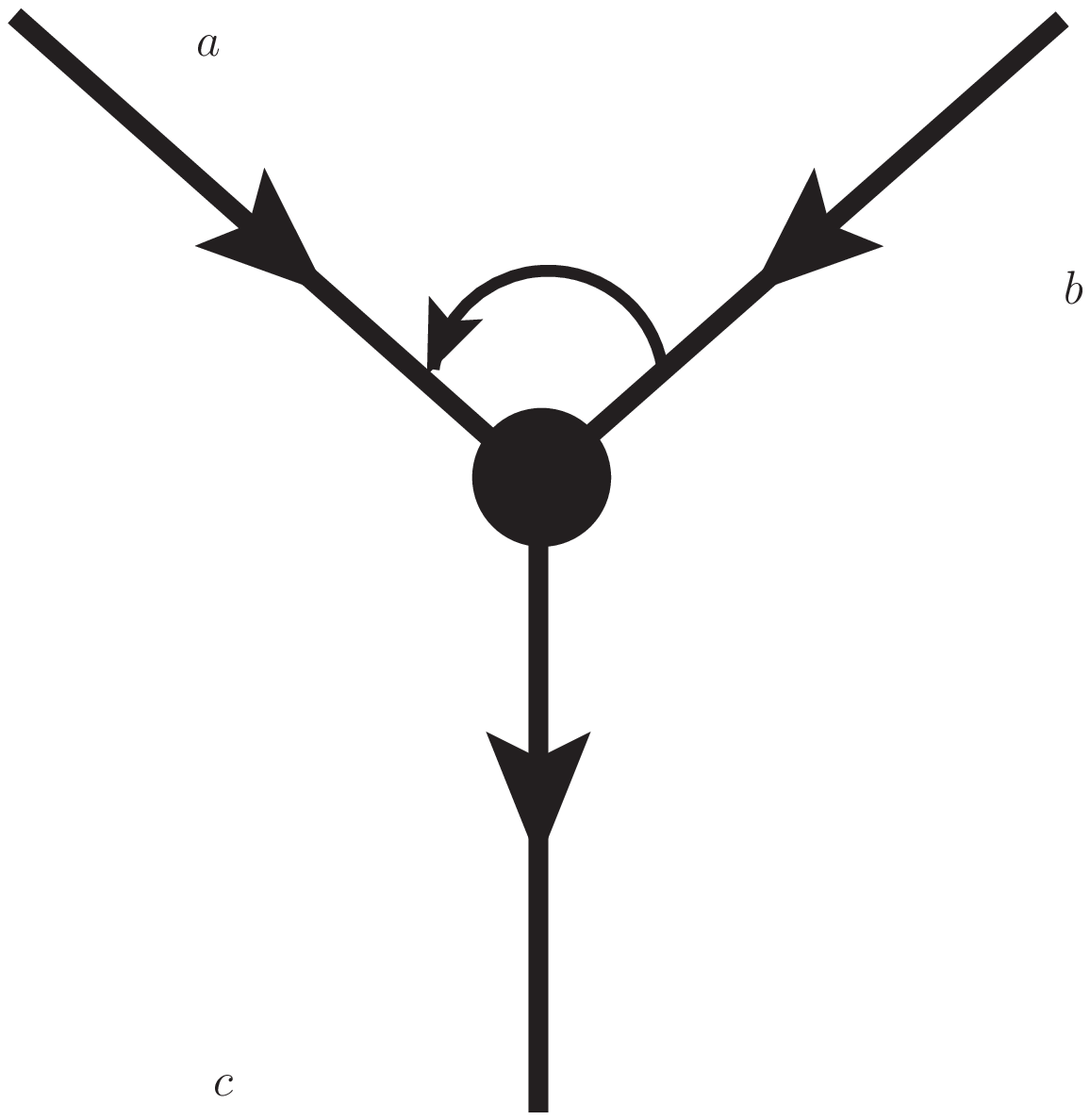}}
\qquad {\rm by}\qquad
\parbox{3.6cm}{\psfrag{a}{$\scriptstyle{a}$}
\psfrag{b}{$\scriptstyle{b}$}\psfrag{c}{$\scriptstyle{c}$}
\psfrag{l}{$\scriptstyle{\frac{a+c-b}{2}}$}
\psfrag{t}{$\scriptstyle{\frac{a+b-c}{2}}$}
\psfrag{r}{$\scriptstyle{\frac{b+c-a}{2}}$}
\includegraphics[width=3.6cm]{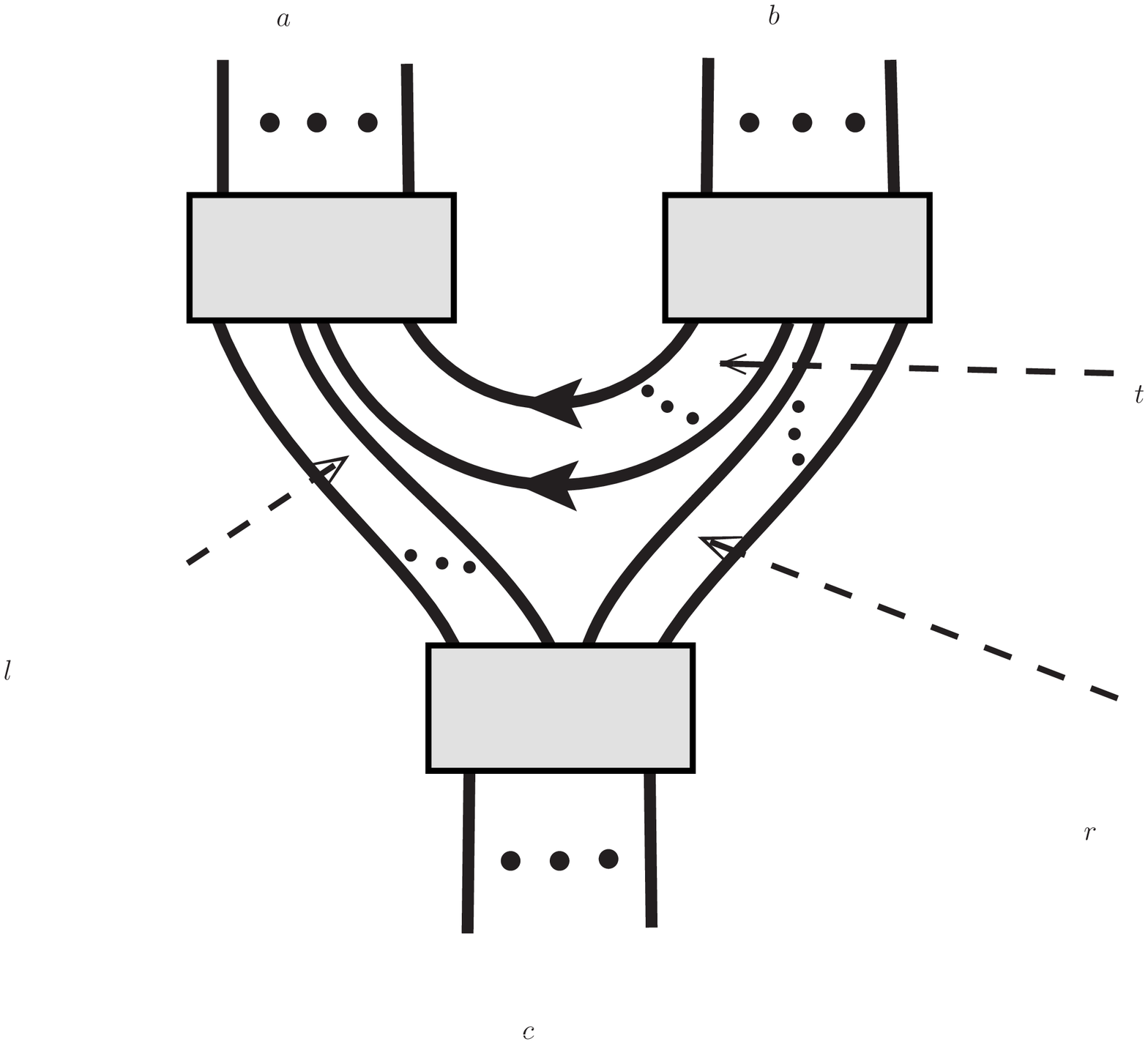}}
\]
\item
evaluating the resulting `microscopic' Feynman diagram
using the rules of FDC, with indices summed over the range $\{1,2\}$.
\end{enumerate}
\end{Definition}

\begin{Remark}
Nothing, in this formalism, relies on how the graph is drawn on the page
or any kind of surface imbedding.
This is the main difference with Penrose's binor calculus
as explained for instance in~\cite{Penrose2,Kauffman}
or~\cite[App. A]{Rovelli}.
\end{Remark}

The precise correspondence between CSN's of \S\ref{introdefsec}
and CG networks is the object of \S\ref{negdimsec}.
We will conclude this all too brief tour of CIT
by a few comments.
It has been said that the main goal of CIT is to prove
finite generation of rings of invariants for reductive groups.
Anyone who would make the effort of reading an original source which either
uses the English Omega Operator and hyperdeterminant formalism,
or the German symbolic method, will see that the goal of CIT is much more
ambitious.
It is to develop {\em explicit} algebraic geometry (not just
{\em effective} algebraic geometry), i.e., nonlinear algebra in the sense
of~\cite{DolotinM}.
Of course one wants an explicit description of a minimal
system of generators for invariant rings, but this is only a means to an end
which is to understand invariants which detect a specific geometrical
event.
The most important perhaps are the multidimensional resultants
which detect when $n$ algebraic hypersurfaces in $\Proj^{n-1}$ have a common
intersection, and multidimensional discriminants which detect when a
hypersurface is singular~\cite{GelfandKZ}.
Apart from their theoretical interest these often
are in practice the most efficient tools for elimination,
indeed faster than Gr\"obner basis methods (see,
e.g.,~\cite{Sturmfels,BresslerGPP}). Note in passing
that one of the first nontrivial examples of Gr\"obner bases
is also due to Gordan~\cite{GordanGrobner}.
A long time after Hilbert is said to have killed CIT, Gordan was still
working on the problem of symbolic forms of resultants.
In one of his last articles~\cite{Gordan2},
he succeeded in finding an explicit symbolic expression
for the resultant of two binary forms $F,G$
of the same order $d$. His formula amounts to a cycle expansion
of the Bezout determinant.
When translated in macroscopic FDC, it is an explicit linear
combination of products of `wheels'
\[
\parbox{5cm}{\psfrag{x}{$x$}
\psfrag{F}{$F$}
\psfrag{G}{$G$}
\psfrag{d}{$\scriptstyle{d}$}
\includegraphics[width=5cm]{Fig105.eps}}\qquad\qquad .
\]
\[
\begin{array}{c}
\ \\
\ 
\end{array}
\]
These essentially are spin networks with external $F$ and $G$ legs.
Finding similar symbolic or microscopic FDC formulae for the
multidimensional case is a very old and very open problem.
It was briefly hinted at in~\cite[Introduction \S VIII]{GelfandKZ}.
The reader who would like to learn more about the CIT of binary forms
must consult the
works of classical masters such as Clebsch, von Gall, Gordan, Stroh, Young, etc.
The introduction via FDC provided in this section
should make the reading much easier.

\section{Clebsch-Gordan networks with external legs and $SL_2$ invariance}
\label{legsection}

In this section, we will slightly generalize the notion of CG networks
$(G,\cO,\ta,\ga)$ by allowing $G$
to have 1-valent vertices.
So now $G$ is any graph with vertices of degree 1 or 3.
The orientation $\cO$ of the edges is smooth in the sense that 3-valent
vertices must have $({\rm indegree},{\rm outdegree})$
equal to $(1,2)$ or $(2,1)$.
The gate signage $\ta$ is an ordering of the gates at each
3-valent vertex and $\ga$ is, as before, a decoration
of the edges by an admissible collection of nonnegative integers.
The 1-valent vertices impose no new constraints as far as admissiblility goes.
We will use the expression `external leg'
of the network indiscriminately for a 1-valent vertex or for the unique
edge incident to it.
A 1-valent vertex $\parbox{0.6cm}{\includegraphics[width=0.6cm]{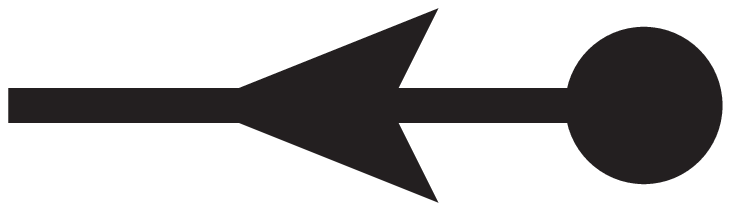}}$
of outdegree 1 is called an entry vertex.
A 1-valent vertex $\parbox{0.6cm}{\includegraphics[width=0.6cm]{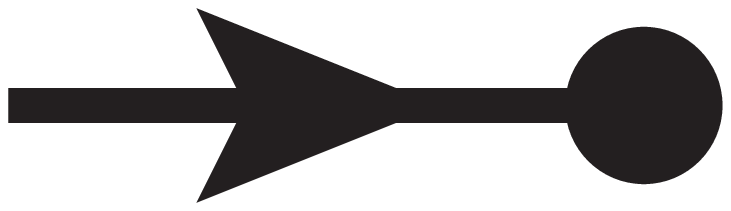}}$
of indegree 1 is called an exit vertex.
We denote by $V_{\rm in}(G)$, $V_{\rm out}(G)$
the set of entry, exit vertices respectively.
We will commit a slight abuse of notation by writing $\ga(v)$
for the decoration of the unique edge incident to a 1-valent vertex $v$.

We will now define as in Def. \ref{CGevaldef}
the Clebsch-Gordan evaluation $\<G,\cO,\ta,\ga\>^{CG}$ of such a network
which will be a tensor living in
\[
\left(
\bigotimes\limits_{v\in V_{\rm in}(G)}
\cH_{\ga(v)}^{\ast}
\right)
\otimes
\left(
\bigotimes\limits_{v\in V_{\rm out}(G)}
\cH_{\ga(v)}
\right)
\]
where $\cH_a={\rm Sym}^a(V^{\ast})$,
$\cH_a^{\ast}$ is the dual ${\rm Sym}^a(V)$,
and $V=\C^2$ as in \S\ref{CGsection}.
The rules for constructing $\<G,\cO,\ta,\ga\>^{CG}$
are the same as before except that a 1-valent vertex gives
\[
\parbox{0.8cm}{\psfrag{a}{$\scriptstyle{a}$}
\psfrag{o}{${\rm or}$}
\includegraphics[width=0.8cm]{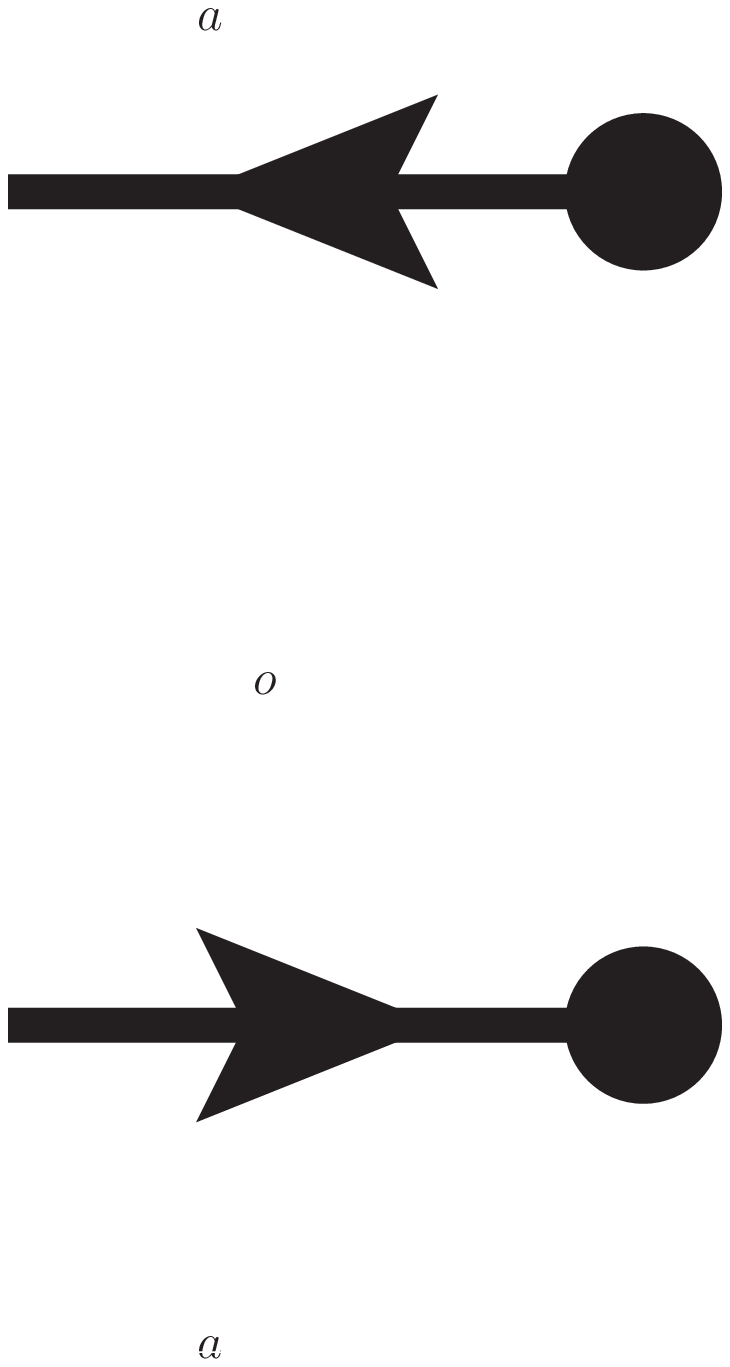}}
\qquad \longrightarrow \qquad
\left.
\parbox{1.4cm}{\psfrag{1}{$\scriptstyle{i_1}$}
\psfrag{2}{$\scriptstyle{i_2}$}\psfrag{a}{$\scriptstyle{i_a}$}
\includegraphics[width=1.4cm]{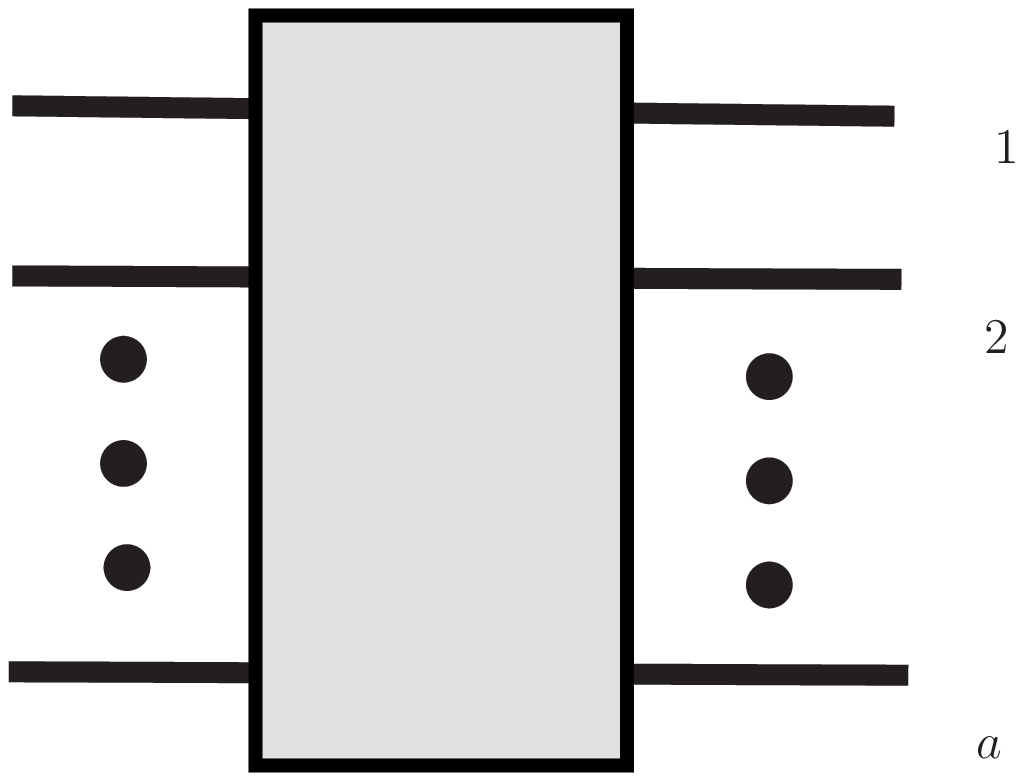}}
\ \right\}\ {\rm external\ indices}\qquad .
\]
For instance identity (\ref{Clebschident})
is
\begin{equation}
\<
\parbox{1.6cm}{\psfrag{p}{$\scriptstyle{p}$}
\psfrag{q}{$\scriptstyle{q}$}\psfrag{m}{$\scriptstyle{m}$}
\psfrag{n}{$\scriptstyle{n}$}
\includegraphics[width=1.6cm]{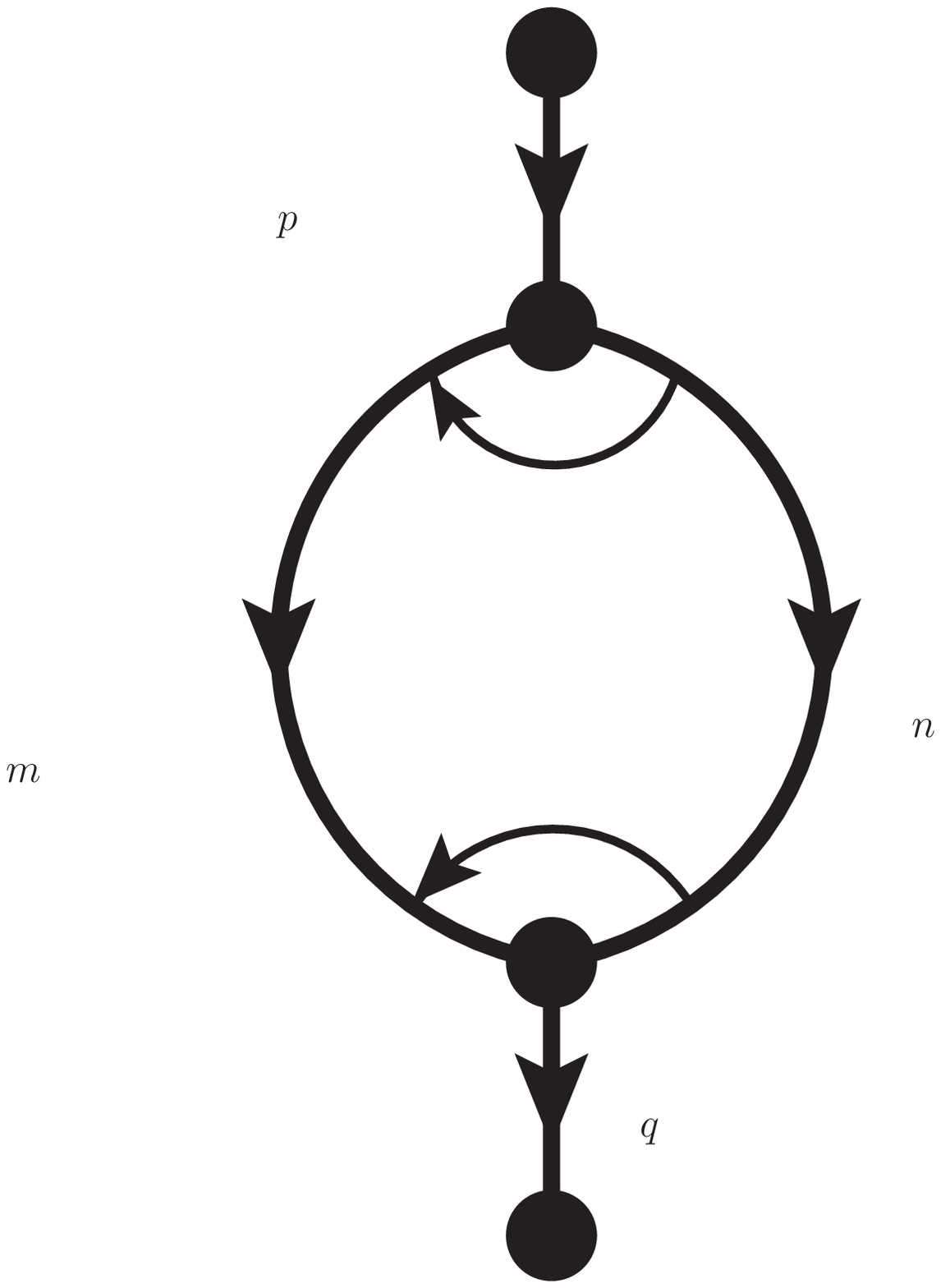}}
\>^{CG}=
\de_{pq}
\frac{k!\ (m+n-k+1)!\ (m-k)!\ (n-k)!}
{m!\ n!\ (m+n-2k+1)!}
\<
\parbox{0.6cm}{\psfrag{p}{$\scriptstyle{p}$}
\includegraphics[width=0.6cm]{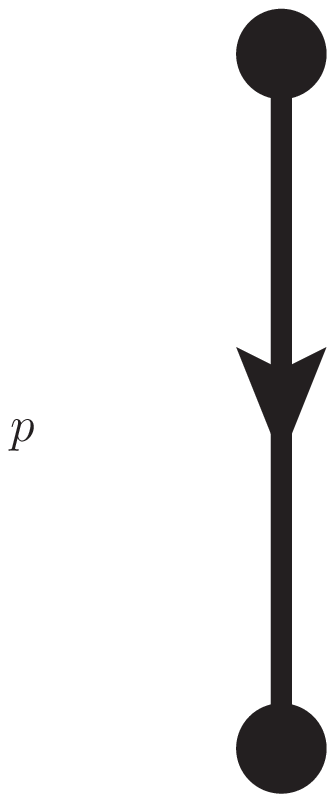}}\ \ 
\>^{CG}
\label{ClebschidentCG}
\end{equation}
and the Gordan series (\ref{Gordanseries}) is
\begin{equation}
\<
\parbox{1.3cm}{\psfrag{m}{$\scriptstyle{m}$}
\psfrag{n}{$\scriptstyle{n}$}
\includegraphics[width=1.3cm]{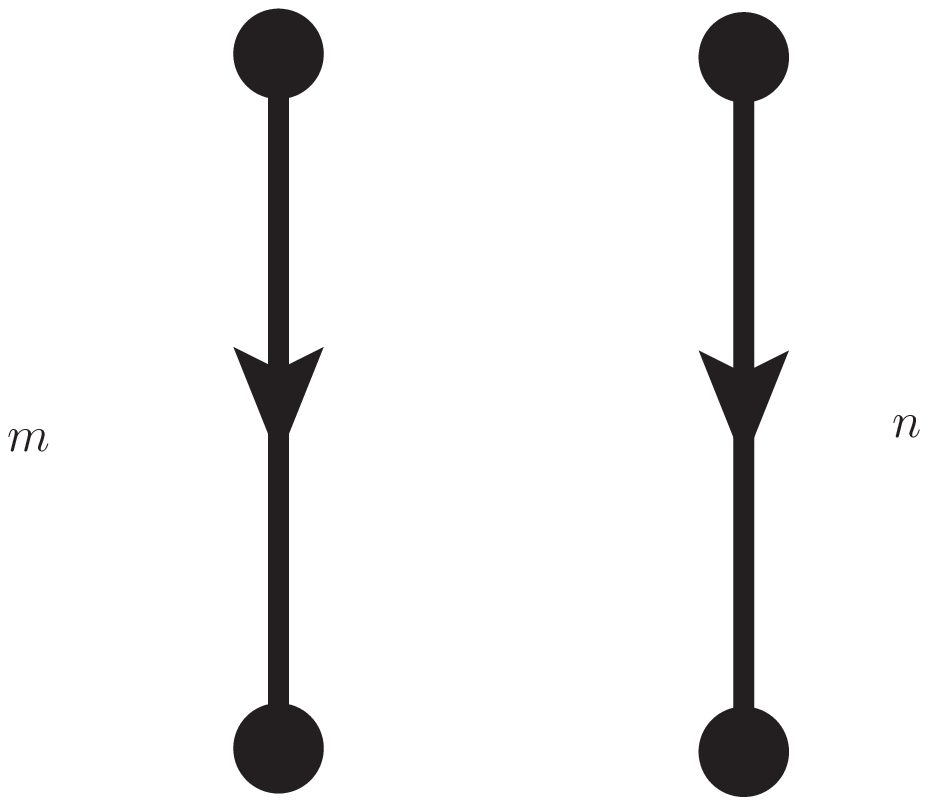}}
\>^{CG} =
\sum\limits_{k=0}^{\min(m,n)}
\frac{
\left(
\begin{array}{c}
m\\
k
\end{array}
\right)
\left(
\begin{array}{c}
n\\
k
\end{array}
\right)
}{
\left(
\begin{array}{c}
m+n-k+1\\
k
\end{array}
\right)
}
\<
\ \ 
\parbox{2cm}{\psfrag{m}{$\scriptstyle{m}$}
\psfrag{n}{$\scriptstyle{n}$}
\psfrag{j}{$\scriptstyle{m+n-2k}$}
\includegraphics[width=1.4cm]{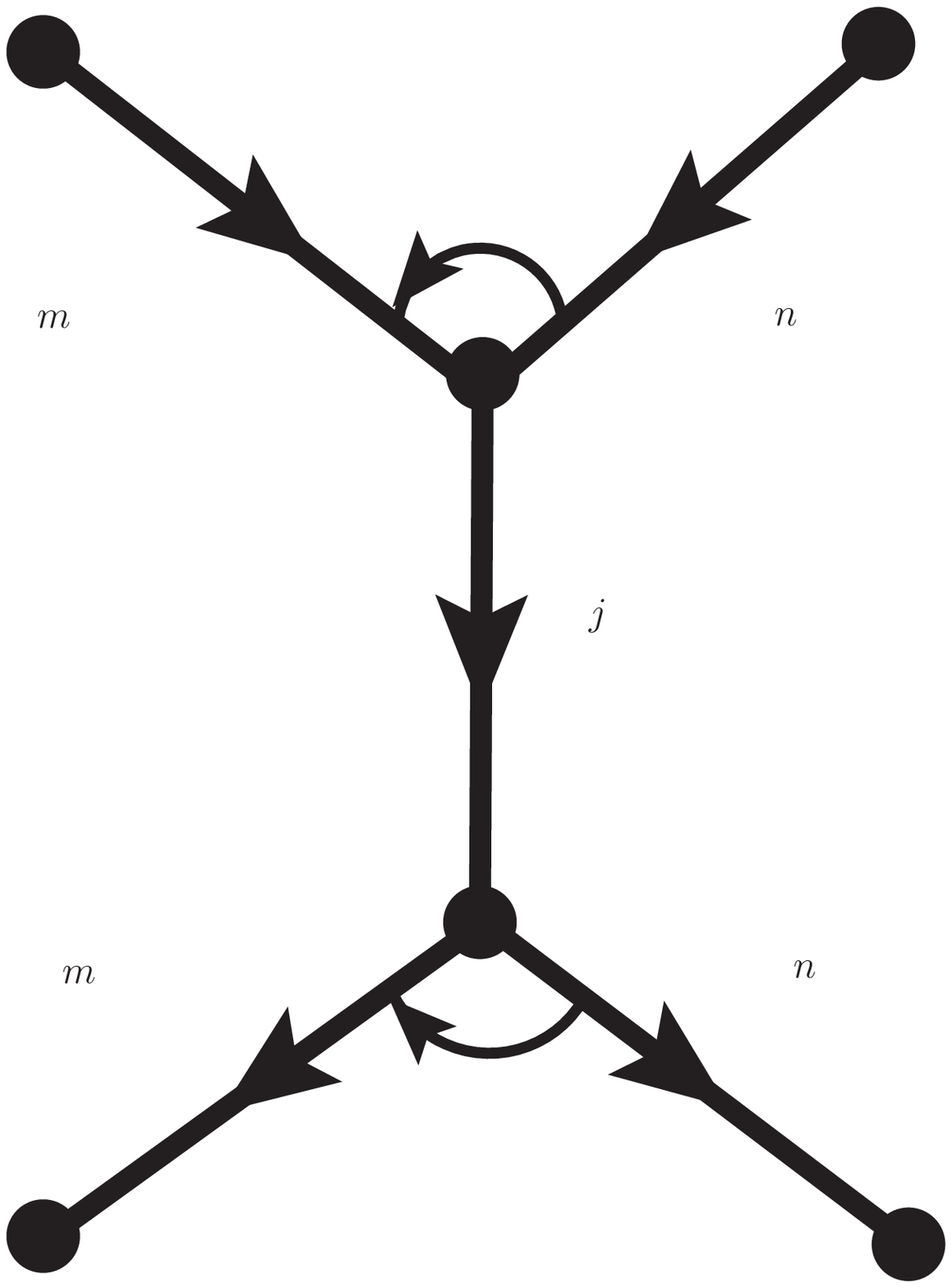}}
\>^{CG}\ .
\label{GordanseriesCG}
\end{equation}
Classical transvectants as in (\ref{transvectdef})
are obtained by plugging the `holes' corrsponding to the $\cH_{\ga(v)}^{\ast}$'s
by blobs 
$\parbox{0.8cm}{\psfrag{F}{$\scriptstyle{F}$}
\psfrag{g}{$\scriptstyle{\underbrace{\ \ \ }_{\ga(v)}}$}
\includegraphics[width=0.8cm]{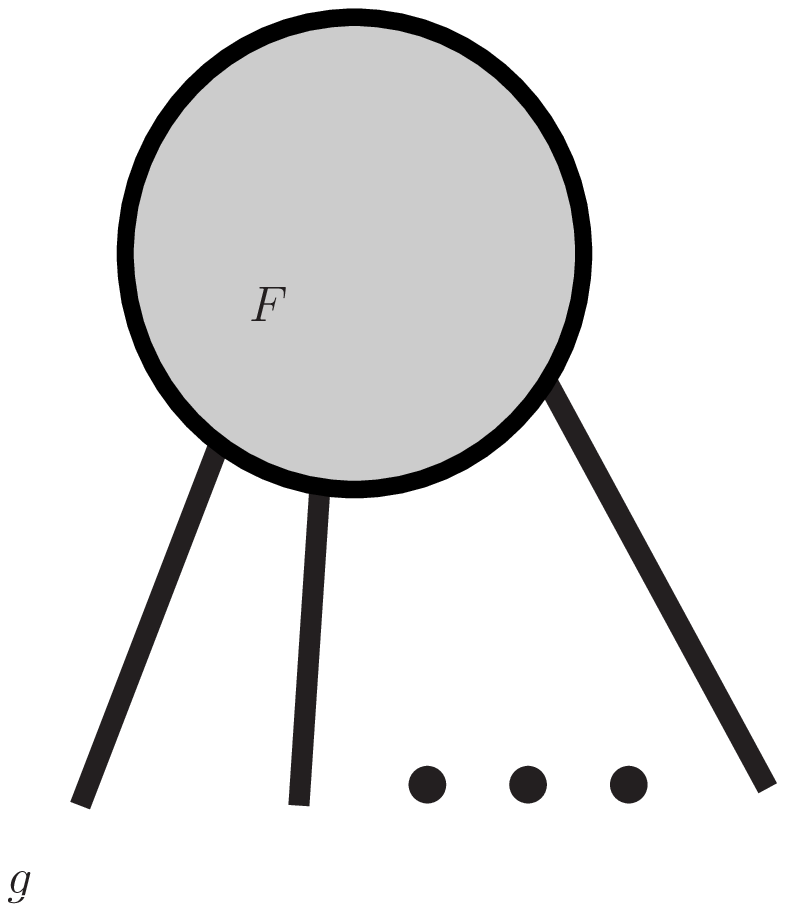}}$
of binary forms, and those corresponding to 

\medskip\noindent
$\cH_{\ga(v)}$'s
by $\parbox{1cm}{\psfrag{x}{$\scriptstyle{x}$}
\psfrag{g}{$\scriptstyle{\underbrace{\ \ \ }_{\ga(v)}}$}
\includegraphics[width=1cm]{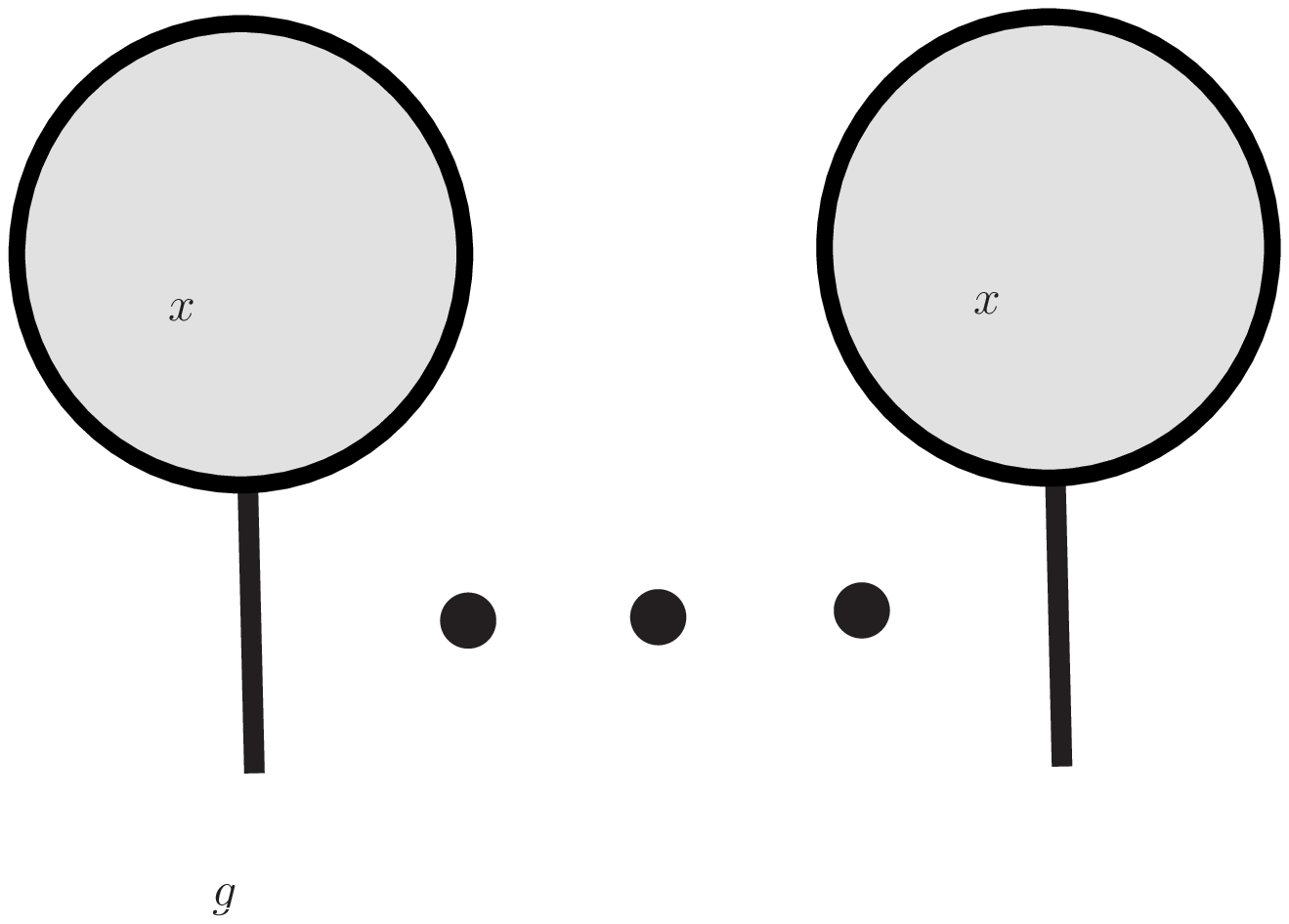}}$
or by other series of binary variables $\uy,\uz,\ldots$
\[
\ 
\]

There are many natural operations one can do with CG networks
$(G,\cO,\ta,\ga)$.
One can cut an edge
\begin{equation}
\parbox{1.2cm}{\psfrag{a}{$\scriptstyle{a}$}
\includegraphics[width=1.2cm]{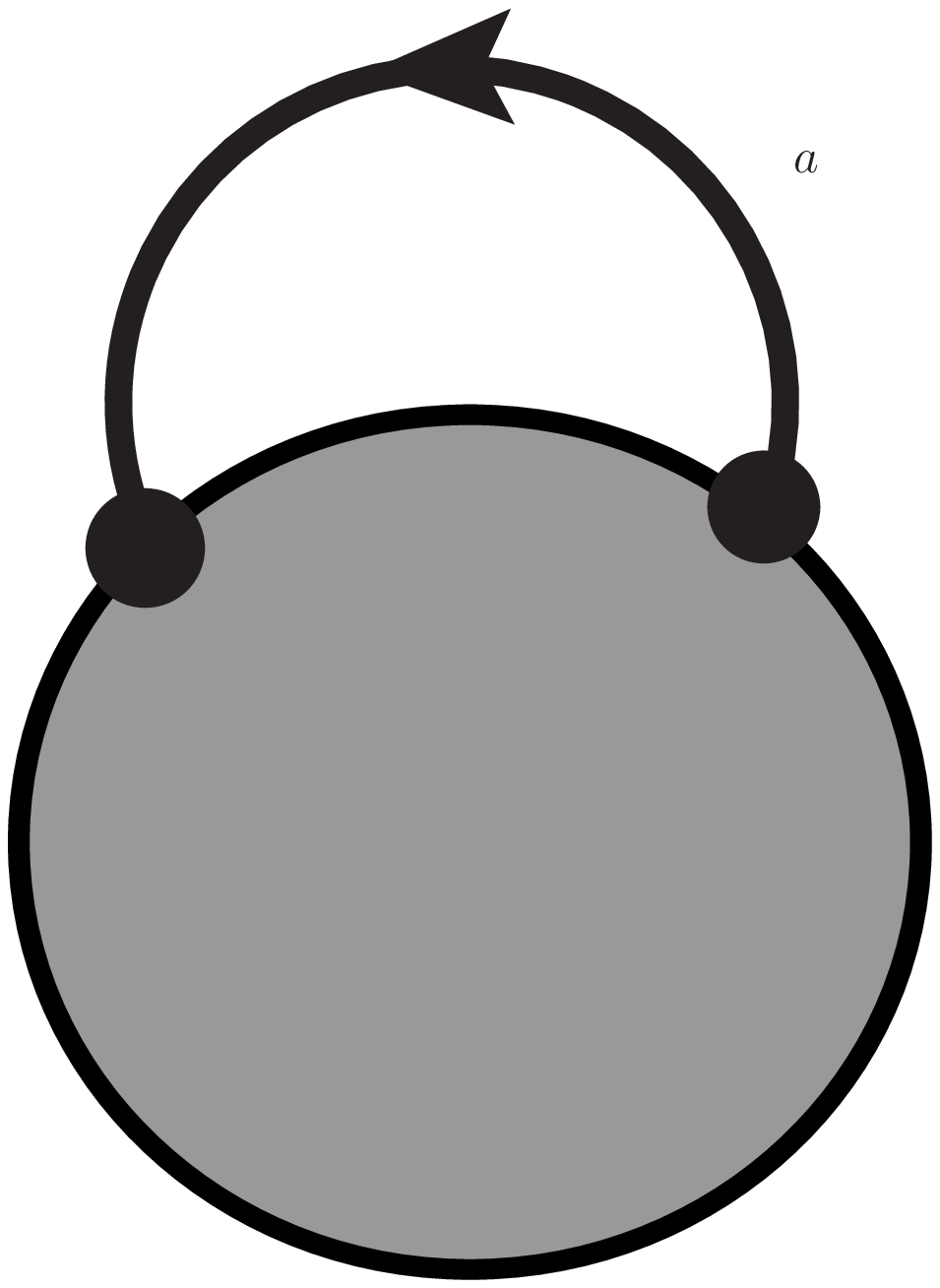}}
\qquad\longrightarrow\qquad
\parbox{1.4cm}{\psfrag{a}{$\scriptstyle{a}$}
\includegraphics[width=1.4cm]{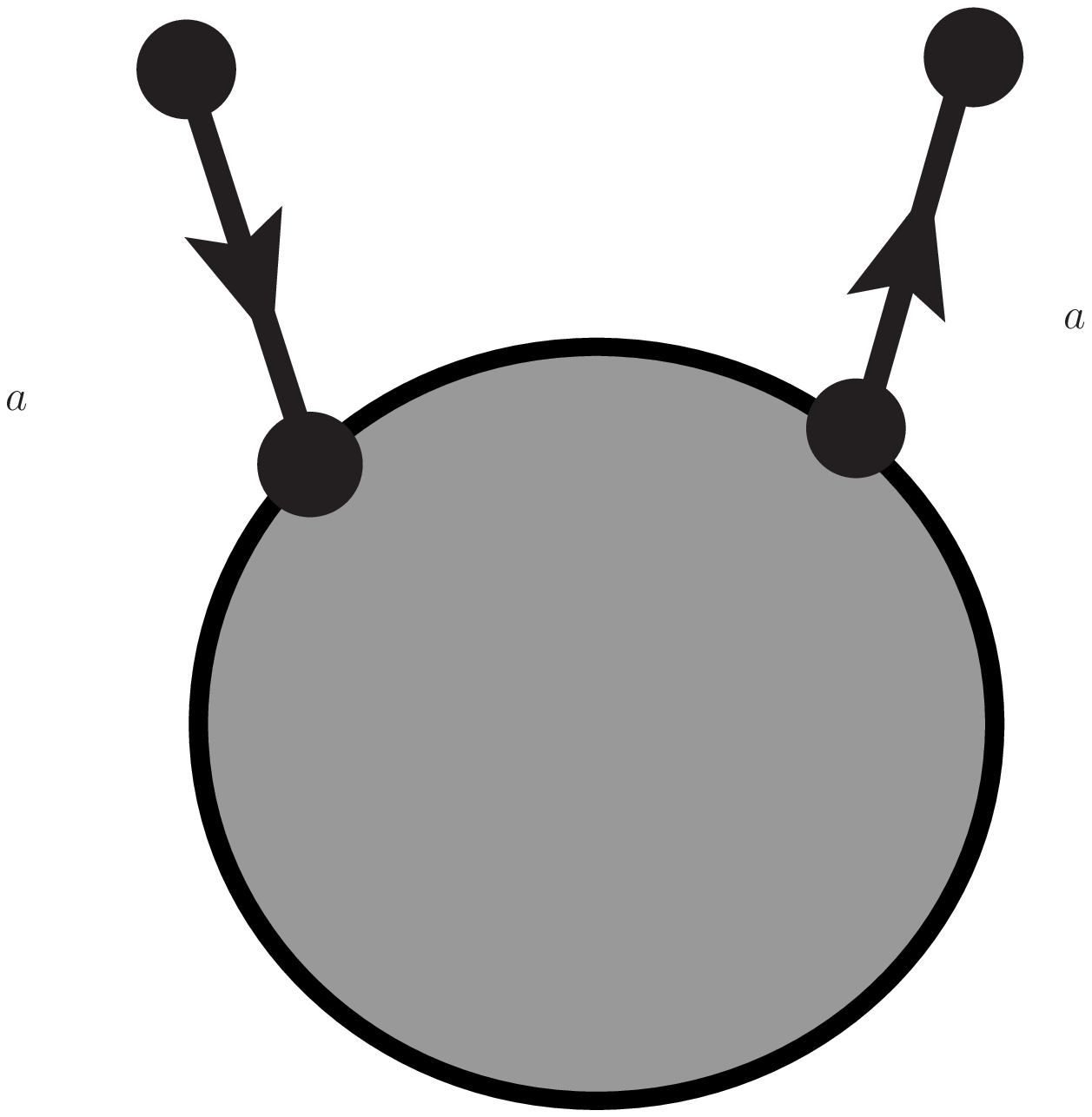}}
\label{cuttingedge}
\end{equation}
and create two new 1-valent vertices.
The shaded areas represent the rest of the graph
which does not necessarily have to be connected.
Conversely, for a pair $v_{\rm in}\in V_{\rm in}(G)$,
$v_{\rm out}\in V_{\rm out}(G)$
with $\ga(v_{\rm in})=\ga(v_{\rm out})$
one can glue them
\[
\parbox{1.4cm}{\psfrag{a}{$\scriptstyle{a}$}
\psfrag{i}{$\scriptstyle{v_{\rm in}}$}
\psfrag{o}{$\scriptstyle{v_{\rm out}}$}
\includegraphics[width=1.4cm]{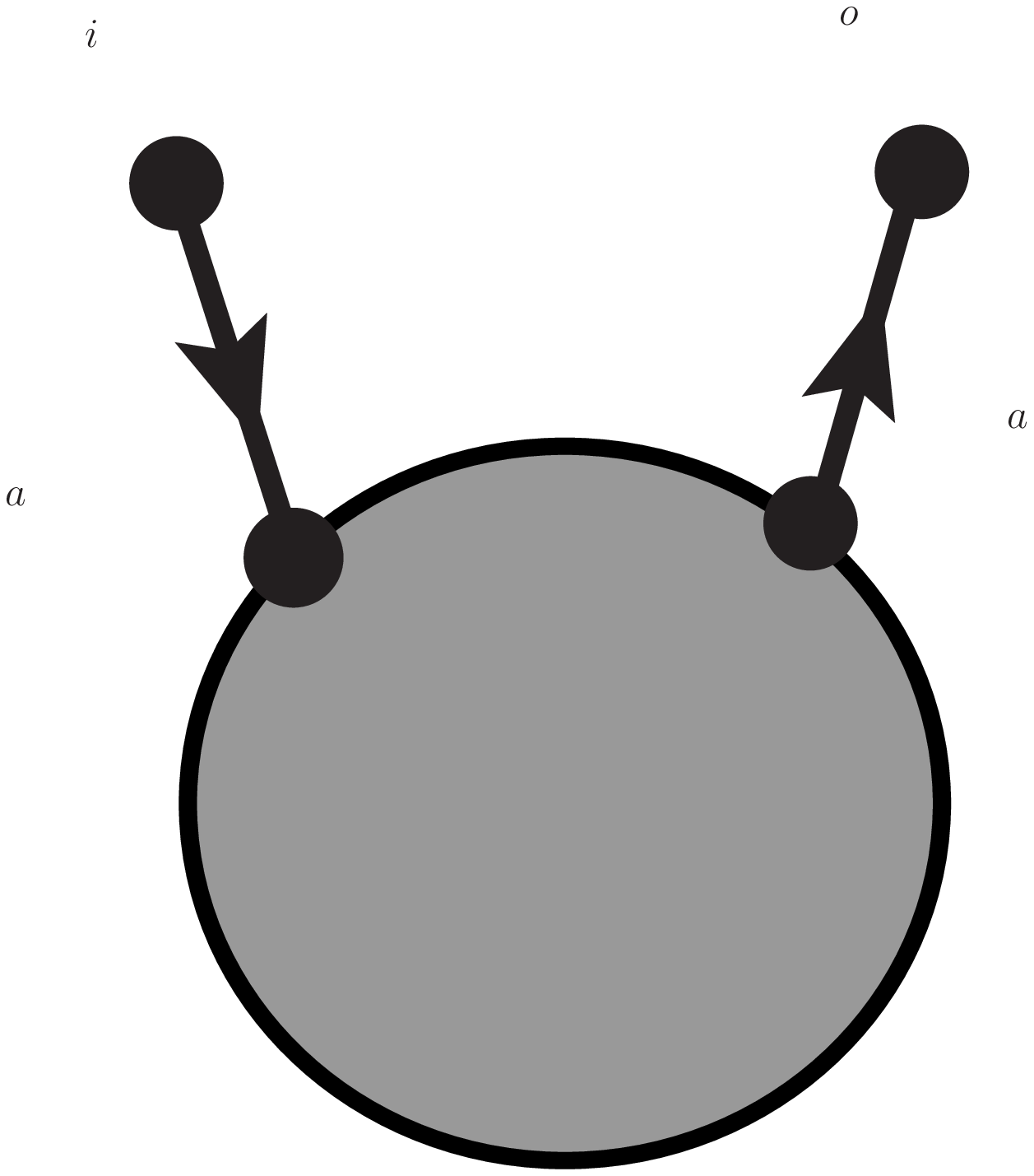}}
\qquad\longrightarrow\qquad
\parbox{1.2cm}{\psfrag{a}{$\scriptstyle{a}$}
\includegraphics[width=1.2cm]{Fig116.eps}}\ \ .
\]

Now a matrix $g\in GL_2(\C)$ naturally acts on the tensor
$\<G,\cO,\ta,\ga\>^{CG}$ following (\ref{xtransf}) and (\ref{Ftransf}).
The resulting tensor $g\cdot\<G,\cO,\ta,\ga\>^{CG}$ 
can be obtained using FDC as before, except that an entry leg contributes
\[
\parbox{1.2cm}{\psfrag{a}{$\scriptstyle{a}$}
\includegraphics[width=1.2cm]{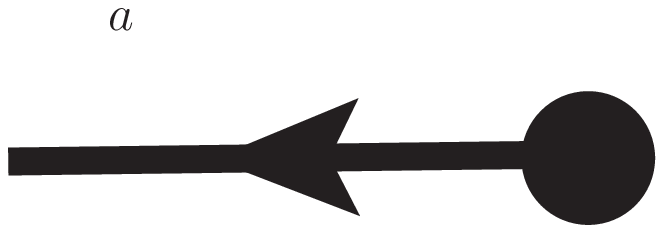}}
\qquad\longrightarrow\qquad
\parbox{2.6cm}{\psfrag{a}{$\scriptstyle{i_a}$}
\psfrag{1}{$\scriptstyle{i_1}$}\psfrag{g}{$\scriptstyle{g}$}
\includegraphics[width=2.6cm]{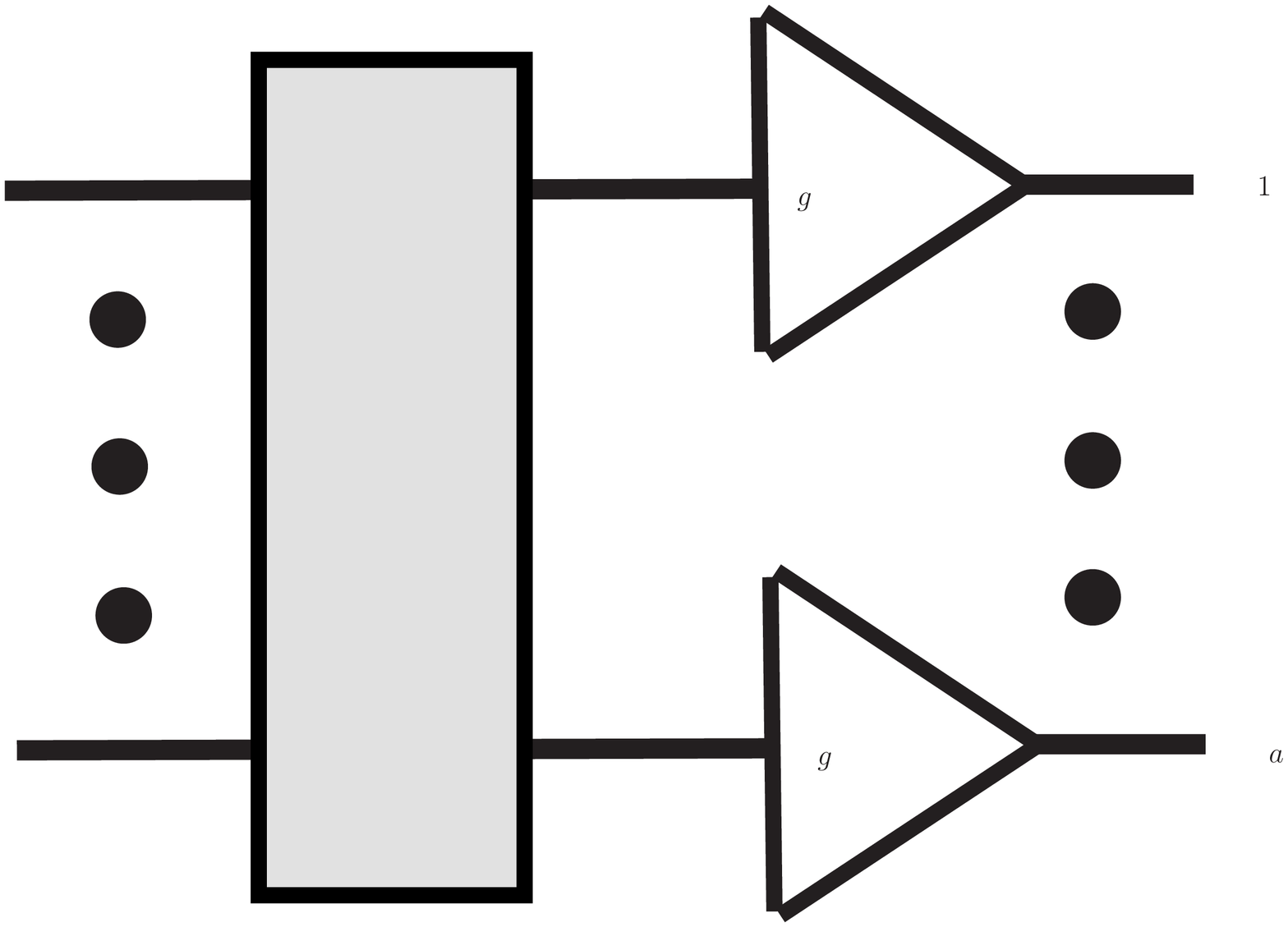}}
\]
and for an exit leg
\[
\parbox{1.2cm}{\psfrag{a}{$\scriptstyle{a}$}
\includegraphics[width=1.2cm]{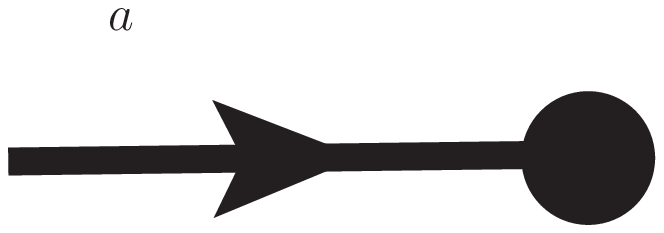}}
\qquad\longrightarrow\qquad
\parbox{3.5cm}{\psfrag{a}{$\scriptstyle{i_a}$}
\psfrag{1}{$\scriptstyle{i_1}$}\psfrag{g}{$\scriptscriptstyle{g^{-1}}$}
\includegraphics[width=3.5cm]{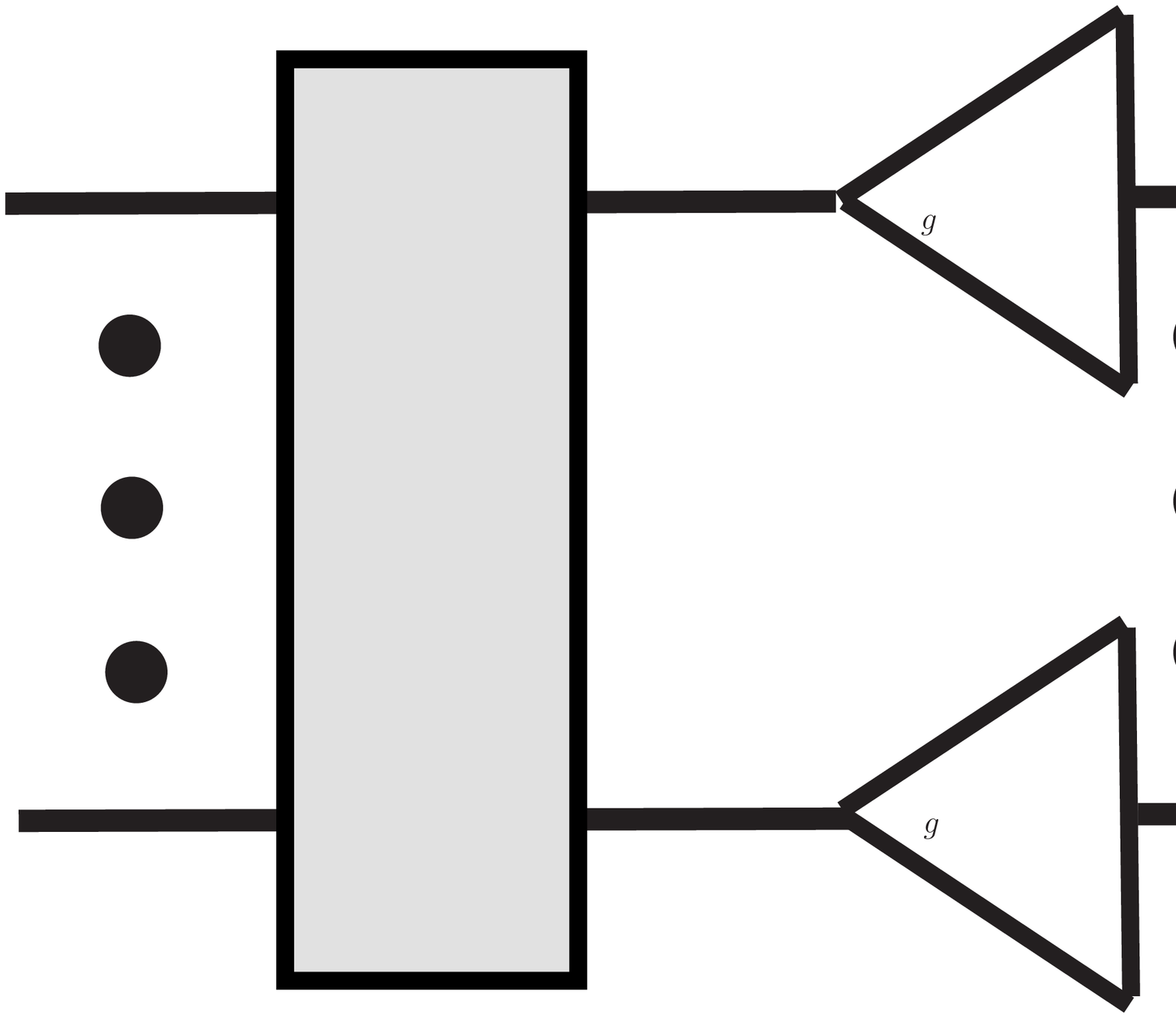}}\ \ .
\]
An important property of CG networks is $SL_2$
invariance.
\begin{Proposition}
For any $g\in SL_2(\C)$, $g\cdot\<G,\cO,\ta,\ga\>^{CG}=\<G,\cO,\ta,\ga\>^{CG}$.
\end{Proposition}
\noindent{\bf Proof:}
Consider the microscopic Feynman diagram for $g\cdot\<G,\cO,\ta,\ga\>^{CG}$.
For any edge between 3-valent vertices (eventually the same in case of
a loop edge) insert the identity
\[
\parbox{0.8cm}{\psfrag{a}{$\scriptstyle{a}$}
\includegraphics[width=0.8cm]{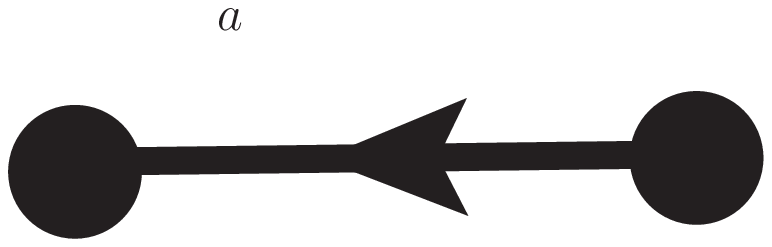}}
\qquad\longrightarrow\qquad
\left.
\parbox{1.8cm}{\includegraphics[width=1.8cm]{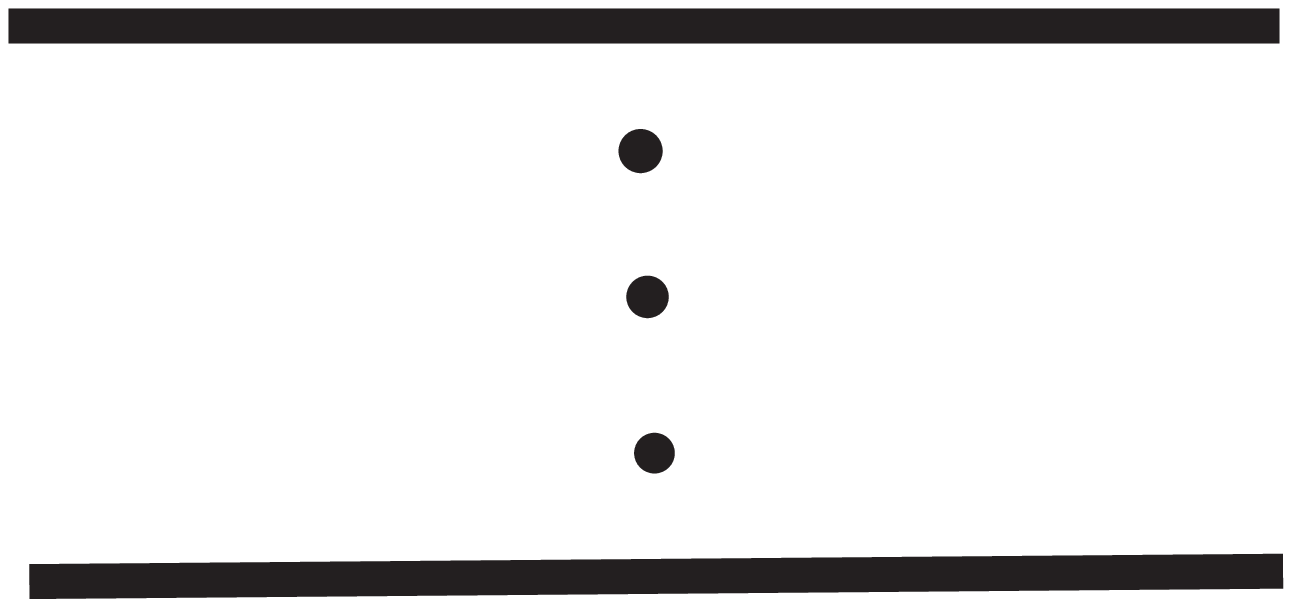}}
\right\} {\scriptstyle a}
=
\left.
\parbox{3cm}{\psfrag{1}{$\scriptstyle{g}$}
\psfrag{2}{$\scriptscriptstyle{g^{-1}}$}
\includegraphics[width=3cm]{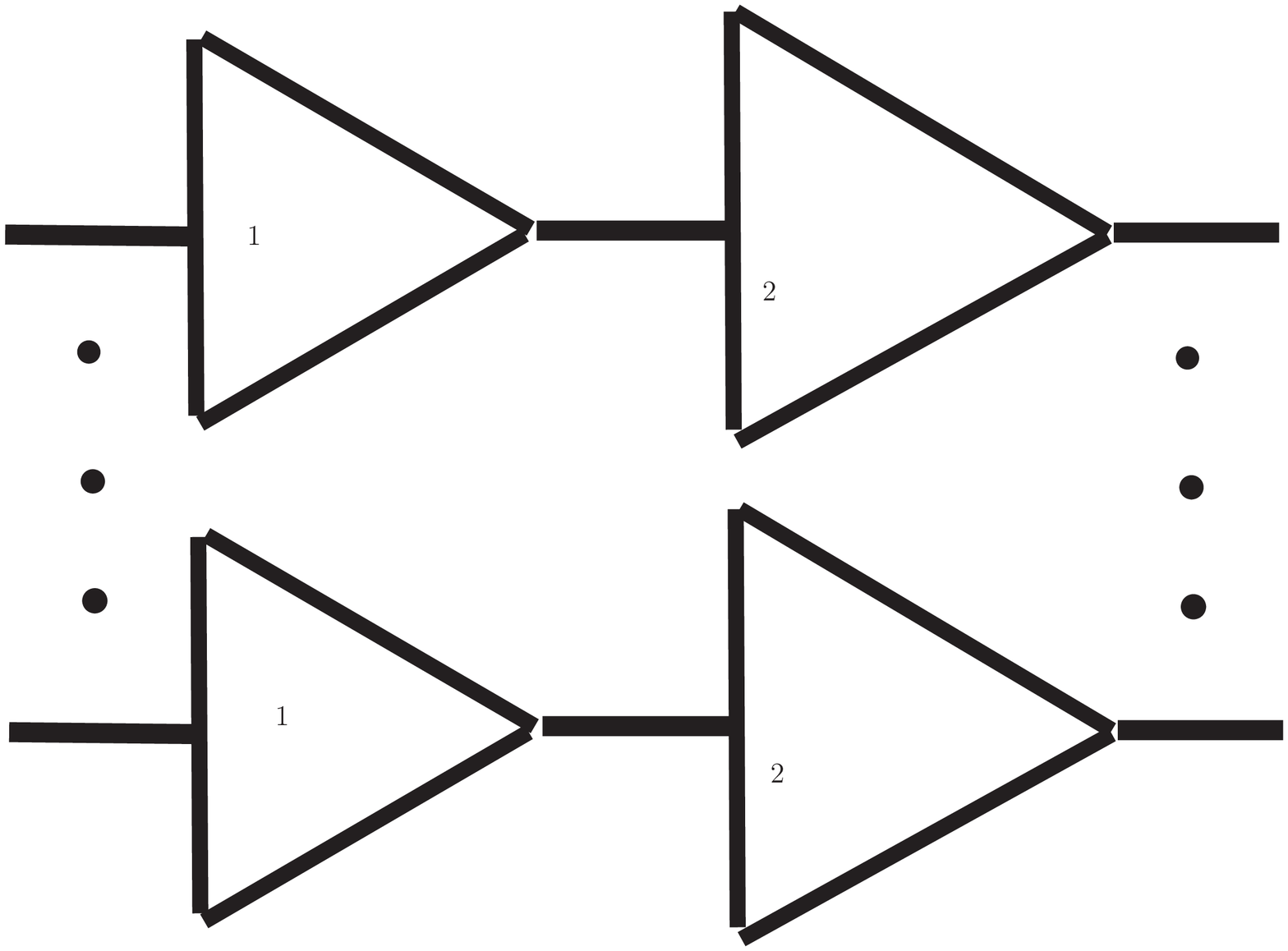}}
\right\} {\scriptstyle a}
\]
with due care for the indicated orientation.
Now the vicinity of every 3-valent vertex becomes
\[
\begin{array}{c}
\parbox{1.6cm}{\psfrag{a}{$\scriptstyle{a}$}
\psfrag{b}{$\scriptstyle{b}$}\psfrag{c}{$\scriptstyle{c}$}
\includegraphics[width=1.6cm]{Fig103.eps}}\\
\ \\
{\rm macroscopic}
\end{array}
\qquad\longrightarrow\qquad
\begin{array}{c}
\parbox{5cm}{\psfrag{a}{$\scriptstyle{a}$}
\psfrag{b}{$\scriptstyle{b}$}\psfrag{c}{$\scriptstyle{c}$}
\psfrag{g}{$\scriptstyle{g}$}\psfrag{h}{$\scriptscriptstyle{g^{-1}}$}
\includegraphics[width=5cm]{Fig126.eps}}\\
\ \\
{\rm microscopic}
\end{array}
\]
or
\[
\parbox{1.6cm}{\psfrag{a}{$\scriptstyle{a}$}
\psfrag{b}{$\scriptstyle{b}$}\psfrag{c}{$\scriptstyle{c}$}
\includegraphics[width=1.6cm]{Fig102.eps}}
\qquad\longrightarrow\qquad
\parbox{5cm}{\psfrag{a}{$\scriptstyle{a}$}
\psfrag{b}{$\scriptstyle{b}$}\psfrag{c}{$\scriptstyle{c}$}
\psfrag{h}{$\scriptstyle{g}$}\psfrag{g}{$\scriptscriptstyle{g^{-1}}$}
\includegraphics[width=5cm]{Fig127.eps}}\ \ .
\]

\medskip
One can easily see that one can push the matrices through the symmetrizers
\begin{equation}
\parbox{2.5cm}{\psfrag{g}{$\scriptstyle{M}$}
\includegraphics[width=2.5cm]{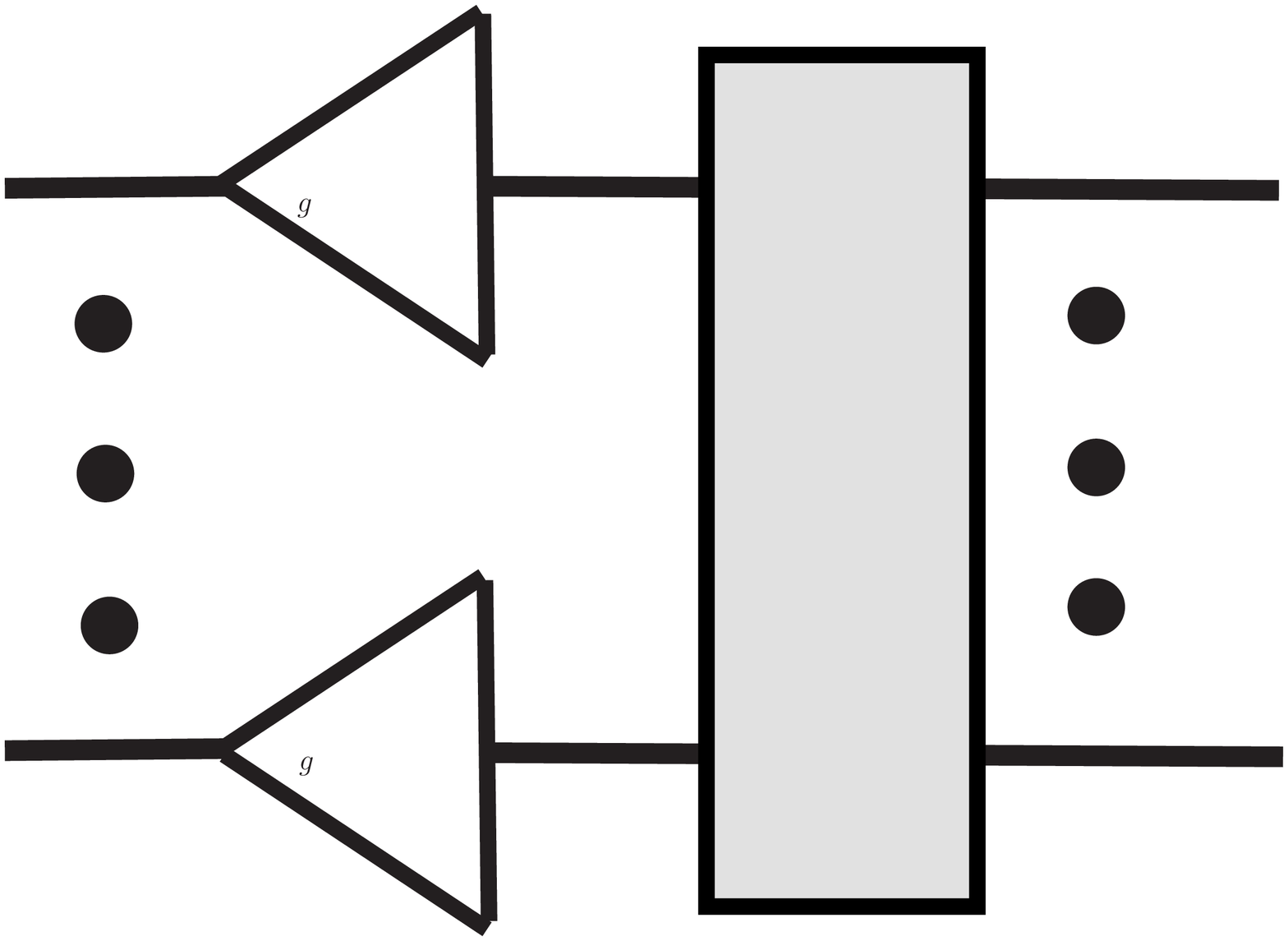}}
=
\parbox{2.5cm}{\psfrag{g}{$\scriptstyle{M}$}
\includegraphics[width=2.5cm]{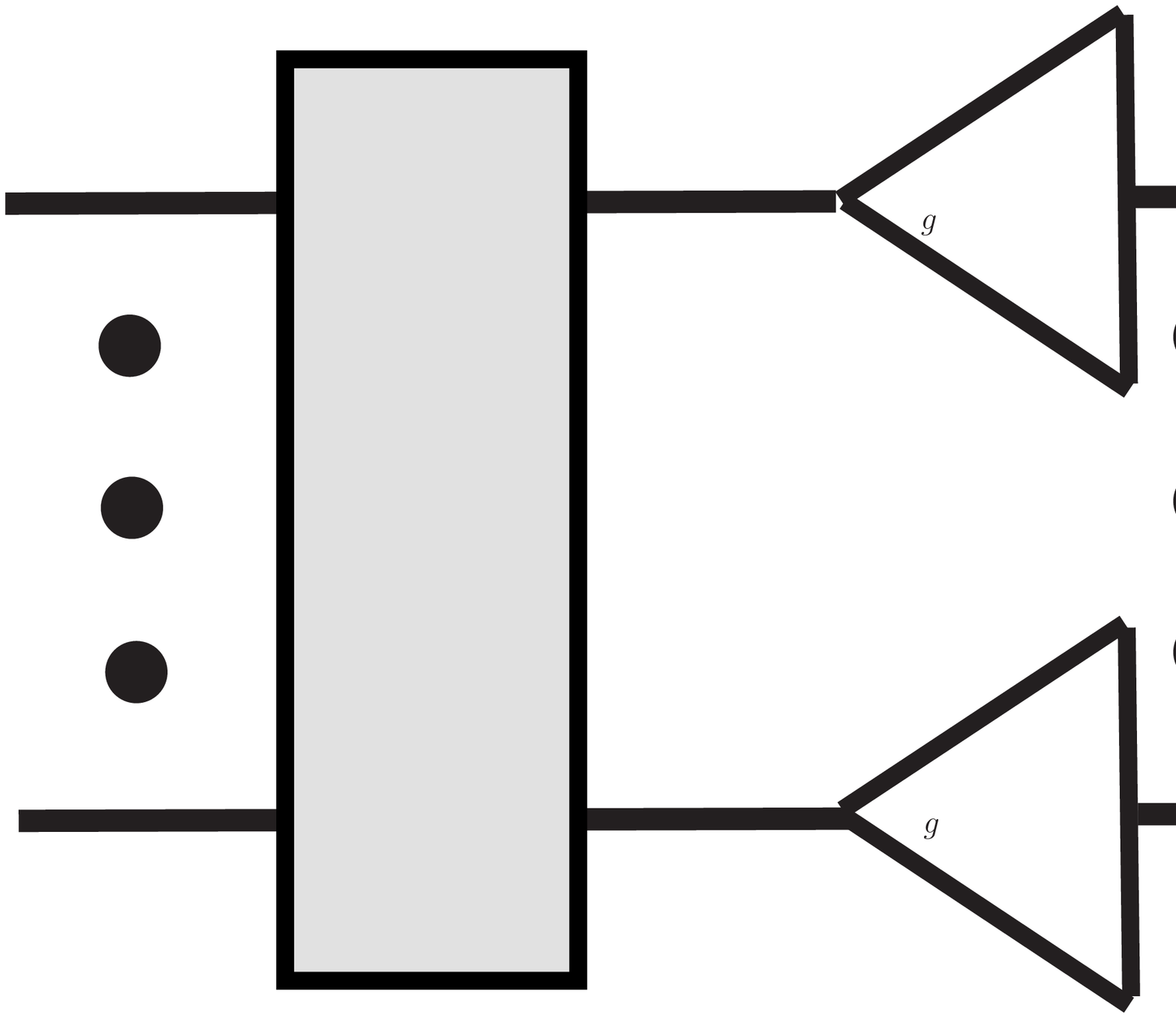}}\ \ .
\label{matrixpush}
\end{equation}
Therefore, using (\ref{epsilonident}), ${\rm det}(g)=1$ and
$g\cdot g^{-1}={\rm Id}$
one obtains the original blown-up vertex from Rule 2) in Def. \ref{CGevaldef}.
\qed

Now identity (\ref{ClebschidentCG})
simply is an explicit instance of Schur's Lemma.

\begin{Proposition}
More generally for CG networks with one entry and one exit leg
we have
\begin{equation}
\<
\parbox{2.5cm}{\psfrag{a}{$\scriptstyle{a}$}
\psfrag{b}{$\scriptstyle{b}$}
\includegraphics[width=2.5cm]{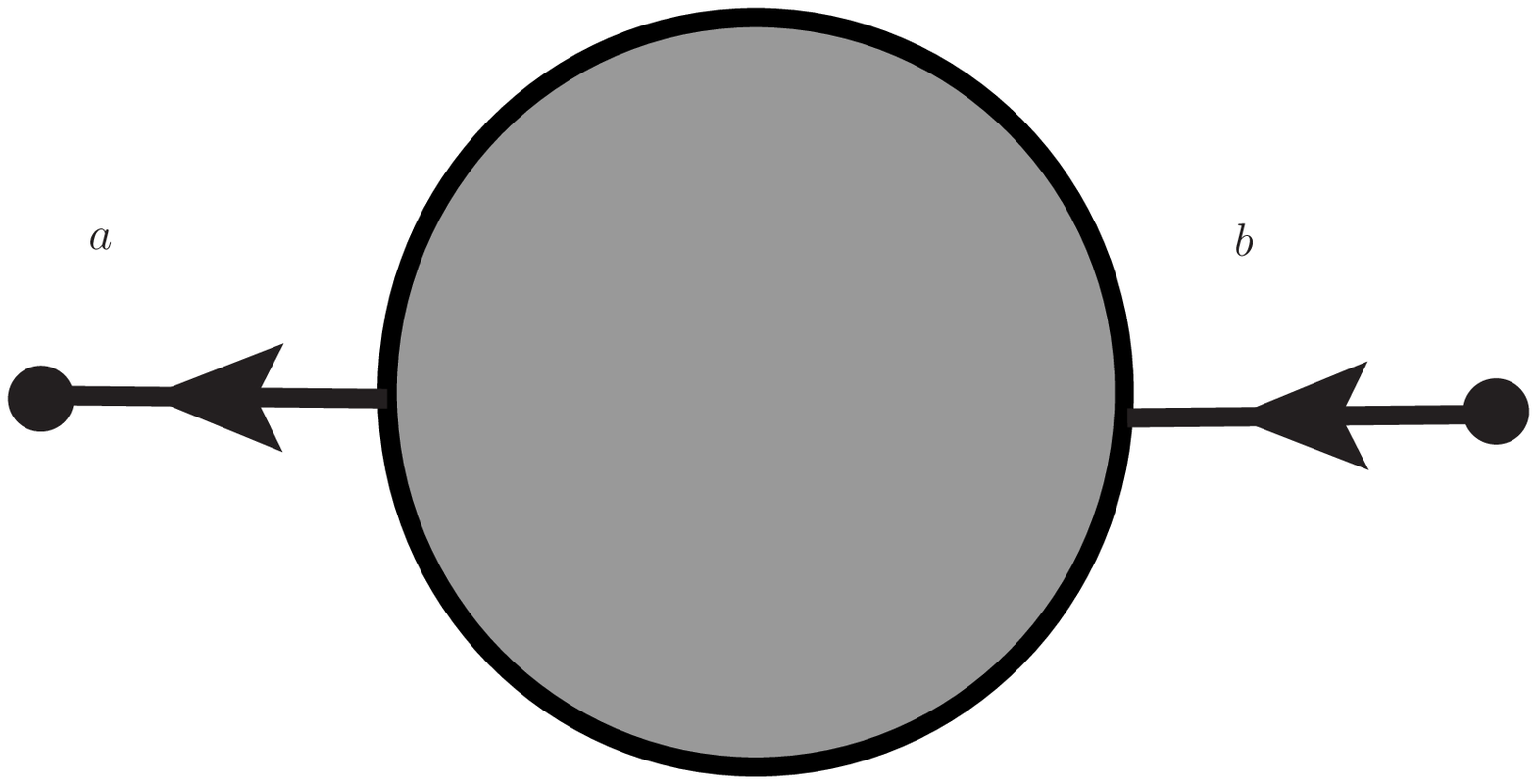}}
\>^{CG}=\frac{\de_{a,b}}{a+1}
\<
\parbox{1.4cm}{\psfrag{a}{$\scriptstyle{a}$}
\includegraphics[width=1.4cm]{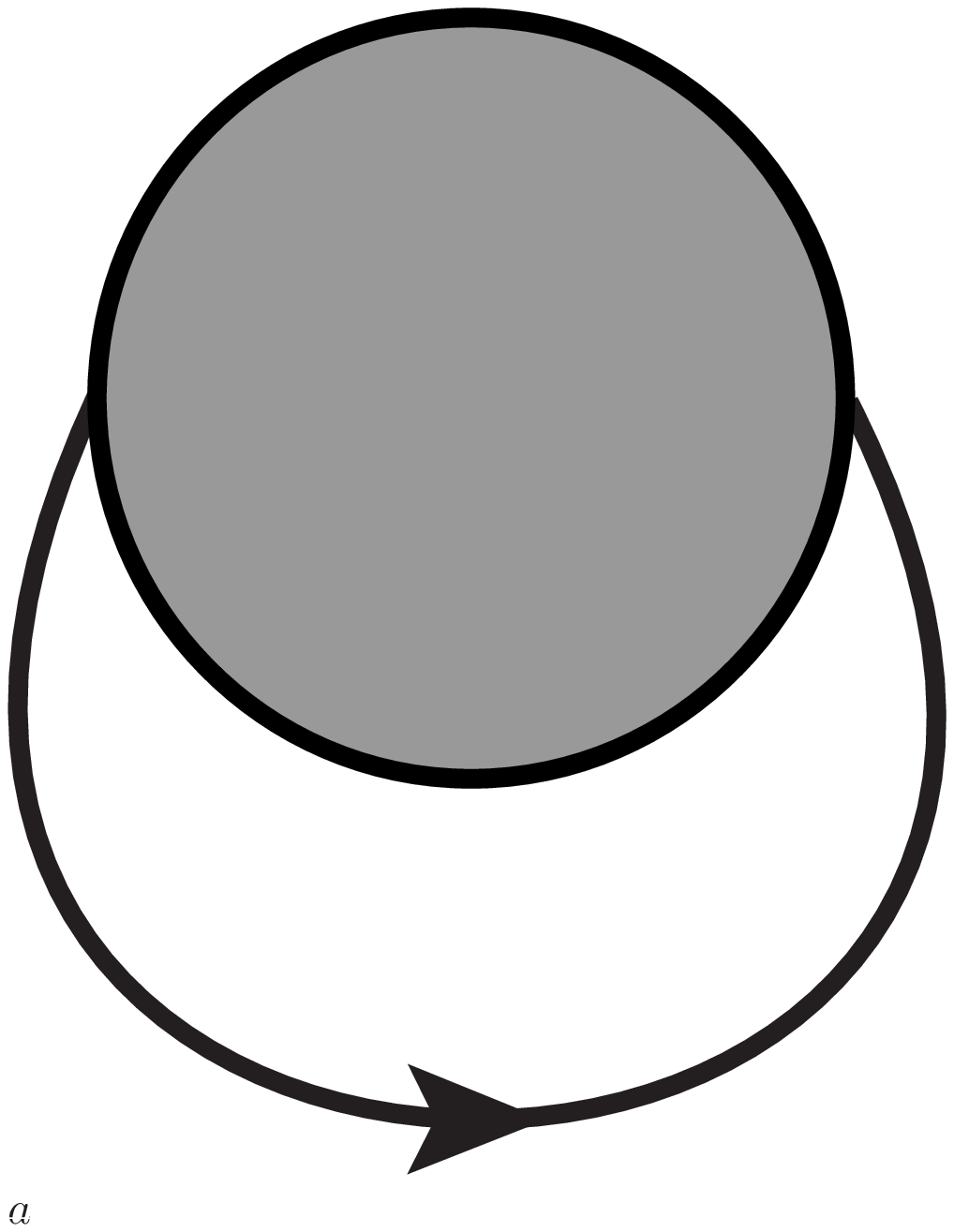}}
\>^{CG} \<
\parbox{1.4cm}{\psfrag{a}{$\scriptstyle{a}$}
\includegraphics[width=1.4cm]{Fig123.eps}}
\>^{CG}\ .
\label{Schur}
\end{equation}
\end{Proposition}
\noindent{\bf Proof:}
Schur's Lemma tells us that the left hand side vanishes unless
$b=a$, i.e., the two irreducible representations are the same.
Besides, if $a=b$, then the left-hand side is a multiple of the identity
with a proportionality constant which can be determined as the ratio
of the traces of these two operators. \qed

In general, a tensor $\<G,\cO,\ta,\ga\>^{CG}$
can be seen as an $SL_2$-equivariant
map from $\bigotimes\limits_{v\in V_{\rm in}(G)}\cH_{\ga(v)}$ to
$\bigotimes\limits_{v\in V_{\rm out}(G)}\cH_{\ga(v)}$.
For further use, we state the following trivial consequence.
\begin{Corollary}\label{oneleg}
If $(G,\cO,\ta,\ga)$ is a CG network with only one external leg $v$,
its evaluation $\<G,\cO,\ta,\ga\>^{CG}$ vanishes unless $\ga(v)=0$.
\end{Corollary}
Classically this amounts to the fact that there are no
nonzero linear invariants for a binary form. This follows from the
First Fundamental Theorem stated earlier.

An important remark is that the Feynman diagram piece used to compute
the contribution of a 3-valent vertex by Rule 2) of Def. \ref{CGevaldef},
is the same regardless of the two possible orientations of the edges
incident to the vertex.
Besides, this contribution is real for any values of the $a+b+c$
indices.
As a result,
the two equivariant maps
\[
\Pi_{CG}:\cH_m\otimes\cH_n\rightarrow \cH_{m+n-2k}\qquad
{\rm given\ by}\qquad
\<
\ \ 
\parbox{1.6cm}{\psfrag{a}{$\scriptstyle{m}$}
\psfrag{b}{$\scriptstyle{n}$}
\psfrag{c}{$\scriptstyle{m+n-2k}$}
\includegraphics[width=1.6cm]{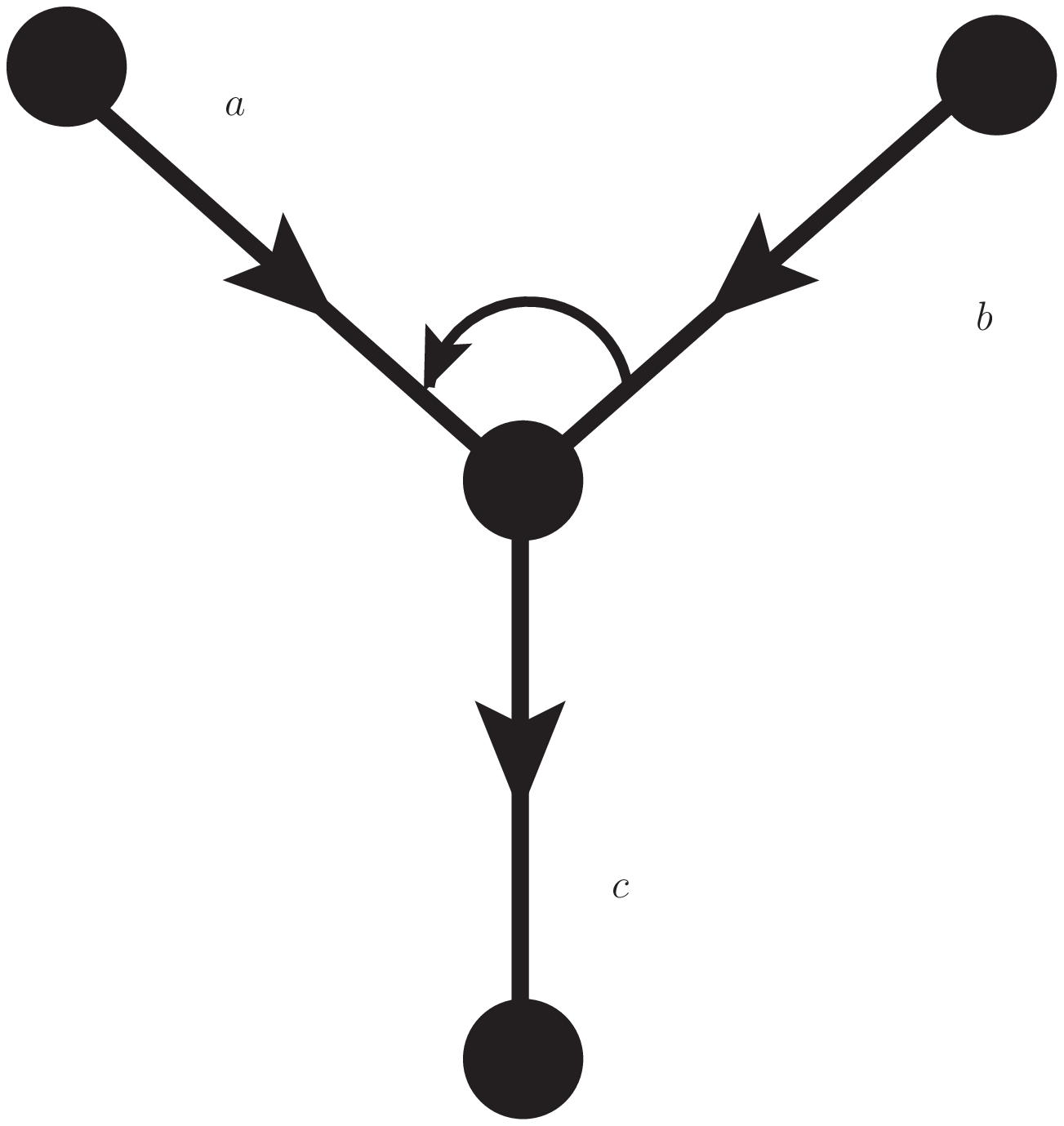}}
\qquad\>^{CG}
\]
and
\[
I_{CG}:\cH_{m+n-2k}\rightarrow \cH_m\otimes\cH_n\qquad
{\rm given\ by}\qquad
\<
\ \ 
\parbox{1.6cm}{\psfrag{a}{$\scriptstyle{m}$}
\psfrag{b}{$\scriptstyle{n}$}
\psfrag{c}{$\scriptstyle{m+n-2k}$}
\includegraphics[width=1.6cm]{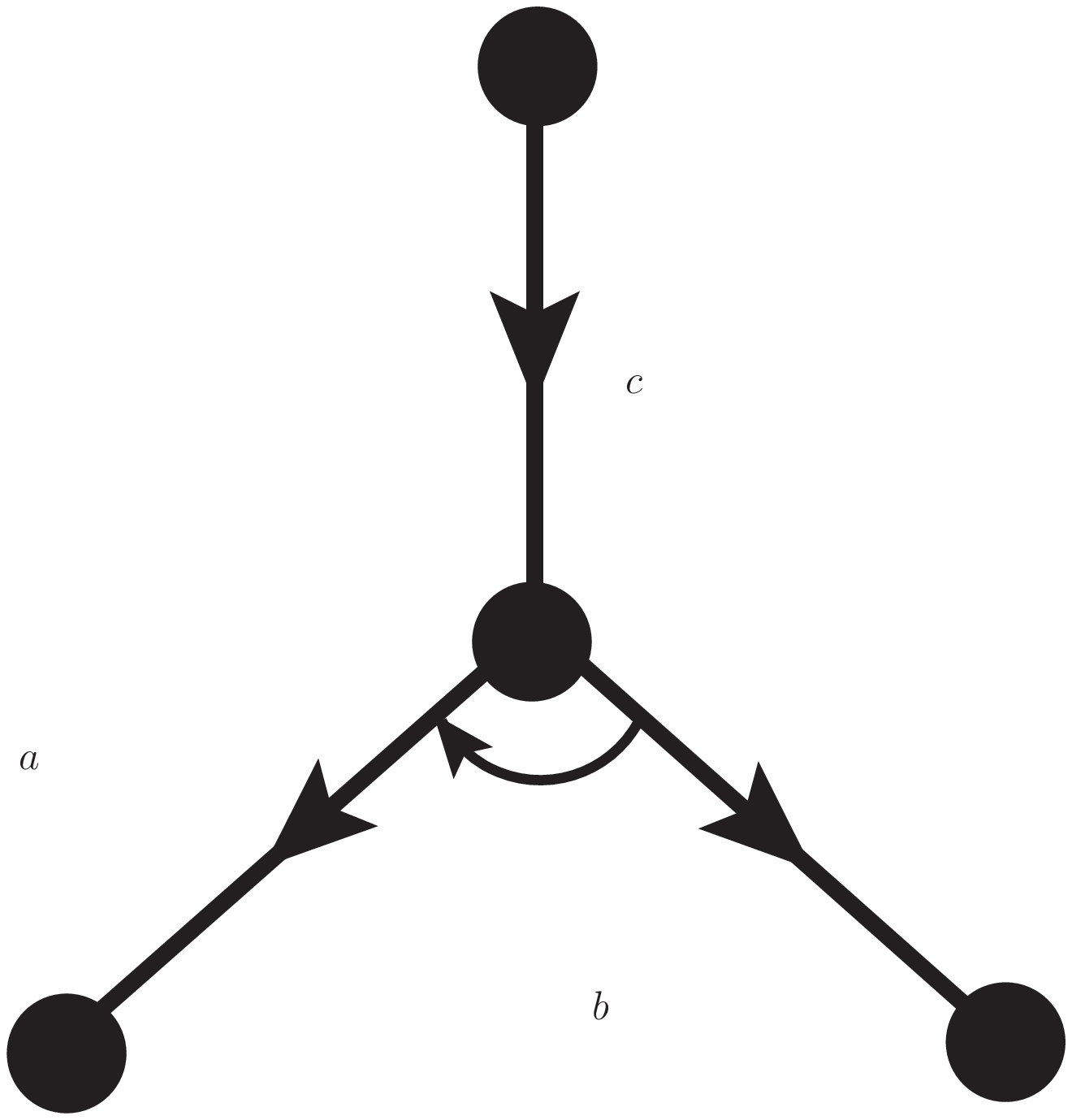}}
\qquad
\>^{CG}
\]
are Hermitian conjugates.
Indeed, for a binary form $F\in\cH_{m+n-2k}$
and a bihomogeneous form $B\in\cH_m\otimes\cH_n$ we have
\[
\langle F|\Pi_{CG}(B)\rangle=
\parbox{4.6cm}{\psfrag{m}{$\scriptstyle{m}$}
\psfrag{n}{$\scriptstyle{n}$}
\psfrag{k}{$\scriptstyle{m+n-2k}$}
\psfrag{B}{$B$}
\psfrag{F}{$\overline{F}$}
\psfrag{i}{$\begin{array}{c}
{\scriptstyle\rm inner}\\
{\scriptstyle\rm product}\\
{\scriptstyle\rm split}
\end{array}$}
\includegraphics[width=4.6cm]{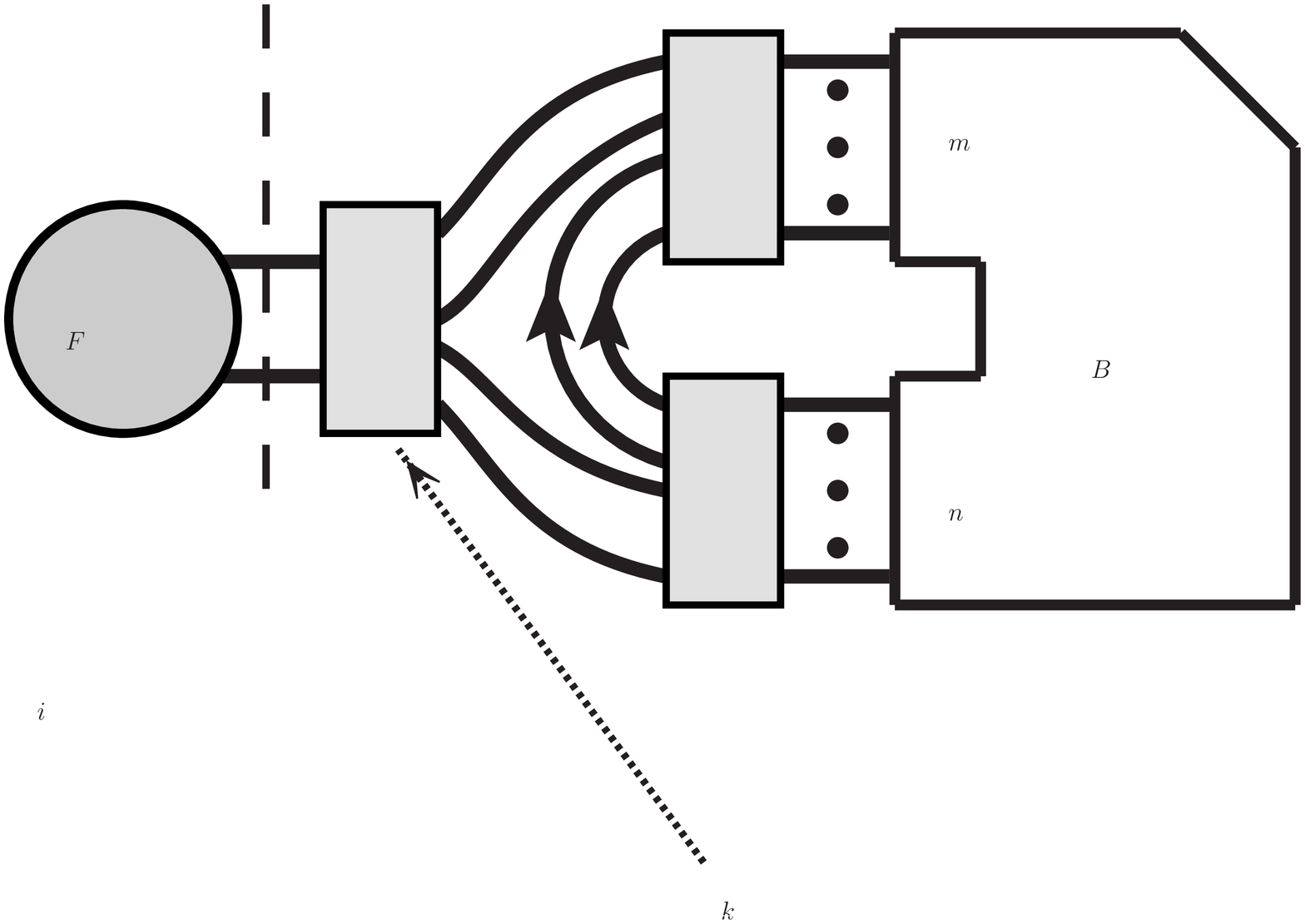}}
\]
\[
\ 
\]
\begin{equation}
=\overline{\left(
\parbox{4.6cm}{\psfrag{m}{$\scriptstyle{m}$}
\psfrag{n}{$\scriptstyle{n}$}
\psfrag{k}{$\scriptstyle{m+n-2k}$}
\psfrag{F}{$F$}
\psfrag{B}{$\overline{B}$}
\psfrag{i}{$\begin{array}{c}
{\scriptstyle\rm inner}\\
{\scriptstyle\rm product}\\
{\scriptstyle\rm split}
\end{array}$}
\includegraphics[width=4.6cm]{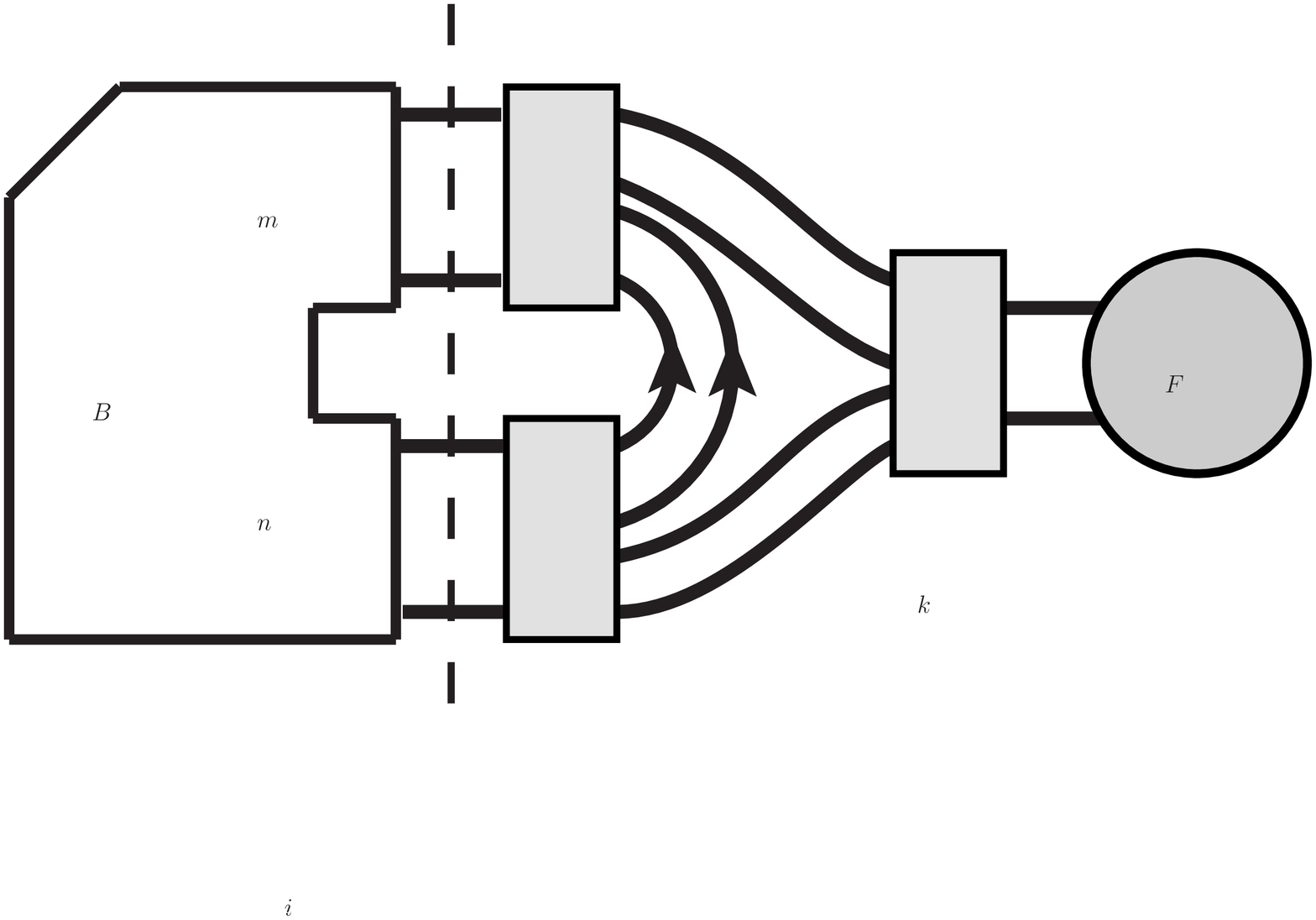}}
\right)}
\label{innerprodB}
\end{equation}
\[
\begin{array}{c}
\ \\
\ 
\end{array}
\]
\[
=\overline{\langle B|I_{CG}(F)\rangle}=\langle I_{CG}(F)| B\rangle\ .
\]
Using identities (\ref{ClebschidentCG}) and
(\ref{GordanseriesCG}) and similar graphical
computations of inner products, one can easily see that the map
$\io_k:\cH_{m+n-2k}\rightarrow\cH_m\otimes\cH_n$
defined by
\[
\io_k=\sqrt{
\frac{m!\ n!\ (m+n-2k+1)!}{k!\ (m+n-k+1)!\ (m-k)!\ (n-k)!}
}\times I_{CG}
\]
and the map $\pi_k:\cH_m\otimes\cH_n\rightarrow\cH_{m+n-2k}$
defined by
\[
\pi_k=\sqrt{
\frac{m!\ n!\ (m+n-2k+1)!}{k!\ (m+n-k+1)!\ (m-k)!\ (n-k)!}
}\times \Pi_{CG}
\]
satisfy the following properties.
\begin{Proposition}\label{piiotaprop}{\ }

\begin{enumerate}
\item
$\io_k$ is an isometric injection.
\item
$\pi_k$ is the Hermitian conjugate of $\io_k$.
\item
$\pi_k\circ\io_k={\rm Id}_{\cH_{m+n-2k}}$.
\item
The images of the $\io_k$, $0\le k\le \min(m,n)$,
are orthogonal.
\item
One has the decomposition of the identity
\[
\sum\limits_{k=0}^{\min(m,n)} \io_k\circ\pi_k={\rm Id}_{\cH_m\otimes\cH_n}\ .
\]
\end{enumerate}
\end{Proposition}
These are the standard normalizations for the maps underlying
the QAMT in the physics literature.
See~\cite[\S7]{AC3}
for a precise dictionary with the CIT approach, and formulae for
6-j and 9-j symbols.

\begin{Remark}\label{looprem}
In some instances one might need to allow trivial components
$
\parbox{0.8cm}{
\psfrag{a}{$\scriptstyle{a}$}
\includegraphics[width=0.8cm]{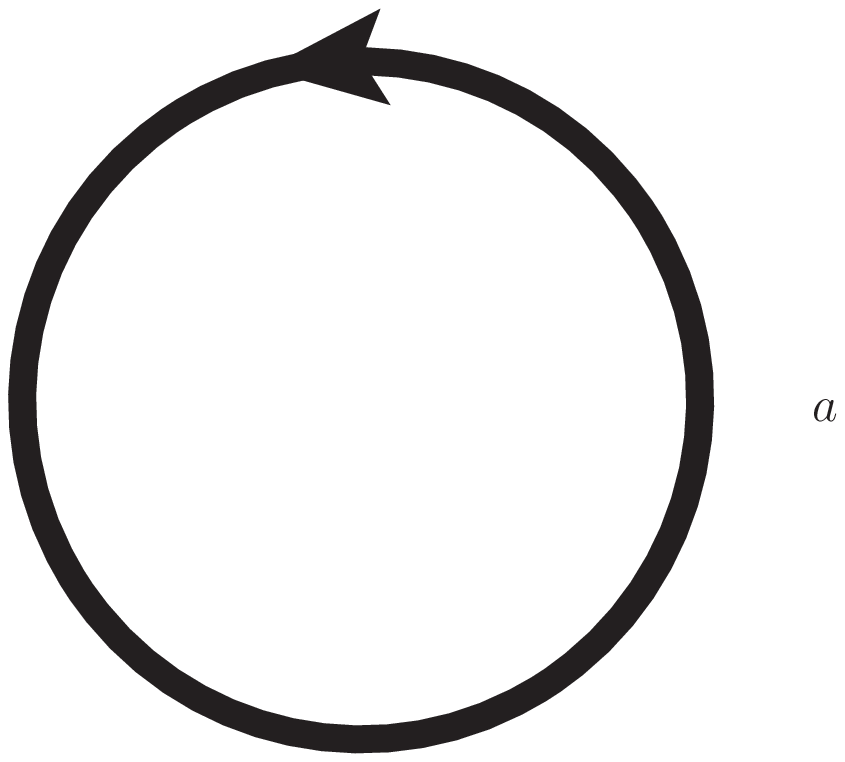}}
$ in a CG network even though the underlying graph $G$ does not
satisfy
the usual definitions of a graph.
One can decide that the evaluation of such a component is
\[
\<
\parbox{0.8cm}{
\psfrag{a}{$\scriptstyle{a}$}
\includegraphics[width=0.8cm]{Fig136.eps}}
\>^{CG}=
\parbox{1.4cm}{
\psfrag{a}{$\scriptstyle{a}$}
\includegraphics[width=1.4cm]{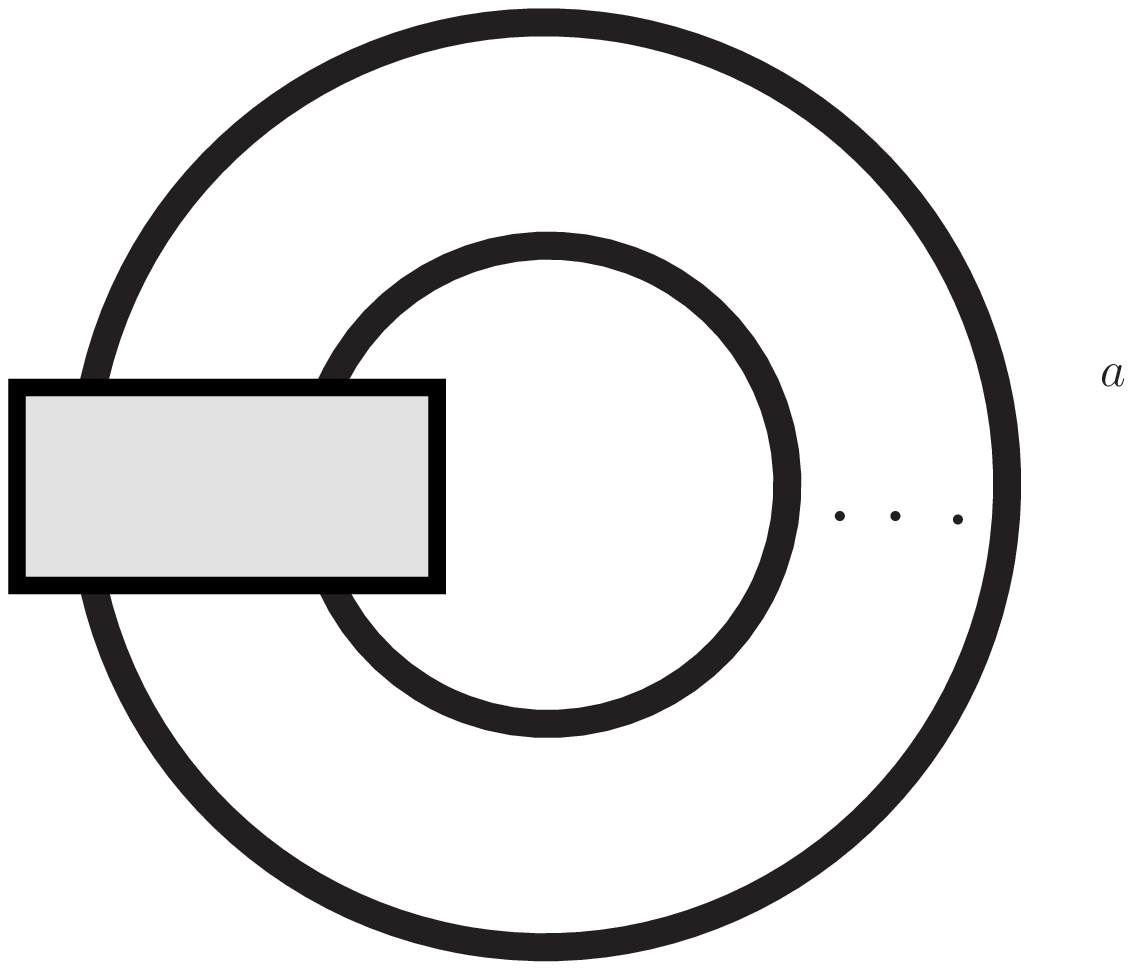}}
= a+1
\]
from (\ref{loopevalCG}). 
\end{Remark}

\section{The negative dimensionality theorem}\label{negdimsec}
We now turn to the precise relationship between spin networks of
\S\ref{introdefsec} and
CG networks.
Let $(\Ga,\ga)$
be a spin network as in \S\ref{introdefsec}.
For more precision let us write $\Ga=(G,R)$
where $G$ is the underlying cubic graph and $R$ is the extra structure
corresponding to the pure rotation system needed in order
to define the imbedding of
$\Ga$ in an orientable compact Riemann surface $\Si$.
Suppose one also has a smooth orientation $\cO$ of the edges
of $G$ and a gate signage $\ta$ as in \S\ref{CGsection}.
Then one has the following correspondence which is reminiscent
of the negative dimensionality theorem of~\cite{ElvangCK} (see 
also~\cite[Ch. 13]{Cvitanovic}).

\begin{Theorem}\label{negdimthm}
One has
\[
\<\Ga,\ga\>^P=\mu\times
\left(
\prod\limits_{e\in E(G)} \ga(e)!
\right)\times \<G,\cO,\ta,\ga\>^{CG}
\]
where $\mu=\pm$ is a global sign which depends on the given combinatorial
data.
\end{Theorem}
\noindent{\bf Proof:}
Because of the factorization property of Lemma \ref{factolemma}
which also holds for CG evaluations of CG networks,
it is enough
to consider connected graphs.
The case of trivial components is dealt with by the comparison
of Lemma \ref{trivcomplemma} and Remark \ref{looprem}.
We now only consider nontrivial connected
cubic graphs which may have loops and multiple edges.
Define $H(G)$ the set of half-edges of $G$ obtained
by cutting the edges of $G$ according to
\[
\parbox{0.8cm}{
\psfrag{e}{$\scriptstyle{e}$}
\includegraphics[width=0.8cm]{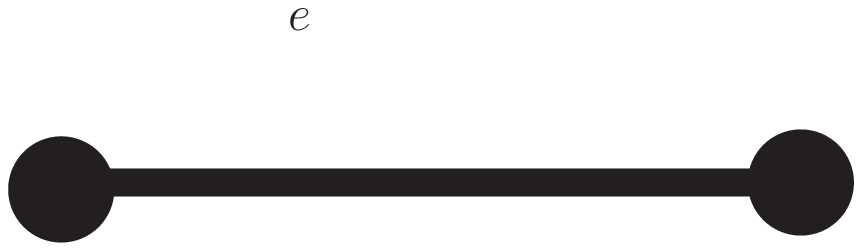}}
\qquad\longrightarrow\qquad
\parbox{1.1cm}{
\psfrag{1}{$\scriptstyle{h_1}$}
\psfrag{2}{$\scriptstyle{h_2}$}
\includegraphics[width=1.1cm]{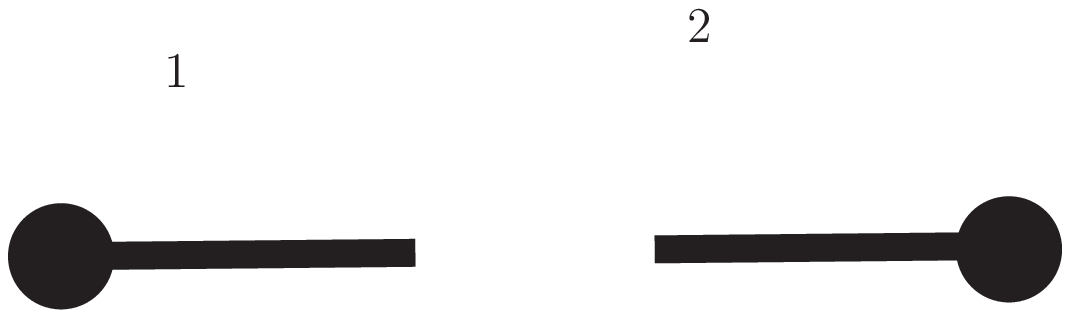}}
\]
\[
\parbox{0.8cm}{
\psfrag{e}{$\scriptstyle{e}$}
\includegraphics[width=0.8cm]{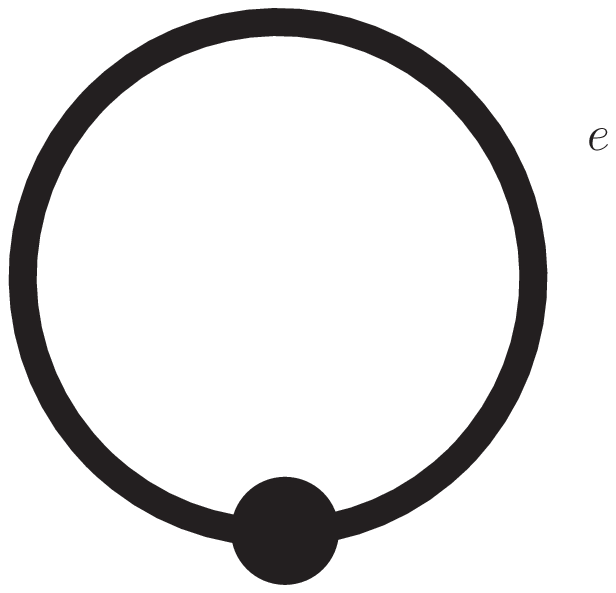}}
\qquad\longrightarrow\qquad
\parbox{1.1cm}{
\psfrag{1}{$\scriptstyle{h_1}$}
\psfrag{2}{$\scriptstyle{h_2}$}
\includegraphics[width=1.1cm]{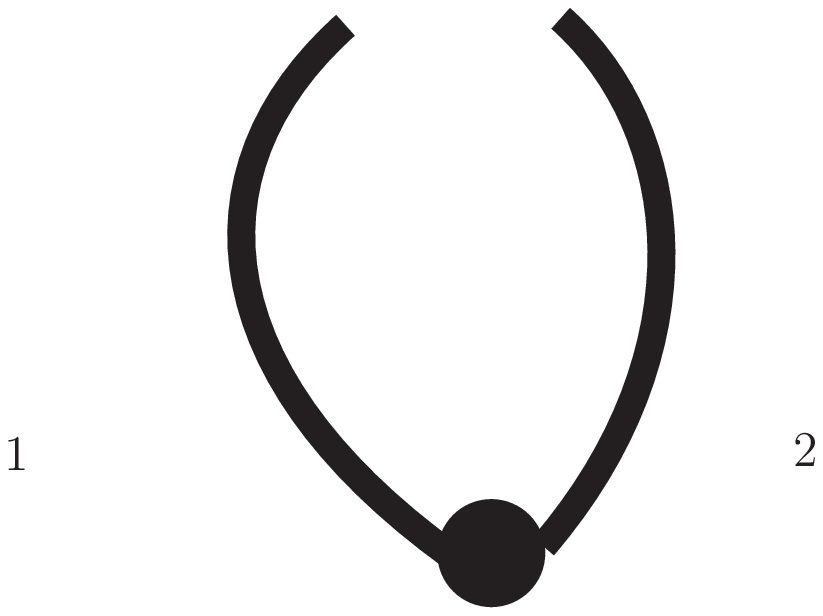}}
\]
where $e\in E(G)$ and $h_1\neq h_2$ belong to $H(G)$.
Two half-edges are called edge-adjacent if they are distinct
and come from the same edge.
Two half-edges are called vertex-adjacent if they are distinct
and attached to the same vertex.

A directed path is a sequence
\[
v_0,h_1,h'_1,v_1,h_2,h'_2,\ldots
v_{n-1},h_n,h'_n,v_n
\]
where $n\ge 0$ and where the $v_i$ are vertices
and the $h_i,h'_i$ are half-edges.
We also impose the condition that $h_i$ is incident to $v_{i-1}$,
for $1\le i\le n$, that $h'_i$ is incident to $v_i$, for $1\le i\le n$,
and that $h_i,h'_i$ are edge-adjacent for $1\le i\le n$.
Finally we enforce a no backtracking clause
$h'_i\neq h_{i+1}$ for any $i$, $1\le i< n$.
\noindent{\bf Example:}
On the graph
\[
\parbox{2cm}{
\psfrag{1}{$\scriptstyle{h_1}$}
\psfrag{2}{$\scriptstyle{h_2}$}
\psfrag{3}{$\scriptstyle{h_3}$}
\psfrag{4}{$\scriptstyle{h_4}$}
\psfrag{5}{$\scriptstyle{h_5}$}
\psfrag{6}{$\scriptstyle{h_6}$}
\psfrag{a}{$\scriptstyle{v_1}$}
\psfrag{b}{$\scriptstyle{v_2}$}
\includegraphics[width=2cm]{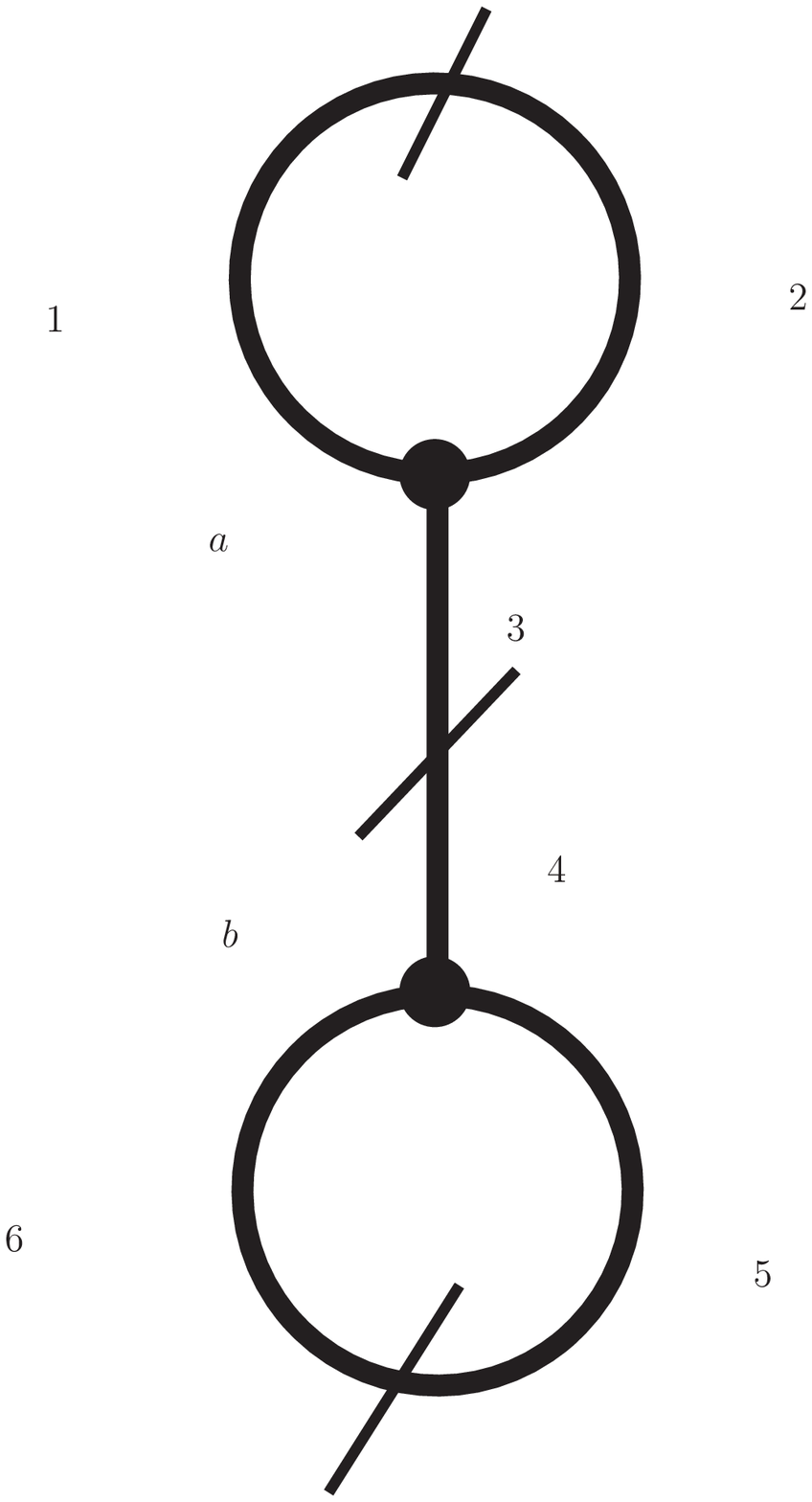}}
\]
the path
\[
\parbox{4cm}{\psfrag{s}{${\rm \scriptstyle{start}}$}
\psfrag{f}{${\rm \scriptstyle{finish}}$}
\includegraphics[width=2cm]{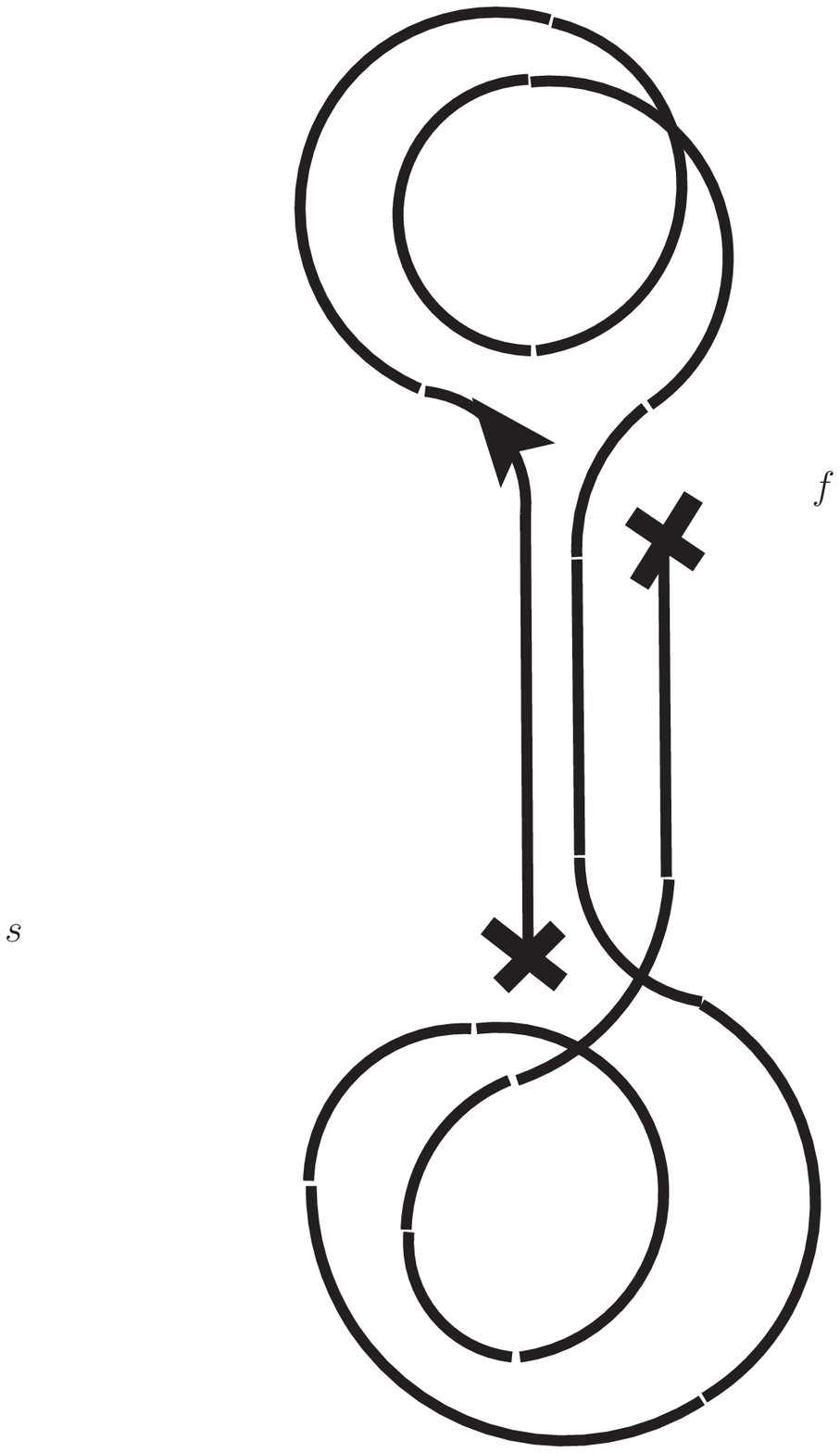}}
\]
is allowed and corresponds to the sequence
\[
v_2,h_4,h_3,v_1,h_1,h_2,v_1,h_1,h_2,v_1,h_3,h_4,v_2,h_5,h_6,v_2,h_5,h_6,
v_2,h_4,h_3,v_1.
\]
An undirected path is an equivalence class of directed paths with respect
to reversal
\[
v_0,h_1,h'_1,\ldots,h'_n,v_n\qquad\longrightarrow\qquad
v_n,h'_n,h_n,\ldots,h_1,v_0\ .
\]
A pointed closed directed path is a directed path with $n\ge 1$
such that $v_n=v_0$
and $h_1\neq h'_n$.
Note that on the same graph as before
\[
\parbox{1.2cm}{
\includegraphics[width=1.2cm]{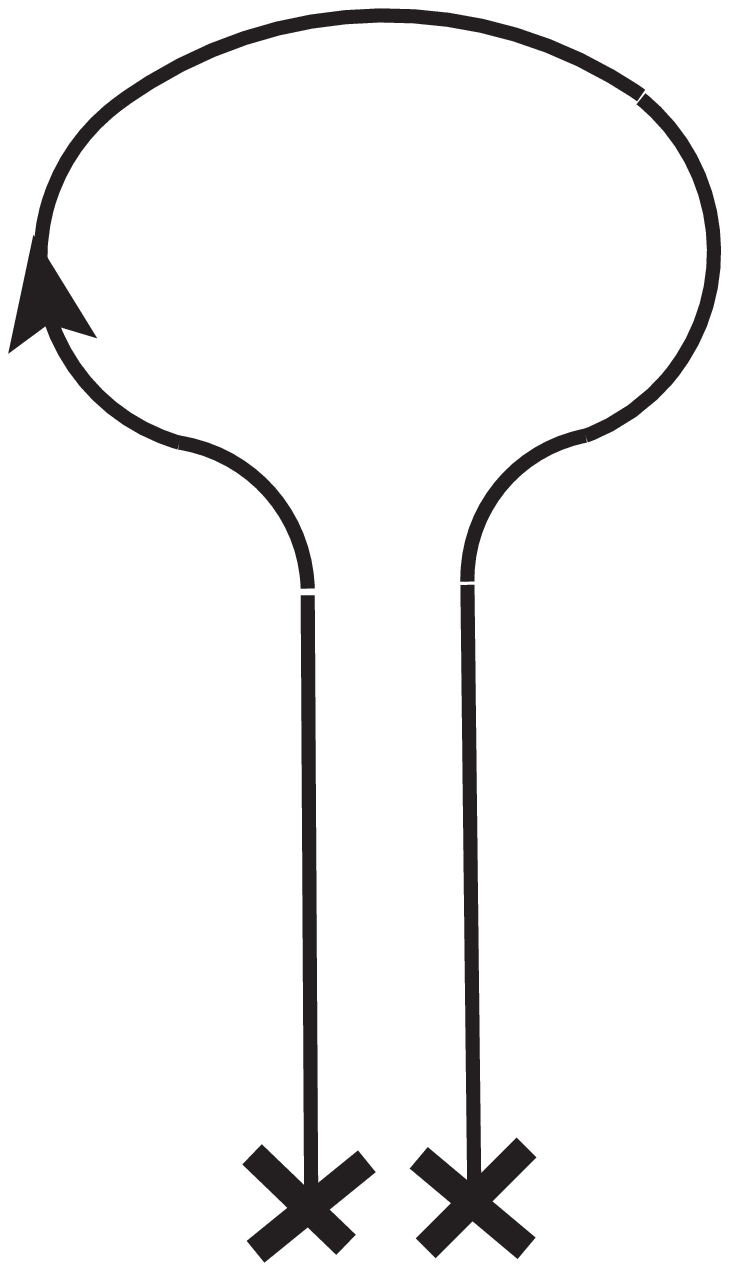}}
\]
violates the last condition and is not allowed.
A closed directed path is an equivalence class
of pointed closed directed paths under cyclic
transformation
\[
v_0,h_1,h'_1,\ldots,h'_n,v_n\ \longrightarrow\ 
v_i,h_{i+1},h'_{i+1},\ldots,h'_n,v_n,h_1,h'_1,v_1,\ldots
h_i,h'_i,v_i
\]
where $v_0$ has been deleted and $v_i$ has been repeated.
A closed path is an equivalence class of pointed directed paths under both
cyclic transformations and reversal.

When expanding the state sum over all connecting permutations in the definition
of $\<\Ga,\ga\>^P$, each continuous curve drawn on $\Si$
defines in an unamiguous manner a closed
path in the previous sense.
Each term contributes $(-1)^{C(\vec{\si})}(-2)^{N(\vec(\si))}$
where $\vec{\si}=(\si_e)_{e\in E(G)}\in\prod_{e\in E(G)}\gS_{\ga(e)}$
is the given configuration
or state, $C(\vec{\si})$ denotes
the total number of strand crossings and $N(\vec{\si})$
is the number of closed curves.
Consider such a curve $\cC$
or rather its correponding closed path
$[v_0,h_1,h'_1,\ldots,h'_n,v_n]$ where brackets mean equivalence class.
To such data one can associate
the number $B(\cC)$ of gates crossed.
Indeed, the half-edges inherit the $\cO$
orientation of their parent edge.
Gates correspond to vertex-adjacent pairs of half-edges with opposite
orientations
\[
\parbox{1.2cm}{
\includegraphics[width=1cm]{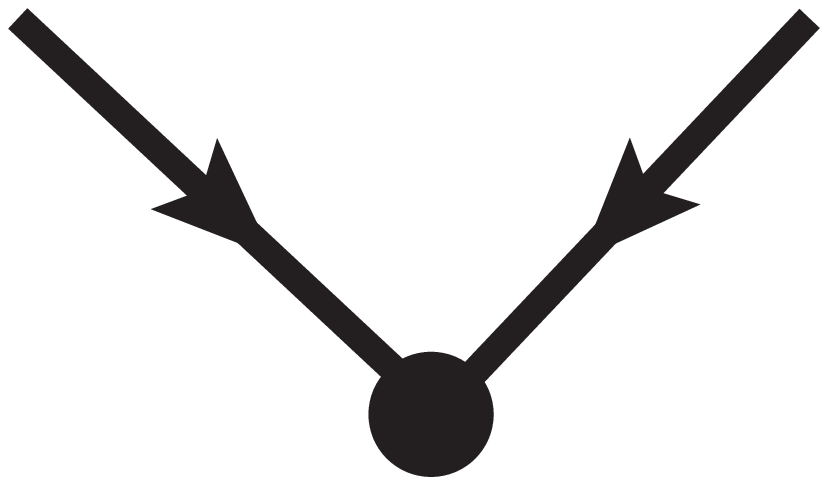}}
\qquad {\rm or}\qquad
\parbox{1.2cm}{
\includegraphics[width=1cm]{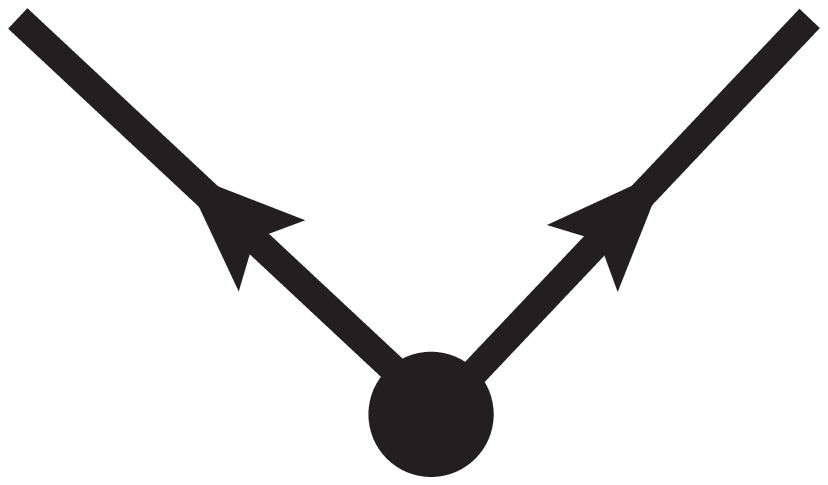}}\ \ .
\]
A gate crossing is an index $i$, $1\le i\le n$,
such that either $i<n$
and $h'_i,h_{i+1}$ form a gate
or $i=n$ and $h'_n,h_1$ form a gate.
We now make the following observation.

\noindent{\bf Key observation:}
The number $B(\cC)$ of gates crossed must be even.

\noindent
Indeed, this is the number of changes of direction (relative to the edge
orientation $\cO$) encountered as one travels along the closed path
starting and ending at the same half-edge say $h_1$.

Now choose a closed directed path structure for the closed path of $\cC$.
This amounts to picking a direction of travel.
One then defines a good gate crossing as a gate crossing of the form
\[
\parbox{2cm}{\psfrag{d}{${\rm \scriptstyle{direction\ of\ travel}}$}
\includegraphics[width=2cm]{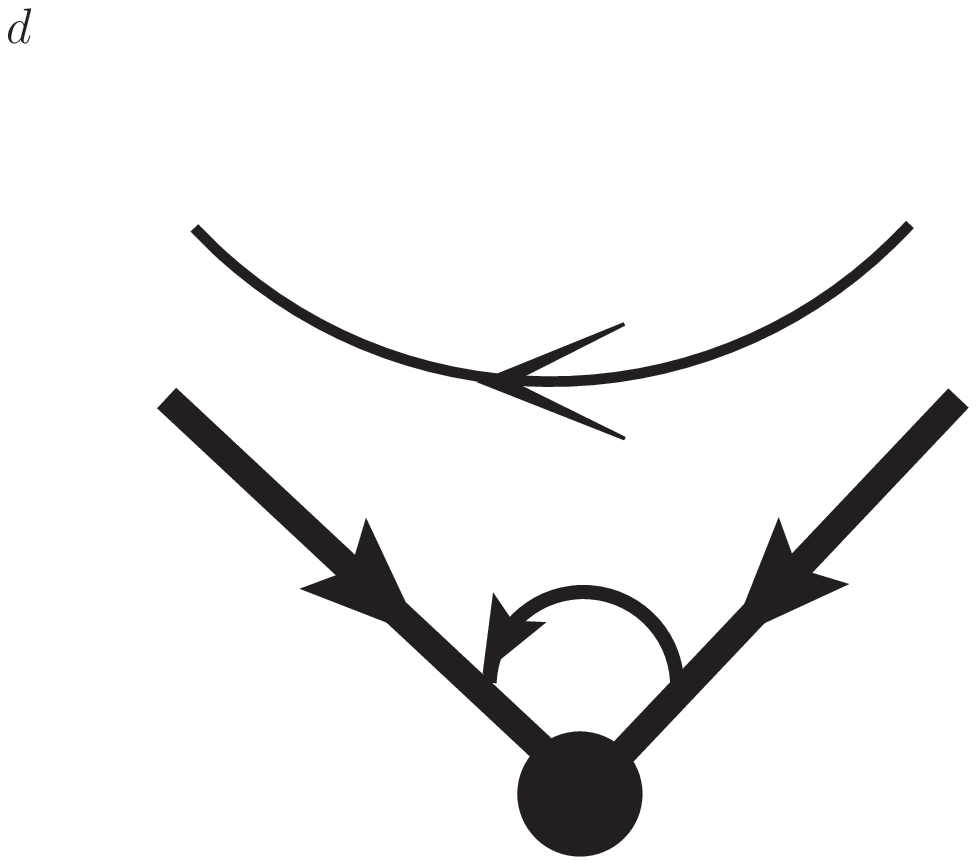}}
\qquad {\rm or}\qquad
\parbox{1.6cm}{\includegraphics[width=1.6cm]{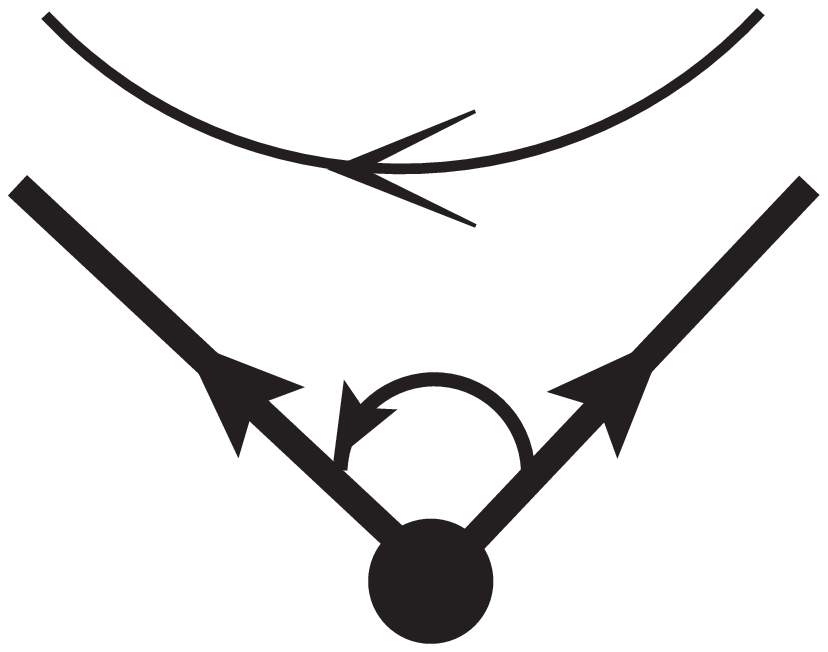}}
\]
which is a notion relative to the gate signage $\ta$.
The other crossings of the form
\[
\parbox{1.6cm}{\includegraphics[width=1.6cm]{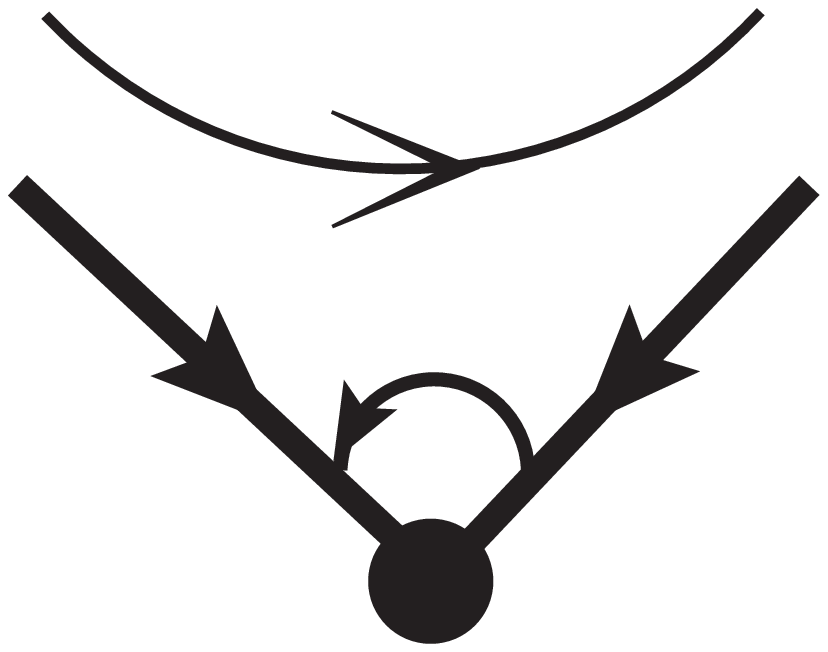}}
\qquad {\rm or}\qquad
\parbox{1.6cm}{\includegraphics[width=1.6cm]{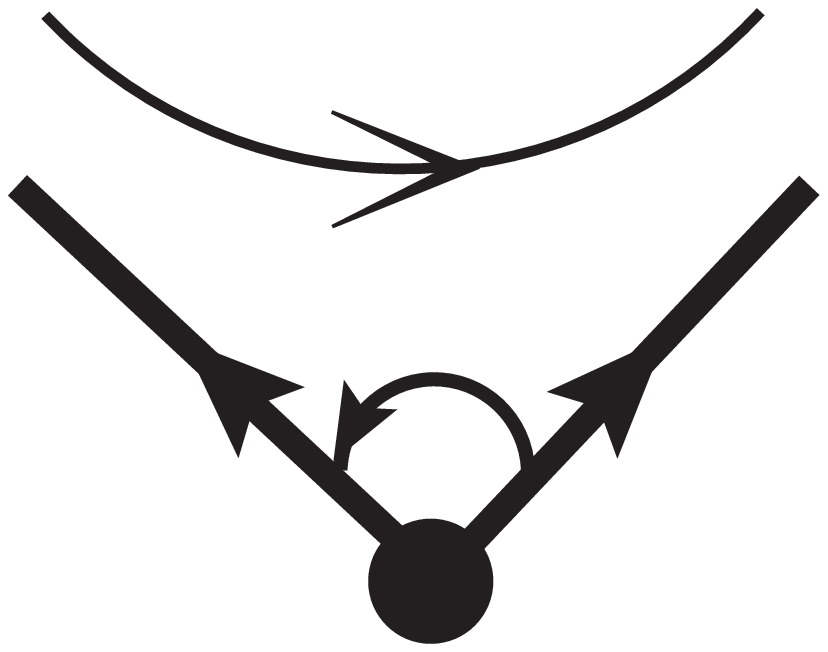}}
\]
are bad gate crossings.
Let $B_{+}(\cC)$ be the number of good gate crossings
and $B_{-}(\cC)$ that of bad gate crossings.
Obviously $B(\cC)=B_{+}(\cC)+B_{-}(\cC)$
and choosing the other direction of travel merely
exchanges $B_{+}(\cC)$ and $B_{-}(\cC)$.
Because of the key observation that $B(\cC)$ is even,
the sign $(-1)^{B_{+}(\cC)}=(-1)^{B_{-}(\cC)}$
is independent of the direction of travel.
Let $B_{+}(\vec{\si})$
be the sum, over all curves $\cC$ present in a given configuration
produced by $\vec{\si}$, of the $B_{+}(\cC)$'s.
The main fact needed to prove the theorem is the following statement.

\noindent{\bf Claim:}
The sign $\tilde{\mu}=(-1)^{C(\vec{\si})+N(\vec{\si})+B_{+}(\vec{\si})}$
is independent of $\vec{\si}$.

\noindent
Since all configurations $\vec{\si}$
are related by sequences of transpositions it is enough
to show the invariance of this sign with respect to any such
transposition of the strands at some edge of the cubic graph.
Of course $(-1)^{C(\vec{\si})}$
changes to the opposite sign under a transposition and we need to show
that so does $(-1)^{N(\vec{\si})+B_{+}(\vec{\si})}$.
For a strand transposition
\[
\parbox{0.4cm}{\includegraphics[width=0.4cm]{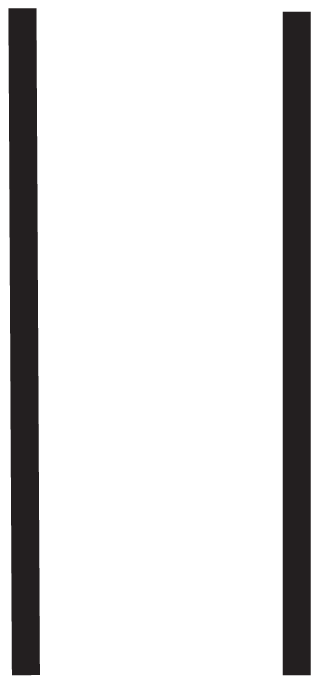}}
\qquad\longrightarrow\qquad
\parbox{0.4cm}{\includegraphics[width=0.4cm]{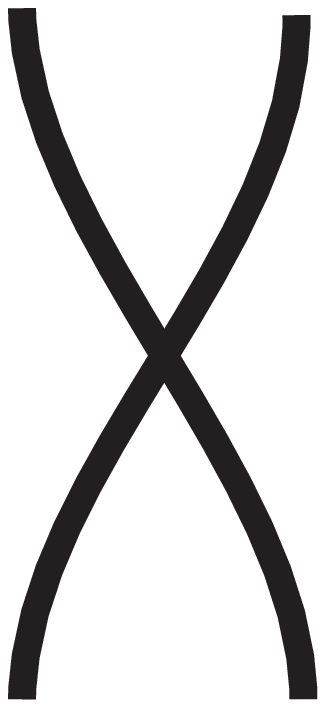}}
\]
depending on the configuration of curves $\cC$
to which these strands belong, there are three cases to consider
(and not two as in~\cite[Ch. 13]{Cvitanovic}):
\[
\begin{array}{c}
\parbox{3.6cm}{\includegraphics[width=3.6cm]{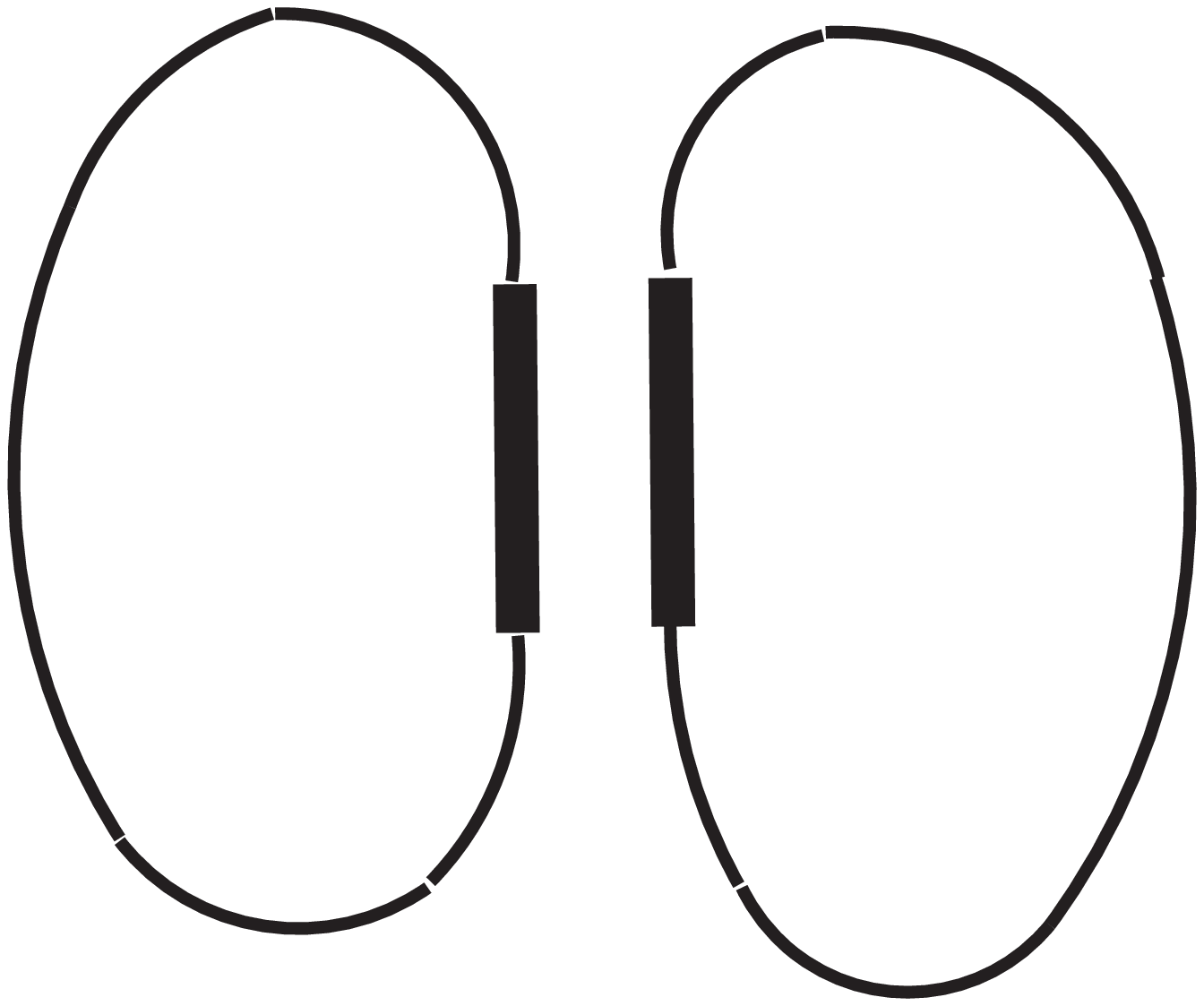}}\\
{\rm case\ I}
\end{array}
\qquad
\begin{array}{c}
\parbox{1.5cm}{\includegraphics[width=1.5cm]{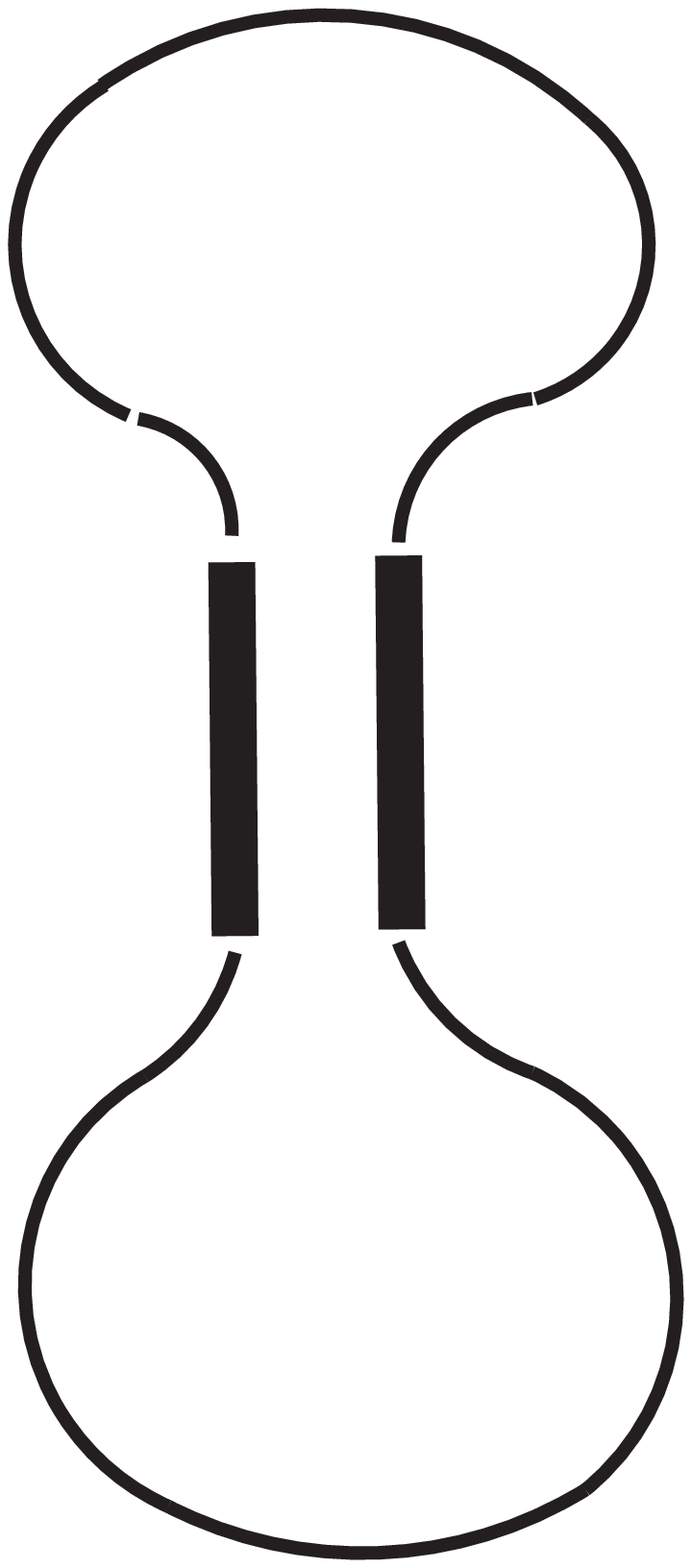}}\\
{\rm case\ II}
\end{array}
\qquad
\begin{array}{c}
\parbox{2cm}{\includegraphics[width=2cm]{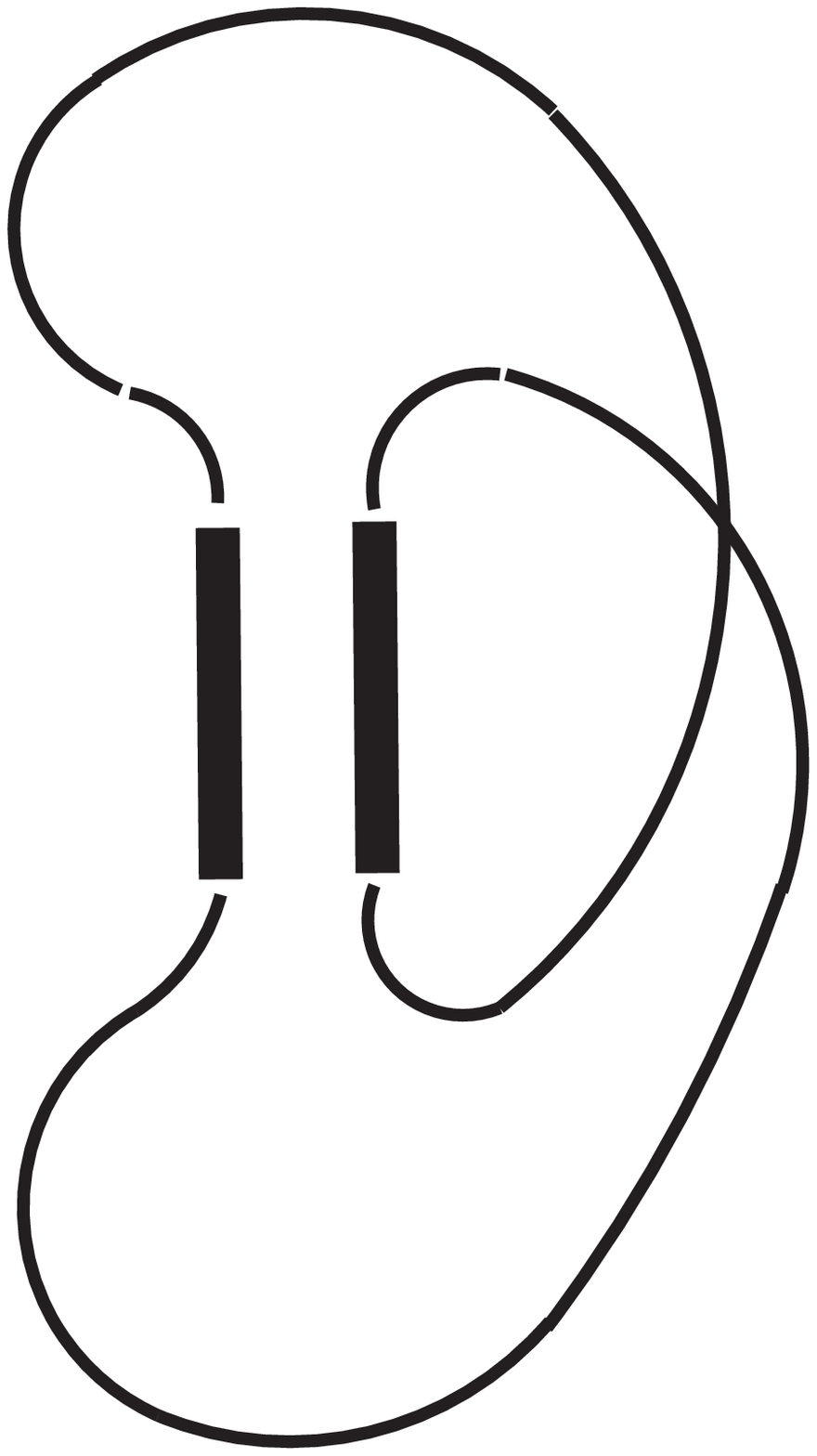}}\\
{\rm case\ III}
\end{array}
\]
where we highlighted the strands along the edge of the transposition.

\noindent{\bf Case I:}
\[
\parbox{3.6cm}{\includegraphics[width=3.6cm]{Fig153.eps}}
\qquad
\stackrel{\rm transposition}{\longrightarrow}\qquad
\parbox{3.6cm}{\includegraphics[width=3.6cm]{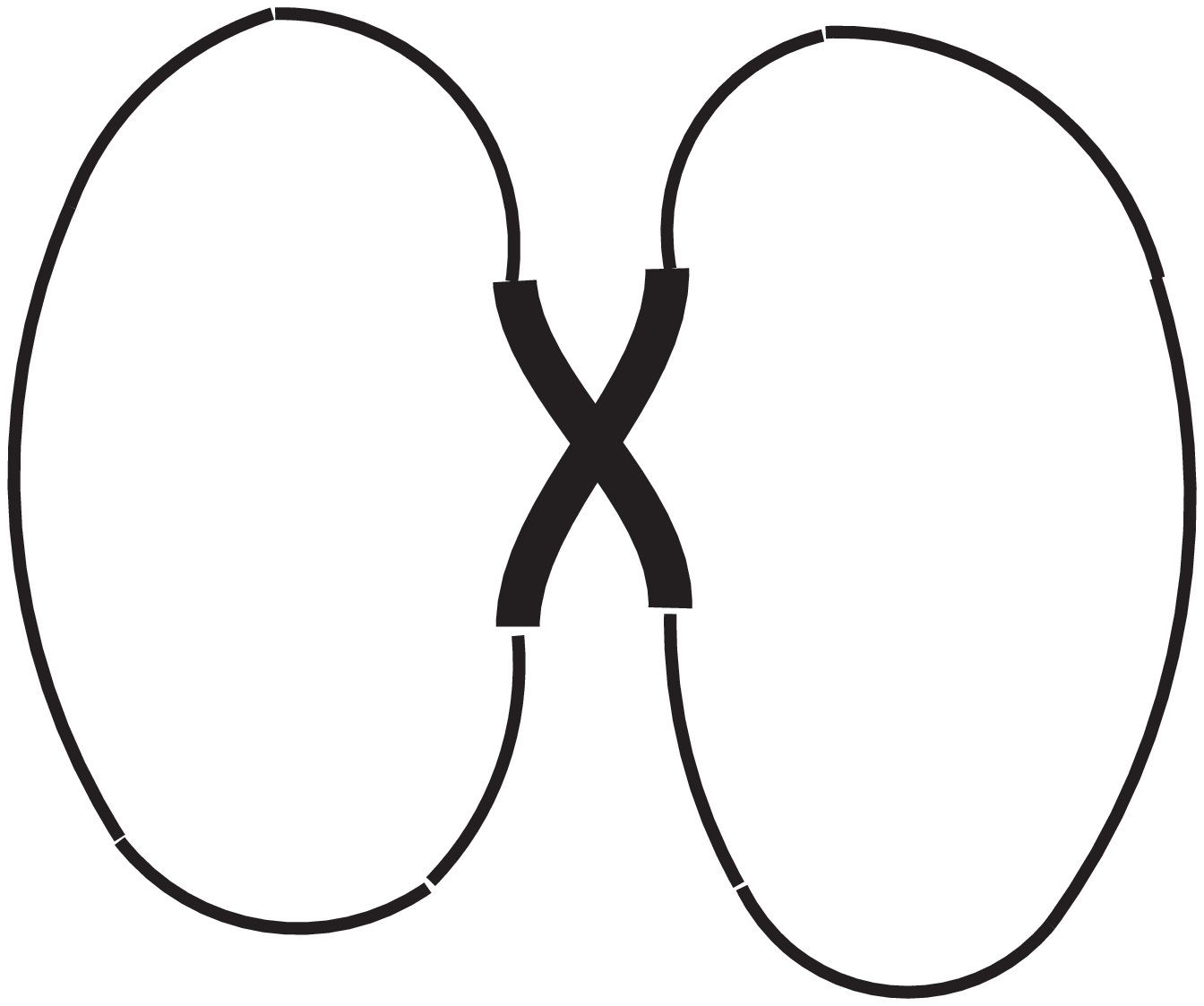}}\ \ .
\]
This fuses two curves or closed paths into one.
Choose a direction of travel on the resulting closed
path and deduce from it a direction of travel
on both initial curves as in
\[
\parbox{3.6cm}{\psfrag{1}{$\cC_1$}
\psfrag{2}{$\cC_2$}
\includegraphics[width=3.6cm]{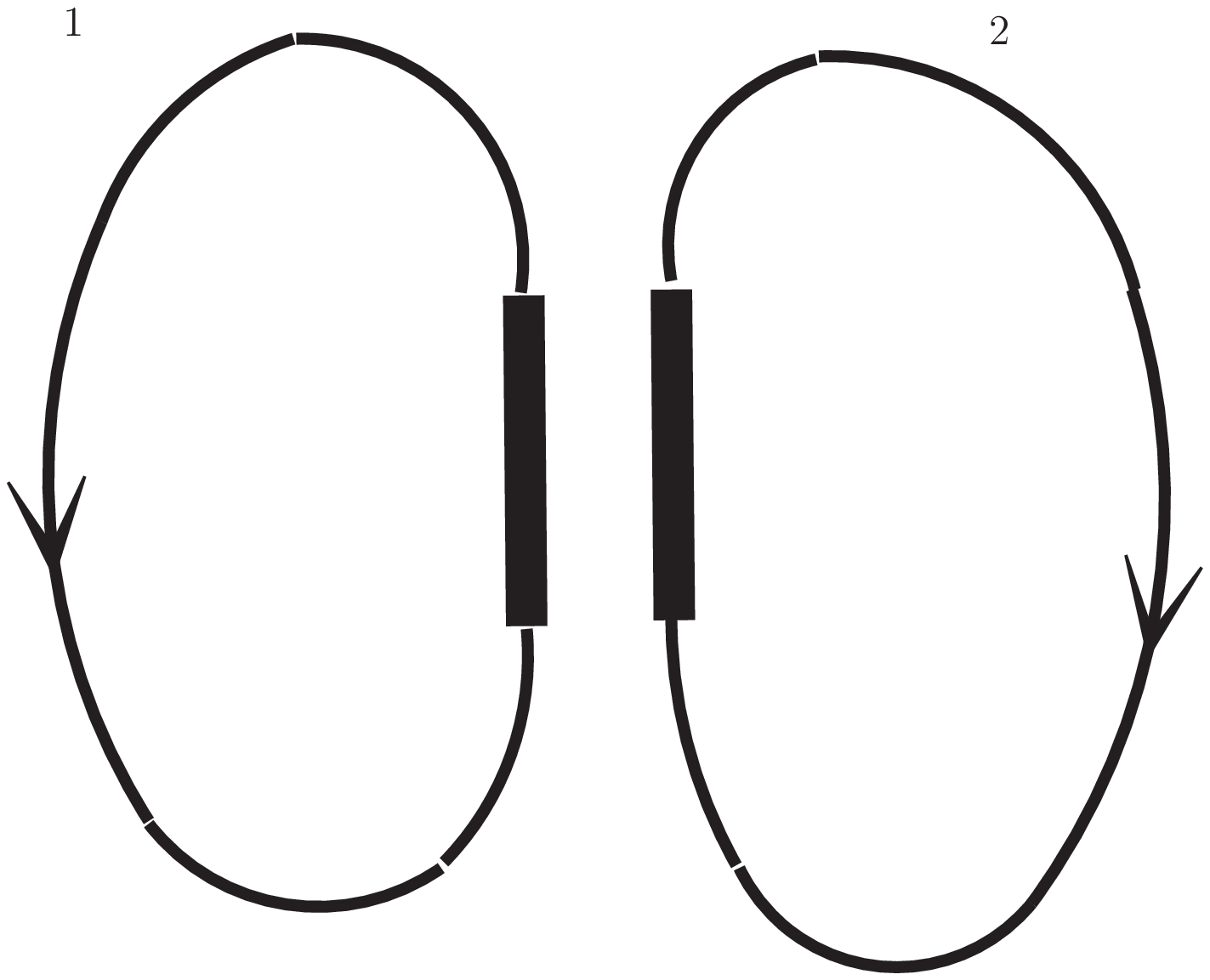}}
\qquad , \qquad
\parbox{3.6cm}{\psfrag{C}{$\cC$}
\includegraphics[width=3.6cm]{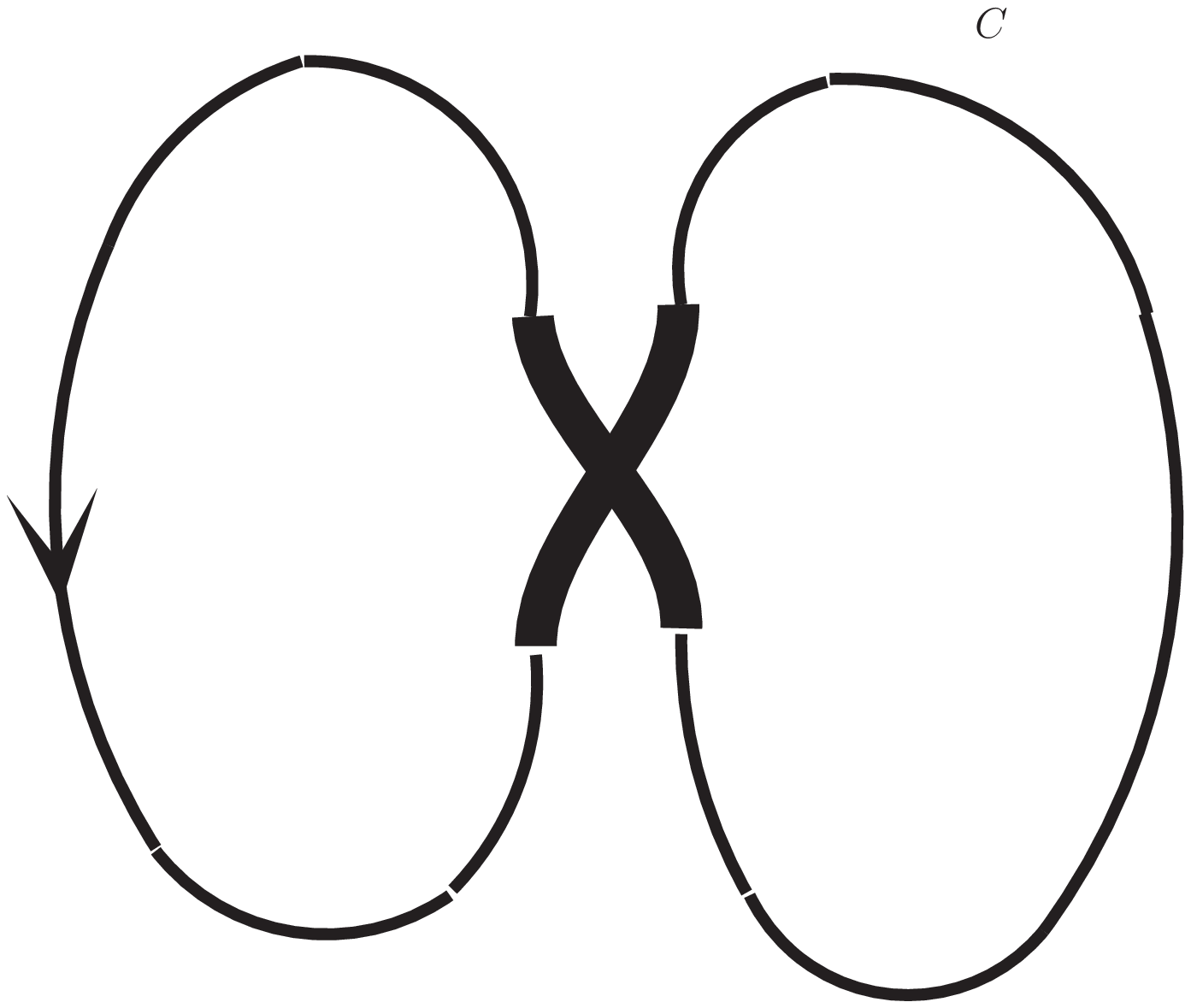}}\ \ .
\]
We have that $N(\vec{\si})$ dropped to $N(\vec{\si})-1$.
On the other hand, one clearly has
$B_{+}(\cC)=B_{+}(\cC_1)+B_{+}(\cC_2)$
therefore $B_{+}(\vec{\si})$ remains unchanged.

\noindent{\bf Case II:}
\[
\parbox{2cm}{\psfrag{1}{$\cC_1$}
\includegraphics[width=2cm]{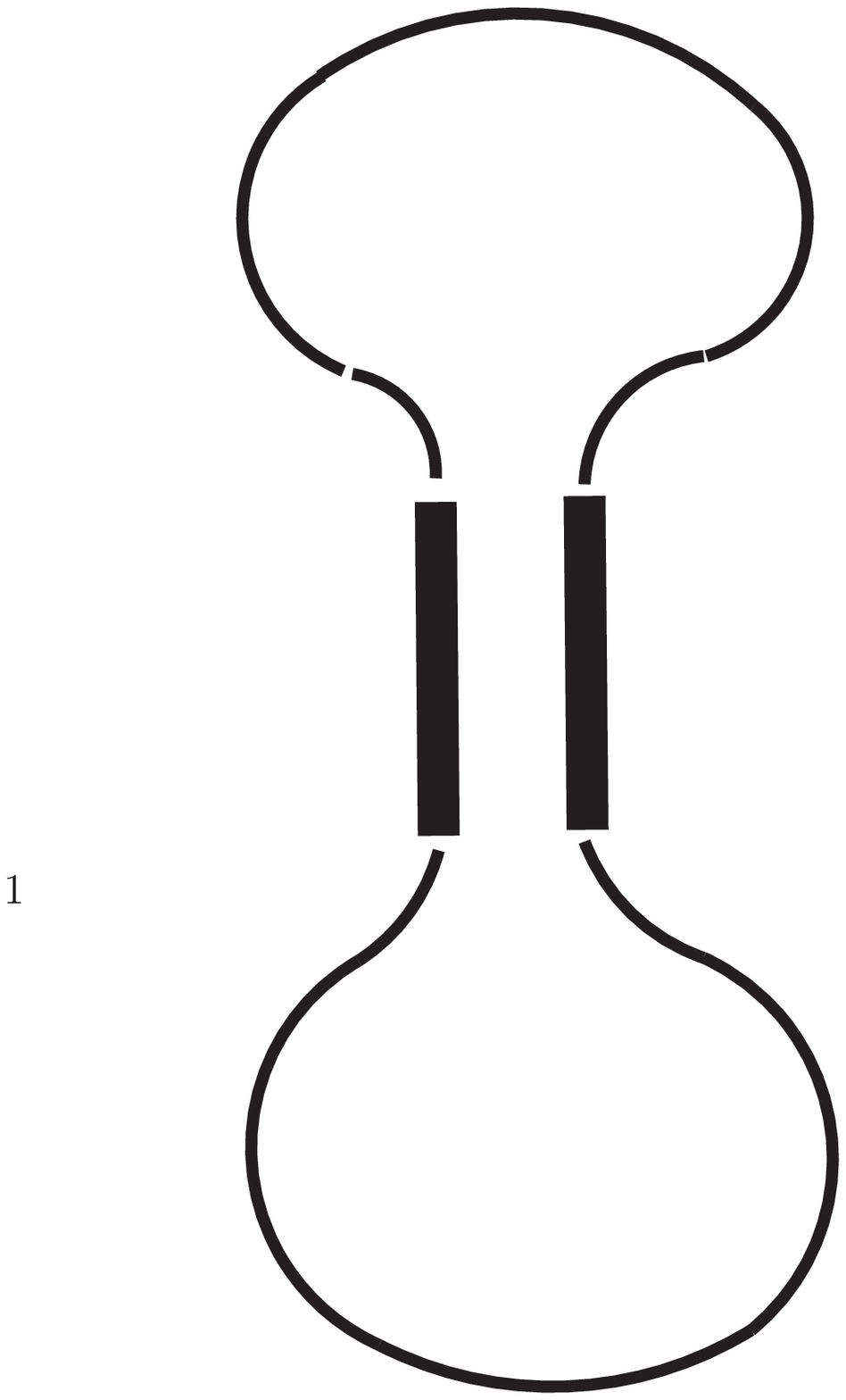}}
\qquad
\stackrel{\rm transposition}{\longrightarrow}\qquad
\parbox{2cm}{\psfrag{2}{$\cC_2$}
\includegraphics[width=2cm]{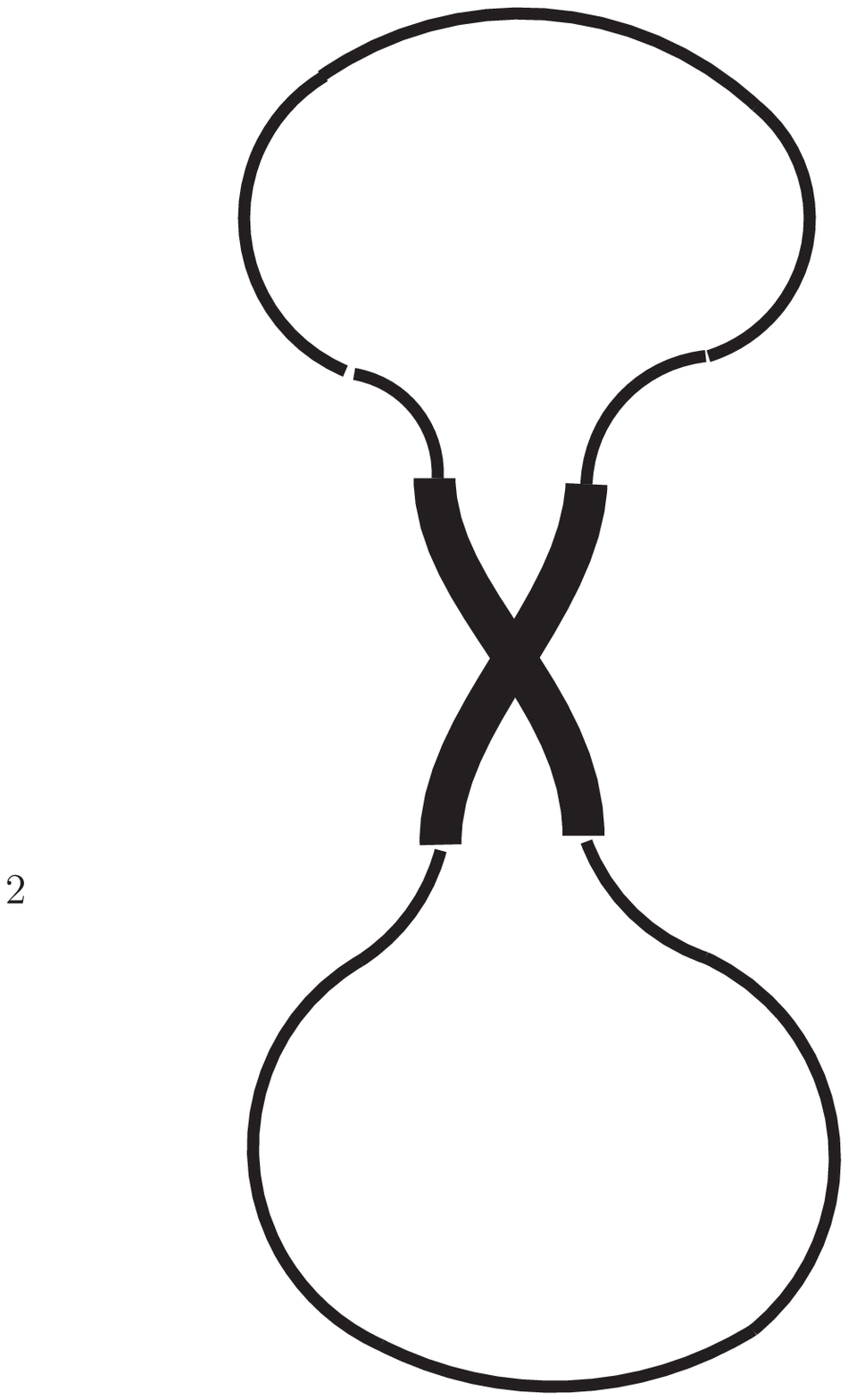}}\ \ .
\]
Choose a direction of travel on $\cC_1$ and $\cC_2$.
These will coincide on one of the two portions
of the path (say the bottom one on the picture)
and differ on the other portion
\[
\parbox{2cm}{\psfrag{1}{$\cC_1$}
\includegraphics[width=2cm]{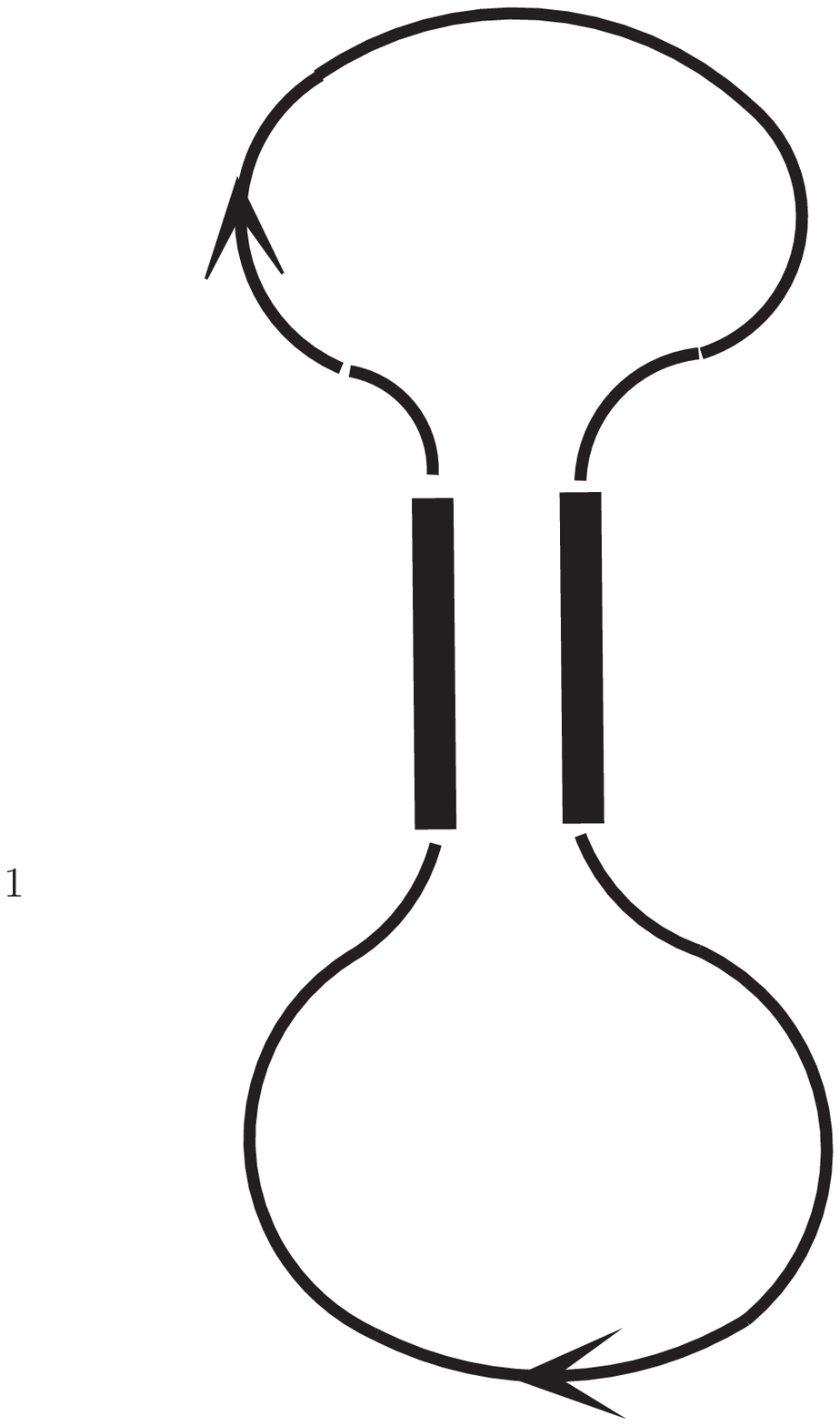}}
\qquad , \qquad
\parbox{2cm}{\psfrag{2}{$\cC_2$}
\includegraphics[width=2cm]{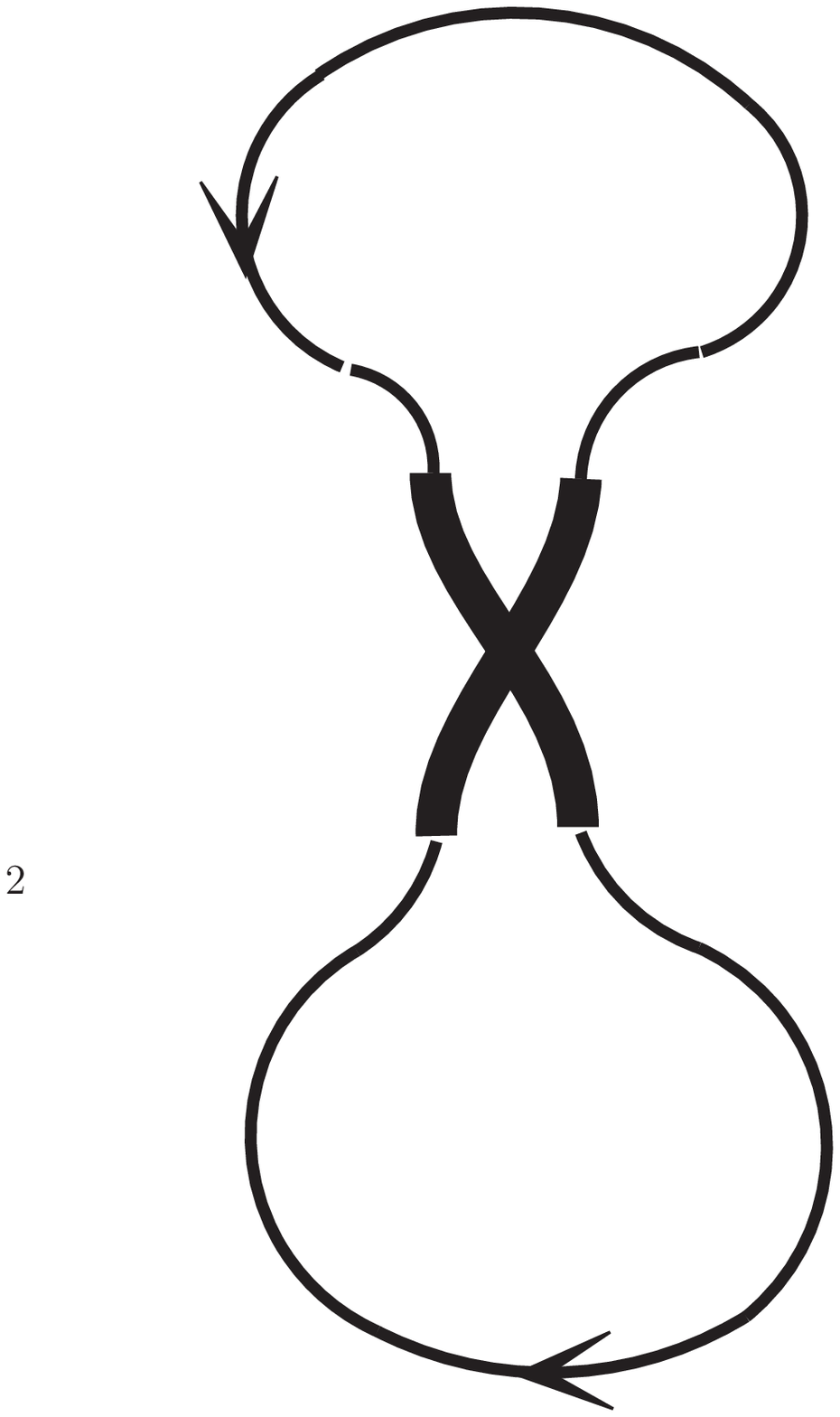}}\ \ .
\]
By the same argument as in the key observation
one easily sees
that the number of gate crossings in the upper portion is odd.
Therefore the number of good and bad crossings on this portion
add up to an odd number and therefore have different parity.
The reversal of direction of travel on the upper portion
will exchange these two numbers
and thus $(-1)^{B_{+}(\cC_2)}=-(-1)^{B_{+}(\cC_1)}$ while
$N(\vec{\si})$ remains unchanged.

\noindent{\bf Case III:}
\[
\parbox{2.6cm}{\psfrag{C}{$\cC$}
\includegraphics[width=2.6cm]{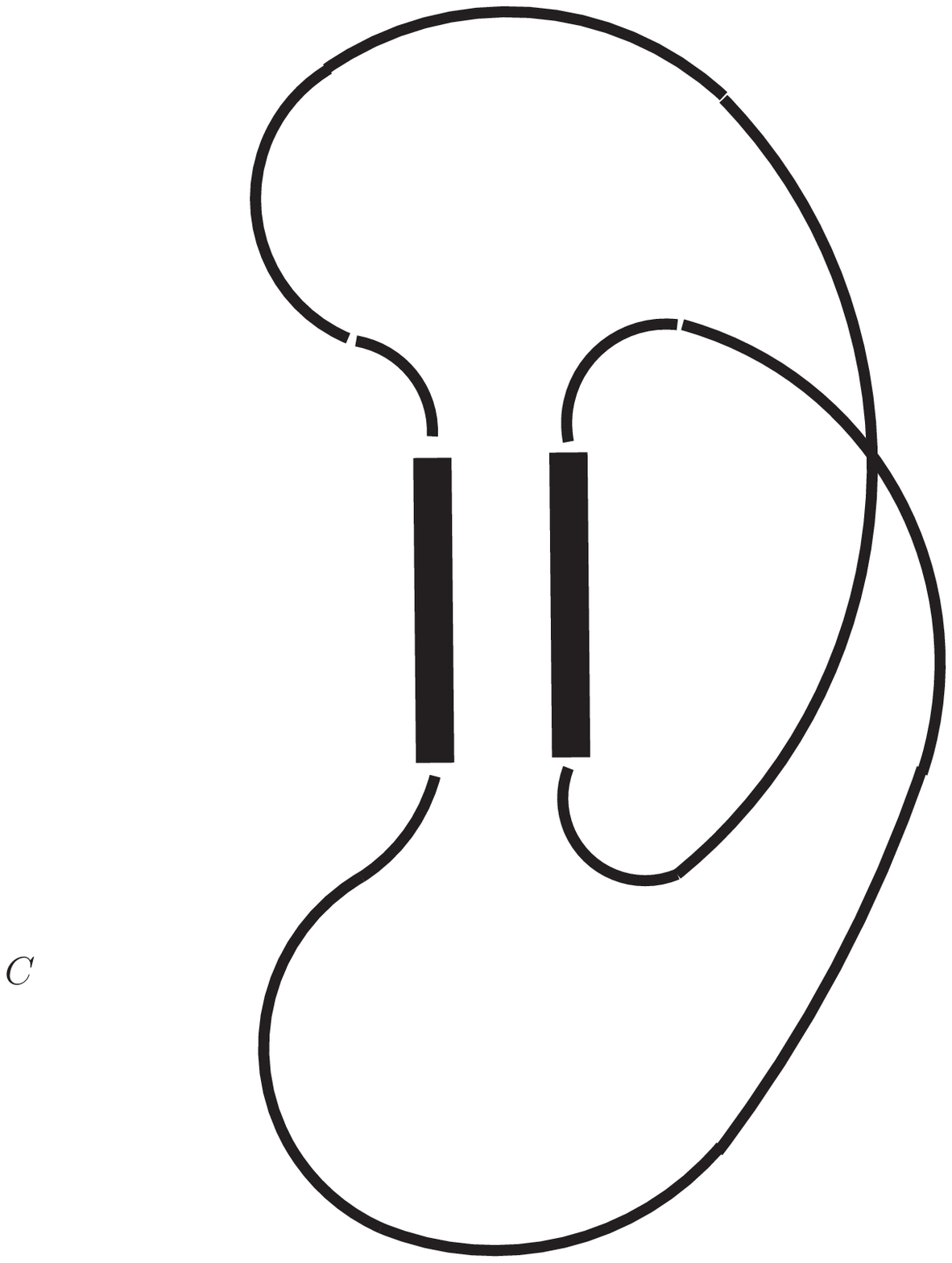}}
\qquad
\stackrel{\rm transposition}{\longrightarrow}\qquad
\parbox{2.6cm}{\psfrag{1}{$\cC_1$}\psfrag{2}{$\cC_2$}
\includegraphics[width=2.6cm]{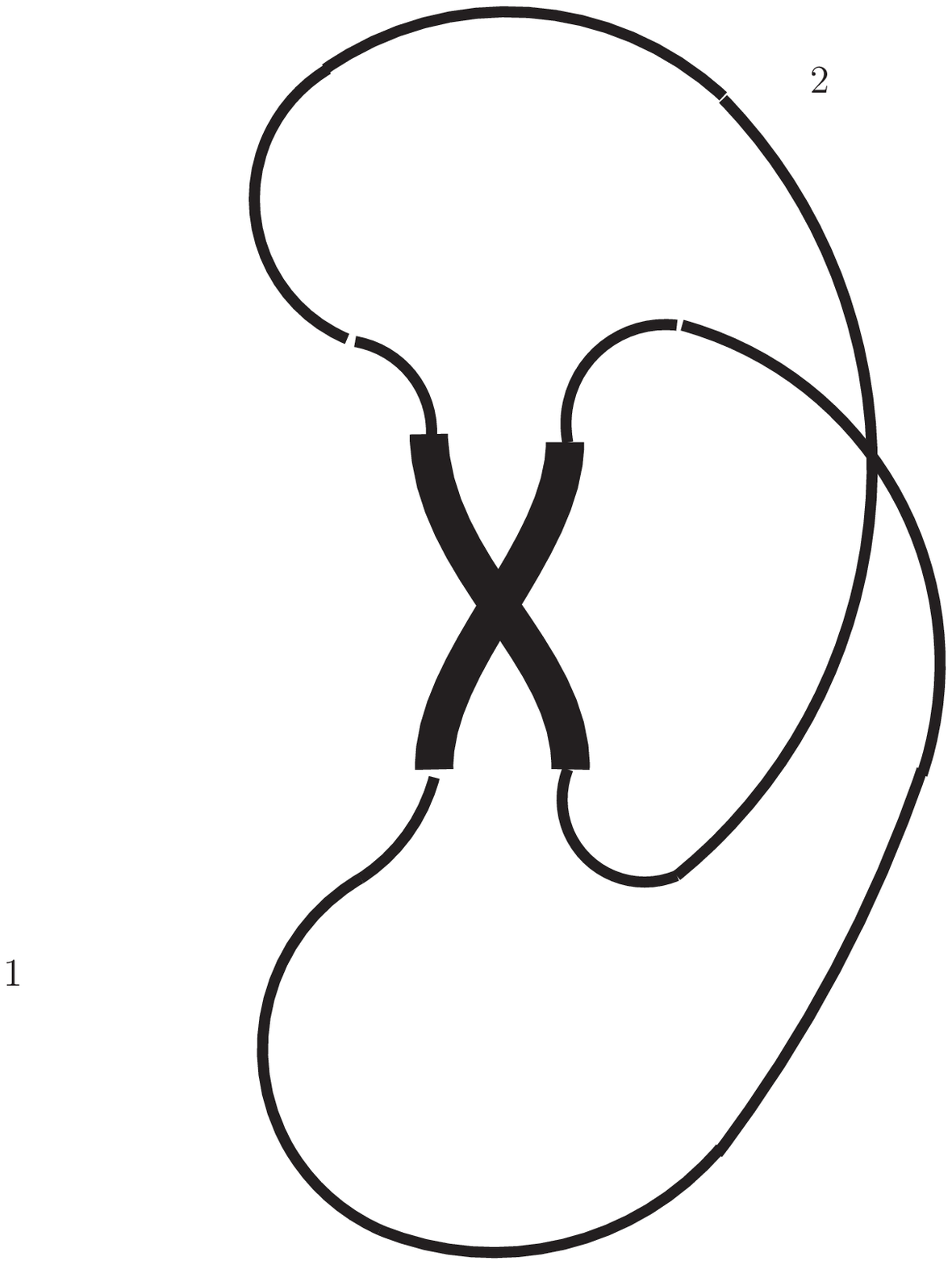}}\ \ .
\]
This is the undoing of case I.
One can again choose a direction of travel on $\cC$
and deduce from it
directions of travel on the two resuting curves $\cC_1$ and $\cC_2$
\[
\parbox{2.6cm}{\psfrag{C}{$\cC$}
\includegraphics[width=2.6cm]{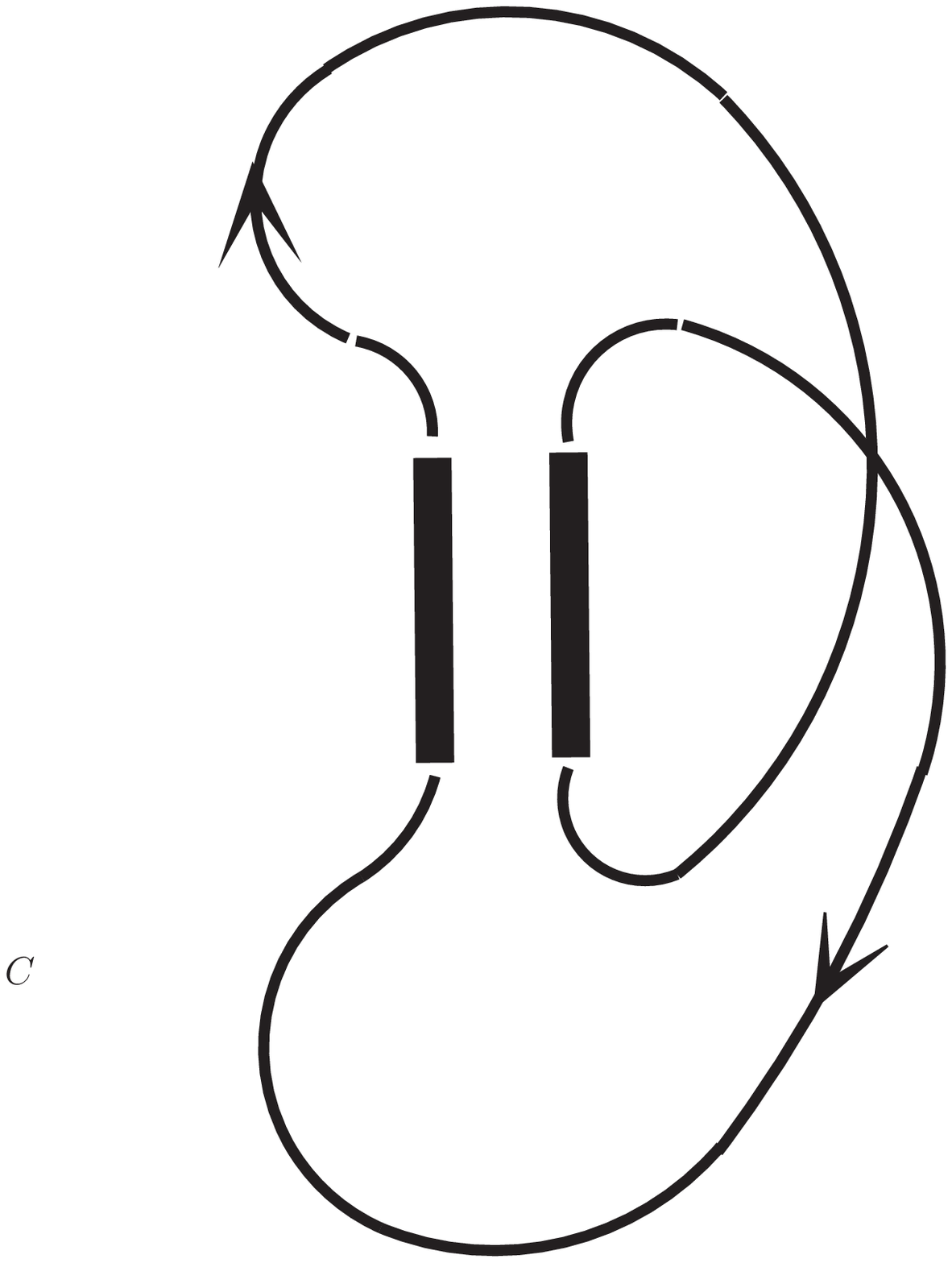}}
\qquad , \qquad
\parbox{2.6cm}{\psfrag{1}{$\cC_1$}\psfrag{2}{$\cC_2$}
\includegraphics[width=2.6cm]{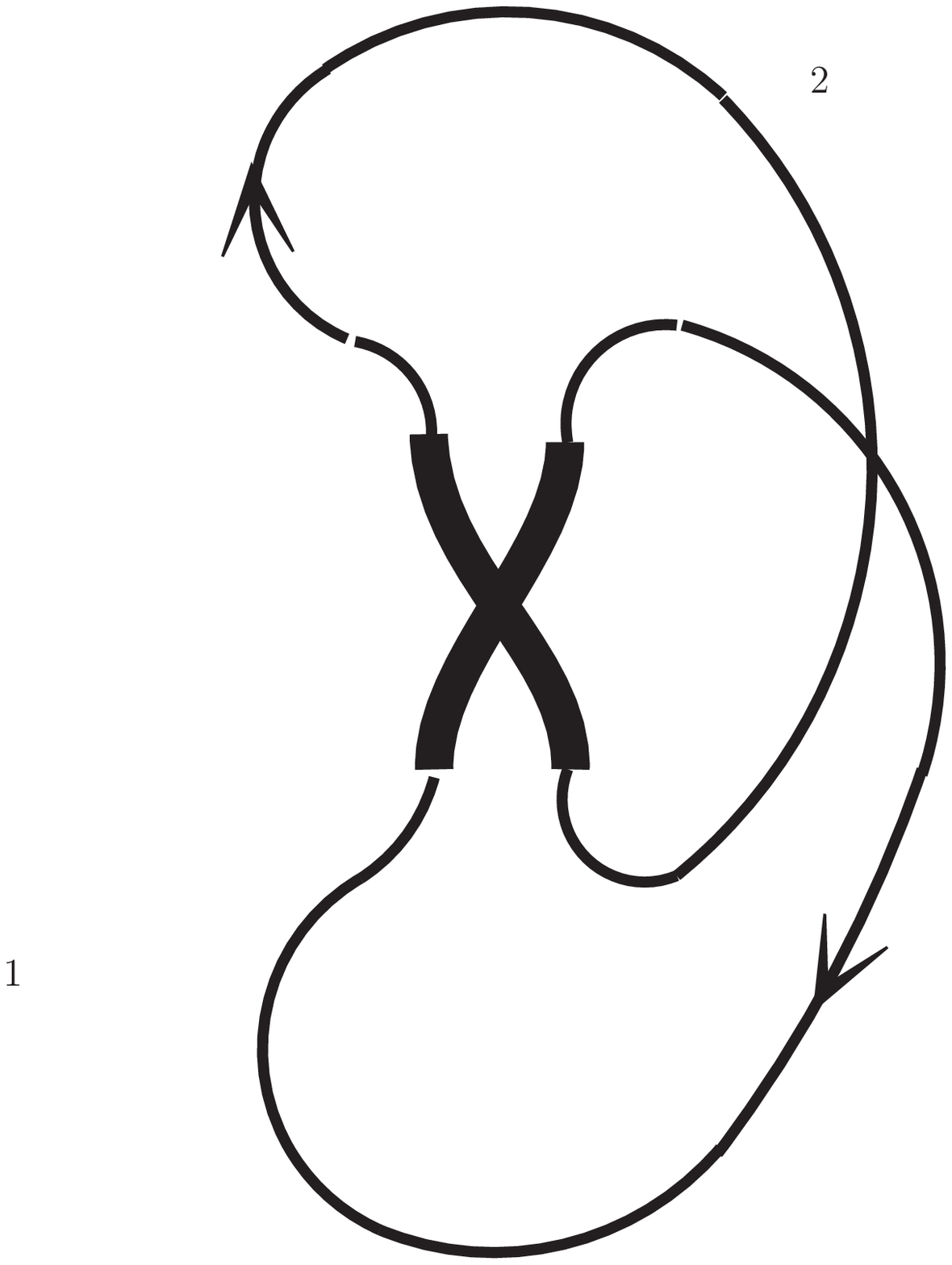}}\ \ .
\]
One then has $N(\vec{\si})\rightarrow N(\vec{\si})+1$
and $(-1)^{B_{+}(\vec{\si})}$ remains unchanged.
This concludes the proof of the claim.

As a result, one has
\[
\<\Ga,\ga\>^P=\sum_{\vec{\si}} (-1)^{C(\vec{\si})} (-2)^{N(\vec{\si})}
\]
\[
= \tilde{\mu} \sum_{\vec{\si}} (-1)^{B_{+}(\vec{\si})} 2^{N(\vec{\si})}\ \ .
\]
Now by the rules of Def. \ref{CGevaldef}
\[
\<G,\cO,\ta,\ga\>^{CG}=\left(
\prod\limits_{e\in E(G)} \ga(e)!
\right)^{-1}\times
\sum_{\si} \prod_{\cC} {\rm tr}\ \cC
\]
where the contribution ${tr}\ \cC$
of each curve is the trace of a product of $2\times 2$ matrices
taken among the commuting matrices $I,\ep,-\ep$.
More precisely
\[
{\rm tr}\ \cC={\rm tr}\left[
\ep^{B_{+}(\cC)} (-\ep)^{B_{-}(\cC)}
\right]
\]
\[
= (-1)^{B_{-}(\cC)}\ {\rm tr}\left[
\ep^{B(\cC)}\right]=(-1)^{B_{+}(\cC)}\ {\rm tr}\left[(-I)^{\frac{B(\cC)}{2}}\right]
\]
since $B(\cC)$ is even by the key observation.
Therefore ${\rm tr}\ \cC=(-1)^{B_{+}(\cC)}\times 2\times (-1)^{\frac{B(\cC)}{2}}$
and
\[
\<\Ga,\ga\>^P=\mu\times
\left(
\prod\limits_{e\in E(G)} \ga(e)!
\right)\times \<G,\cO,\ta,\ga\>^{CG}
\]
with $\mu=\tilde{\mu}\times (-1)^{\frac{k}{2}}$
where $k$ is the sum over curves $\cC$ of the number of
gate crossings $B(\cC)$.
This is the same as the total number of epsilon arrows
in the microscopic Feynman diagram used to evaluate the CG network.
\qed

\begin{Remark}
A byproduct of this proof is that the number $k$
of $\ep$'s is even.
Another way to see this is to write $k=k_\pi+k_\io$
where
$k_\pi$ is the total number of epsilons coming from vertices
\[
\parbox{1.4cm}{\psfrag{a}{$\scriptstyle{a}$}
\psfrag{b}{$\scriptstyle{b}$}\psfrag{c}{$\scriptstyle{c}$}
\includegraphics[width=1.4cm]{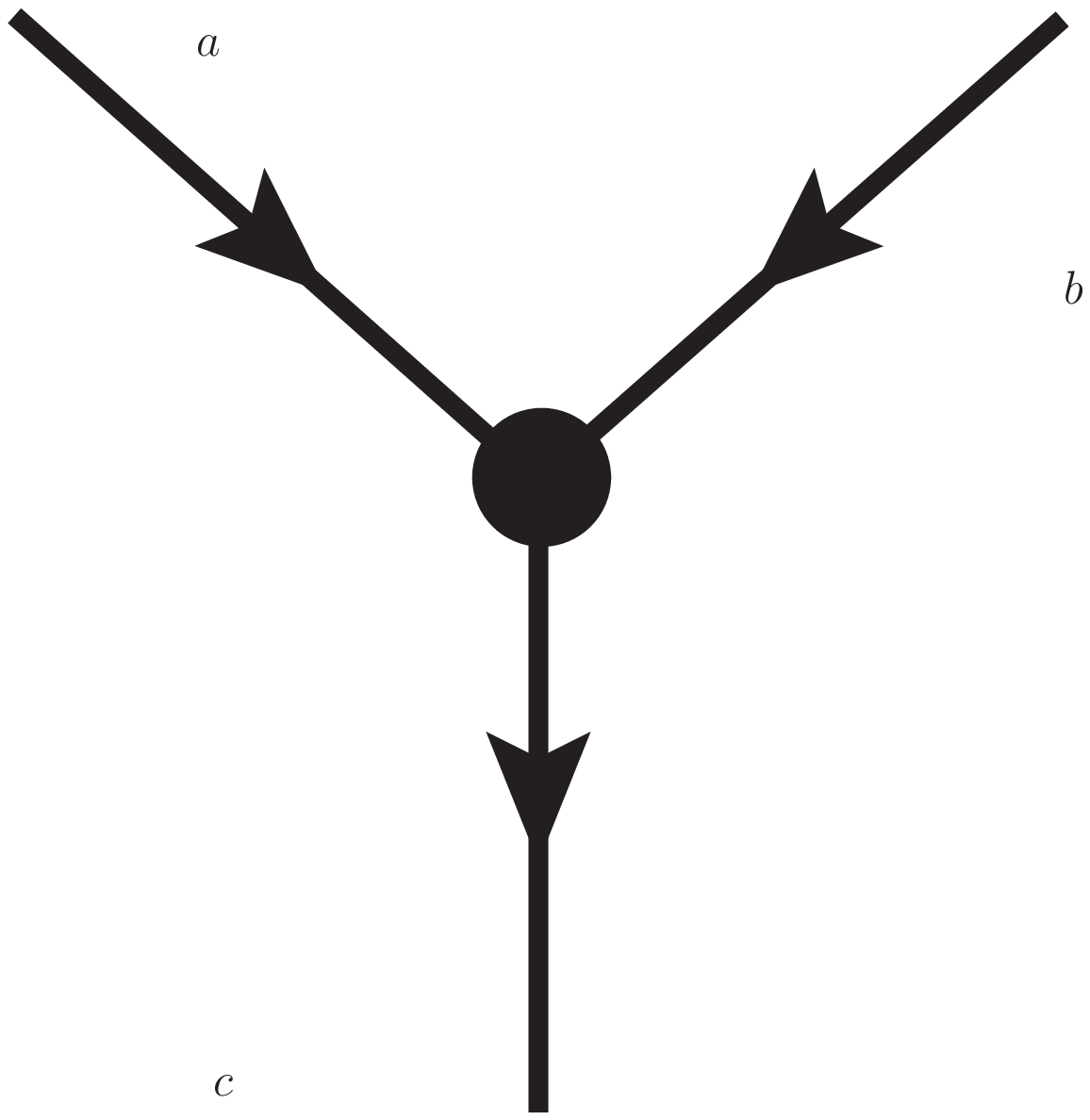}}
\]
and $k_\io$ is the total number of epsilons coming
from vertices
\[
\parbox{1.4cm}{\psfrag{a}{$\scriptstyle{a}$}
\psfrag{b}{$\scriptstyle{b}$}\psfrag{c}{$\scriptstyle{c}$}
\includegraphics[width=1.4cm]{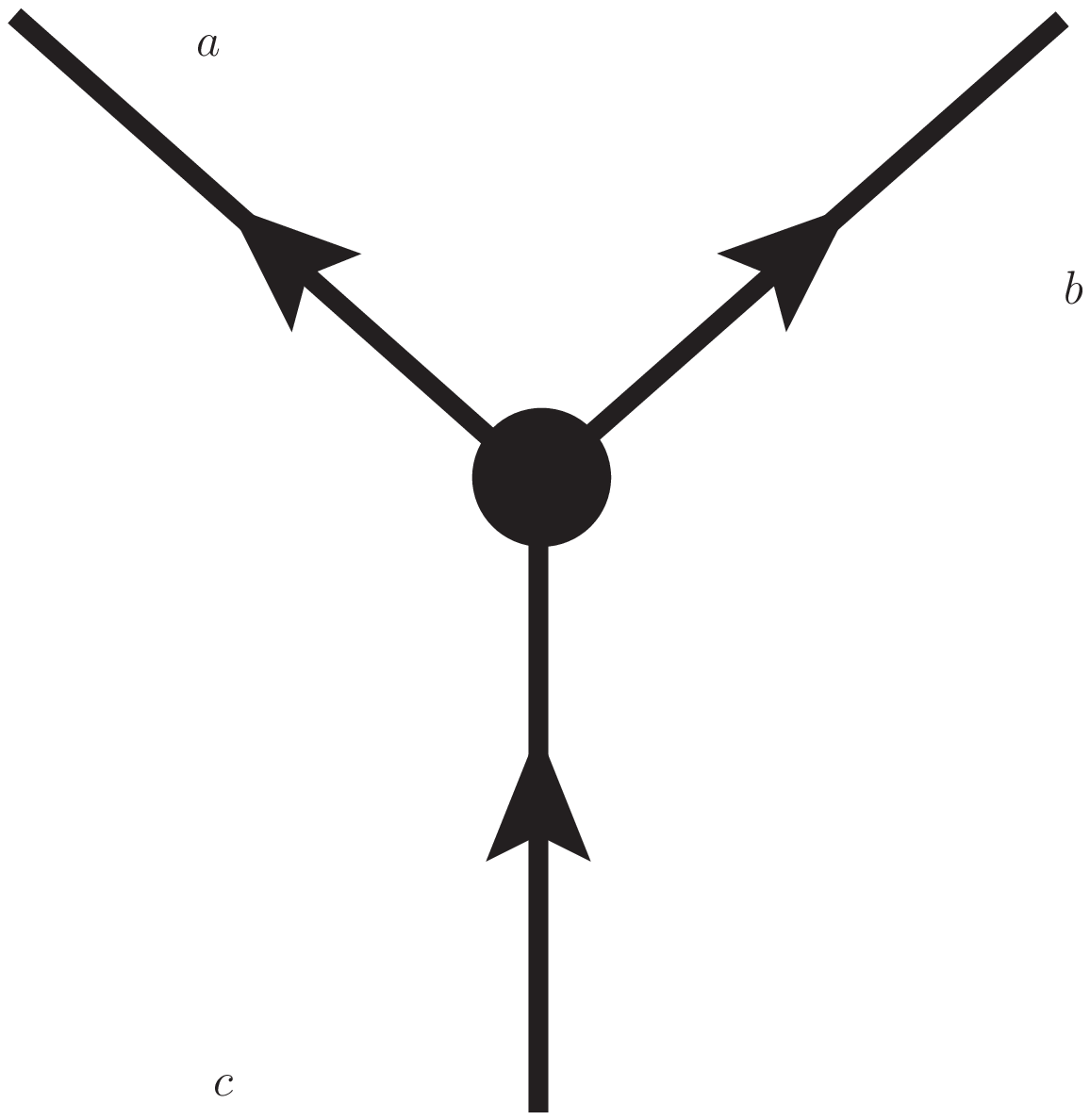}}\ \ .
\]
By using the $\frac{a+b-c}{2}$ formula to count the $\ep$'s
at these
vertices one easily finds that $k_\pi-k_\io=0$.
Indeed, an edge decoration is counted positively at the vertex of destination and negatively
at the vertex of origin, with respect to the orientation $\cO$.
\end{Remark}

\begin{Remark}
One can use the same parity arguments on the number of gate crossings
in order to prove
Corollary \ref{oneleg} in a purely diagrammatic way, without the invocation
of Schur's Lemma and the irreducibility of the $\cH_d$'s as $SU_2$
representations.
For a graph with only one external leg
\[
\parbox{2cm}{\psfrag{a}{$\scriptstyle{a}$}
\psfrag{1}{$\begin{array}{c}
{\scriptstyle\rm 1-valent}\\
{\scriptstyle\rm vertex}
\end{array}$}
\includegraphics[width=2cm]{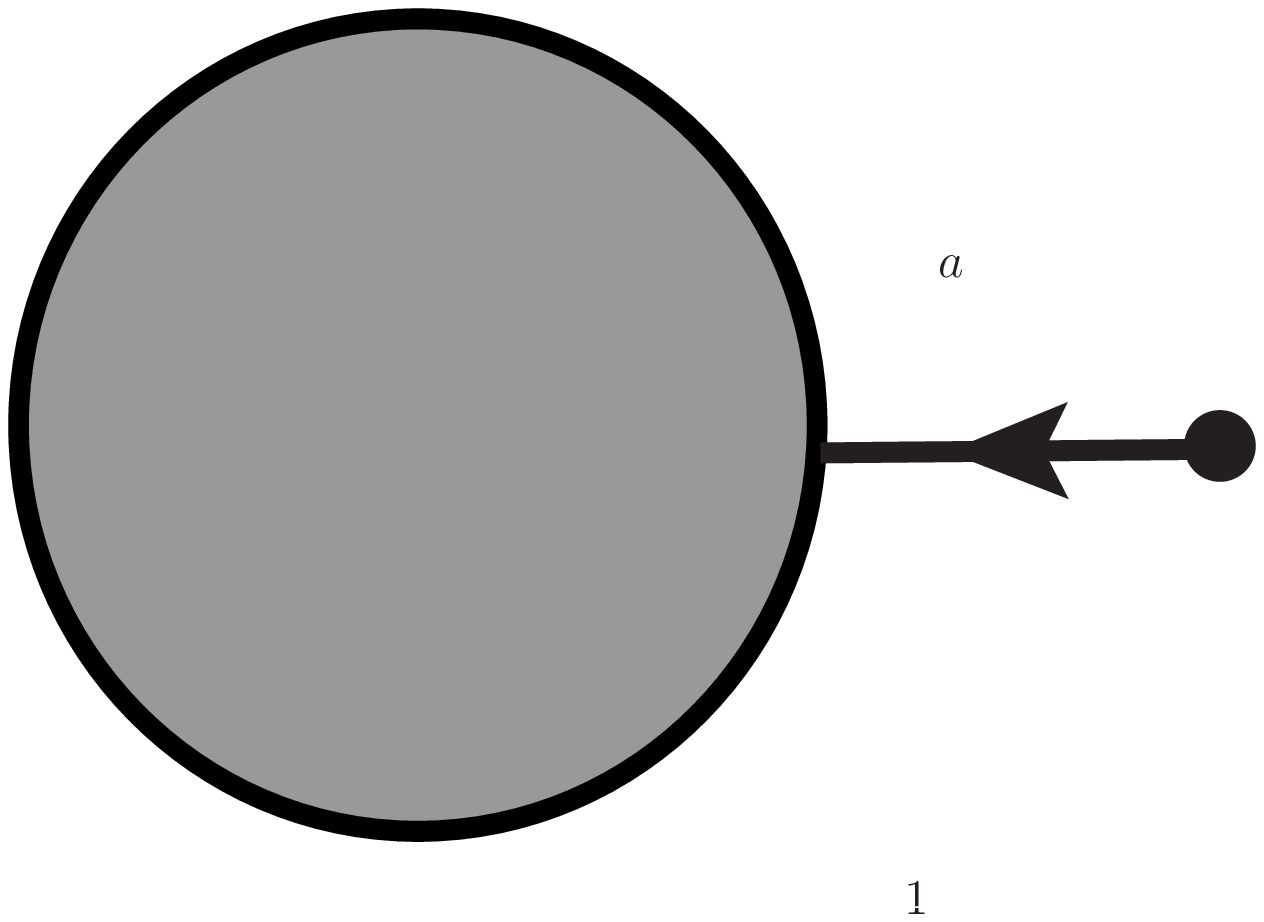}}
\qquad
\stackrel{{\rm FDC}}{\longrightarrow}
\qquad
\parbox{3cm}{\psfrag{1}{$\scriptstyle{i_1}$}
\psfrag{a}{$\scriptstyle{i_a}$}
\includegraphics[width=3cm]{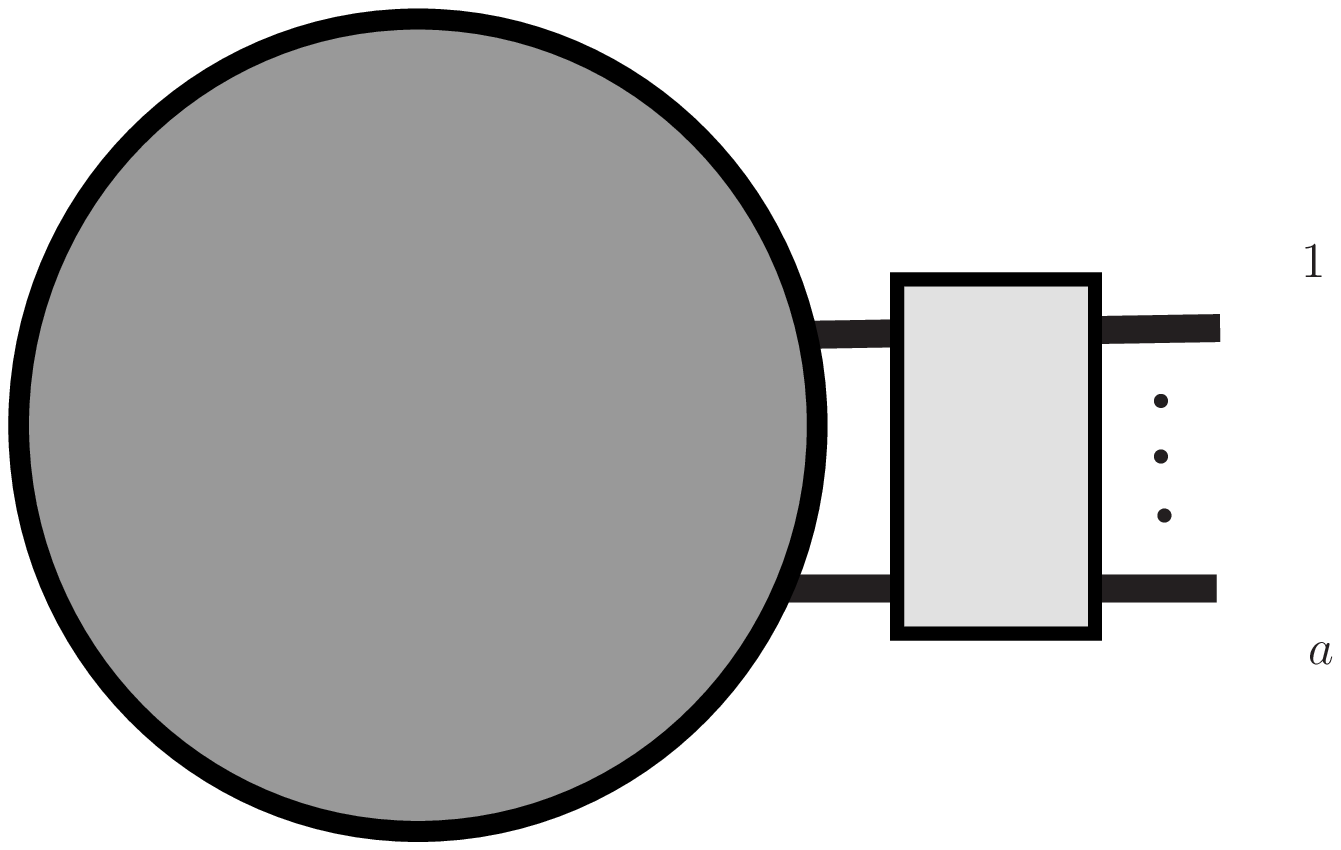}}
\]
\[
\ 
\]
any strand attached to an external index which goes in must come out.
Besides the number of gates crossed, i.e., the number of $\ep$ arrows
along the way
must be odd, as in Case II above.
Therefore, the contribution $\pm \ep^{2p+1}=\pm\ep$
is antisymmetric and is killed by the symmetrization over the indices
$i_1,\ldots,i_a$.
Note that one can also show Schur's Lemma in the same manner since
a strand that comes in one way must come out the other way. 
\end{Remark}

We will later show in \S\ref{smoothsec}
that there always exist smooth orientations $\cO$
and, consequently, that any CSN is amenable to the evaluation of a CG network.
In order to make the dictionary between CSN's and CG networks more
satisfactory, it would be desirable to have a canonical way to associate
a smooth orientation to a surface imbedding
$\Ga\stackrel{R}{\rightarrow} \Si$. Perhaps the methods of~\cite{Mohar}
could help towards that goal.
By contrast, a possible canonical choice for the gate signage $\ta$
is trivially obtained by orienting
the small curved arrows counterclockwise along
the faces of the map as in
\[
\parbox{4cm}{\includegraphics[width=4cm]{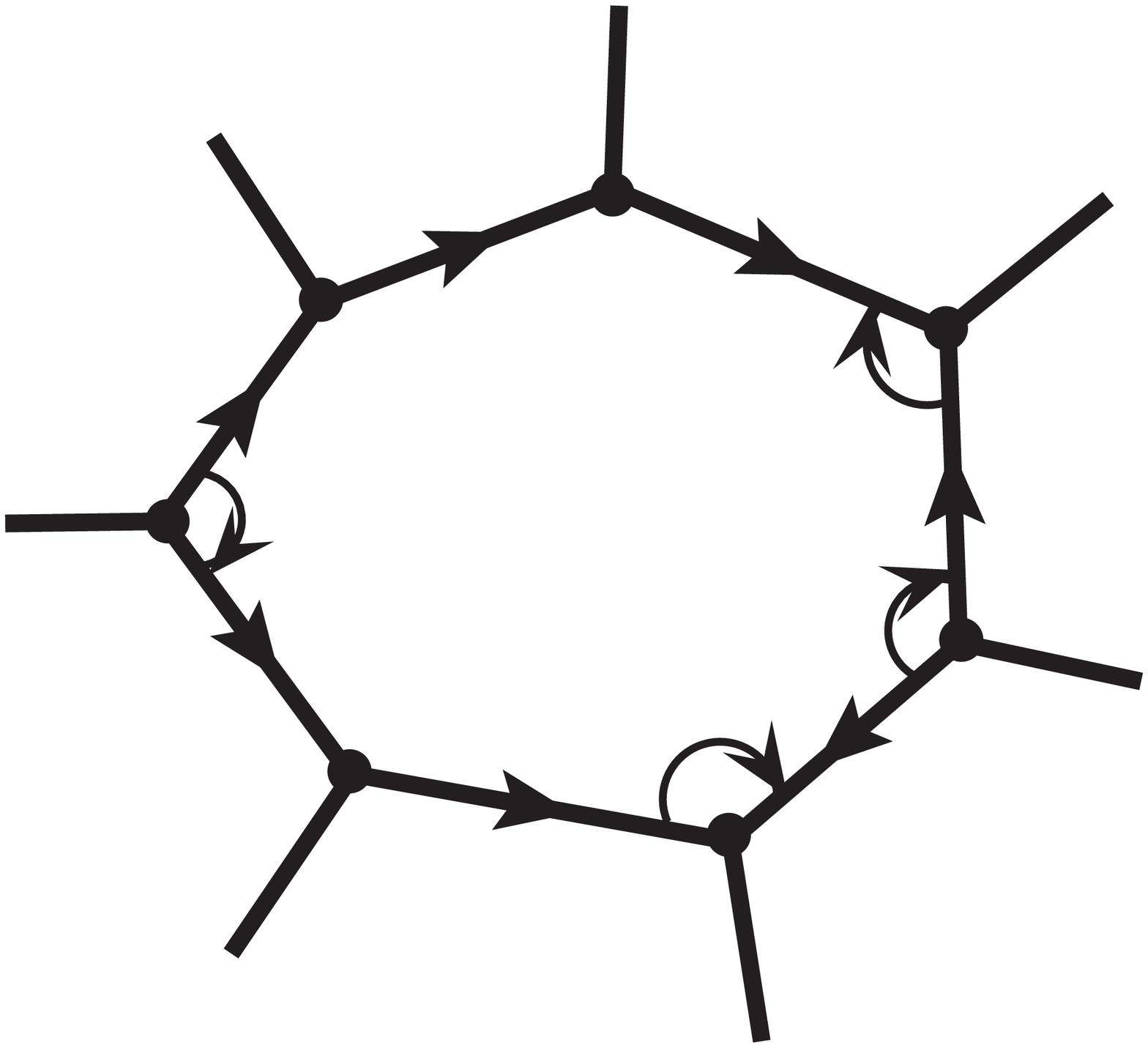}}\ \ .
\]
It would also be nice to be able to compute the sign $\mu$ for such
canonical prescriptions.

\begin{Remark}
Note that an alternative to the proof of Theorem \ref{negdimthm} is that sketched
in~\cite[\S2.2]{Roberts1}. Up to a harmless global
prefactor, the setting of~\cite[Def. 3]{Roberts1}
corresponds to a CG network evaluated by FDC as in Def. \ref{CGevaldef}
but using
$\ep$ arrows and no Kronecker delta's. One can connect Roberts
definition to ours by pushing the edge epsilons into the vertices,
in the direction
indicated by the orientation $\cO$, and eliminating those
which can be eliminated by the relation $\ep^2=-{\rm Id}$.
\end{Remark}

We now define an alternative evaluation for CG networks
which is closely related to the unitary evaluation of CSN's.
Given a CG network $(G,\cO,\ta,\ga)$,
possibly with external legs, we let the $\pi$\&$\io$
evaluation $\<G,\cO,\ta,\ga\>^{\pi\io}$
be defined by the same rules as in \S\ref{legsection}
except that the contribution of a vertex
\[
\parbox{1.3cm}{\psfrag{a}{$\scriptstyle{a}$}
\psfrag{b}{$\scriptstyle{b}$}
\psfrag{c}{$\scriptstyle{c}$}
\includegraphics[width=1.3cm]{Fig102.eps}}
\qquad {\rm or}\qquad
\parbox{1.3cm}{\psfrag{a}{$\scriptstyle{a}$}
\psfrag{b}{$\scriptstyle{b}$}
\psfrag{c}{$\scriptstyle{c}$}
\includegraphics[width=1.3cm]{Fig103.eps}}
\]
gets multiplied by the factor
\[
\sqrt{\frac{a!\ b!\ (c+1)!}{\left(\frac{a+b+c}{2}+1\right)!\ 
\left(\frac{a+b-c}{2}\right)!\ 
\left(\frac{a+c-b}{2}\right)!\ 
\left(\frac{b+c-a}{2}\right)!}}
\]
dictated by the considerations preceeding Proposition \ref{piiotaprop}.
For such a vertex $v$ we define ${\rm dim}(v)=c+1$.
Now an immediate corollary of Theorem \ref{negdimthm} is as follows.

\begin{Corollary}\label{negdimU}
With the same hypotheses and notations as in Theorem \ref{negdimthm}
one has
\[
\<\Ga,\ga\>^U=\mu \times
\left(\prod_{v\in V(G)} \frac{1}{\sqrt{{\rm dim}(v)}}\right)
\times
\<G,\cO,\ta,\ga\>^{\pi\io}\ \ .
\]
\end{Corollary}

\section{The existence of smooth orientations}\label{smoothsec}

The following proposition garantees that one can always use the CG formalism
in order to evaluate a spin network.

\begin{Proposition}\label{smoothprop}
For any spin network $(\Ga,\ga)$ without trivial components,
the underlying graph
has a smooth orientation.
\end{Proposition}

We can of course assume $\Ga=(G,R)$
where $G$ is a connected cubic graph, possibly containing loops
and multiple edges. For good measure, we will provide several proofs.

\bigskip
\noindent{\bf 1st proof for the case of bridgeless graphs:}
(Indicated to us by Bill Jackson and Gordon Royle)
If $g$ is bridgeless or 2-edge-connected, one can use Petersen's famous
CIT-inspired 1-factor Theorem~\cite{Petersen} (see also~\cite{Sabidussi}).
Namely, there exists a 1-factor of $G$, i.e., a set of edges $E_1\subset E(G)$
such that every vertex $v\in V(G)$
has degree 1 in the spanning subgraph given by $E_1$.
The edge complement is therefore a collection
of cycles.
For each such cycle, choose a coherent orientation, i.e.,
an orientation which follows a direction of travel along
the cycle. Finally, take any orientation of the edges in $E_1$. The obtained
orientation $\cO$ is smooth.
\qed

\bigskip
\noindent{\bf 2nd proof in the general case:}
Choose a spanning tree $T$ in the connected graph $G$.
Choose $v_0$ among the leafs of the tree, i.e., vertices of degree 1
relatively to $T$.
Choose $v_0$ as a root for $T$, and orient all the edges of $T$ towards
the root $v_0$.
For example:
\[
\parbox{6cm}{\psfrag{G}{$G$}
\includegraphics[width=4.2cm]{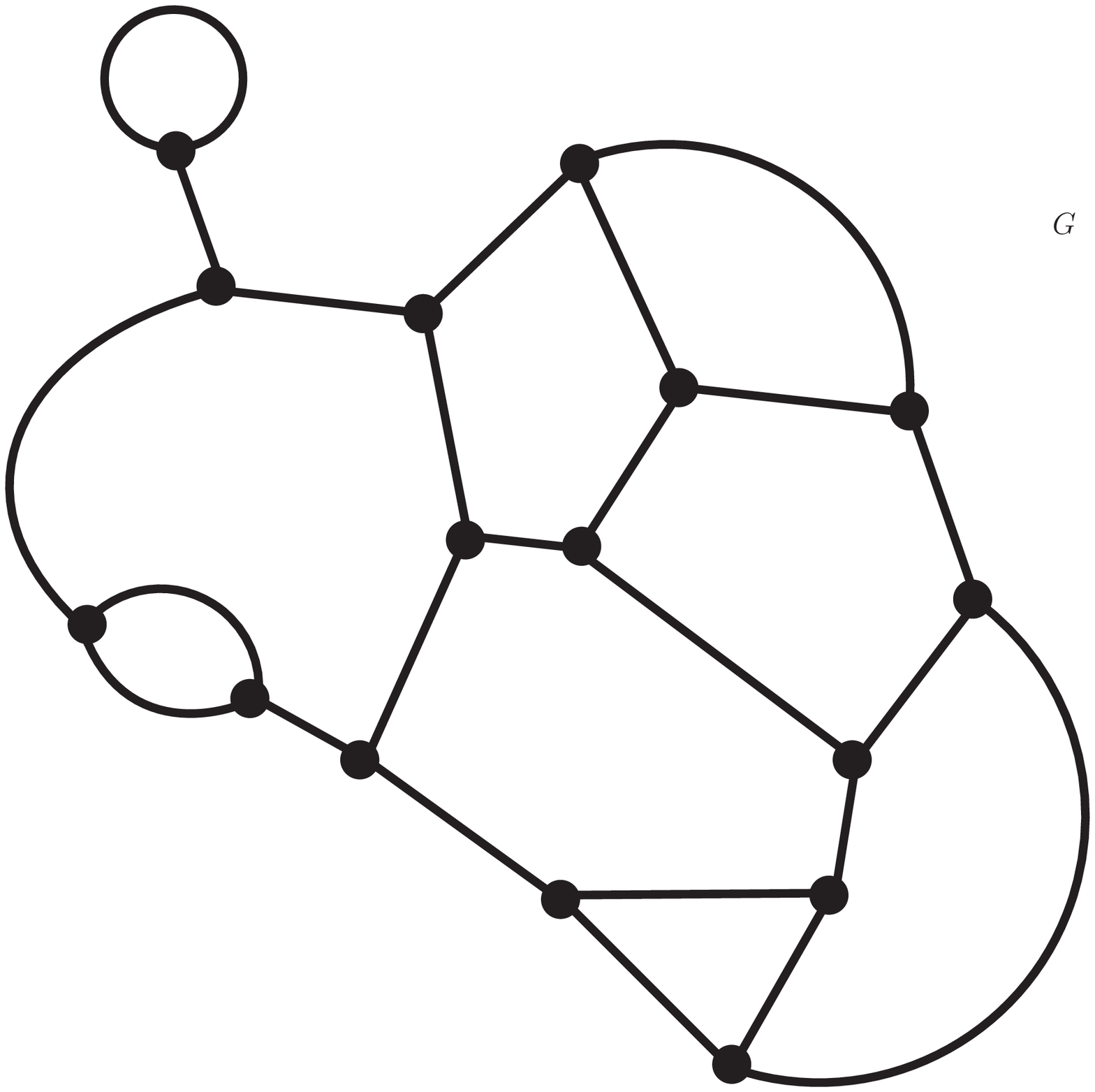}}
\qquad\longrightarrow\qquad
\parbox{6cm}{\psfrag{v}{$v_0$}
\includegraphics[width=4.2cm]{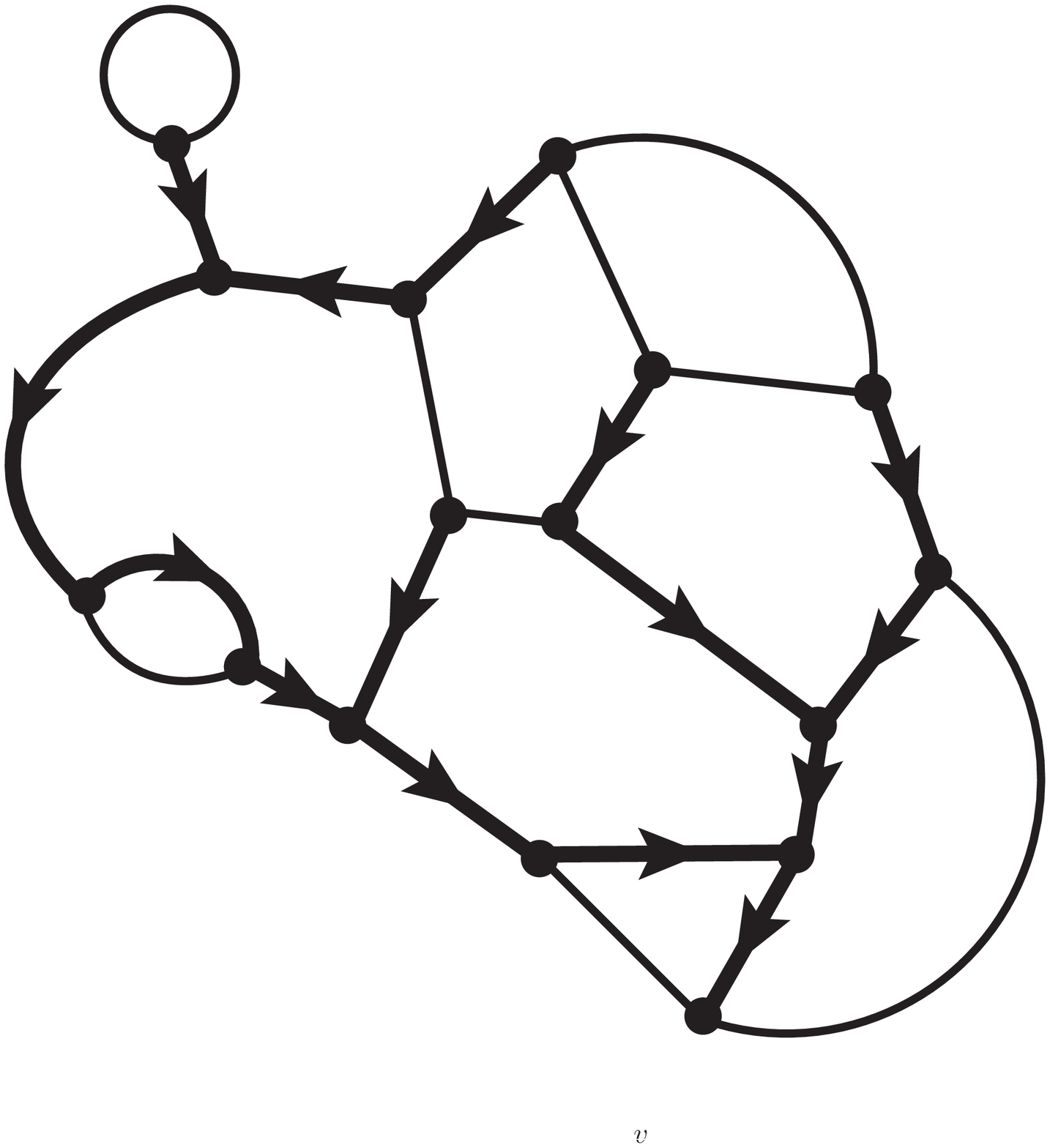}}
\]
where the edges of $T$ are indicated by thicker lines.
Because we have chosen a leaf to be the root,
the vertices of degree 3 in the tree are of the form $\parbox{0.7cm}{
\includegraphics[width=0.7cm]{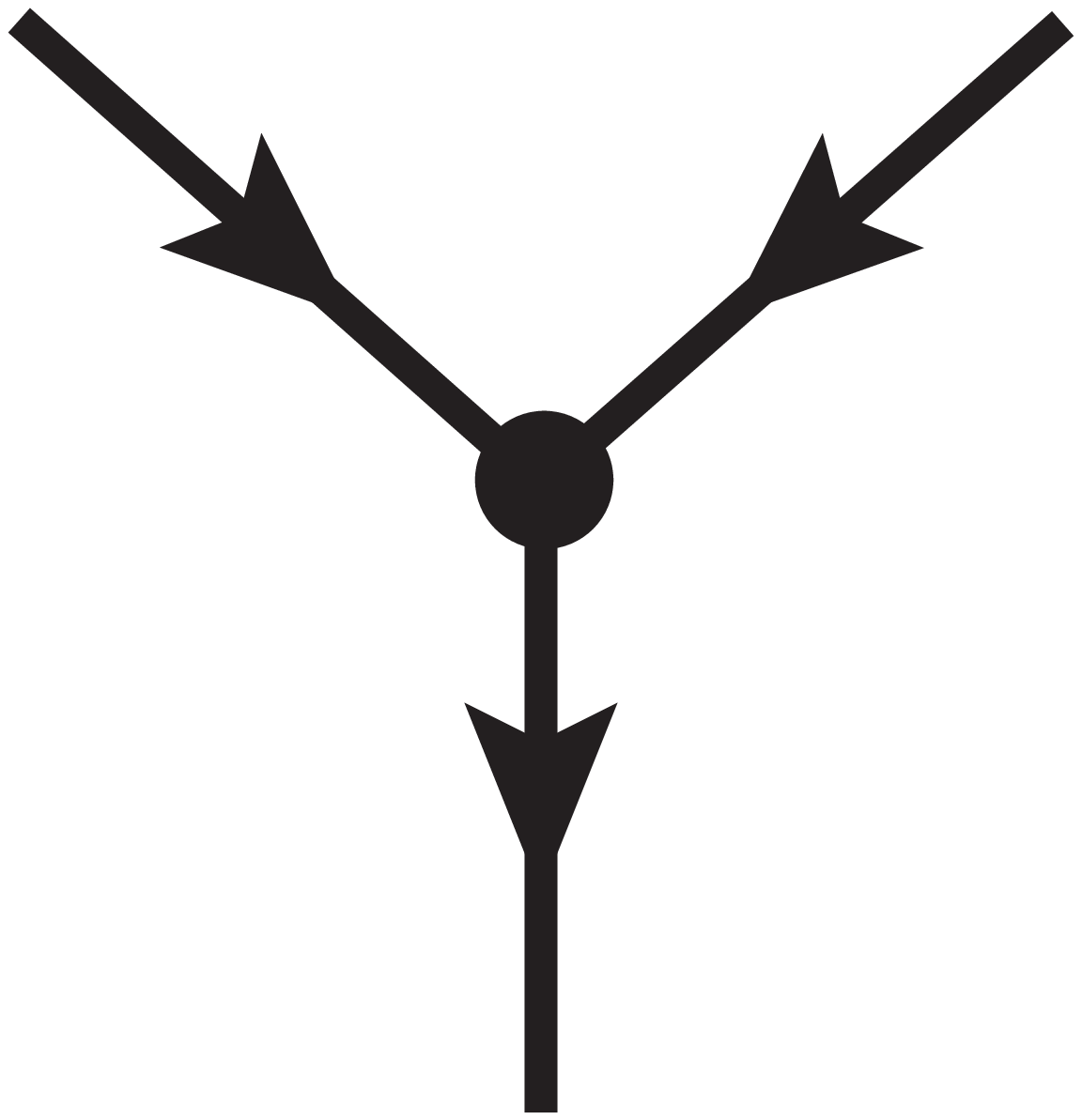}}$
which is acceptable.
The vertices of degree 2 in $T$ must be of the form $\parbox{0.7cm}{
\includegraphics[width=0.7cm]{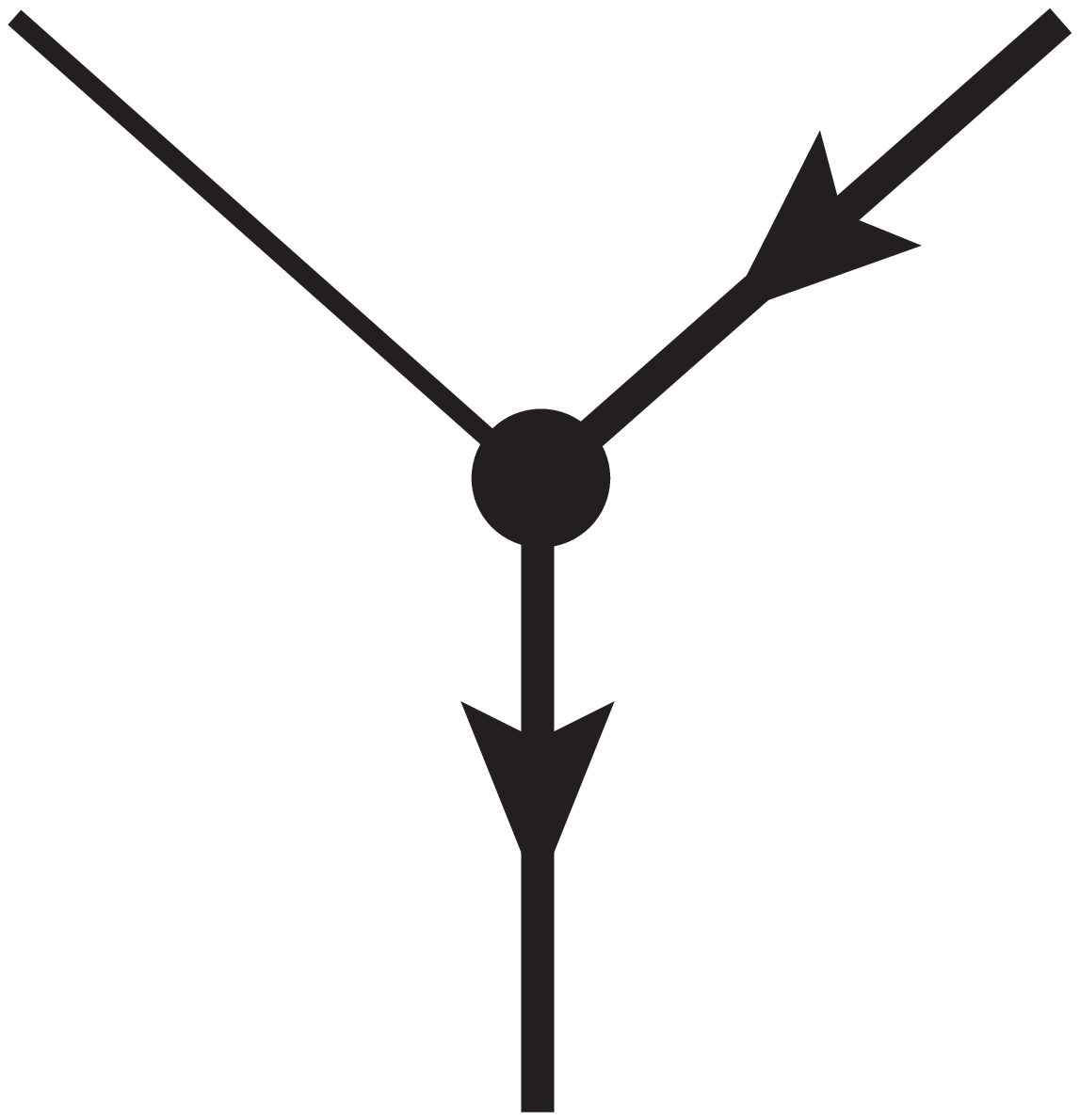}}$.
Clearly, the remaining edges can be organized
into a collection of paths and cycles (which can be reduced to a loop).
For the cycles which bounce around vertices of degree 1 in $T$ and
then close upon themselves, choose a coherent orientation following
a direction of travel along the cycle.
For a path which must start and end at two vertices of degree 2 in $T$,
and, in between, may bounce around vertices of degree 1 in $T$,
again choose a coherent orientation along the path.
The resulting orientation $\cO$ is easily seen to be smooth.\qed

For the previous example a possible outcome of this procedure is
\[
\parbox{5cm}{
\includegraphics[width=5cm]{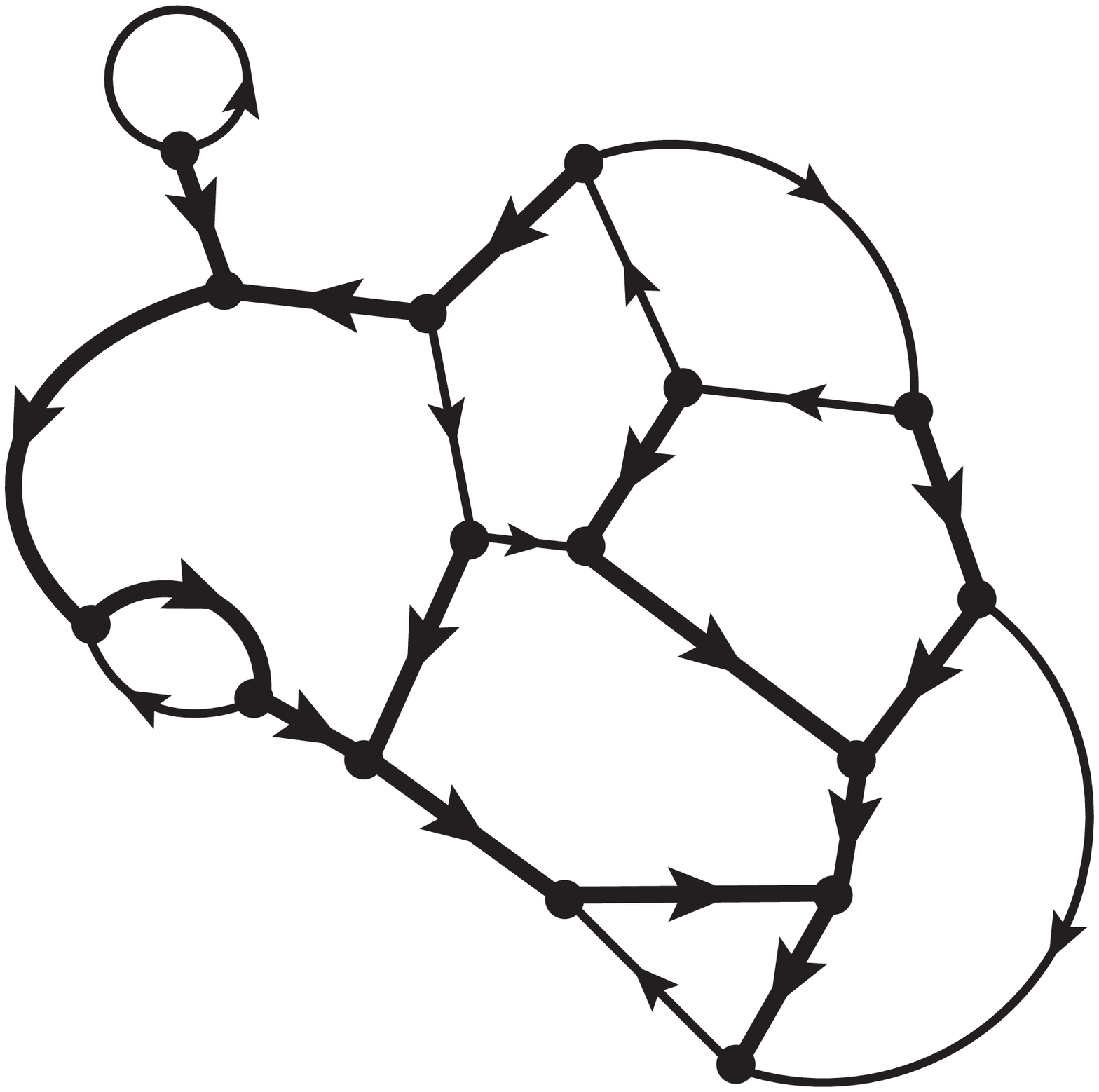}}\qquad .
\]

\noindent{\bf 3rd proof in the general case:} (Indicated to us by Bill Jackson)
We will show that all graphs $G$, and not only cubic graphs, have an
orientation such that, for any vertex, the absolute value of the difference
between indegree and outdegree is most one.
The proof is by induction on the number of edges.
If some edge $e$ is incident
to a vertex of odd degree, then by the induction hypothesis
$G-e$ (i.e., the graph with edge $e$ removed)
has such an orientation $\cO'$.
If $e$ is a loop, then any orientation of $e$ will do.
Otherwise, $e$ is incident to two distinct vertices
$v_1,v_2$
\[
\parbox{5cm}{\psfrag{e}{$\scriptstyle{e}$}
\psfrag{1}{$\scriptstyle{v_1}$}\psfrag{2}{$\scriptstyle{v_2}$}
\psfrag{a}{$\scriptstyle{a}$}\psfrag{b}{$\scriptstyle{b}$}
\psfrag{c}{$\scriptstyle{c}$}\psfrag{d}{$\scriptstyle{d}$}
\includegraphics[width=5cm]{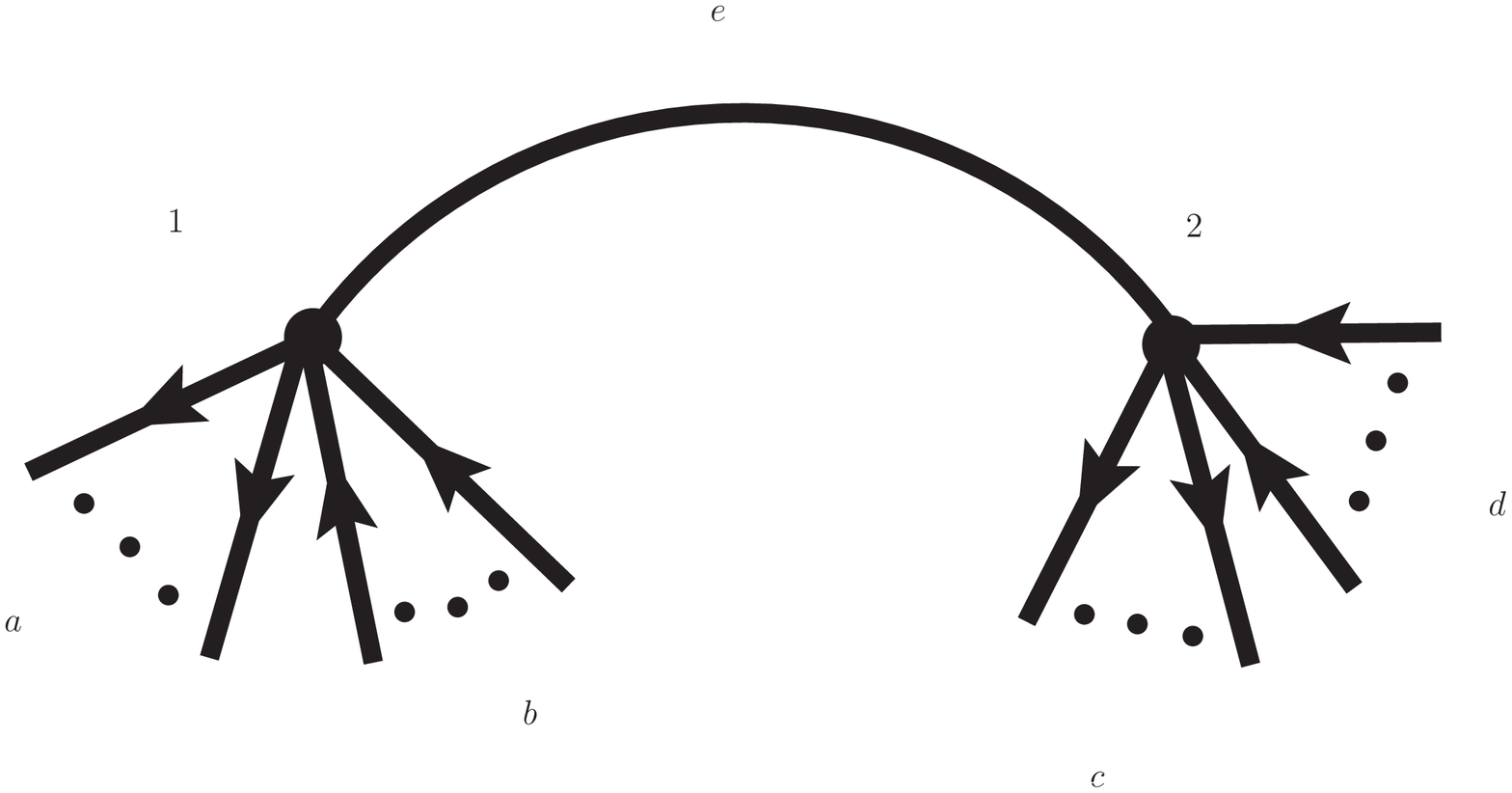}}
\]
such that say $v_1$
is of odd degree in $G$. The orientations indicated
are the ones provided by $\cO'$.
Since $|a-b|\le 1$ and $a+b+1$ is odd, we must have $a=b$.
If $c=d$, any orientation of $e$ will do.
If $c=d-1$,
then orient $e$ from $v_2$ to $v_1$.
If $c=d+1$, then orient $e$ from $v_1$ to $v_2$.
The resulting orientation is then acceptable.

If no edge is incident to a vertex of odd degree, then each connected
component of $G$ is Eulerian and one can direct
the edges around
an Euler tour of each connected component.
The resulting orientation will have indegree$=$outdegree
at each vertex.
\qed

\begin{Remark}
Note that one can find a blend of Proof 1 and Proof 3
in~\cite[Prop 3.4]{Mohar},
concerning the notion of `NS orientations'.
These are synonymous with smooth orientations, as well as orientations
without sources and sinks,
in the case of cubic graphs.
\end{Remark}

\section{The bridge reduction}
Let $(\Ga,\ga)$
be a connected CSN which is not reduced to a trivial component.
Let us also assume there is a cut-edge or bridge $e_0$, i.e.,
an edge which when removed makes the graph disconnected.
\begin{Lemma}
If $\ga(e_0)\neq 0$ then $\<\Ga,\ga\>^P=0$.
\end{Lemma}
\noindent{\bf Proof:}
Choose a smooth orientation $\cO$ on the underlying graph
$G$ by Proposition\ref{smoothprop}.
Choose a gate signage $\ta$.
By Theorem \ref{negdimthm}, we only need to show that
$\<G,\cO,\ta,\ga\>^{CG}=0$. Let $i_1,\ldots,i_d$, with $d=\ga(e_0)\neq 0$,
be the indices summed over at the strands of the
Feynman diagram used to compute
$\<G,\cO,\ta,\ga\>^{CG}$.
By Corollary \ref{oneleg}, the contribution of the CG networks
with one external leg on either side of $e_0$
must vanish, for any choice of the indices $i_1,\ldots,i_d$
in $\{1,2\}$.
\qed

We now assume $\ga(e_0)=0$.
This forces the decorations at the endpoints of $e_0$
to be of the form
\[
\parbox{1.4cm}{
\psfrag{a}{$\scriptstyle{a}$}
\psfrag{b}{$\scriptstyle{b}$}
\psfrag{0}{$\scriptstyle{0}$}
\psfrag{e}{$\scriptstyle{e_0}$}
\includegraphics[width=1.4cm]{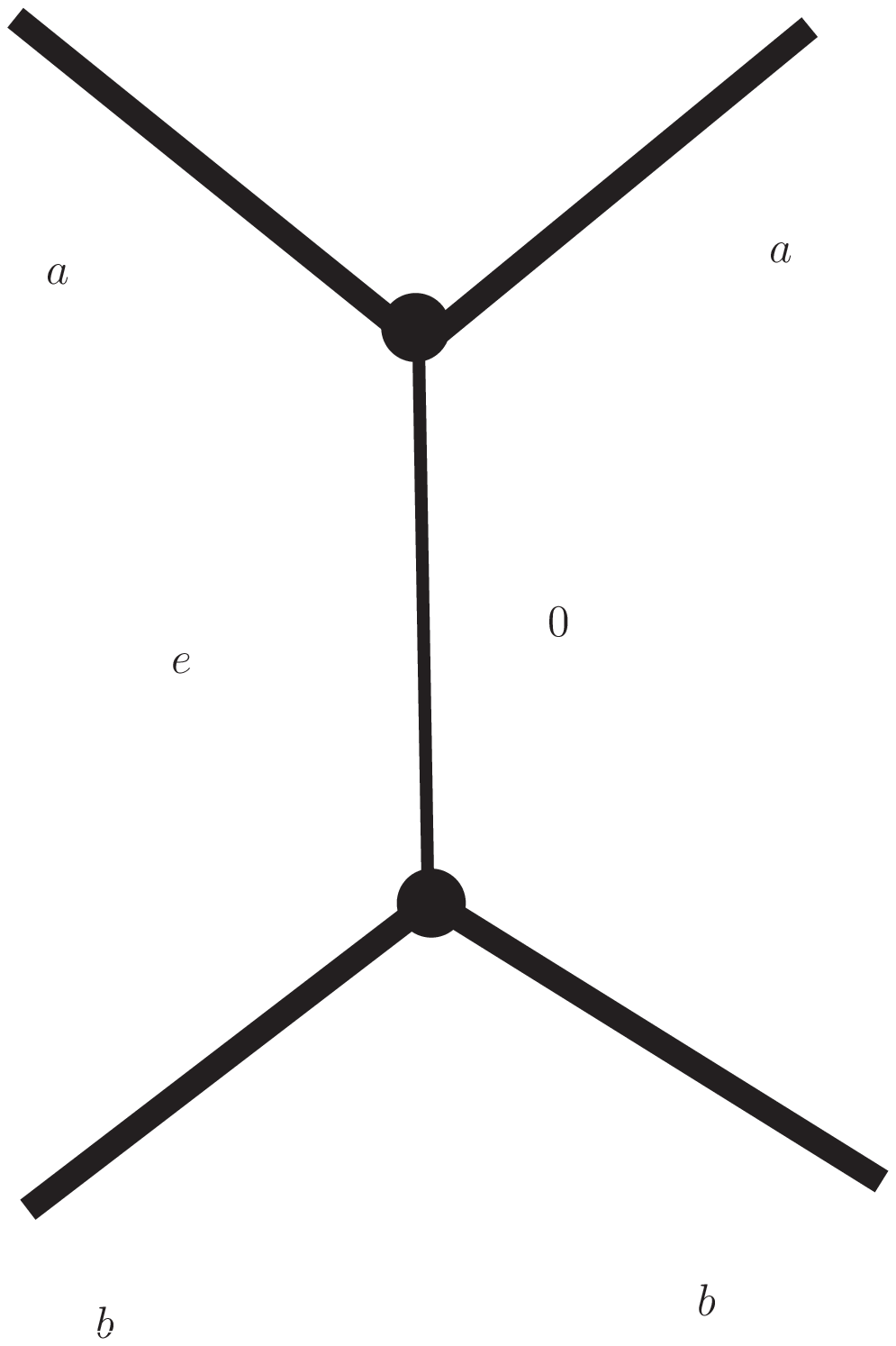}}
\]
by the admissibility condition.
Let $\Ga_1$, $\Ga_2$
be the connected ribbon graphs one obtains by removing
$e_0$
and erasing the vertices $v_1,v_2$ that $e_0$ was incident to
\[
\parbox{3.5cm}{
\psfrag{a}{$\scriptstyle{a_1}$}
\psfrag{b}{$\scriptstyle{a_2}$}
\psfrag{0}{$\scriptstyle{0}$}
\psfrag{1}{$\scriptstyle{v_1}$}
\psfrag{2}{$\scriptstyle{v_2}$}
\psfrag{G}{$\Ga$}
\includegraphics[width=3.5cm]{Fig173.eps}}
\qquad\longrightarrow\qquad
\parbox{3.5cm}{
\psfrag{a}{$\scriptstyle{a_1}$}
\psfrag{b}{$\scriptstyle{a_2}$}
\psfrag{1}{$\Ga_1$}
\psfrag{2}{$\Ga_2$}
\includegraphics[width=3.5cm]{Fig174.eps}}\qquad .
\]
Of course, $\Ga_1,\Ga_2$ inherit decorations $\ga_1,\ga_2$
coming from $\ga$.
\begin{Proposition}\label{bridgeprop}
\[
\<\Ga,\ga\>^U=\<\Ga_1,\ga_1\>^U\ \<\Ga_2,\ga_2\>^U\ 
\de(v_1)\de(v_2)
\]
where $\de(v)$ for a vertex of the form
\[
\parbox{1.5cm}{
\psfrag{a}{$\scriptstyle{a}$}
\psfrag{v}{$\scriptstyle{v}$}
\psfrag{0}{$\scriptstyle{0}$}
\includegraphics[width=1.5cm]{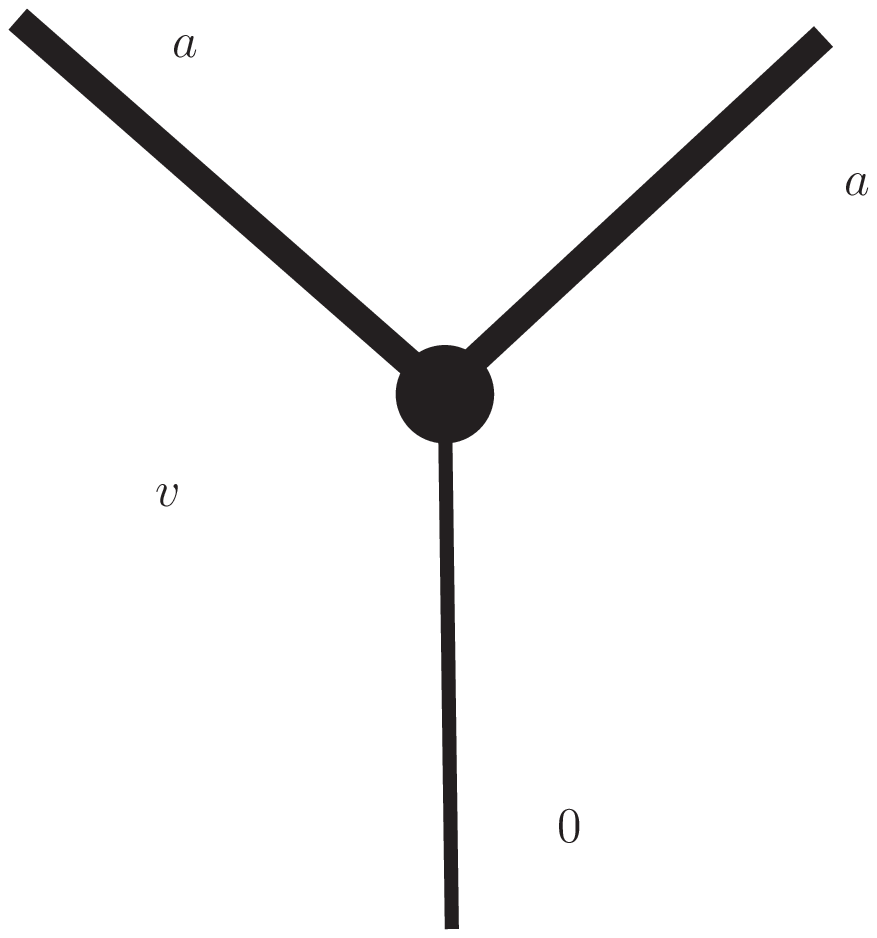}}
\]
is $\de(v)=\frac{1}{a!\sqrt{a+1}}$ if $v$ is a loop
vertex $\parbox{0.7cm}{
\psfrag{a}{$\scriptstyle{a}$}
\includegraphics[width=0.7cm]{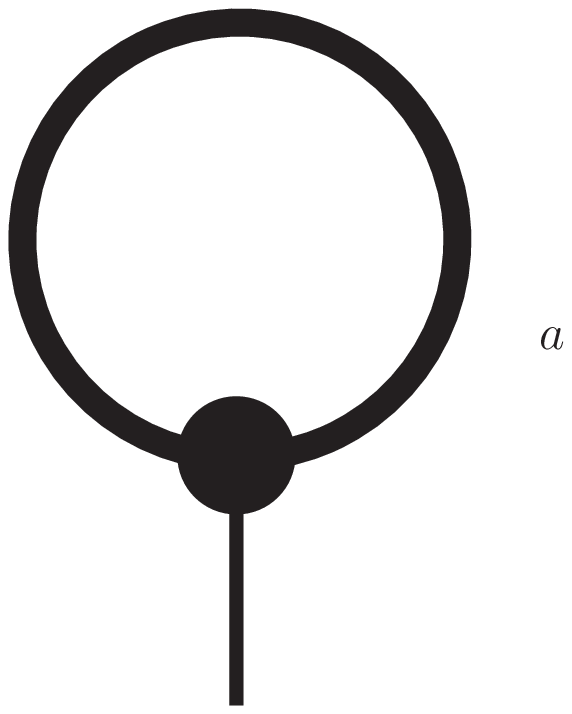}}$
and $\de(v)=\frac{1}{\sqrt{a+1}}$ otherwise.
\end{Proposition}
\noindent{\bf Proof:}
First let us cut the edge $e_0$ of the ribbon graph $\Ga$
by introducing two 1-valent vertices $v'_1$ and $v'_2$:
\[
\parbox{1.8cm}{
\psfrag{e}{$\scriptstyle{e_0}$}
\includegraphics[width=1.8cm]{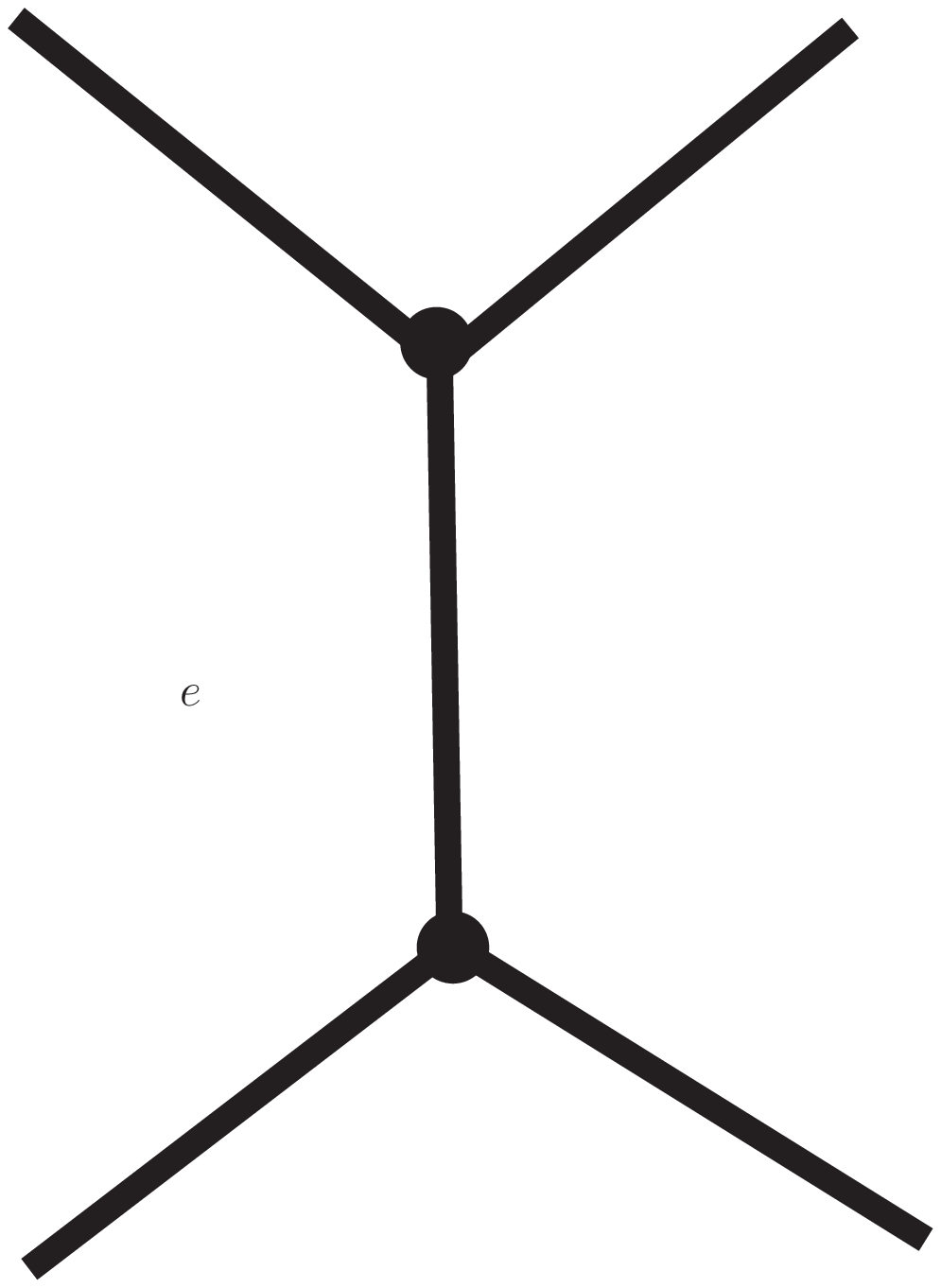}}
\qquad\longrightarrow\qquad
\parbox{1.8cm}{
\psfrag{1}{$\scriptstyle{v_1}$}
\psfrag{2}{$\scriptstyle{v_2}$}
\psfrag{3}{$\scriptstyle{v'_1}$}
\psfrag{4}{$\scriptstyle{v'_2}$}
\includegraphics[width=1.8cm]{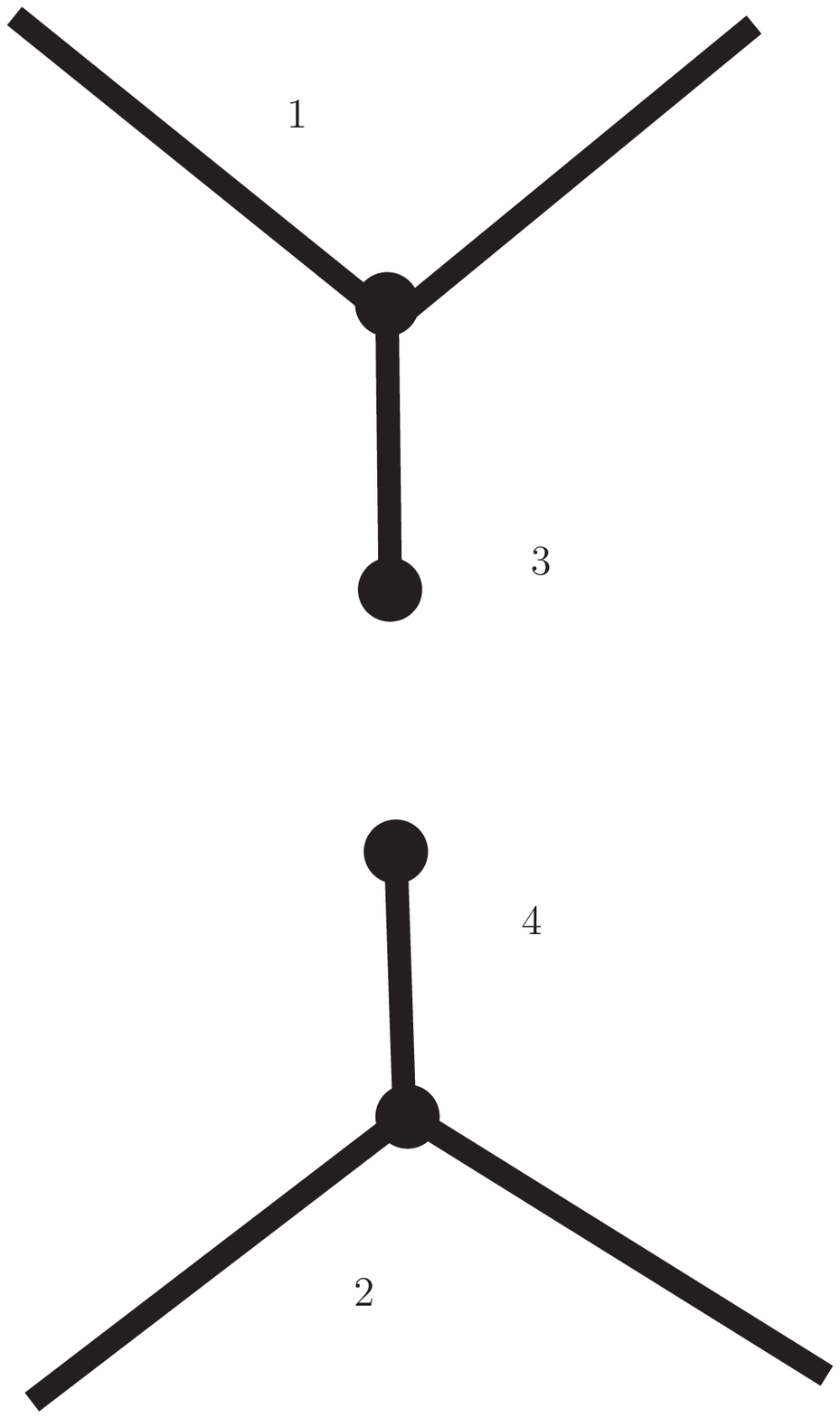}}\qquad .
\]
One therefore obtains two (nonregular) ribbon graphs $\hat{\Ga}_1$ and
$\hat{\Ga}_2$
together with their imbeddings in surfaces $\Si_1$ and $\Si_2$ respectively:
\[
\parbox{3.4cm}{
\psfrag{1}{$\scriptstyle{v_i}$}
\psfrag{2}{$\scriptstyle{v'_i}$}
\psfrag{S}{$\Si_i$}
\psfrag{G}{$\scriptstyle{\Ga_i}$}
\includegraphics[width=3.4cm]{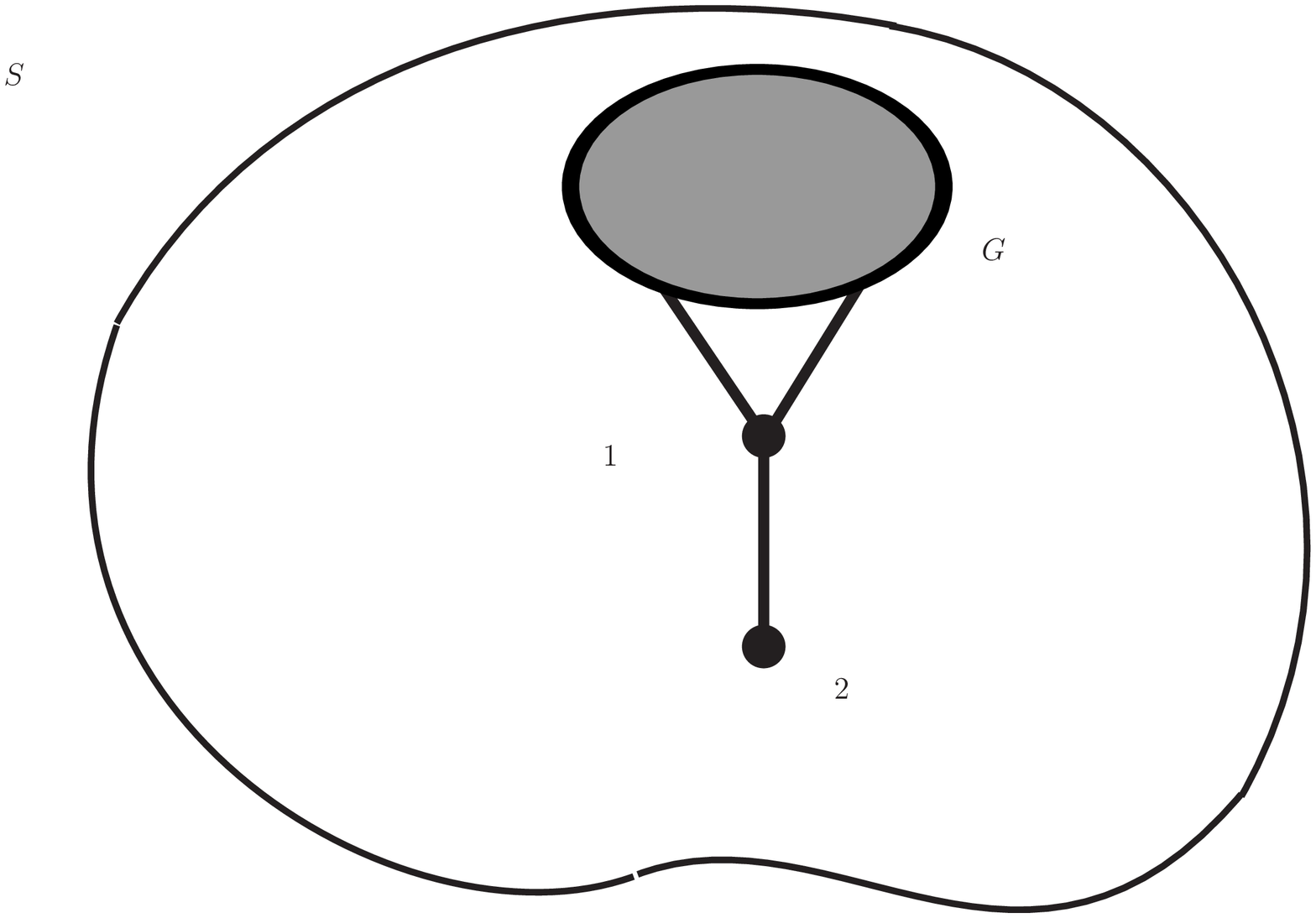}}\qquad .
\]
The shaded area represents the rest of the graph $\hat{\Ga}_i$
which can wrap around the surface $\Si_i$.
Then cut a small disc $D_i$ around $v'_i$.
\[
\parbox{3.4cm}{
\psfrag{1}{$\scriptstyle{v_i}$}
\psfrag{2}{$\scriptstyle{v'_i}$}
\psfrag{S}{$\Si_i$}
\psfrag{D}{$\scriptstyle{D_i}$}
\includegraphics[width=3.4cm]{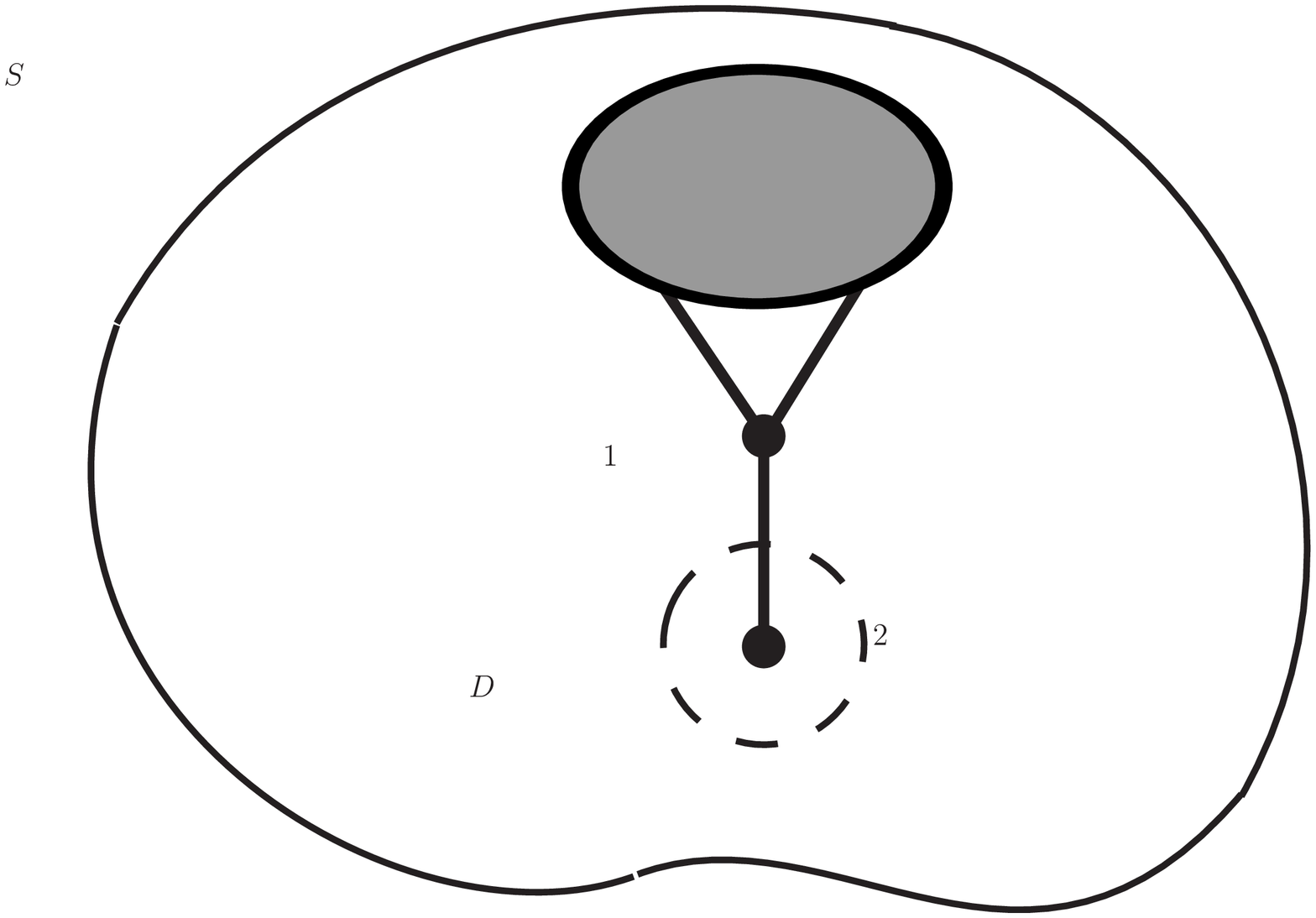}}
\qquad\longrightarrow\qquad
\parbox{3.4cm}{
\psfrag{v}{$\scriptstyle{v_i}$}
\psfrag{S}{$\Si_i$}
\psfrag{D}{$\scriptstyle{\partial D_i}$}
\includegraphics[width=3.4cm]{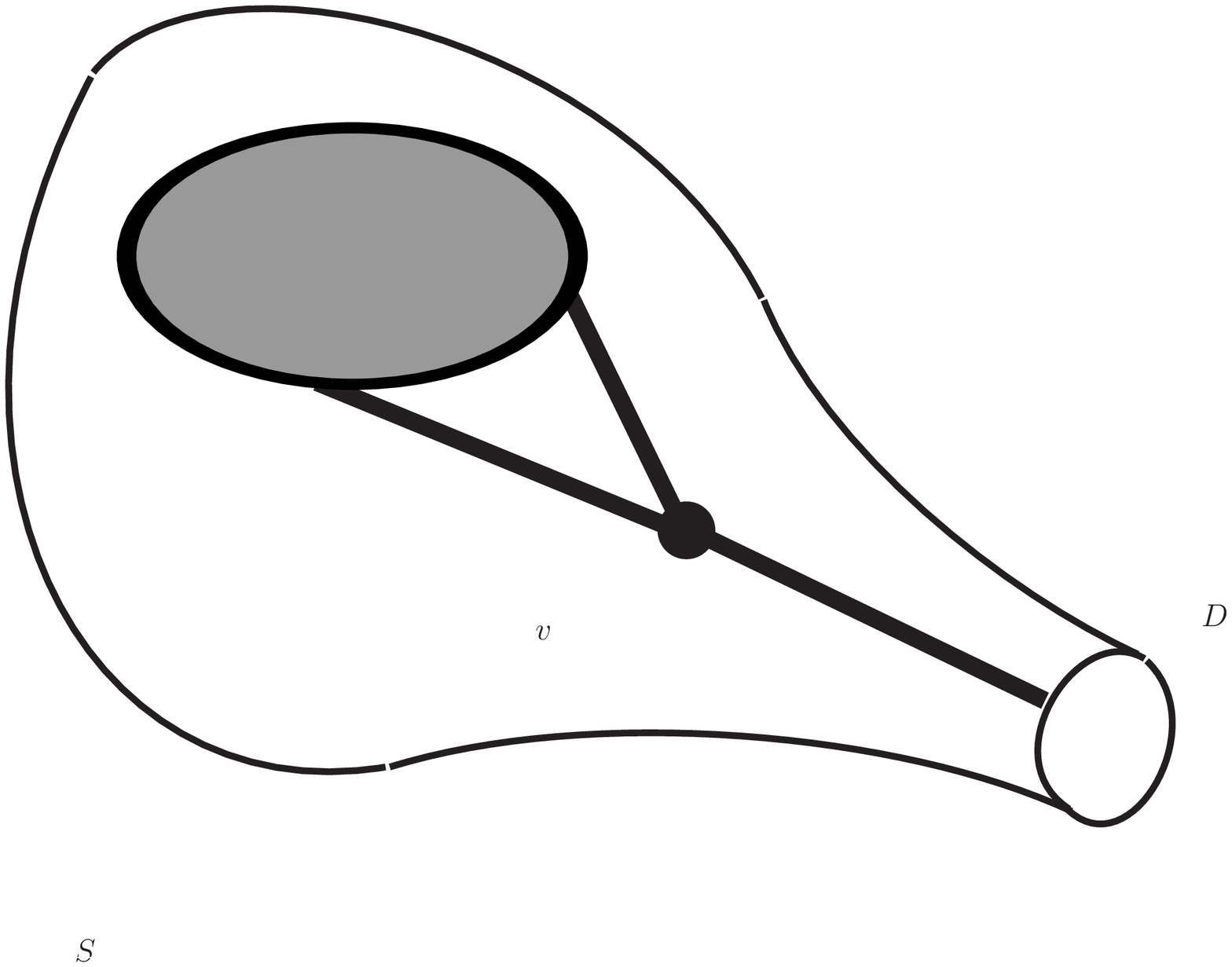}}\qquad .
\]
Clearly the imbedding $\Ga\rightarrow\Si$
can be obtained by gluing $\Si_1$
and $\Si_2$ along the boundaries $\partial D_1$, $\partial D_2$
of the removed discs.
Since the decorations of $v_i$ are of the form $\parbox{1cm}{
\psfrag{a}{$\scriptstyle{a_i}$}
\psfrag{v}{$\scriptstyle{v_i}$}
\psfrag{0}{$\scriptstyle{0}$}
\includegraphics[width=1cm]{Fig175.eps}}$,
there are no strands joining the $\hat{\Ga}_1$ and $\hat{\Ga}_2$
parts when applying the rules of Penrose evaluation in \S\ref{introdefsec}.
Also the application of Rule 3) to $v_i$ produces the same
strand structure as if there was no vertex:
\[
\parbox{1.5cm}{
\psfrag{a}{$\scriptstyle{a_i}$}
\psfrag{v}{$\scriptstyle{v_i}$}
\psfrag{0}{$\scriptstyle{0}$}
\includegraphics[width=1.5cm]{Fig175.eps}}
\ \longrightarrow\ 
\parbox{2.6cm}{
\psfrag{a}{$\scriptstyle{a_i}$}
\includegraphics[width=2.6cm]{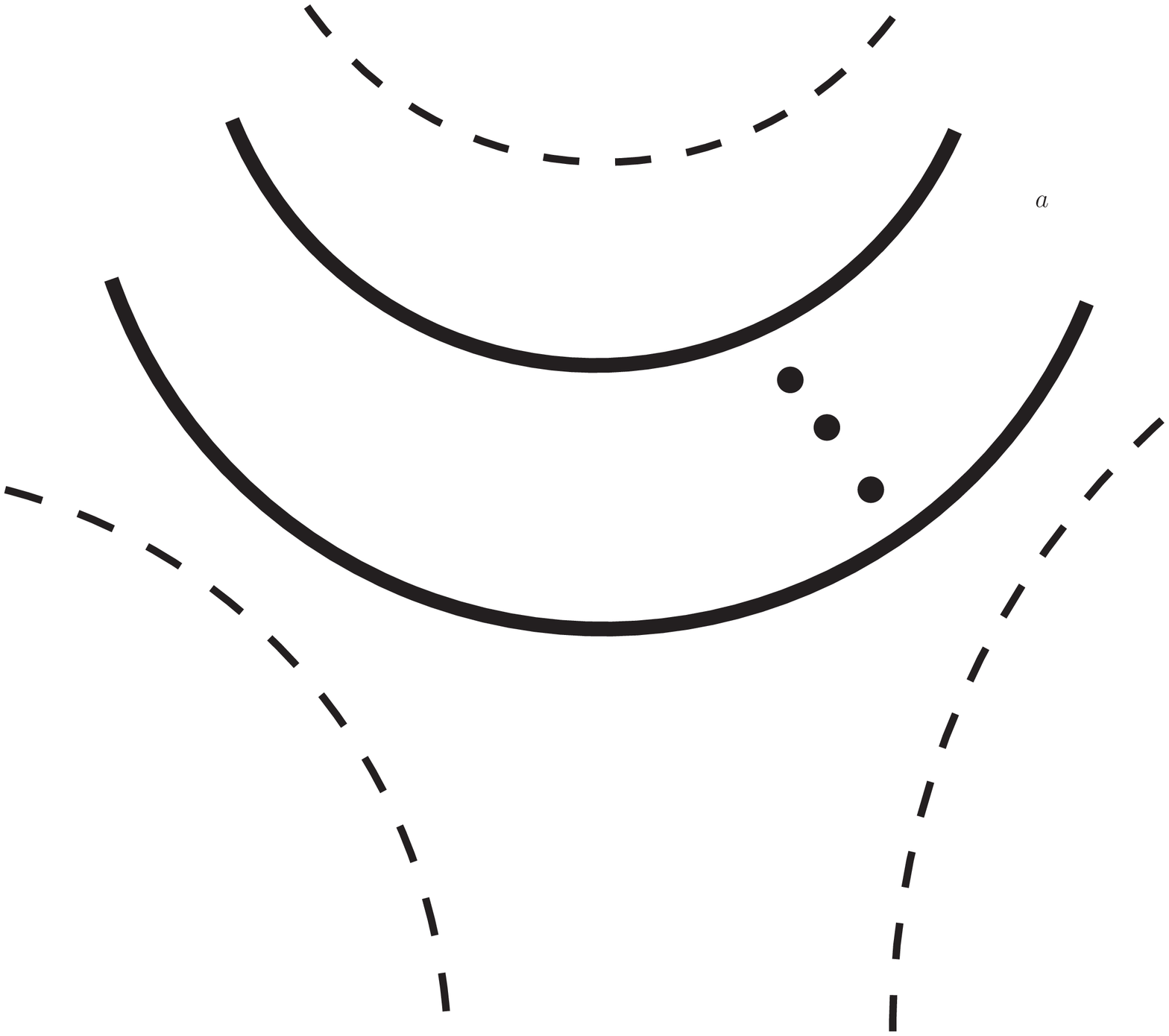}}
\ \longleftarrow\ 
\parbox{1cm}{
\psfrag{a}{$\scriptstyle{a_i}$}
\includegraphics[width=1cm]{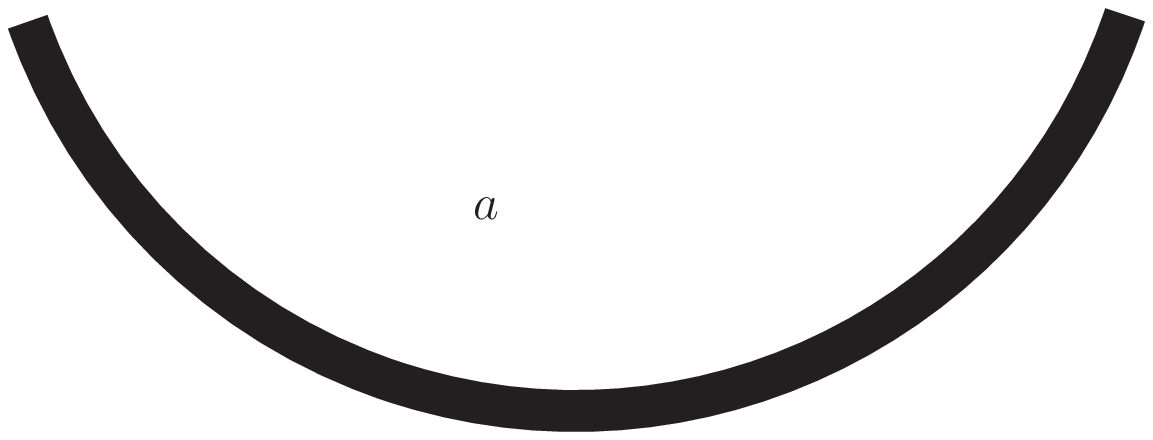}}\qquad .
\]
As a result of the previous considerations, the Penrose evaluation
$\<\Ga,\ga\>^P$ factors into $\<\Ga_1,\ga_1\>^P\times
\<\Ga_2,\ga_2\>^P$
except for one subtlety of the vertex erasure and its effect on the Penrose
bars of Rule 2).
If $v_i$ is not a loop vertex, its contribution is
\[
\parbox{1.5cm}{
\psfrag{a}{$\scriptstyle{a_i}$}
\psfrag{v}{$\scriptstyle{v_i}$}
\psfrag{0}{$\scriptstyle{0}$}
\includegraphics[width=1.5cm]{Fig175.eps}}
\qquad\longrightarrow\qquad
\parbox{2.6cm}{
\psfrag{a}{$\scriptstyle{a_i}$}
\includegraphics[width=2.6cm]{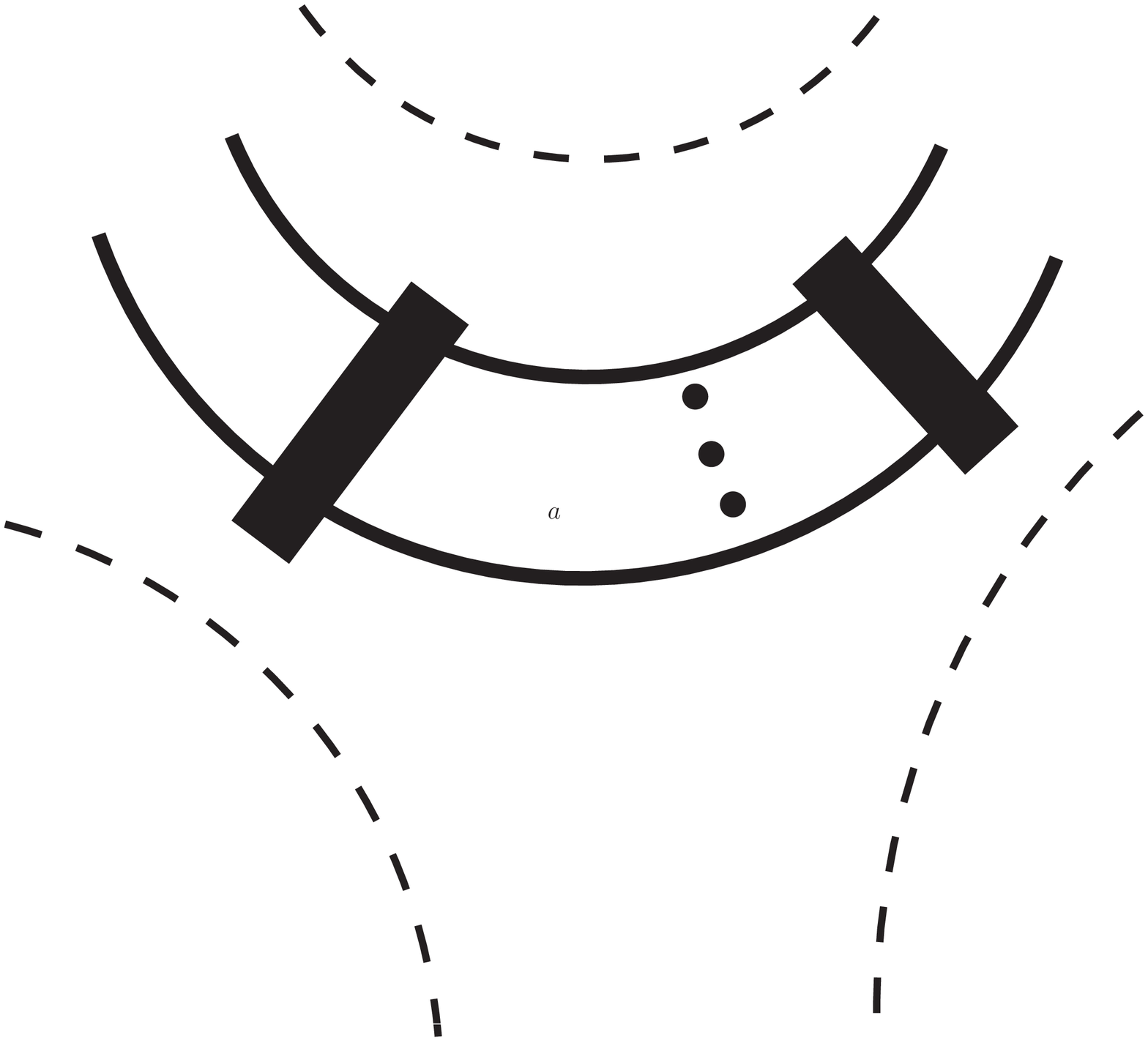}}
\]
namely, $a_i!$ times the contribution
\[
\parbox{1cm}{
\psfrag{a}{$\scriptstyle{a_i}$}
\includegraphics[width=1cm]{Fig184.eps}}
\qquad\longrightarrow\qquad
\parbox{2.6cm}{
\psfrag{a}{$\scriptstyle{a_i}$}
\includegraphics[width=2.6cm]{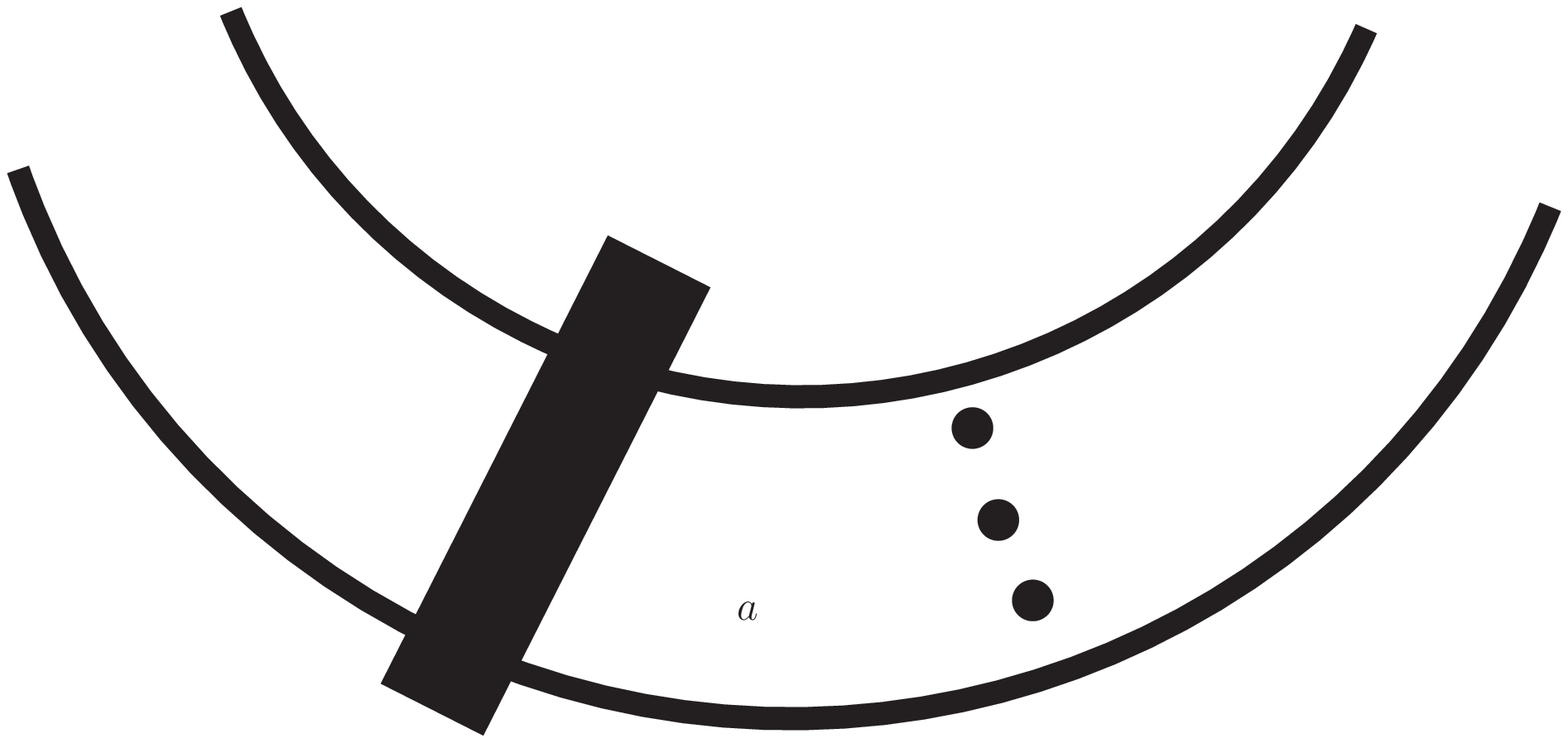}}
\]
one would get after erasing $v_i$ and fusing
the two edges carrying the decoration $a_i$.
If $v_i$ is a loop vertex, then both calculations produce the same result
\[
\parbox{1cm}{
\psfrag{a}{$\scriptstyle{a_i}$}
\psfrag{v}{$\scriptstyle{v_i}$}
\psfrag{G}{$\hat{\Ga}_i$}
\includegraphics[width=1cm]{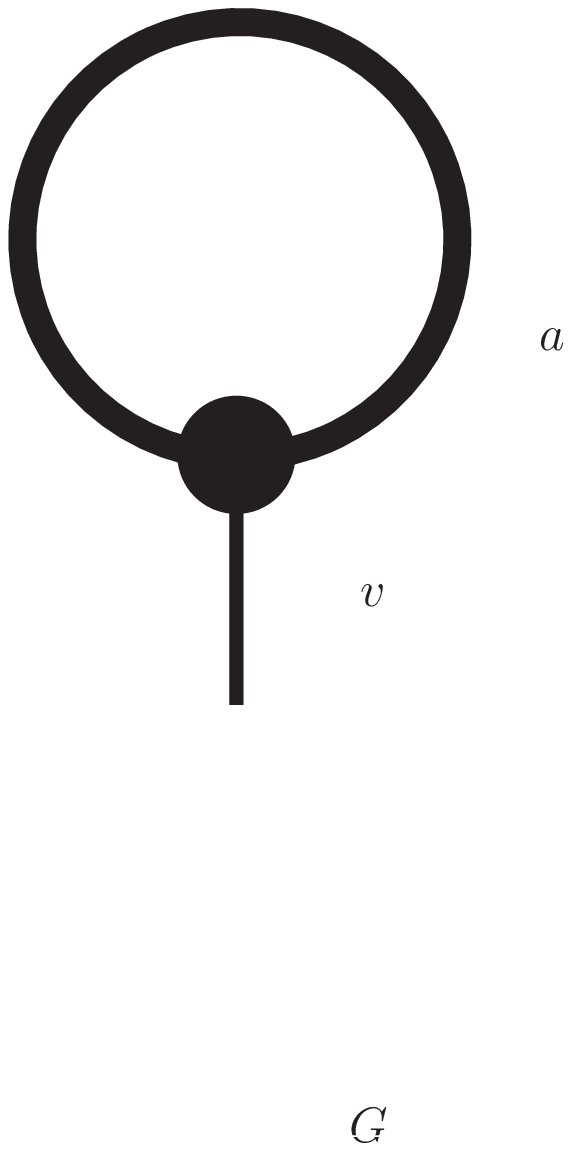}}
\ \longrightarrow\ 
\parbox{2.6cm}{
\psfrag{a}{$\scriptstyle{a_i}$}
\includegraphics[width=2.6cm]{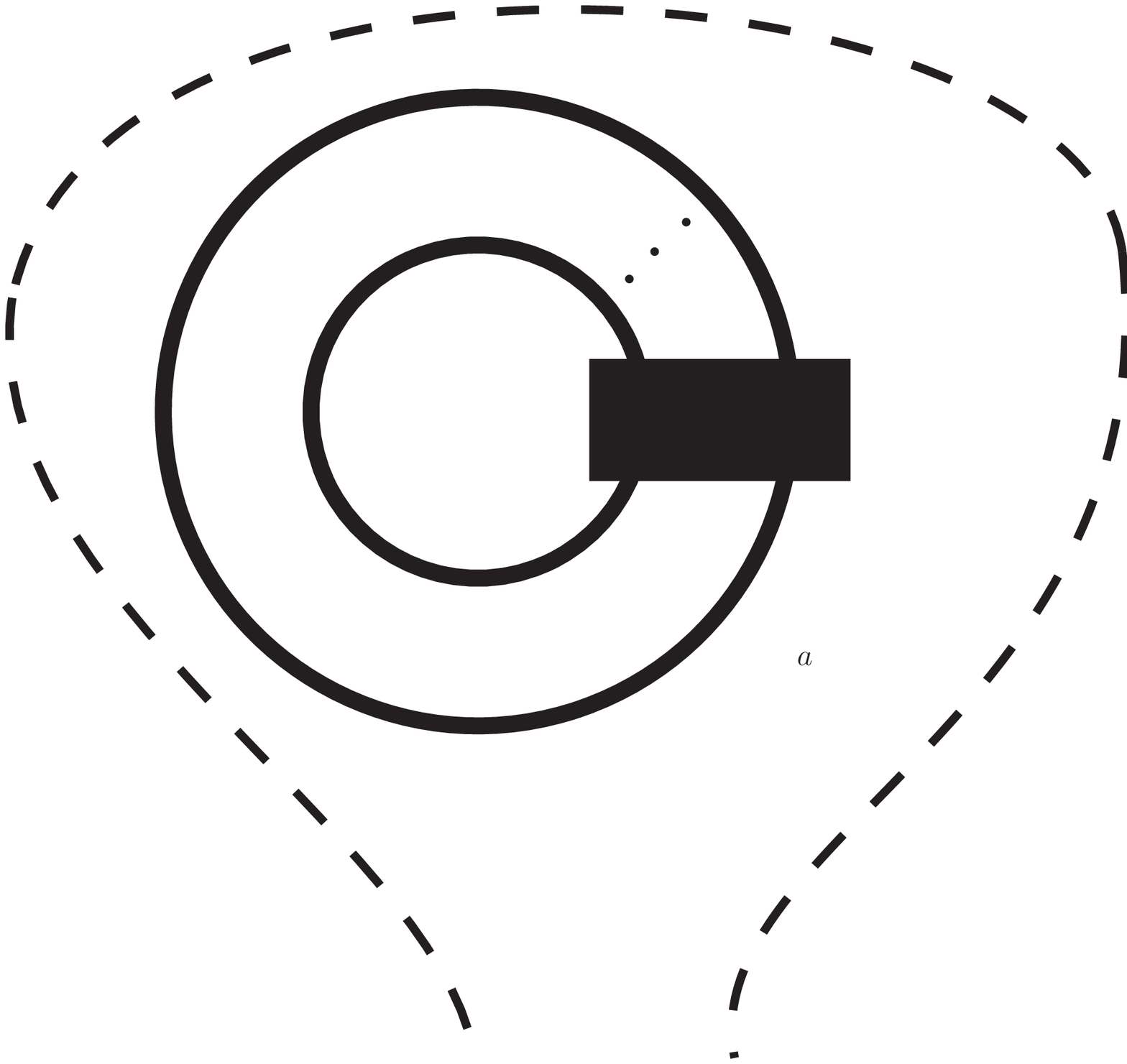}}
\ \longleftarrow\ 
\parbox{1cm}{
\psfrag{a}{$\scriptstyle{a_i}$}
\psfrag{G}{$\Ga_i$}
\includegraphics[width=1cm]{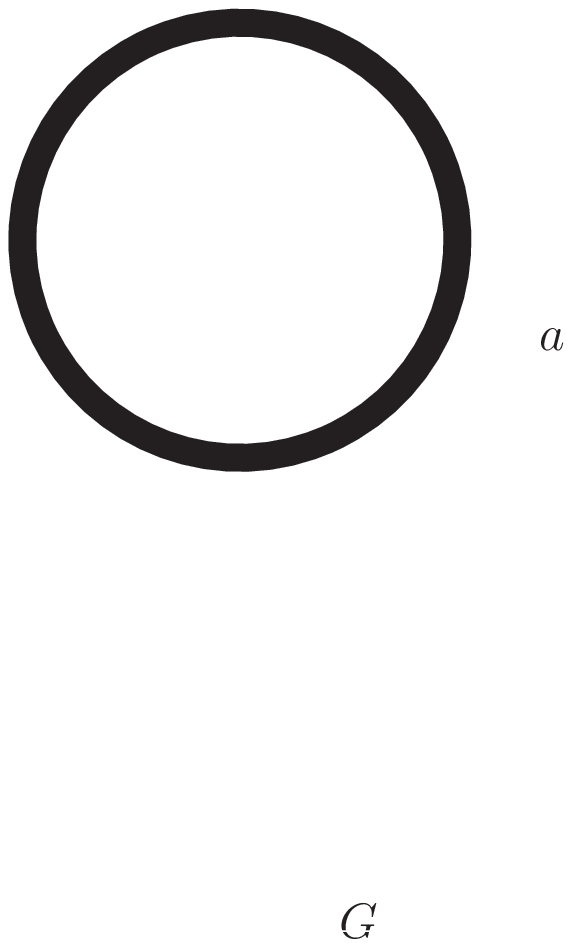}}\qquad .
\]
From the definitions (\ref{Sdef}) and (\ref{Udef}),
the formula relating the Penrose and the unitary evaluations is
\[
\<\Ga,\ga\>^U=\<\Ga,\ga\>^P\times
\prod_{v\in V(\Ga)}
\left\{
\left(\frac{a_v+b_v+c_v}{2}+1\right)!\right.
\]
\[
\left.
\left(\frac{a_v+b_v-c_v}{2}\right)!
\left(\frac{a_v+c_v-b_v}{2}\right)!
\left(\frac{b_v+c_v-a_v}{2}\right)!
\right\}^{-\frac{1}{2}}\ \ .
\]
Since the factors $[(a_i+1)!\ a_i!]^{-\frac{1}{2}}$
corresponding to $v_1,v_2$ are missing in the $\<\Ga_i,\ga_i\>^U$,
the proposition follows.
\qed

\begin{Remark}\label{dumbellrk}
From the proposition and the evaluation of trivial components
in Lemma \ref{trivcomplemma}, the evaluation (\ref{dumbell}) of the dumbell
graph follows.
\end{Remark}

\begin{Proposition}
Theorem \ref{mainthm} implies Theorem \ref{thmwithloops}.
\end{Proposition}

\noindent{\bf Proof:}
By induction on the number of vertices.
By the factorization property of Lemma \ref{factolemma},
it is enough to consider graphs which are connected.
If $(\Ga,\ga)$ has no loops, then $|\<\Ga,\ga\>^U|\le 1$
by Theorem \ref{mainthm} which is assumed to hold.
Otherwise, if there is a loop, then let $e_0$ be the bridge
which connects it to the rest of the graph.
Using the notations and setting of Proposition \ref{bridgeprop},
let $v_1$ be the vertex in the given loop.
Now if the other vertex $v_2$ is also a loop vertex, then by connectedness
$(\Ga,\ga)$ must be a dumbell graph and, for $n\ge 1$, we have
\[
\<
\ 
\parbox{2.5cm}{\psfrag{a}{$\scriptstyle{na}$}
\psfrag{b}{$\scriptstyle{nb}$}\psfrag{c}{$\scriptstyle{nc}$}
\includegraphics[width=2.5cm]{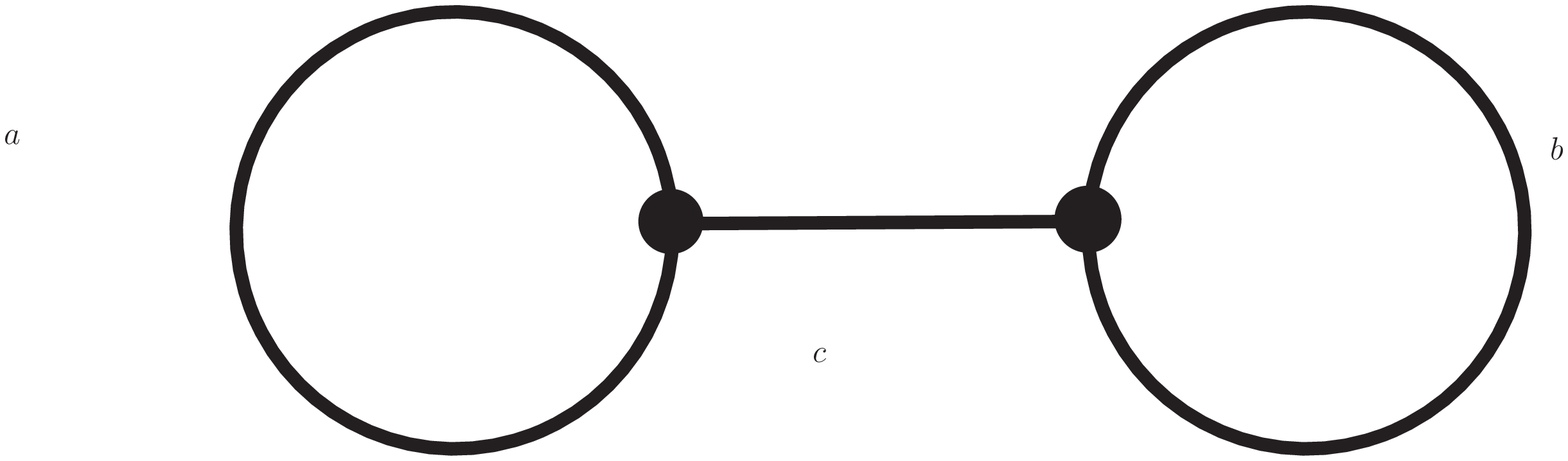}}
\ \ \>^U=\de_{c,0} \sqrt{(na+1)(nb+1)}\ \ .
\]
This provides the desired polynomial bound in $n$.
Else, if $v_2$ is not a loop vertex
then
by Proposition \ref{bridgeprop}
and Lemma \ref{trivcomplemma}
\[
\left|\<\Ga,n\ga\>^U\right|=
\sqrt{n a_1+1}\times\frac{1}{\sqrt{n a_2+1}}\times
\left|\<\Ga_2,n\ga_2\>^U\right|
\]
and $\left|\<\Ga_2,n\ga_2\>^U\right|$
is polynomially bounded by the induction hypothesis. Therefore so
is $\left|\<\Ga,n\ga\>^U\right|$.
\qed

\section{Estimates on special spin networks}

\subsection{The tetrahedron or 6-j symbol}
In this well-known case one has a precise $n\rightarrow \infty$
asymptotic: the Ponzano-Regge formula.
This will provide both upper and lower bounds
on $\sqrt[n]{|\<\Ga,n\ga\>^U|}$ as $n\rightarrow\infty$.
However, we will derive an upper bound using
soft analysis: a simple Cauchy-Schwarz inequality.
This will provide a gentle introduction
to the proof of our main theorem in \S\ref{mainthmsec}.

\subsubsection{Preparation}
Consider the spin network $(\Ga,\ga)$
given by
\[
\parbox{2cm}{
\psfrag{a}{$\scriptstyle{a}$}
\psfrag{b}{$\scriptstyle{b}$}
\psfrag{c}{$\scriptstyle{c}$}
\psfrag{d}{$\scriptstyle{d}$}
\psfrag{e}{$\scriptstyle{e}$}
\psfrag{f}{$\scriptstyle{f}$}
\includegraphics[width=2cm]{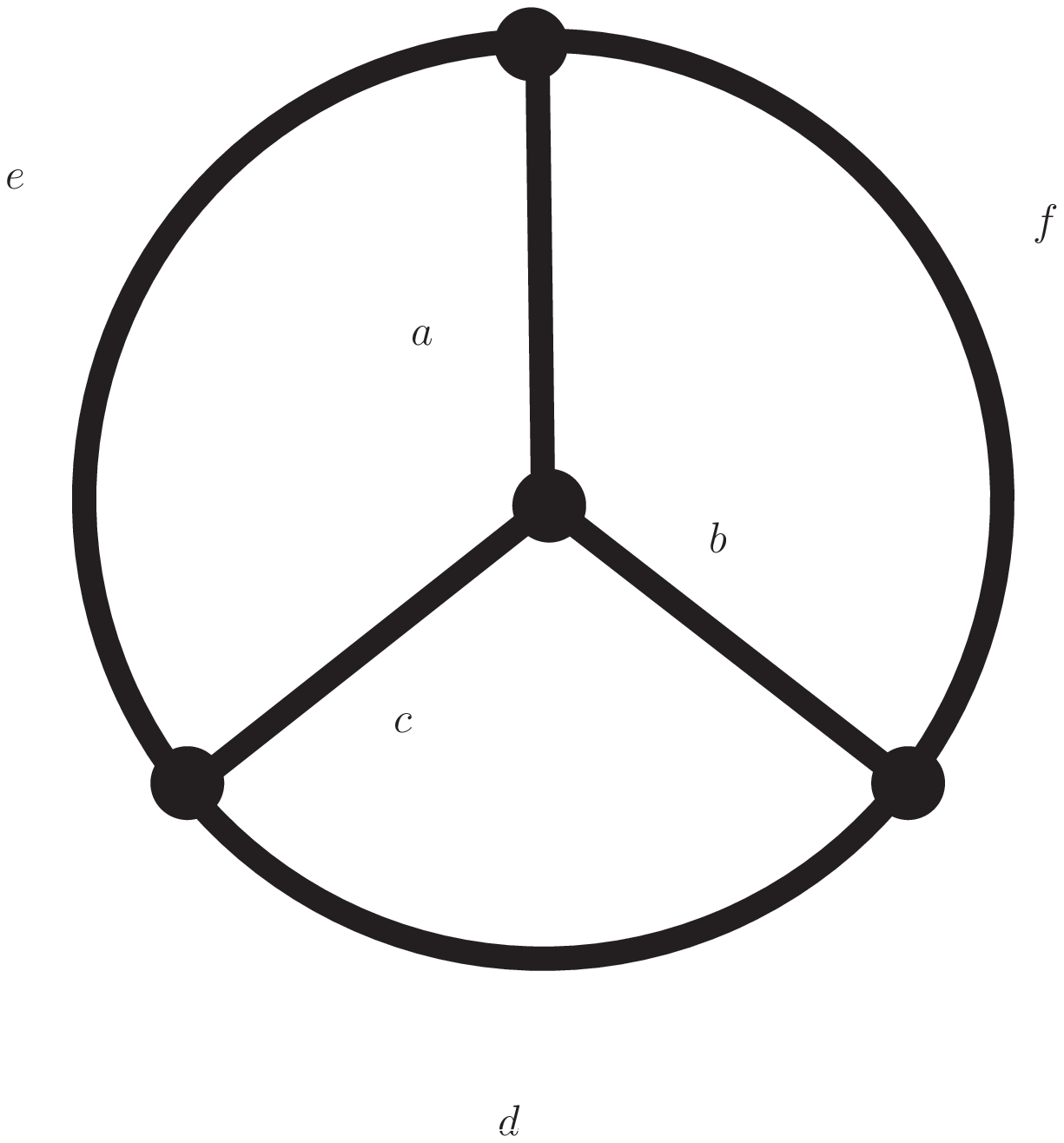}}\qquad .
\]
The cyclic orientations at the vertices are counterclockwise,
and the imbedding is planar.
By Corollary \ref{negdimU},
\[
\left|\<\Ga,\ga\>^U\right|
=\frac{1}{(e+1)\sqrt{(c+1)(f+1)}}\times
\left|
\<
\parbox{2.5cm}{
\psfrag{a}{$\scriptstyle{a}$}
\psfrag{b}{$\scriptstyle{b}$}
\psfrag{c}{$\scriptstyle{c}$}
\psfrag{d}{$\scriptstyle{d}$}
\psfrag{e}{$\scriptstyle{e}$}
\psfrag{f}{$\scriptstyle{f}$}
\includegraphics[width=2.5cm]{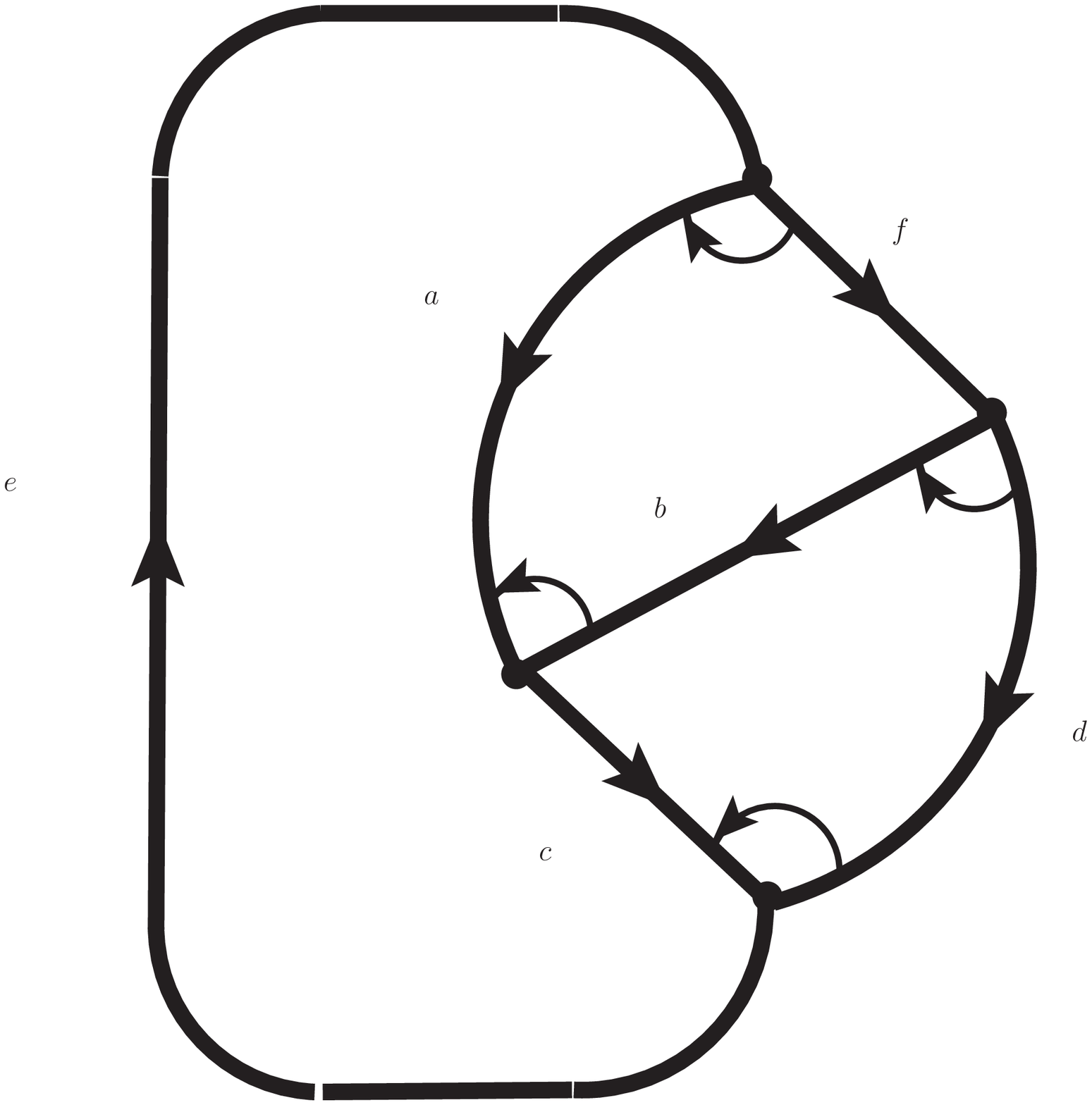}}
\>^{\pi\io}\right|\ \ .
\]
As recalled in~\cite[\S7.6]{AC3}
the corresponding standard Wigner 6-j symbol is
\[
\left\{
\begin{array}{ccc}
\frac{a}{2} & \frac{b}{2} & \frac{c}{2}\\
\frac{d}{2} & \frac{e}{2} & \frac{f}{2}
\end{array}
\right\}
=\frac{(-1)^{a+b+d+e}}{\sqrt{(c+1)(f+1)}}\times\al
\]
where $\alpha$ is the proportionality constant in the explicit
instance of Schur's Lemma
\[
\<
\parbox{1.7cm}{
\psfrag{a}{$\scriptstyle{a}$}
\psfrag{b}{$\scriptstyle{b}$}
\psfrag{c}{$\scriptstyle{c}$}
\psfrag{d}{$\scriptstyle{d}$}
\psfrag{e}{$\scriptstyle{e}$}
\psfrag{f}{$\scriptstyle{f}$}
\includegraphics[width=1.7cm]{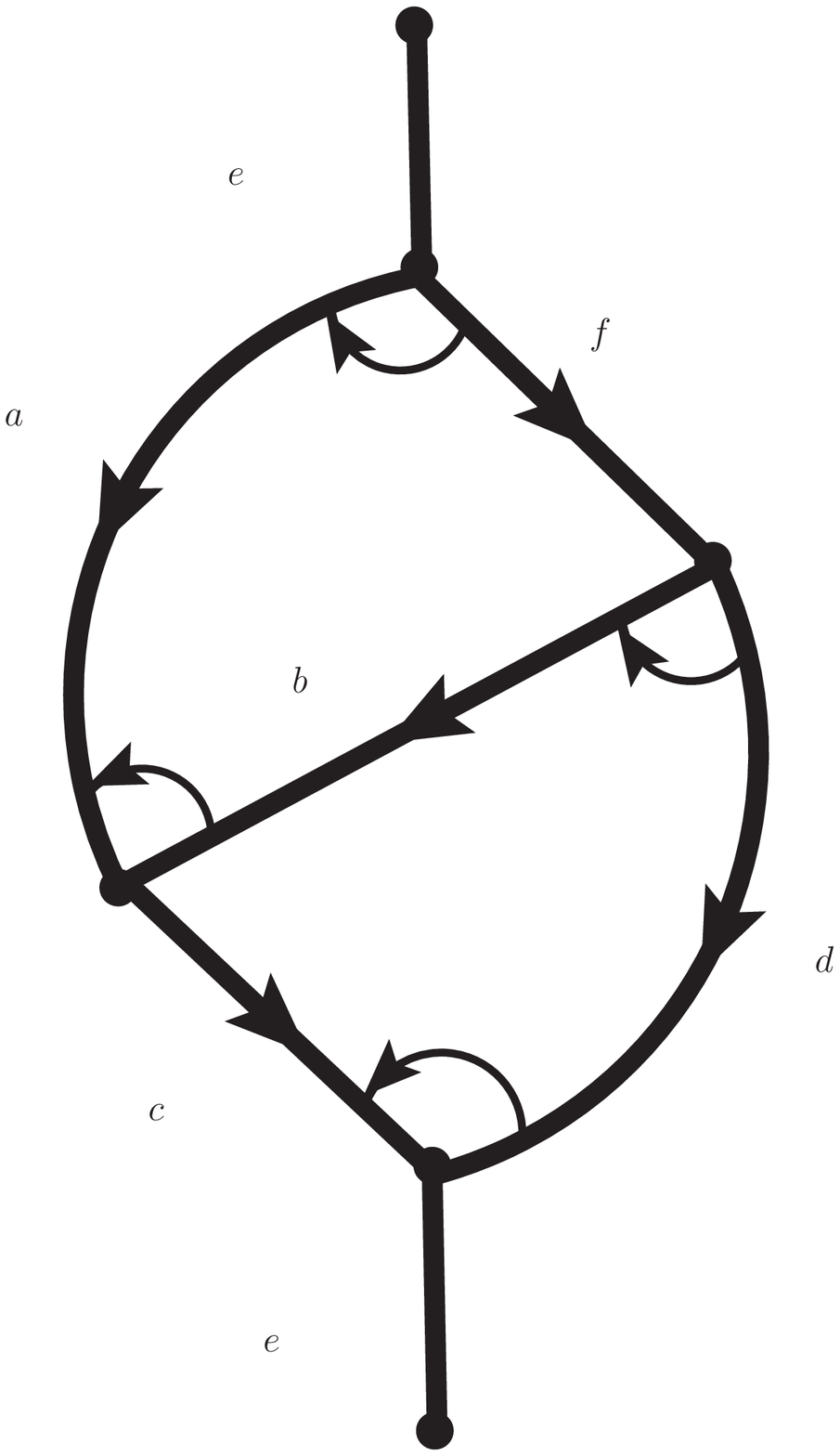}}
\>^{\pi\io}=\alpha\times\<
\parbox{0.7cm}{\psfrag{e}{$\scriptstyle{e}$}
\includegraphics[width=0.7cm]{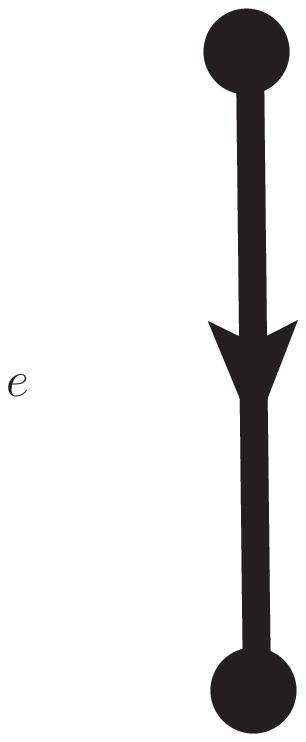}}
\>^{\pi\io}\qquad .
\]
By taking the trace in $\cH_e^{\ast}\otimes\cH_e$,
one has
\[
\alpha =\frac{1}{e+1}\times
\<
\parbox{2.5cm}{
\psfrag{a}{$\scriptstyle{a}$}
\psfrag{b}{$\scriptstyle{b}$}
\psfrag{c}{$\scriptstyle{c}$}
\psfrag{d}{$\scriptstyle{d}$}
\psfrag{e}{$\scriptstyle{e}$}
\psfrag{f}{$\scriptstyle{f}$}
\includegraphics[width=2.5cm]{Fig197.eps}}
\>^{\pi\io}
\]
and therefore
\begin{equation}
\left|\<\Ga,\ga\>^U\right|
=\left|\left\{
\begin{array}{ccc}
\frac{a}{2} & \frac{b}{2} & \frac{c}{2}\\
\frac{d}{2} & \frac{e}{2} & \frac{f}{2}
\end{array}
\right\}\right|\ \ .
\label{6jrelation}
\end{equation}

\subsubsection{The upper bound}\label{6jupsec}
Notice that by introducing a splitting of the graph
and the associated Hilbert space
$\cH=\cH_e^{\ast}\otimes\cH_a\otimes\cH_b\otimes\cH_c$
one has
\[
\<
\parbox{4cm}{
\psfrag{a}{$\scriptstyle{a}$}
\psfrag{b}{$\scriptstyle{b}$}
\psfrag{c}{$\scriptstyle{c}$}
\psfrag{d}{$\scriptstyle{d}$}
\psfrag{e}{$\scriptstyle{e}$}
\psfrag{f}{$\scriptstyle{f}$}
\psfrag{A}{$A$}
\psfrag{B}{$B$}
\includegraphics[width=4cm]{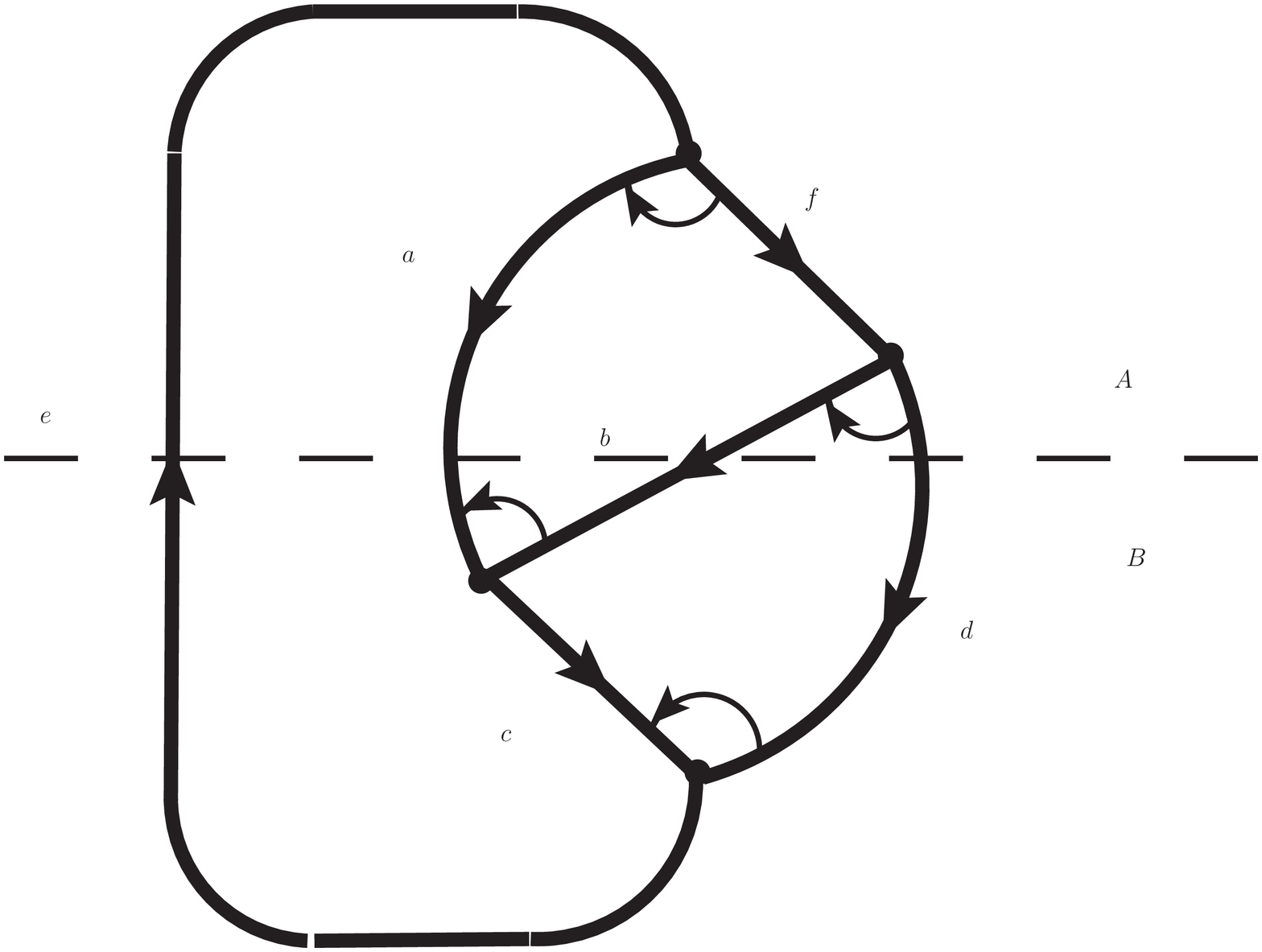}}
\>^{\pi\io}=\<B|A\>_{\cH}
\]
where
\[
A=\<
\parbox{3cm}{
\psfrag{a}{$\scriptstyle{a}$}
\psfrag{b}{$\scriptstyle{b}$}
\psfrag{d}{$\scriptstyle{d}$}
\psfrag{e}{$\scriptstyle{e}$}
\psfrag{f}{$\scriptstyle{f}$}
\includegraphics[width=3cm]{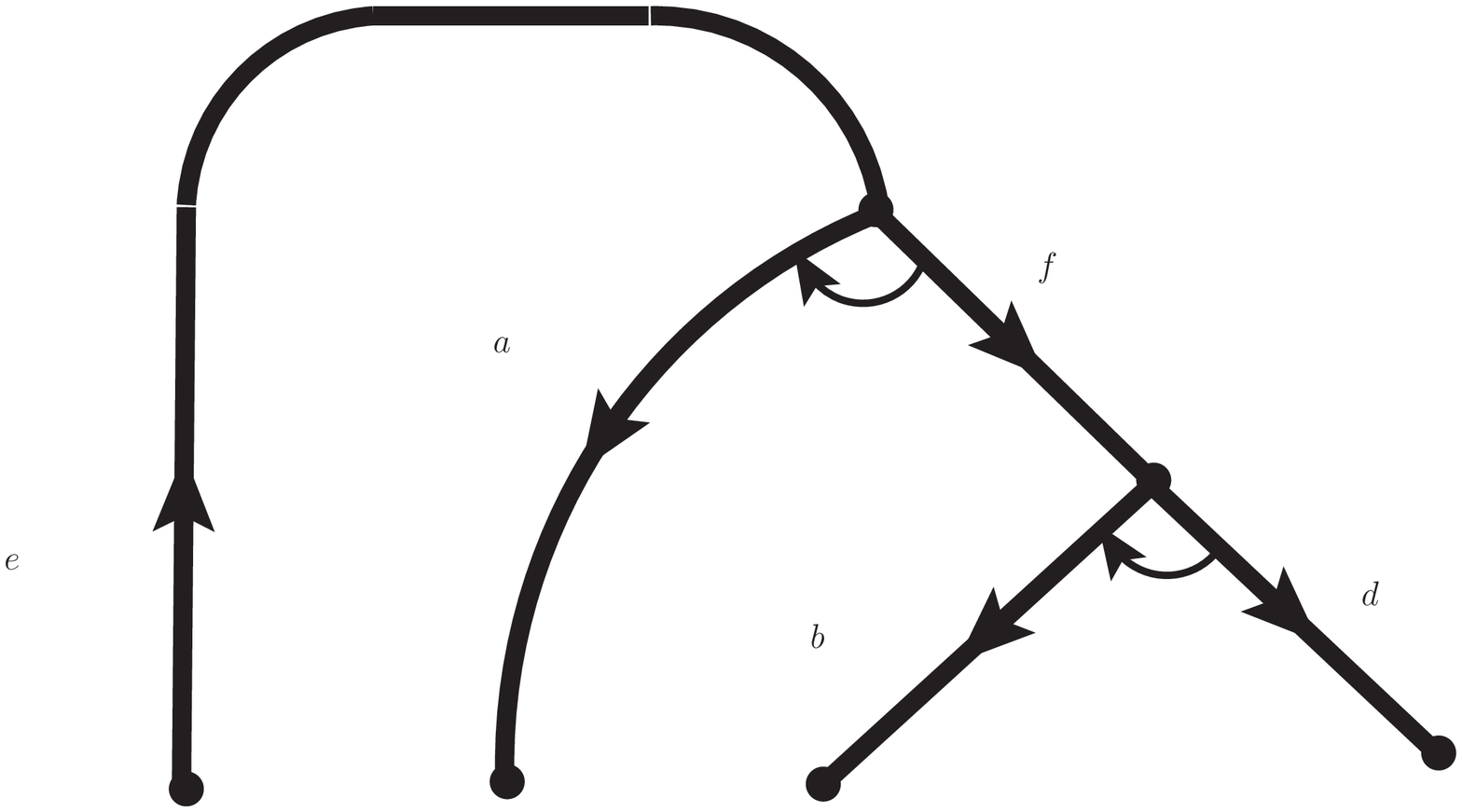}}
\>^{\pi\io}\in\cH
\]
and
\[
B=\<
\parbox{3cm}{
\psfrag{a}{$\scriptstyle{a}$}
\psfrag{b}{$\scriptstyle{b}$}
\psfrag{d}{$\scriptstyle{d}$}
\psfrag{e}{$\scriptstyle{e}$}
\psfrag{c}{$\scriptstyle{c}$}
\includegraphics[width=3cm]{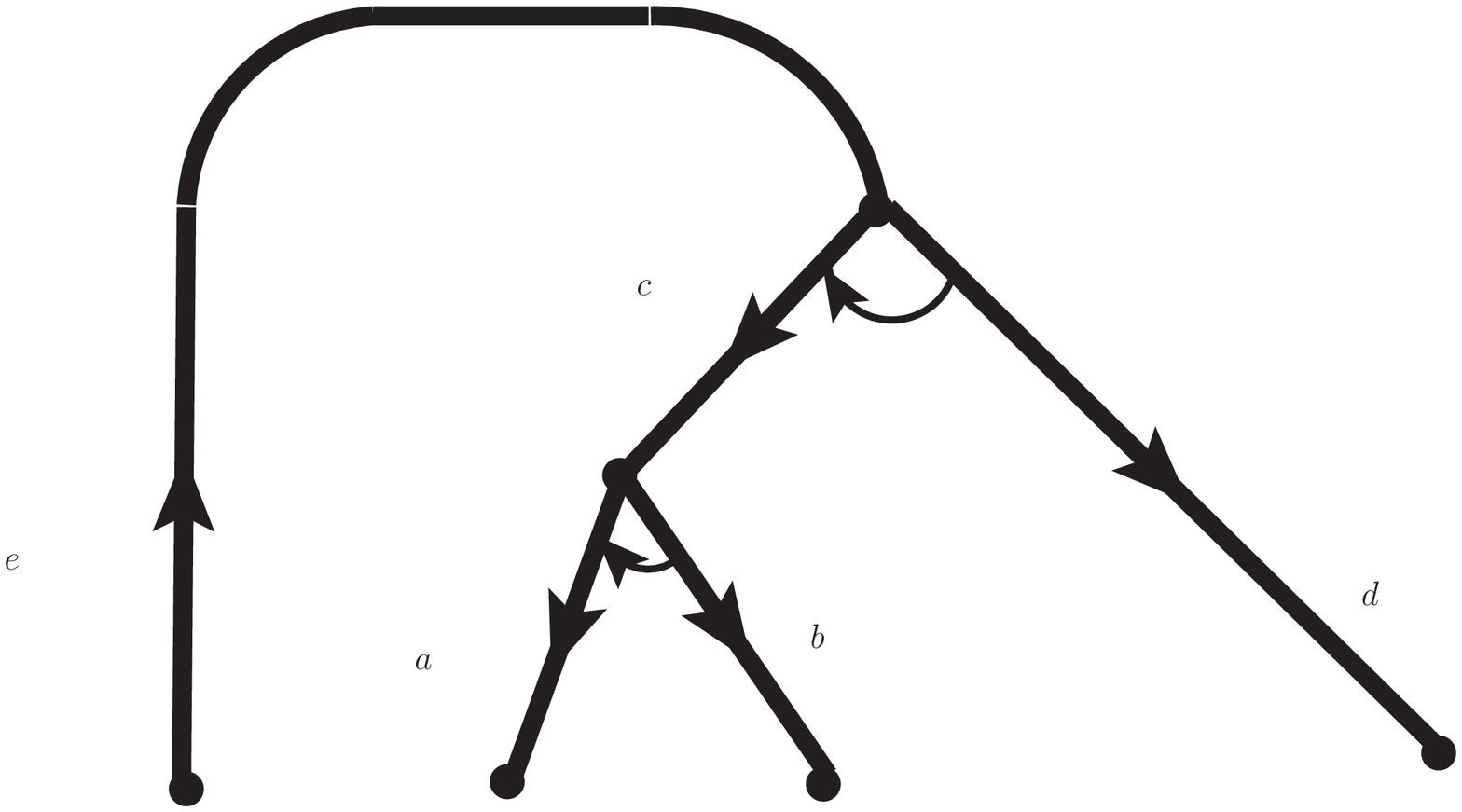}}
\>^{\pi\io}\in\cH
\]
and the inner product is the natural one
as in (\ref{innerprod}) and (\ref{innerprodB}), using FDC.
Note the reversal of orientations in $B$.

By the Cauchy-Schwarz inequality
\[
\left|\<\Ga,\ga\>^U\right|\le
\frac{1}{(e+1)\sqrt{(c+1)(f+1)}}\times
\sqrt{\<A|A\>_{\cH}\<B|B\>_{\cH}}\ \ .
\]
However,
\[
\<A|A\>_{\cH}=\<
\parbox{2.5cm}{
\psfrag{a}{$\scriptstyle{a}$}
\psfrag{b}{$\scriptstyle{b}$}
\psfrag{c}{$\scriptstyle{c}$}
\psfrag{d}{$\scriptstyle{d}$}
\psfrag{e}{$\scriptstyle{e}$}
\psfrag{f}{$\scriptstyle{f}$}
\includegraphics[width=2.5cm]{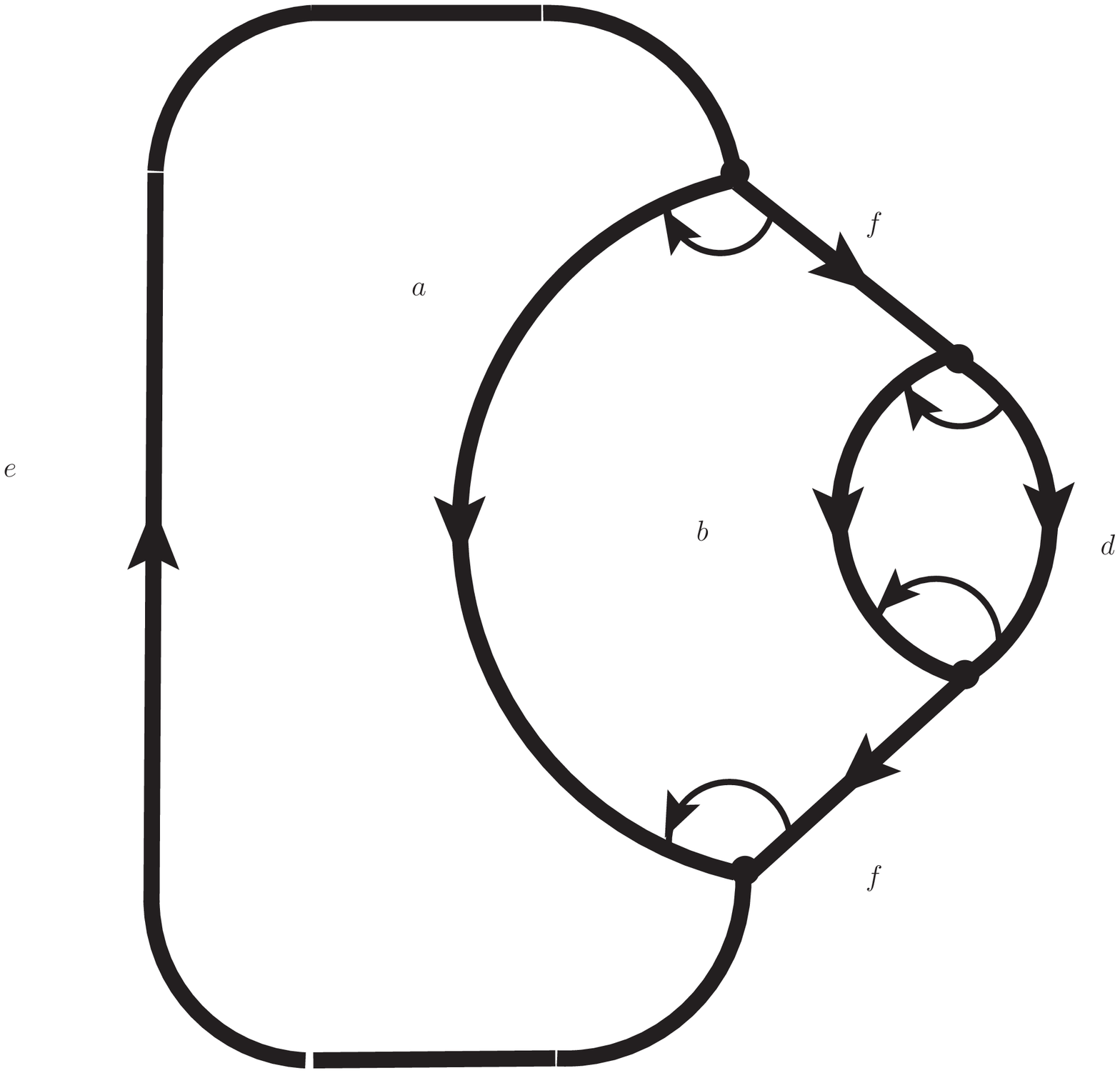}}
\>^{\pi\io}
\]
\[
=\<
\parbox{1.8cm}{
\psfrag{a}{$\scriptstyle{a}$}
\psfrag{e}{$\scriptstyle{e}$}
\psfrag{f}{$\scriptstyle{f}$}
\includegraphics[width=1.8cm]{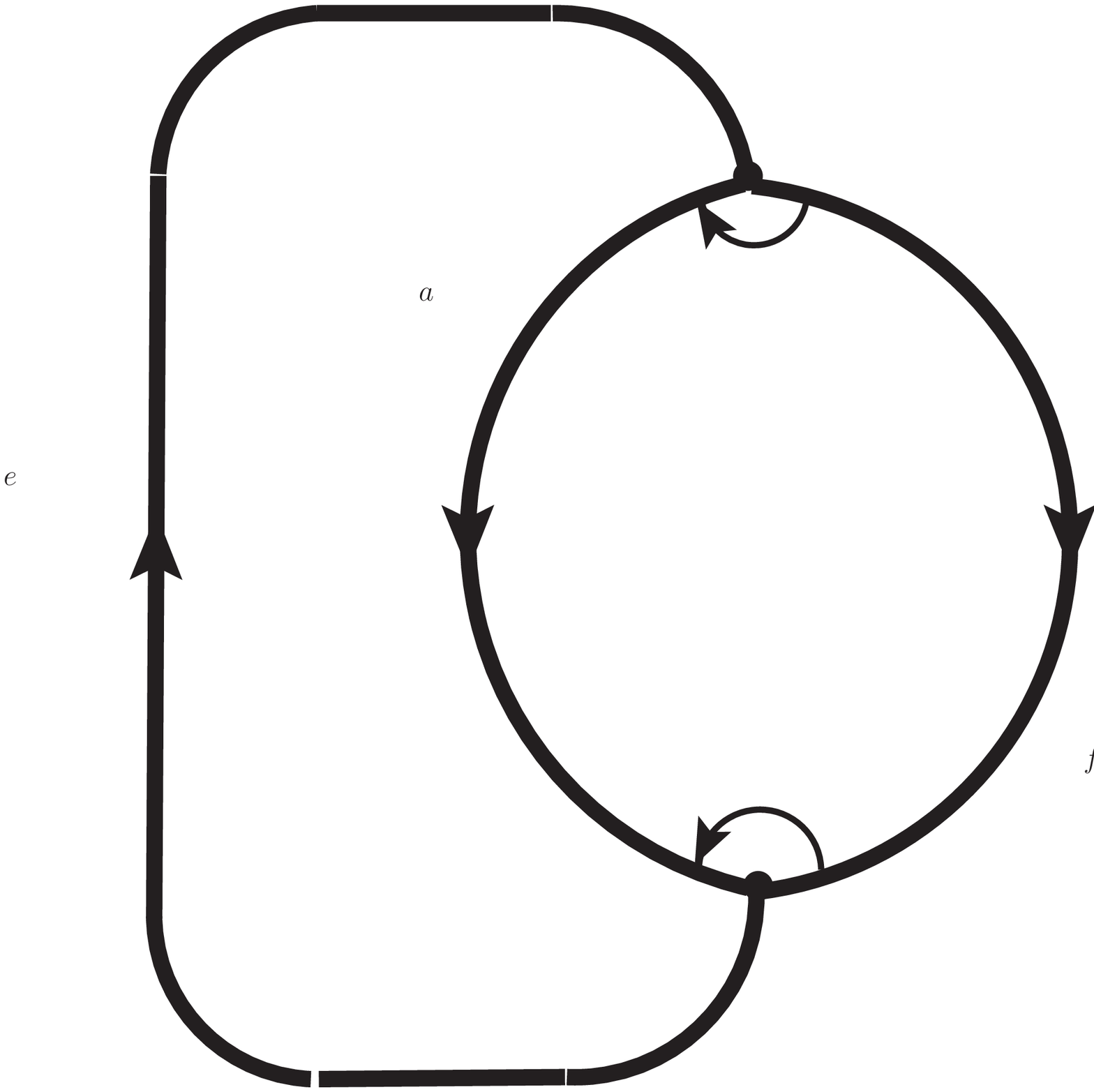}}
\>^{\pi\io}=\<
\parbox{1.2cm}{
\psfrag{e}{$\scriptstyle{e}$}
\includegraphics[width=1.2cm]{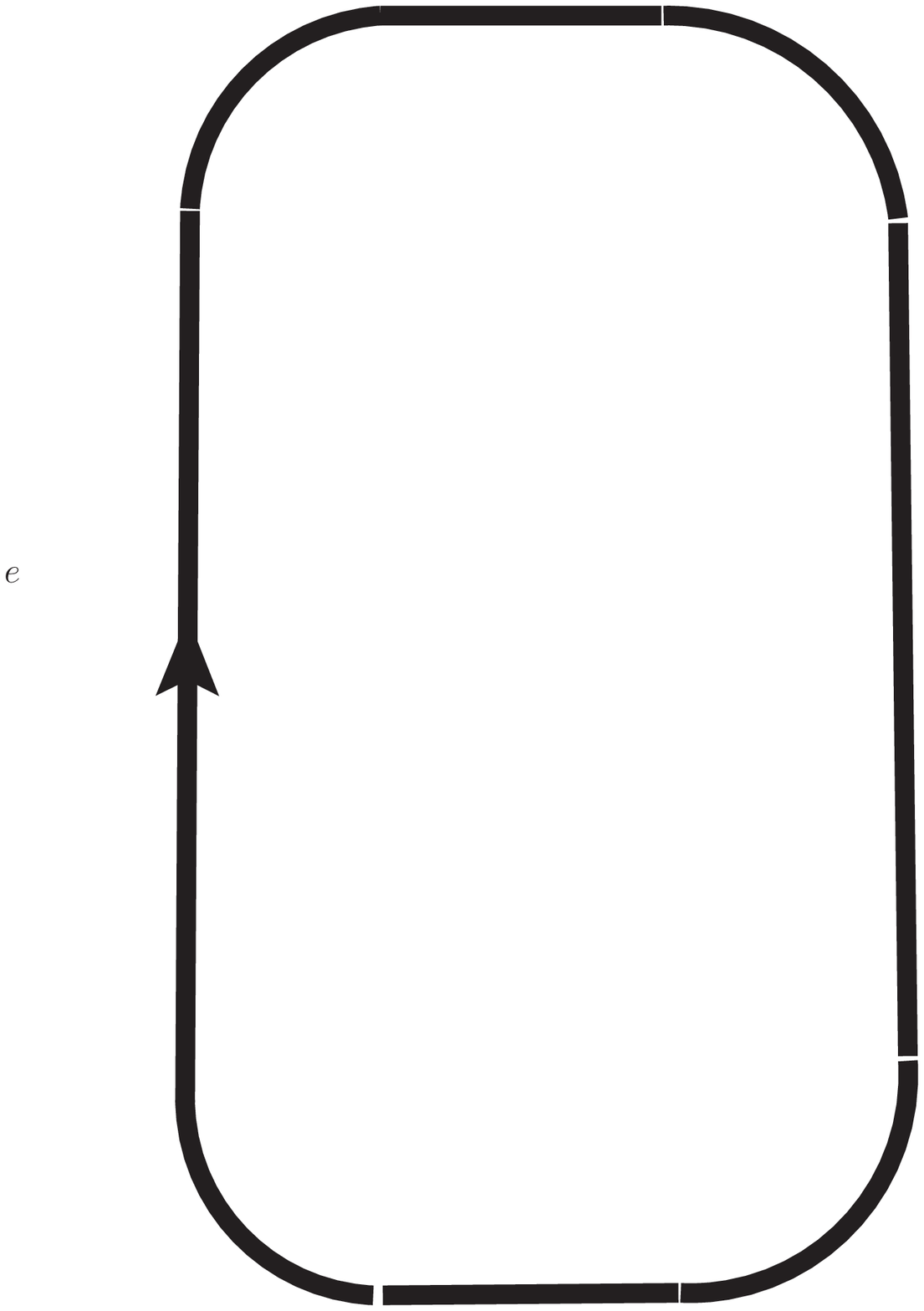}}
\ \ \>^{\pi\io}=e+1
\]
by using part 3) of Proposition \ref{piiotaprop} twice.
Likewise, $\<B|B\>_{\cH}=e+1$.
Therefore
\begin{equation}
\left|\<\Ga,\ga\>^U\right|
=\left|\left\{
\begin{array}{ccc}
\frac{a}{2} & \frac{b}{2} & \frac{c}{2}\\
\frac{d}{2} & \frac{e}{2} & \frac{f}{2}
\end{array}
\right\}\right|
\le\frac{1}{\sqrt{(c+1)(f+1)}}\le 1
\label{6jupperbd}
\end{equation}
which explicitly shows that the tetrahedral
spin network satisfies the statement in Theorem \ref{mainthm}.
Of course, one can use the symmetries of the
the 6-j symbol in order to obtain analogous bounds by
$\frac{1}{\sqrt{(a+1)(d+1)}}$, etc.
Similar bounds for the 6-j symbol are well-known
(see, e.g.,~\cite[App. D]{FriedelL2}).
Also note that similar considerations of inner products
in tensor spaces were used in~\cite[\S3.1]{AC4} in
order to prove
the nonvanishing of some combinations of CG networks (one can also
prove~\cite[Lemma 2.3]{AC3} in the same way).
A thematically similar use of the Cauchy-Schwarz inequality for 3-manifold invariants
or chromatic polynomials of planar graphs can also be found in~\cite[Thm. 2.2]{Garoufalidis}
and~\cite[\S7]{FendleyK}.

\subsubsection{The lower bound}
Lower bounds seem much more difficult to obtain.
In particular, we have not been able to find one for the 6-j
symbol by trying to quantitatively analyse how
far the Cauchy-Schwarz inequality is from an equality.
One difficulty towards this goal is that
Wigner symbols can have accidental zeros which are quite poorly
understood (see~\cite{RaynalJRR} and references therein).
We will rely instead on the Ponzano-Regge asymptotic
formula.
It was conjectured in~\cite{PonzanoR} based on numerical evidence
and some very clever consistency checks, and it was rigorously proved in~\cite{Roberts1}.
Other work related to this asymptotic formula can be found
in~\cite{SchultenG,BarrettS,FriedelL,Charles,Gurau,GvdV}.
Note that~\cite{FriedelL,BarrettS}
consider the square of the 6-j symbol which
is also suitable for the needs of this section.

We will consider the case where the original decorations are uniform
$\ga\equiv 2$ and then the rescaled network $\ga\rightarrow n\ga$.
\begin{Lemma}\label{PRlemma}
\[
\limsup_{n\rightarrow\infty}
\left|
\<
\parbox{1.2cm}{\includegraphics[width=1.2cm]{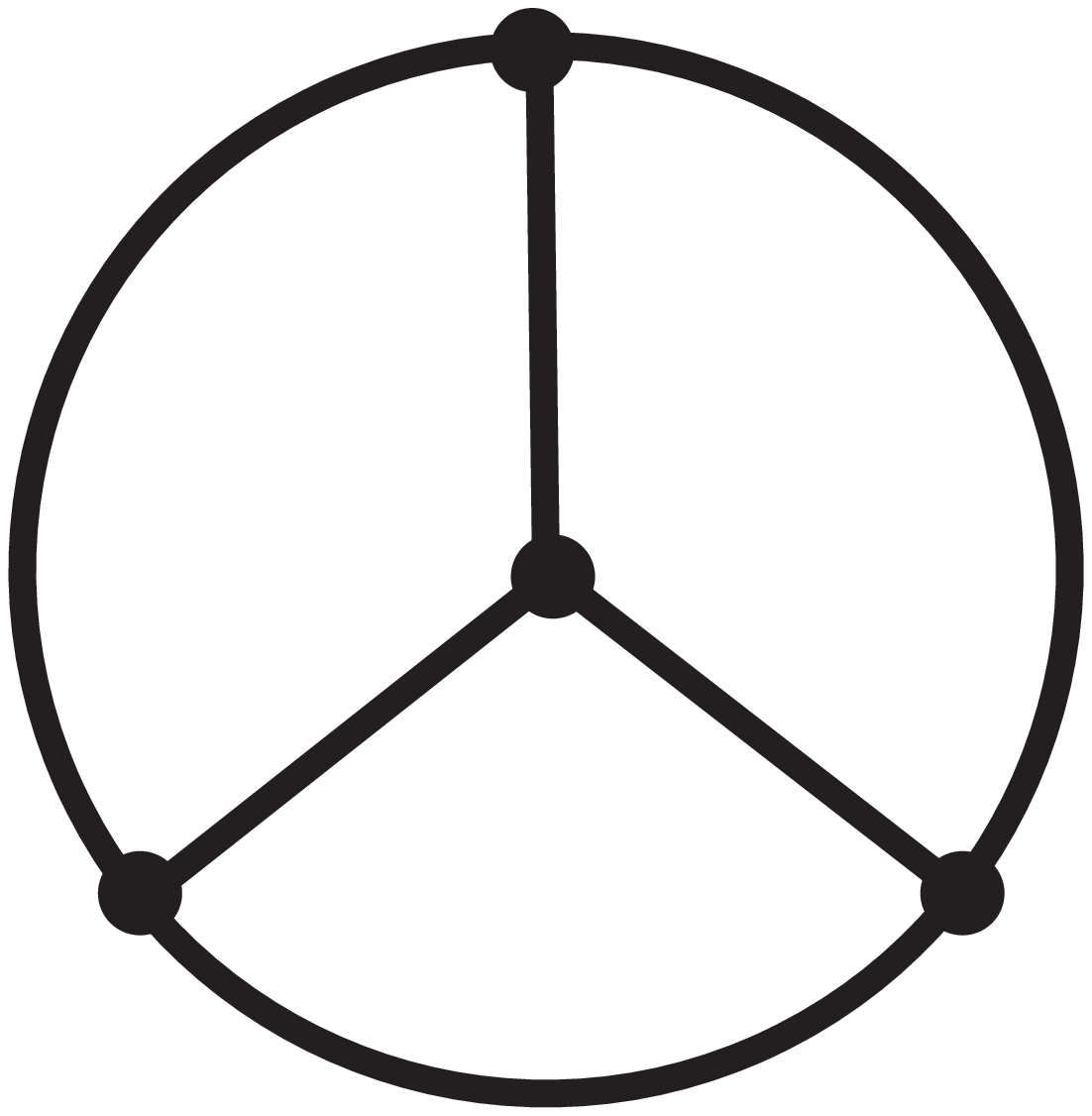}}
,n\ga\>^U
\right|^{\frac{1}{n}}=1\ .
\]
\end{Lemma}

\noindent{\bf Proof:}
By (\ref{6jrelation})
\[
\left|
\<
\parbox{1.2cm}{\includegraphics[width=1.2cm]{Fig207.eps}}
,2n\>^U
\right|=
\left|\left\{
\begin{array}{ccc}
n & n & n\\
n & n & n
\end{array}
\right\}\right|\ \ .
\]
By the Ponzano-Regge formula for the case of a regular (Euclidean)
tetrahedron (see, e.g.,~\cite[\S10]{GvdV}), we have:
\[
\left\{
\begin{array}{ccc}
n & n & n\\
n & n & n
\end{array}
\right\}
=\frac{-1}{2^{\frac{1}{4}}\pi^{\frac{1}{2}}n^{\frac{3}{2}}}
\cos\left(6\left(n+\frac{1}{2}\right)\om-\frac{\pi}{4}\right)
+O\left(n^{-\frac{5}{2}}\right)
\]
when $n\rightarrow\infty$, and where $\om=\arccos\left(\frac{1}{3}\right)$.
Therefore
\[
\limsup_{n\rightarrow\infty}
\left|
\<
\parbox{1.2cm}{\includegraphics[width=1.2cm]{Fig207.eps}}
,n\ga\>^U
\right|^{\frac{1}{n}}=
\limsup_{n\rightarrow\infty}
\left|
\cos\left(6\left(n+\frac{1}{2}\right)\om-\frac{\pi}{4}\right)+\frac{\al_n}{n}
\right|^{\frac{1}{n}}
\]
for some bounded sequence $(\al_n)_{n\ge 1}$.
Trivially, the limsup is $\le 1$.
On the other hand, $\frac{6\om}{2\pi}-1=0.1754\ldots
\in(0,\frac{1}{3})$, hence
one can extract a subsequence for
which the angle is in $[-\frac{\pi}{3},\frac{\pi}{3}]$ ${\rm mod}\ 2\pi$
and the cosine is $\ge\frac{1}{2}$. As a result the limsup is $\ge 1$.
\qed

\noindent
This proves the statement in Conjecture
\ref{mainconj} for the tetrahedron graph, which was altready
known~\cite{GvdV}.

\subsection{The generalized drum}
As in~\cite{GvdV}, we consider the graph ${\rm Drum}_s$:
\[
\parbox{2cm}{\includegraphics[width=2cm]{Fig31.eps}}
\]
with $s$ parallel edges between the two circles.

\subsubsection{The upper bound}
We consider symmetric decorations $\ga$ of the form
\[
\parbox{4cm}{\psfrag{1}{$\scriptstyle{a_1}$}
\psfrag{2}{$\scriptstyle{a_2}$}\psfrag{3}{$\scriptstyle{a_s}$}
\psfrag{4}{$\scriptstyle{b_1}$}\psfrag{5}{$\scriptstyle{b_2}$}
\psfrag{6}{$\scriptstyle{b_3}$}\psfrag{7}{$\scriptstyle{b_s}$}
\includegraphics[width=4cm]{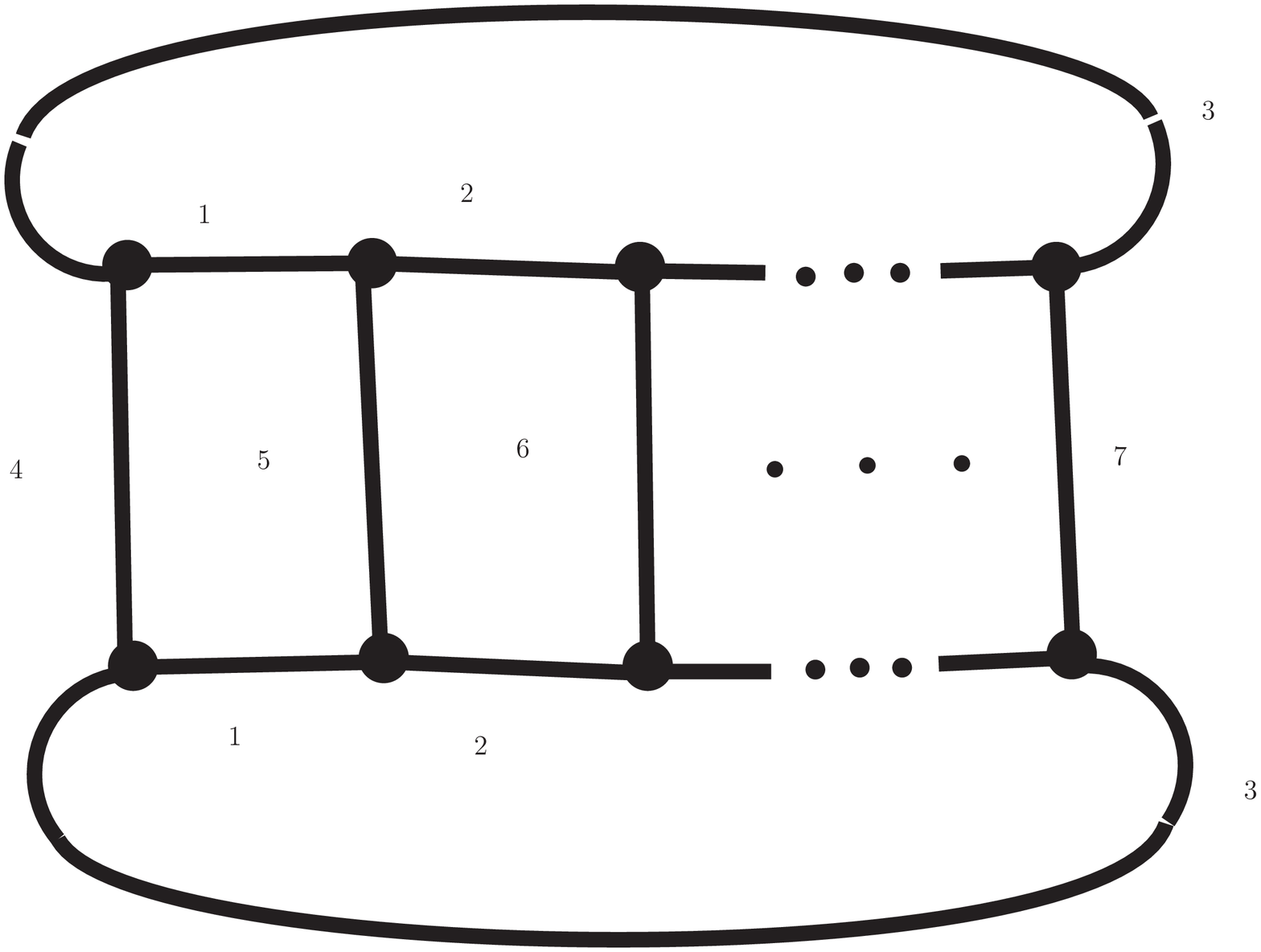}}\qquad .
\]
\begin{Lemma}\label{drumupperbd}
For any $s\ge 1$, and admissible decorations $\ga$ as above
\[
\left|\<{\rm Drum}_s,\ga\>^U\right|\le
\frac{[\min(a_1,\ldots,a_s)+1]^2}{(a_1+1)\cdots(a_s+1)}\ \ .
\]
In particular, the conclusion of Theorem \ref{mainthm} holds
in this case as soon as $s\ge 2$.
\end{Lemma}

\noindent{\bf Proof:}
By Corollary \ref{negdimU} and the indicated choice of smooth orientation
\begin{equation}
\left|\<{\rm Drum}_s,\ga\>^U\right|=
\prod\limits_{i=1}^{s}\frac{1}{a_i+1}
\times
\left|
\<
\parbox{4cm}{\psfrag{1}{$\scriptstyle{a_1}$}
\psfrag{2}{$\scriptstyle{a_2}$}\psfrag{3}{$\scriptstyle{a_s}$}
\psfrag{4}{$\scriptstyle{b_1}$}\psfrag{5}{$\scriptstyle{b_2}$}
\psfrag{6}{$\scriptstyle{b_3}$}\psfrag{7}{$\scriptstyle{b_s}$}
\includegraphics[width=4cm]{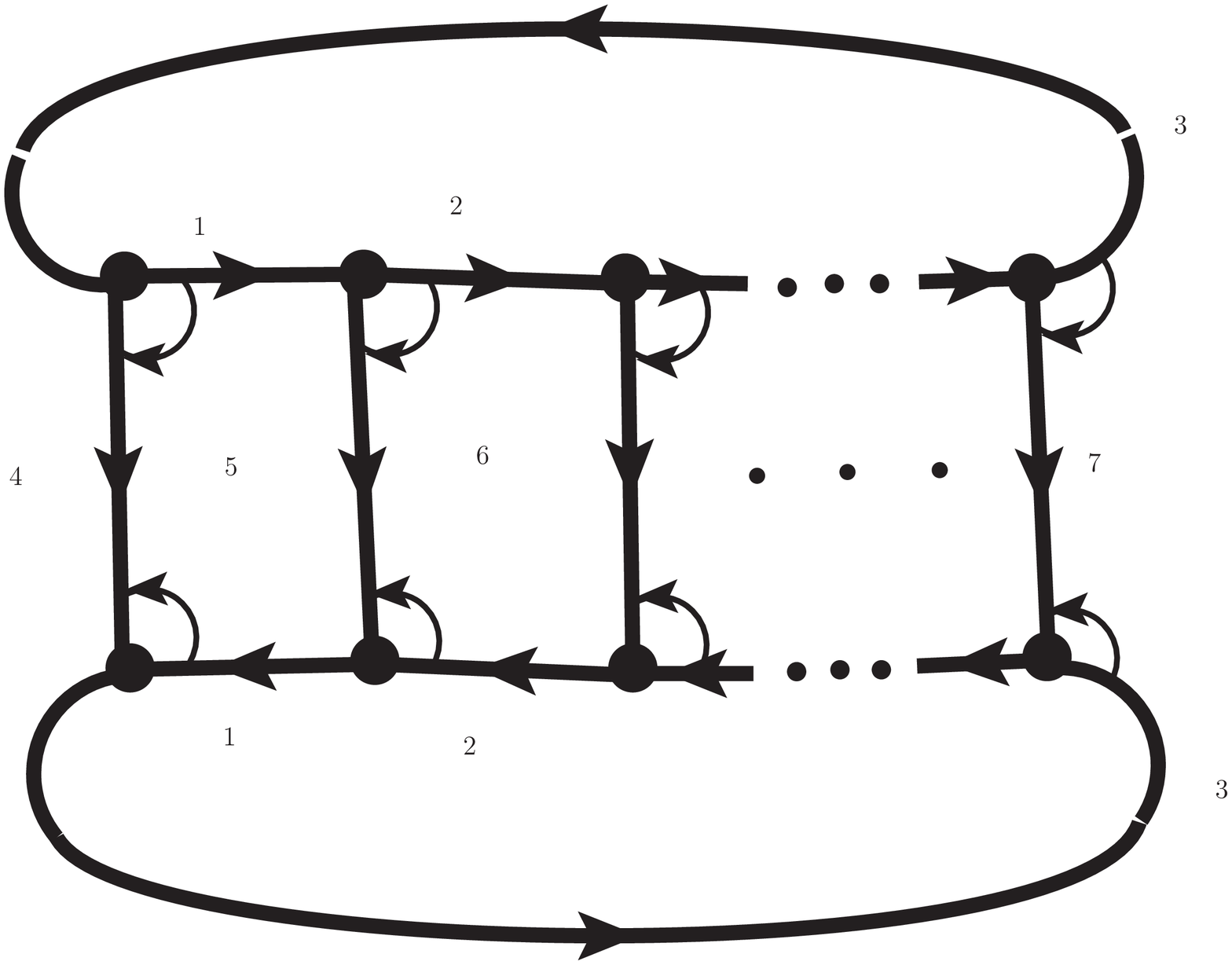}}
\>^{\pi\io}\right|\ \ .
\label{drumtopiio}
\end{equation}
We now perform for every $i$, $1\le i\le s$, the following
FDC calculation
\[
\<
\parbox{1cm}{\psfrag{a}{$\scriptstyle{a_i}$}
\includegraphics[width=1cm]{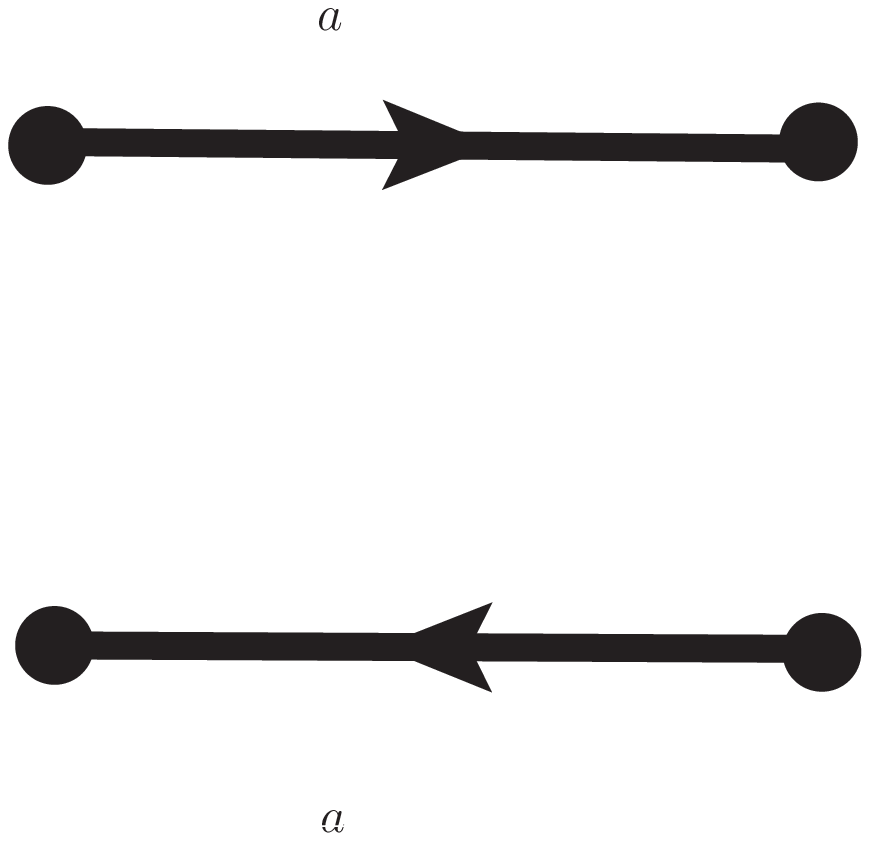}}
\>^{CG}=
\begin{array}{c}
{\scriptstyle{a_i}}\left\{
\parbox{1.5cm}{\includegraphics[width=1.5cm]{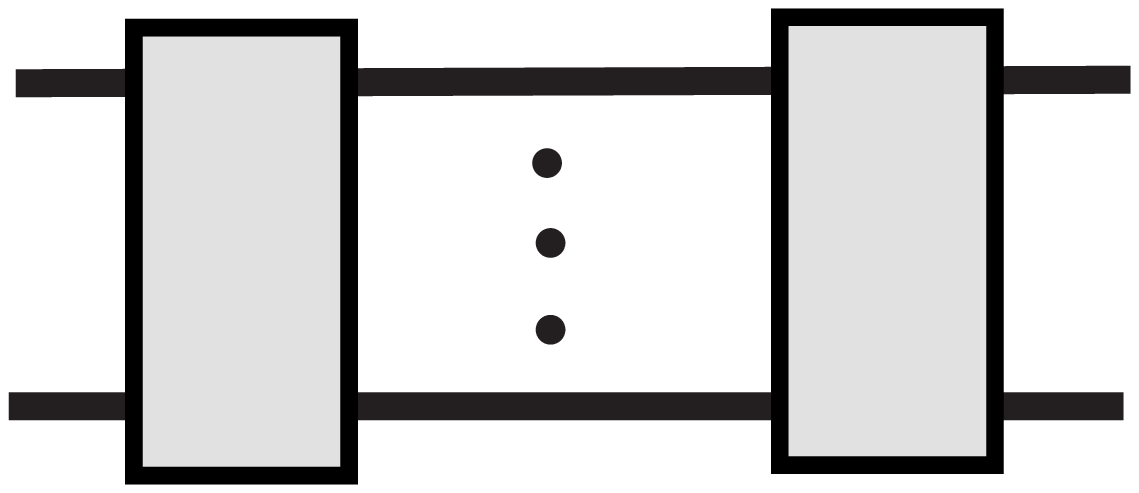}}\right.\\
\ \\
{\scriptstyle{a_i}}\left\{
\parbox{1.5cm}{\includegraphics[width=1.5cm]{Fig211.eps}}\right.
\end{array}
=\qquad
\parbox{2.5cm}{\includegraphics[width=2.5cm]{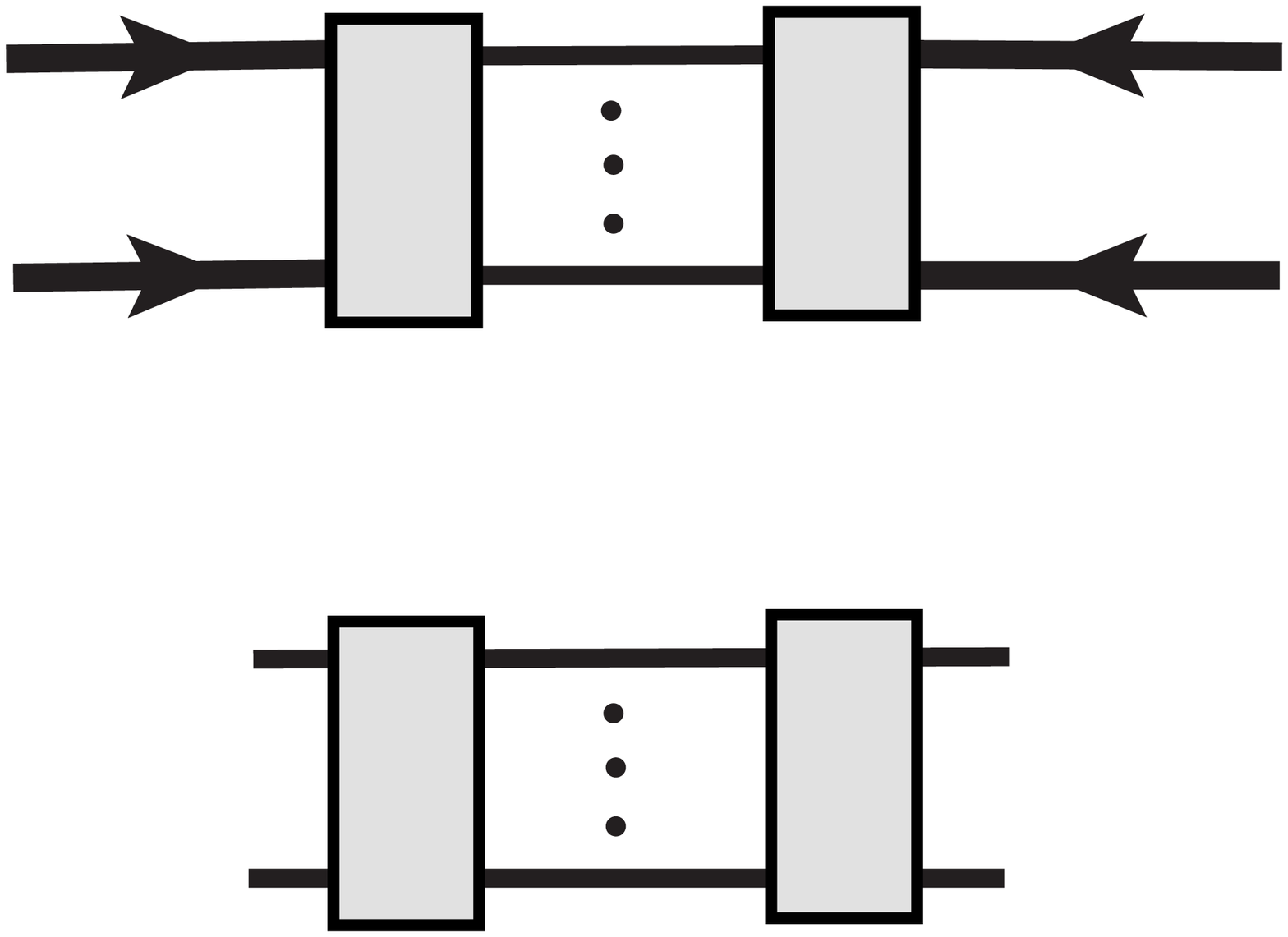}}\ \ .
\]
Then, we insert the identity (\ref{Gordanseries})
\[
\<
\parbox{1cm}{\psfrag{a}{$\scriptstyle{a_i}$}
\includegraphics[width=1cm]{Fig210.eps}}
\>^{CG}=
\sum\limits_{k_i=0}^{a_i}
\frac{
\left(
\begin{array}{c}
a_i\\
k_i
\end{array}
\right)^2
}{
\left(
\begin{array}{c}
2 a_i-k_i+1\\
k_i
\end{array}
\right)}
\ \ 
\parbox{3.6cm}{\psfrag{k}{$\scriptstyle{k_i}$}
\includegraphics[width=3.6cm]{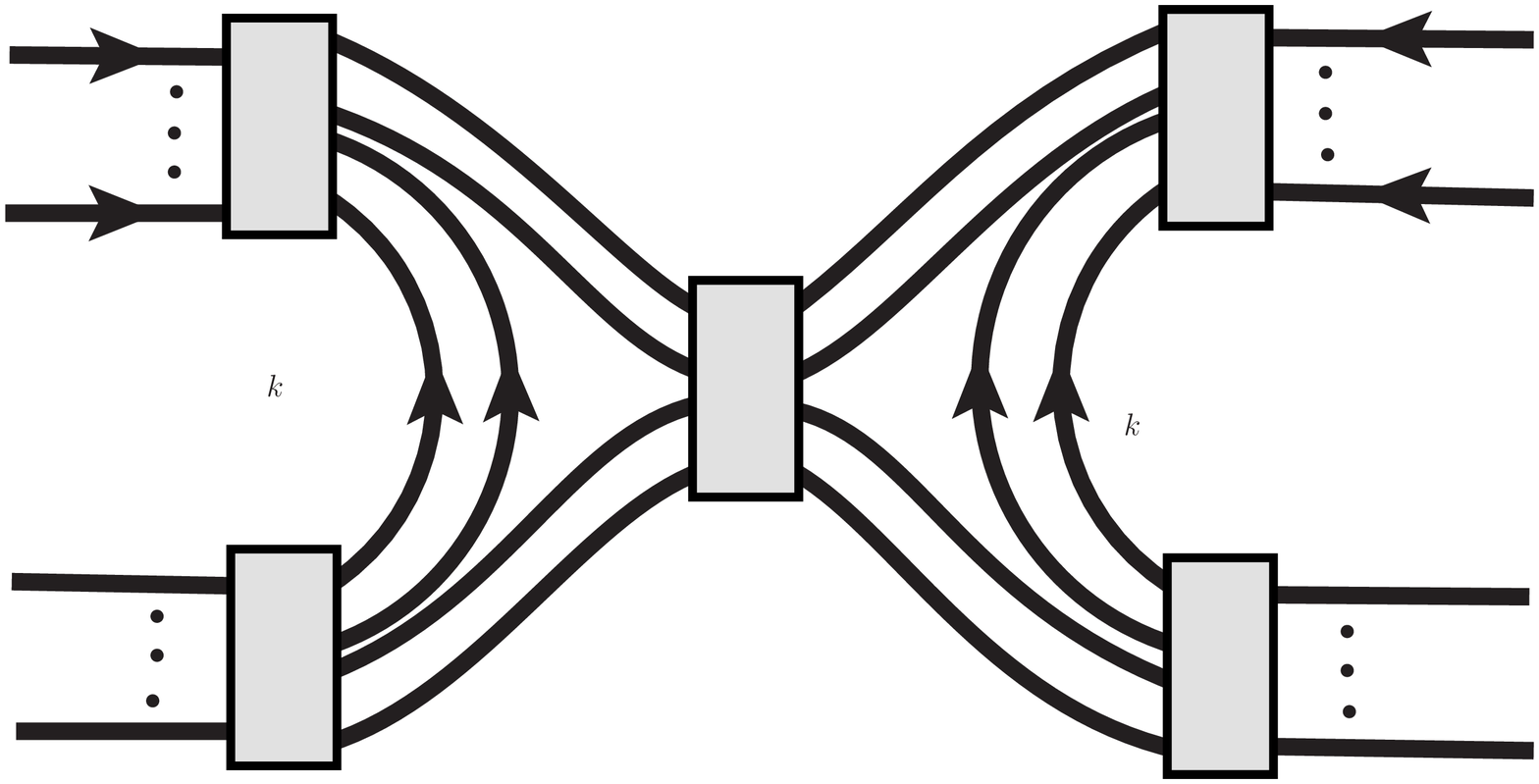}}
\]
\[
=\sum\limits_{k_i=0}^{a_i}
\frac{
\left(
\begin{array}{c}
a_i\\
k_i
\end{array}
\right)^2
}{
\left(
\begin{array}{c}
2 a_i-k_i+1\\
k_i
\end{array}
\right)}
\qquad
\parbox{4.6cm}{\psfrag{k}{$\scriptstyle{k_i}$}
\psfrag{c}{$\scriptstyle{c_i}$}
\includegraphics[width=4.6cm]{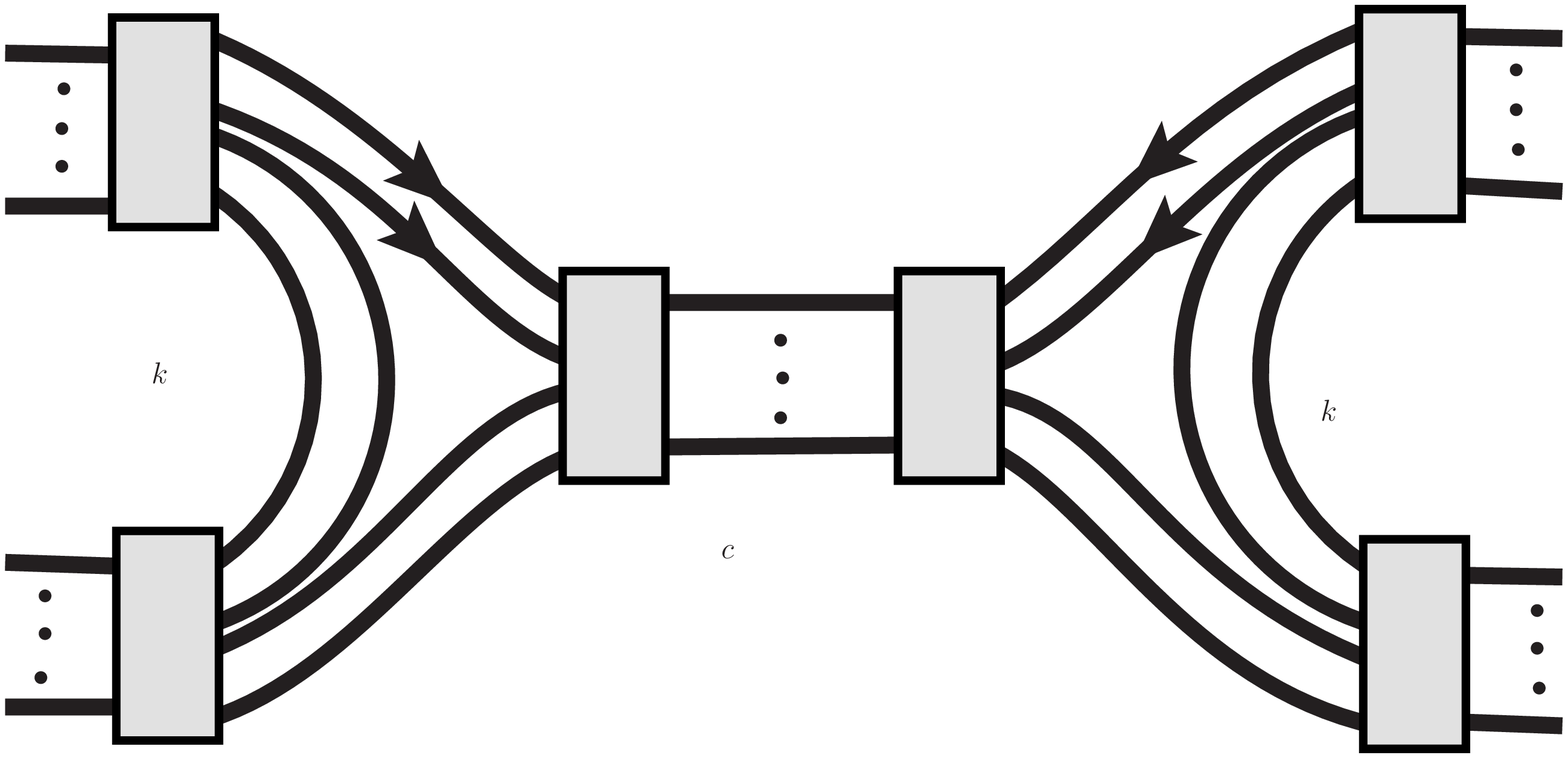}}
\]
by pushing the outside $\ep$ arrows through the facing
symmetrizers as in (\ref{matrixpush}). We also let $c_i=2a_i-k_i$.
Thus
\[
\<
\parbox{1cm}{\psfrag{a}{$\scriptstyle{a_i}$}
\includegraphics[width=1cm]{Fig210.eps}}
\>^{\pi\io}=\<
\parbox{1cm}{\psfrag{a}{$\scriptstyle{a_i}$}
\includegraphics[width=1cm]{Fig210.eps}}
\>^{CG}=
\sum\limits_{{c_i=0}\atop{c_i\ {\rm even}}}^{2a_i}
\frac{c_i+1}{a_i+1}
\<
\ \ \parbox{2.2cm}{\psfrag{a}{$\scriptstyle{a_i}$}
\psfrag{c}{$\scriptstyle{c_i}$}
\includegraphics[width=2.2cm]{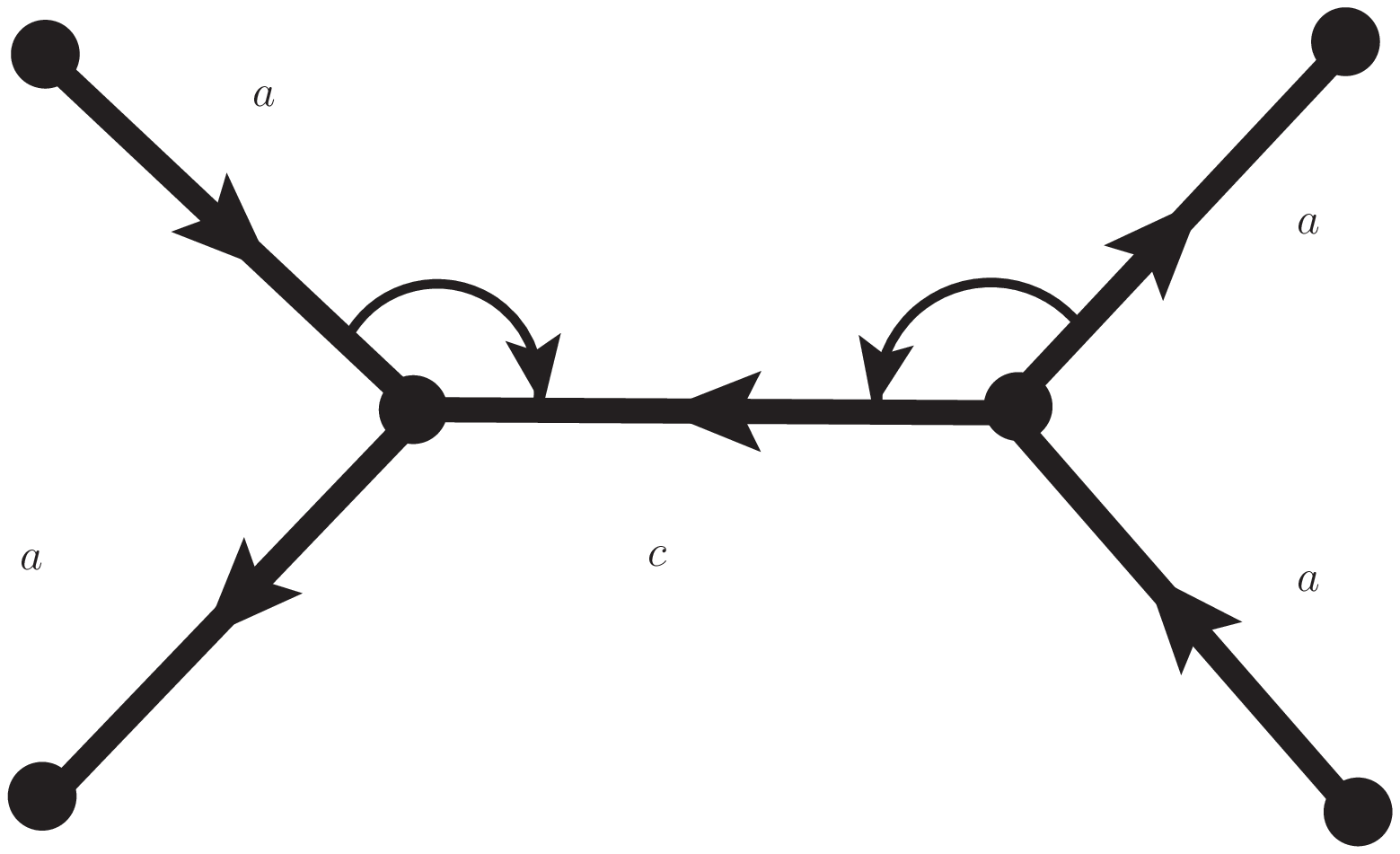}}
\ \ \>^{\pi\io}\ \ .
\]
We now insert this in (\ref{drumtopiio}) and get
\[
\left|\<{\rm Drum}_s,\ga\>^U\right|=
\prod\limits_{i=1}^{s}\frac{1}{(a_i+1)^2}
\times
\left|
\sum\limits_{{c_1,\ldots,c_s=0}\atop{c_i\ {\rm even}}}^{2a_i}
\prod\limits_{i=1}^{s}(c_i+1)\times\right.
\]
\[
\left.
\<
\parbox{7cm}{\psfrag{1}{$\scriptstyle{a_{i-1}}$}
\psfrag{2}{$\scriptstyle{a_i}$}
\psfrag{3}{$\scriptstyle{a_{i+1}}$}
\psfrag{4}{$\scriptstyle{b_i}$}
\psfrag{5}{$\scriptstyle{b_{i+1}}$}
\psfrag{6}{$\scriptstyle{c_{i-1}}$}
\psfrag{7}{$\scriptstyle{c_i}$}
\psfrag{8}{$\scriptstyle{c_{i+1}}$}
\includegraphics[width=7cm]{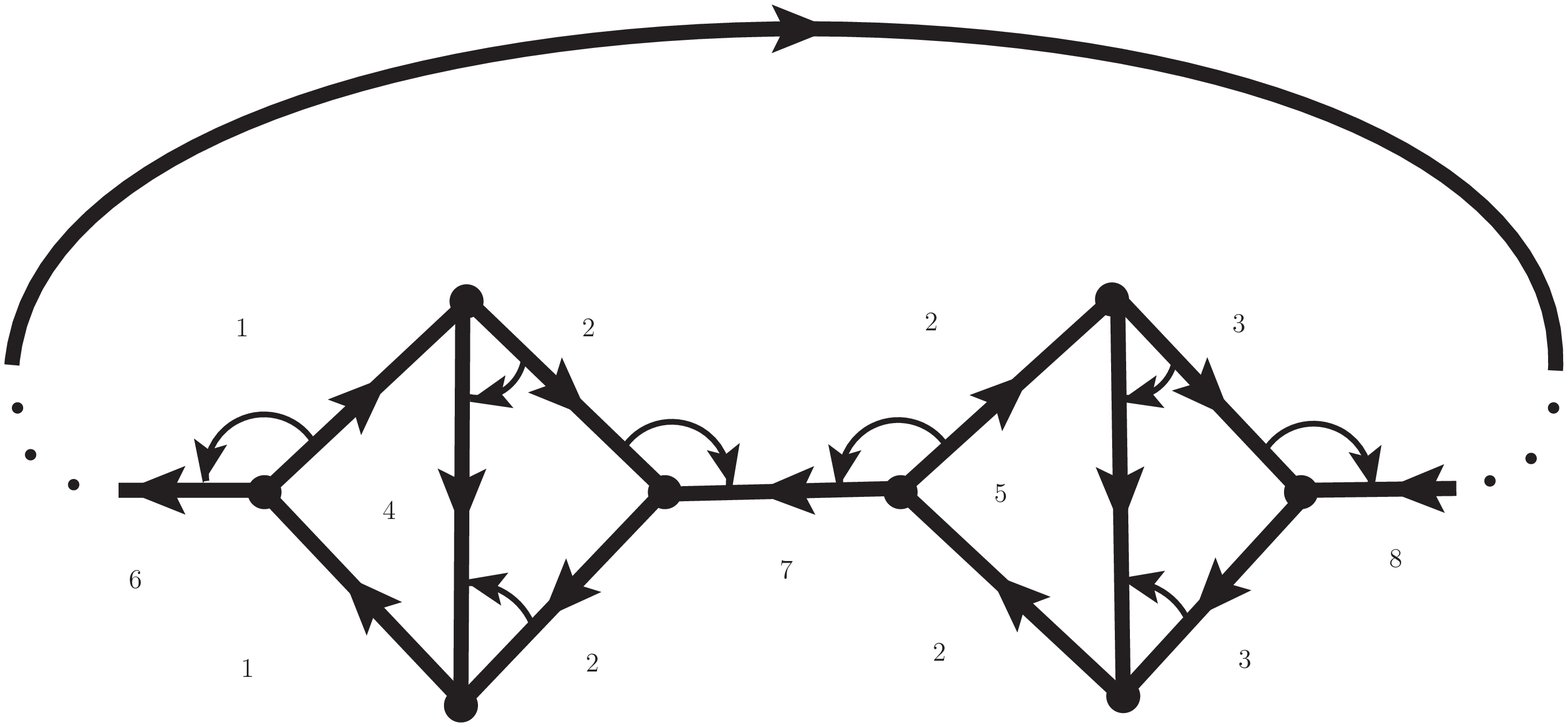}}
\>^{\pi\io}\right|\ \ .
\]
By Schur's Lemma, namely identity (\ref{Schur}),
the $c_i$ have to be the same and one gets
\[
\left|\<{\rm Drum}_s,\ga\>^U\right|=
\prod\limits_{i=1}^{s}\frac{1}{(a_i+1)^2}
\times
\left|
\sum\limits_{{c=0}\atop{c\ {\rm even}}}^{2\min(a_1,\ldots,a_s)}
(c+1)\times\right.
\]
\[
\left.
\times\prod\limits_{i=1}^{s}
\<
\parbox{3.6cm}{\psfrag{1}{$\scriptstyle{a_{i-1}}$}\psfrag{2}{$\scriptstyle{a_i}$}
\psfrag{b}{$\scriptstyle{b_i}$}\psfrag{c}{$\scriptstyle{c}$}
\includegraphics[width=3.6cm]{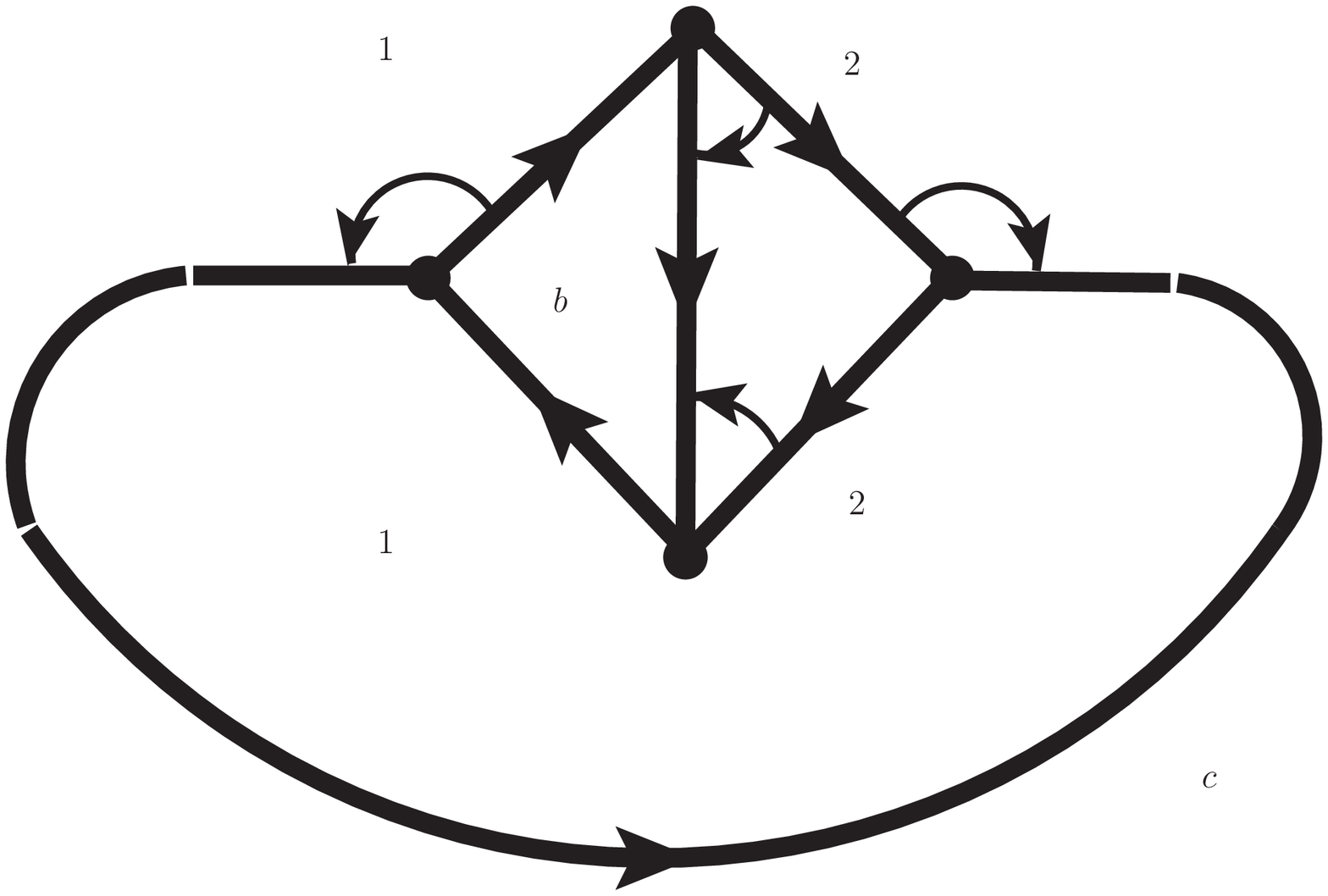}}
\>^{\pi\io}\right|
\]
with the convention $a_0=a_s$. Therefore, by Corollary \ref{negdimU}
\[
\left|\<{\rm Drum}_s,\ga\>^U\right|\le
\prod\limits_{i=1}^{s}\frac{1}{(a_i+1)^2}
\times
\sum\limits_{{c=0}\atop{c\ {\rm even}}}^{2\min(a_1,\ldots,a_s)} (c+1)
\times\prod\limits_{i=1}^{s}
\]
\[
\left[
(a_{i-1}+1)(a_i+1)
\left|\<\ \ 
\parbox{2.5cm}{\psfrag{1}{$\scriptstyle{a_{i-1}}$}\psfrag{2}{$\scriptstyle{a_i}$}
\psfrag{b}{$\scriptstyle{b_i}$}\psfrag{c}{$\scriptstyle{c}$}
\includegraphics[width=2.5cm]{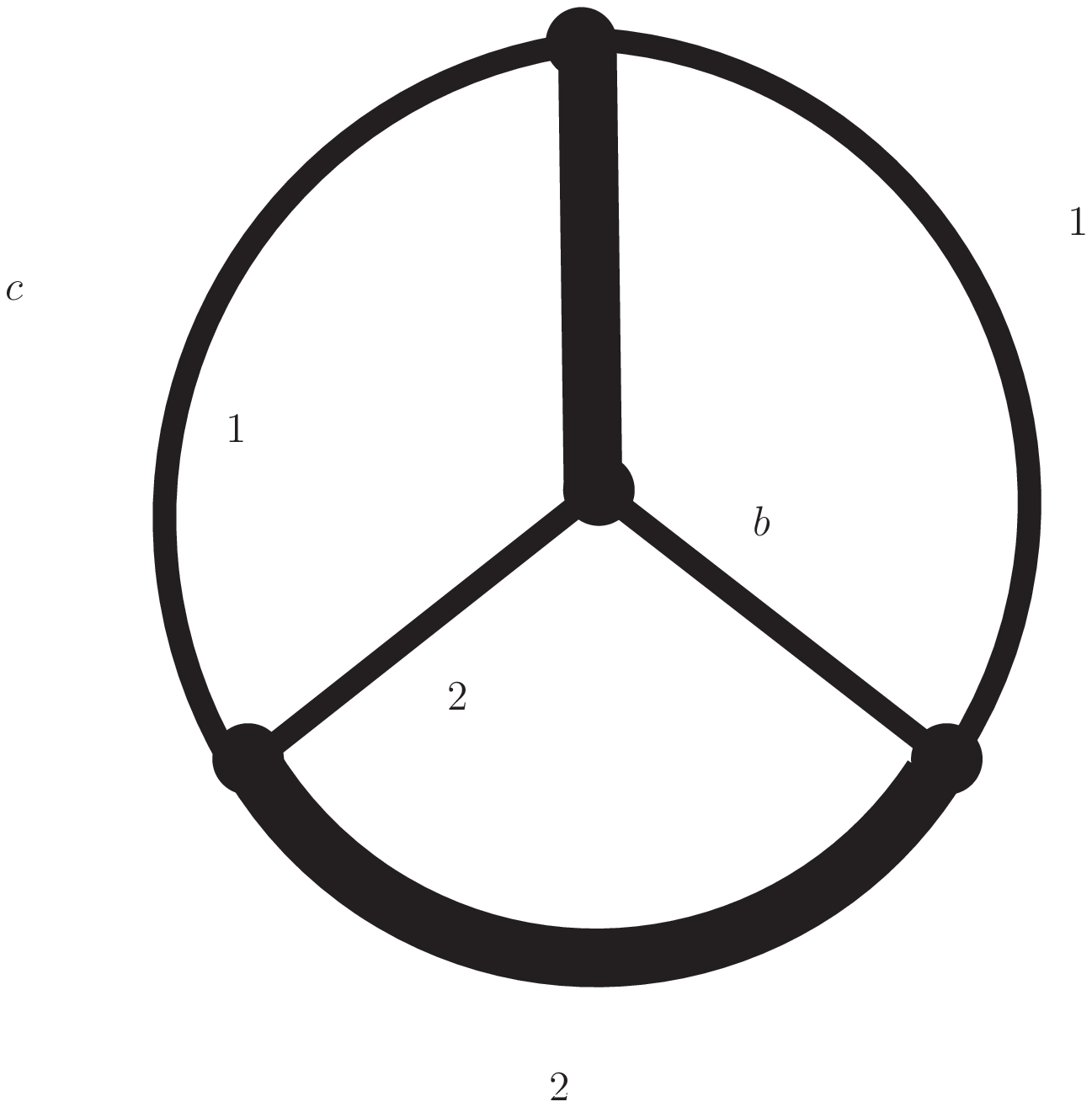}}
\ \ \>^U\right|
\right]
\]
\[
\le \sum\limits_{{c=0}\atop{c\ {\rm even}}}^{2\min(a_1,\ldots,a_s)} (c+1)
\times\prod\limits_{i=1}^{s}
\frac{1}{\sqrt{(a_{i-1}+1)(a_i+1)}}
=\frac{[\min(a_1,\ldots,a_s)+1]^2}{(a_1+1)\cdots(a_s+1)}
\]
by (\ref{6jupperbd}) with the bound corresponding to the highlighted
pair of opposite edges.
\qed

\subsubsection{The lower bound}
We now consider the uniform decoration by $2n$ for ${\rm Drum}_s$,
$s\ge 2$.
\begin{Proposition}
For $s\ge 2$
\[
\limsup_{n\rightarrow\infty}
\left|\<{\rm Drum}_s,2n\>^U\right|^{\frac{1}{n}}=1\ .
\]
\end{Proposition}

\noindent{\bf Proof:}
The limsup is $\le 1$ by Lemma \ref{drumupperbd}.
By (\ref{drumtopiio}) and the same argument as in \S\ref{6jupsec} we have
\[
\left|\<{\rm Drum}_s,2n\>^U\right|=\frac{1}{(2n+1)^s}\times
\<A|A\>_{\cH}
\]
where $\cH=\cH_{2n}^{\otimes s}$
and
\[
A=\<
\underbrace{\parbox{4cm}{\includegraphics[width=4cm]{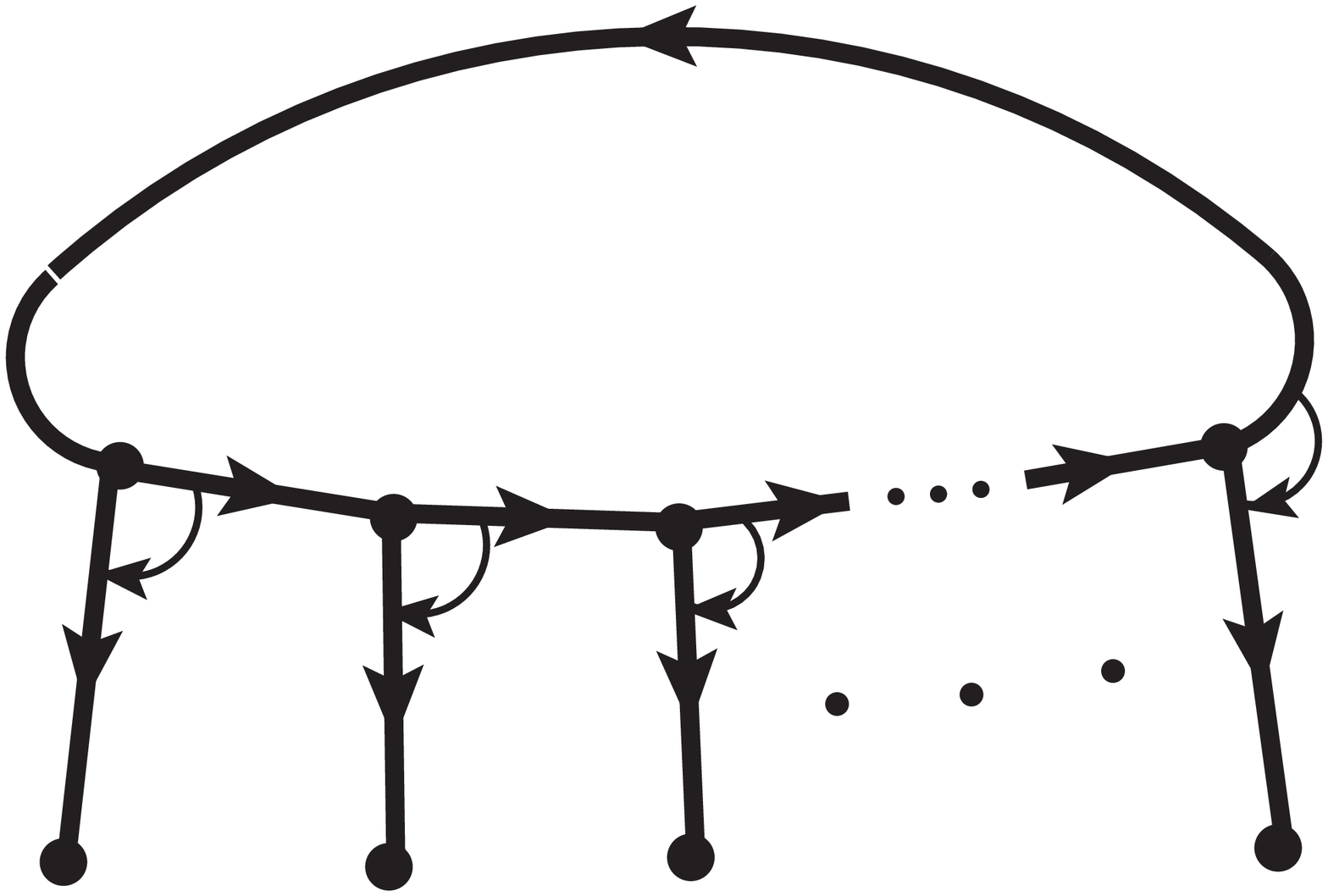}}}_{s}
\>^{\pi\io}
\in\cH
\]
with decorations by $2n$ everywhere.
By the Cauchy-Schwarz inequality
\[
\<A|A\>_{\cH}\ge
\frac{\<A|B\>_{\cH}^2}{\<B|B\>_{\cH}}
\]
for any nonzero $B\in\cH$.
We will obtain adequate lower bounds by making a judicious ansatz
for $B$. For this, we need to distinguish two cases.

\noindent{\bf Case 1 where $s$ is even:}
Take $B$ to be given by the microscopic FDC
formula
\[
B=
\overbrace{
\parbox{1.7cm}{\psfrag{2}{$\scriptstyle{2n}$}
\includegraphics[width=1.7cm]{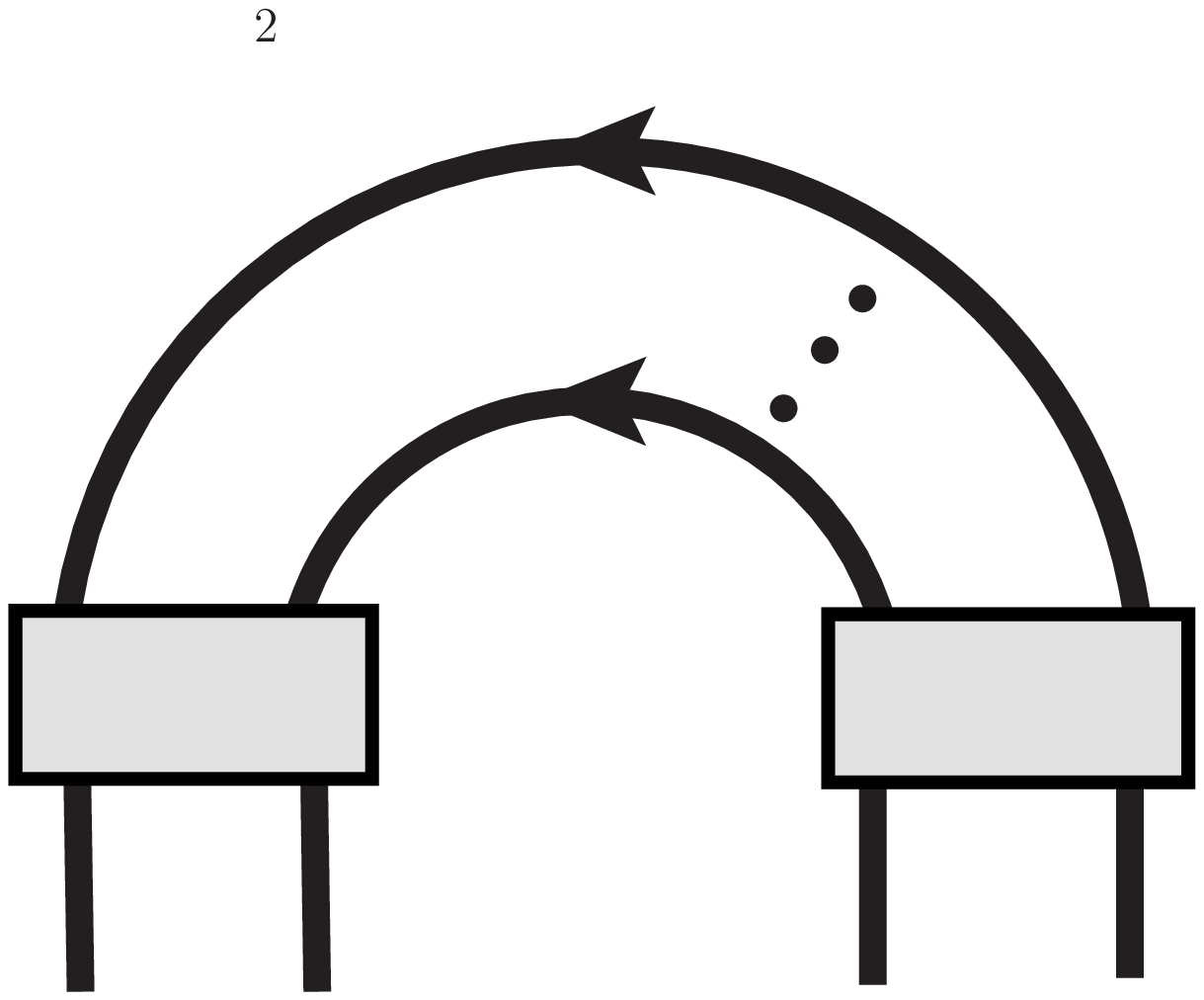}}\qquad .
\ 
\parbox{1.7cm}{\includegraphics[width=1.7cm]{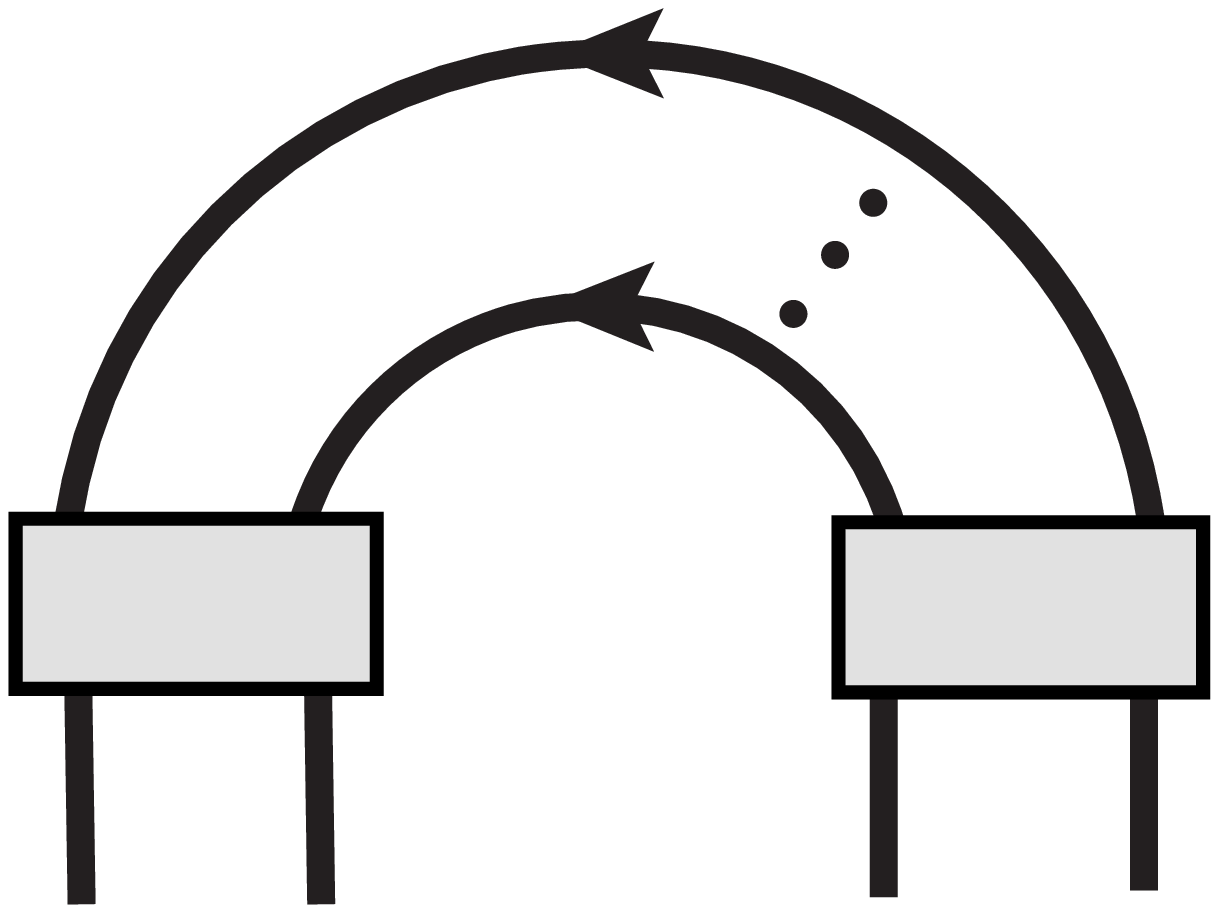}}
\ \cdots\ 
\parbox{1.7cm}{\includegraphics[width=1.7cm]{Fig220bis.eps}}
}^{\frac{s}{2}}
\]
Since by pushing the $\ep$ arrows as indicated
\[
\parbox{5.5cm}{\psfrag{2}{$\scriptstyle{2n}$}
\psfrag{p}{$\scriptstyle{\rm push}$}
\psfrag{i}{$\begin{array}{c}
{\scriptstyle{\rm inner}}\\
{\scriptstyle{\rm product}}
\end{array}$}
\includegraphics[width=5cm]{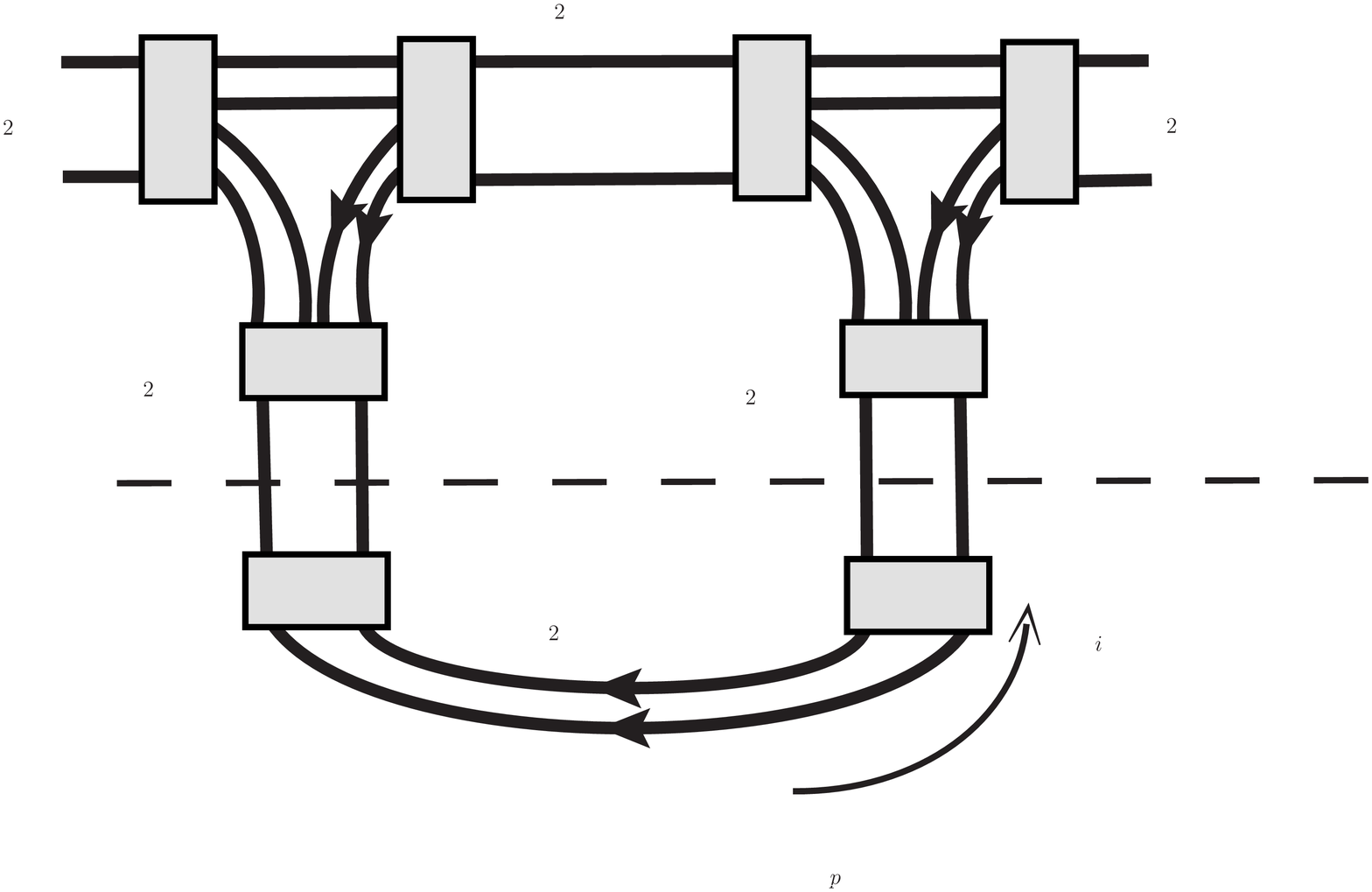}}
=(-1)^n\times
\parbox{5.5cm}{\psfrag{2}{$\scriptstyle{2n}$}
\includegraphics[width=5cm]{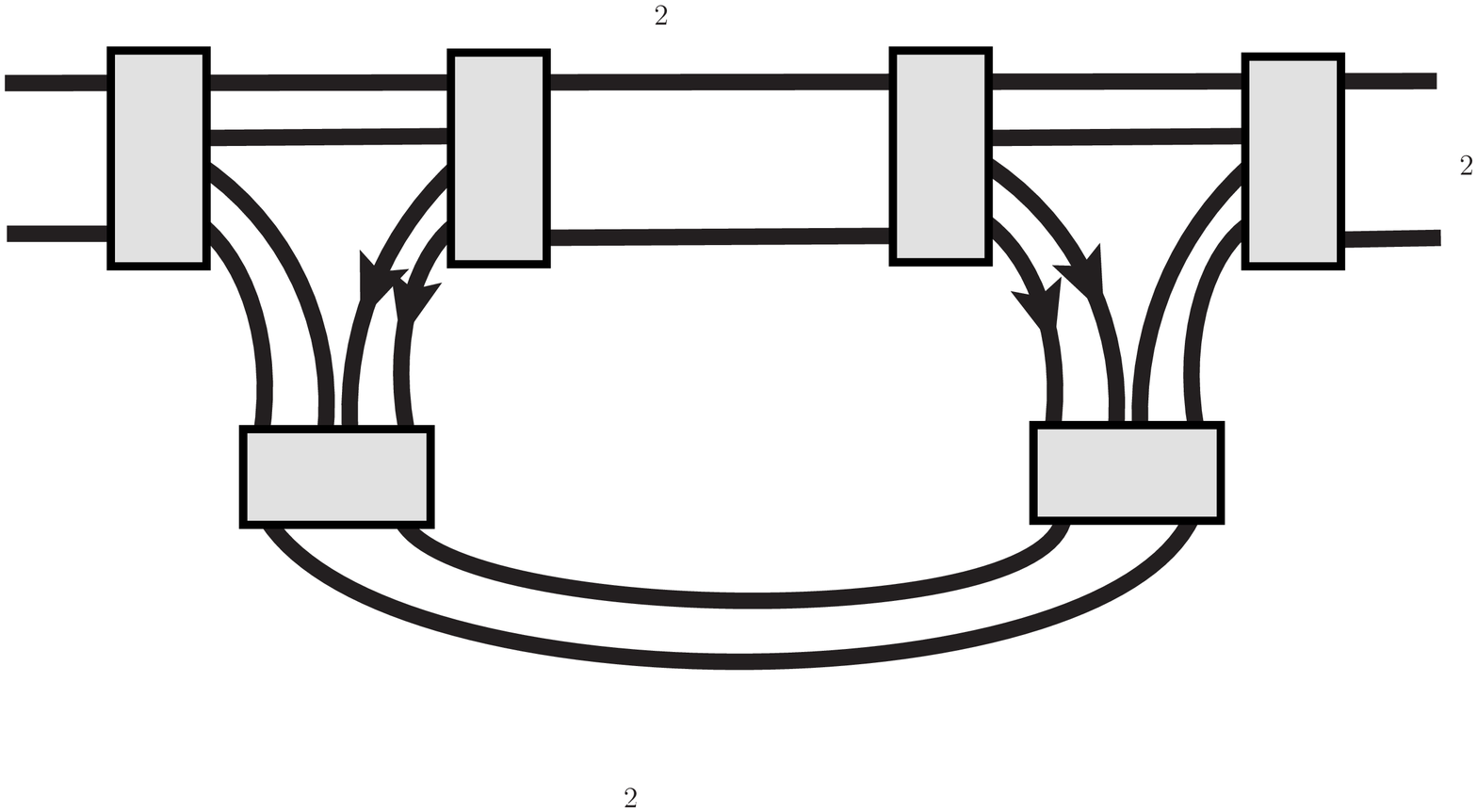}}
\]
one has
\[
\<A|B\>_{\cH}=
(-1)^{\frac{ns}{2}}\times
\<
\parbox{5cm}{\includegraphics[width=5cm]{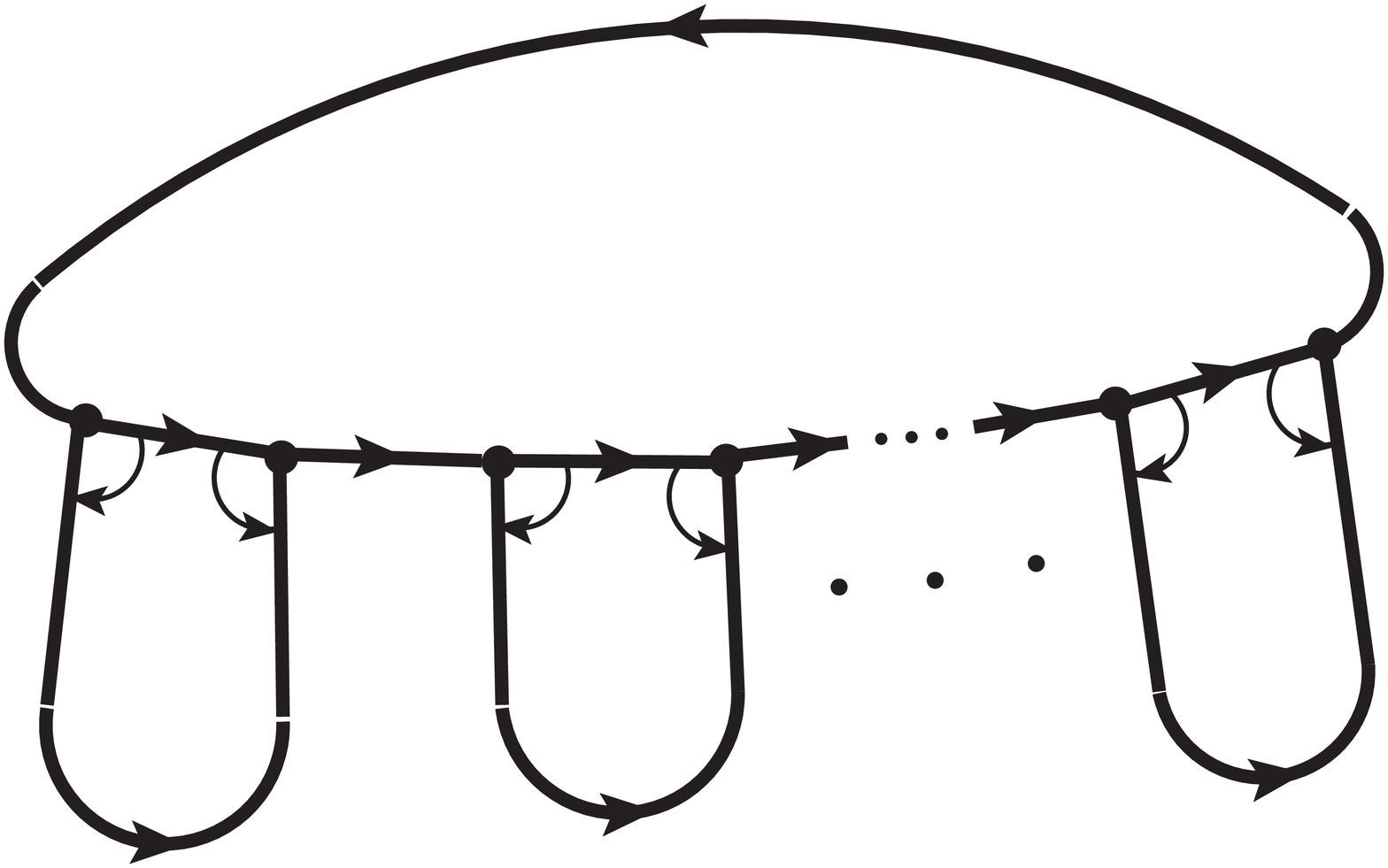}}
\>^{\pi\io}
\]
with $2n$ everywhere.
By Part 3) of Proposition (\ref{piiotaprop})
\[
|\<A|B\>_{\cH}|=\left|
\<
\parbox{1.5cm}{\psfrag{2}{$\scriptstyle{2n}$}
\includegraphics[width=1.5cm]{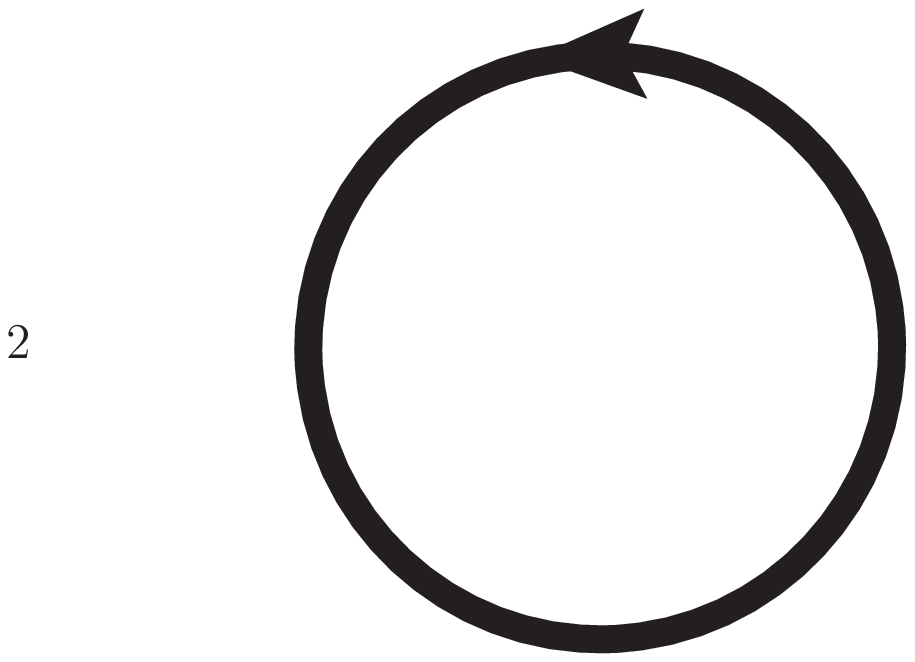}}
\>^{\pi\io}\right|
=2n+1\ .
\]
On the other hand,
\[
\<B|B\>_{\cH}=\left[
\parbox{1.5cm}{\psfrag{2}{$\scriptstyle{2n}$}
\includegraphics[width=1.5cm]{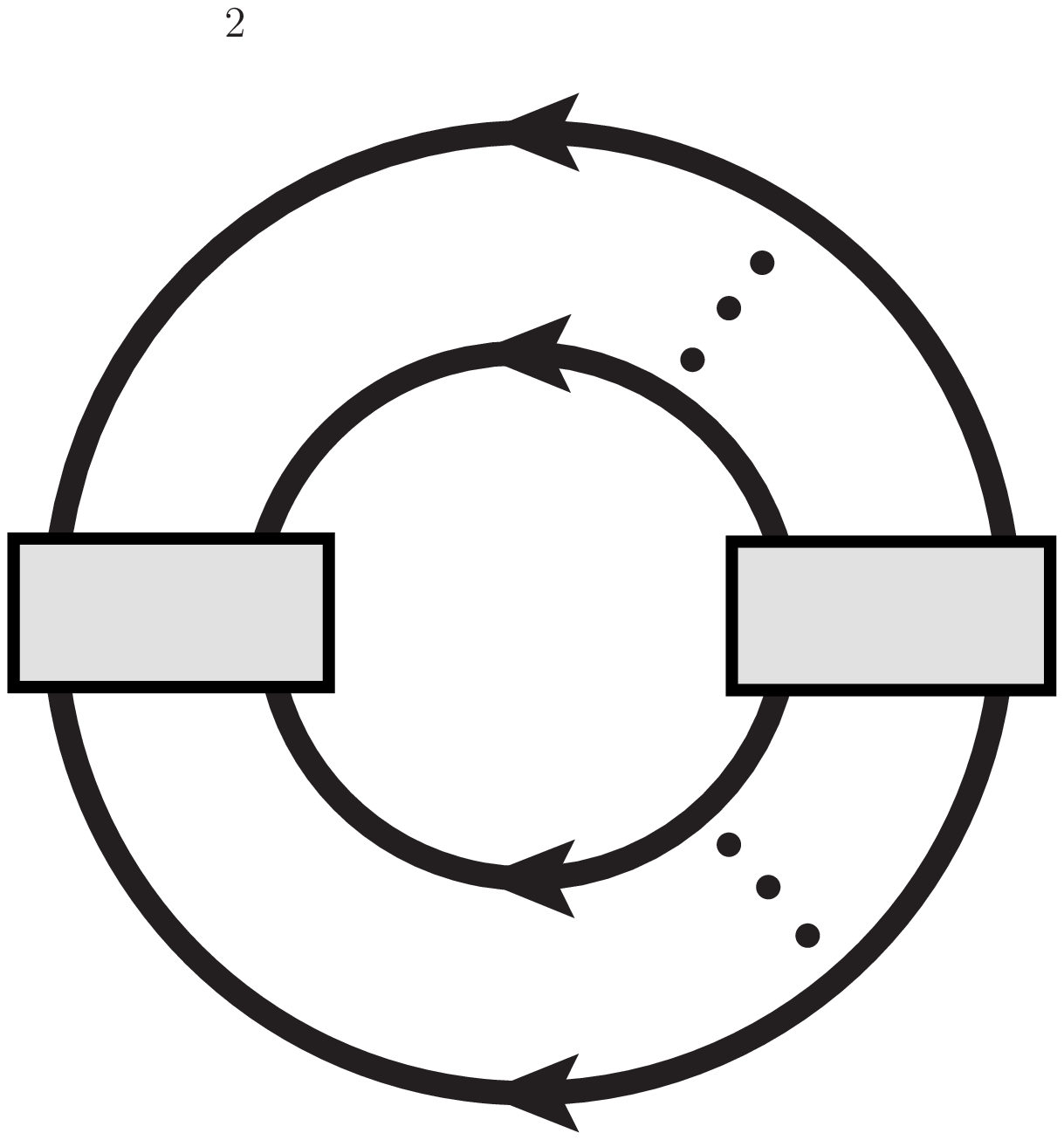}}
\right]^{\frac{s}{2}}
=(2n+1)^{\frac{s}{2}}>0\ .
\]
Finally,
\[
\left|\<{\rm Drum}_s,2n\>^U\right|\ge
\frac{1}{(2n+1)^s}\times\frac{(2n+1)^2}{(2n+1)^{\frac{s}{2}}}=(2n+1)^{-\frac{3s}{2}+2}
\]
and $\limsup_{n\rightarrow\infty}
\left|\<{\rm Drum}_s,2n\>^U\right|^{\frac{1}{n}}=1$ follows.

\noindent{\bf Case 2 where $s$ is odd:}
Since one must have $s\ge 3$, we will take
\[
B=
\parbox{3.5cm}{\psfrag{2}{$\scriptstyle{2n}$}
\includegraphics[width=3.5cm]{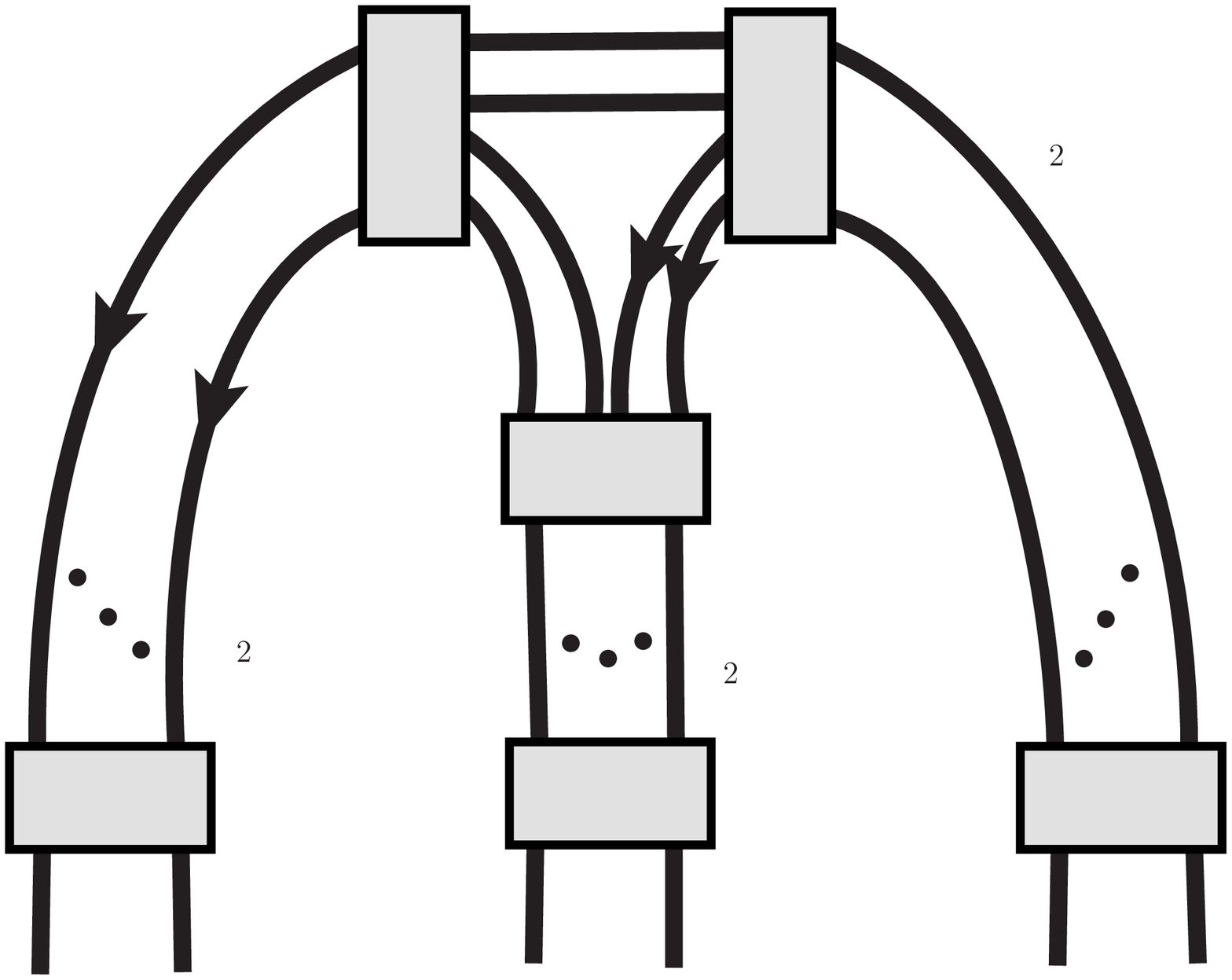}}
\ \ 
\parbox{1.5cm}{\psfrag{2}{$\scriptstyle{2n}$}
\includegraphics[width=1.5cm]{Fig220.eps}}
\ \ \cdots\ \ 
\parbox{1.5cm}{\psfrag{2}{$\scriptstyle{2n}$}
\includegraphics[width=1.5cm]{Fig220.eps}}\qquad .
\]
Then, by pushing the epsilons, we have
\[
\parbox{12cm}{\psfrag{2}{$\scriptstyle{2n}$}
\psfrag{p}{$\scriptstyle{\rm push}$}
\psfrag{u}{${\rm use}\ \ep=\ep^{\rm T}\ep^2$}
\psfrag{t}{${{\scriptstyle{\rm then}}\atop{\scriptstyle{\rm push}}}$}
\psfrag{i}{${{\scriptstyle{\rm inner}}\atop{\scriptstyle{\rm product}}}$}
\includegraphics[width=7cm]{Fig227.eps}}
\]
\[
\ 
\]
\[
=\qquad
\parbox{10cm}{\psfrag{2}{$\scriptstyle{2n}$}
\includegraphics[width=8cm]{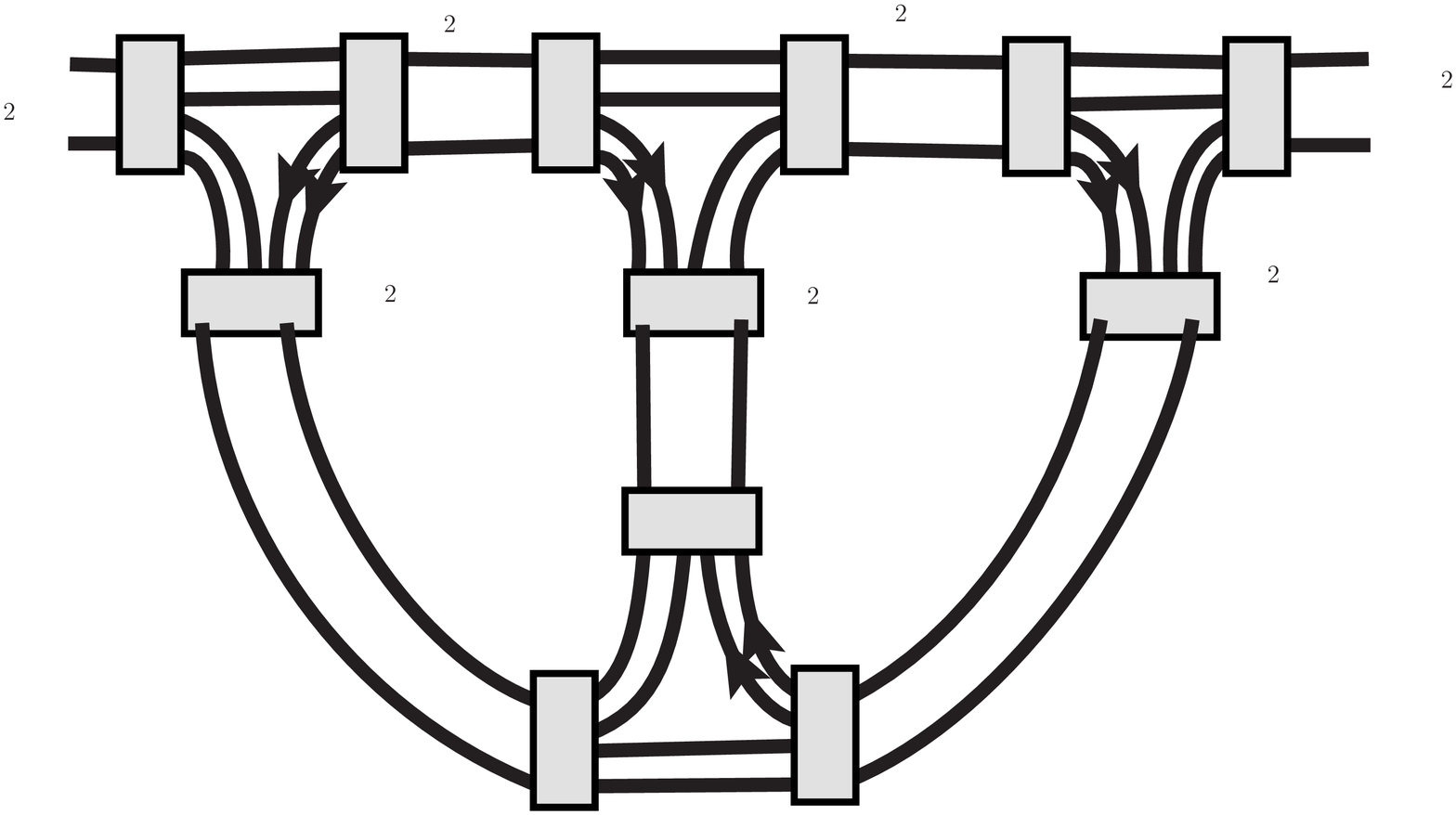}}\ \ .
\]
Hence
\[
\<A|B\>_{\cH}=(-1)^{\frac{n(s-3)}{2}}
\<
\parbox{6cm}{\includegraphics[width=6cm]{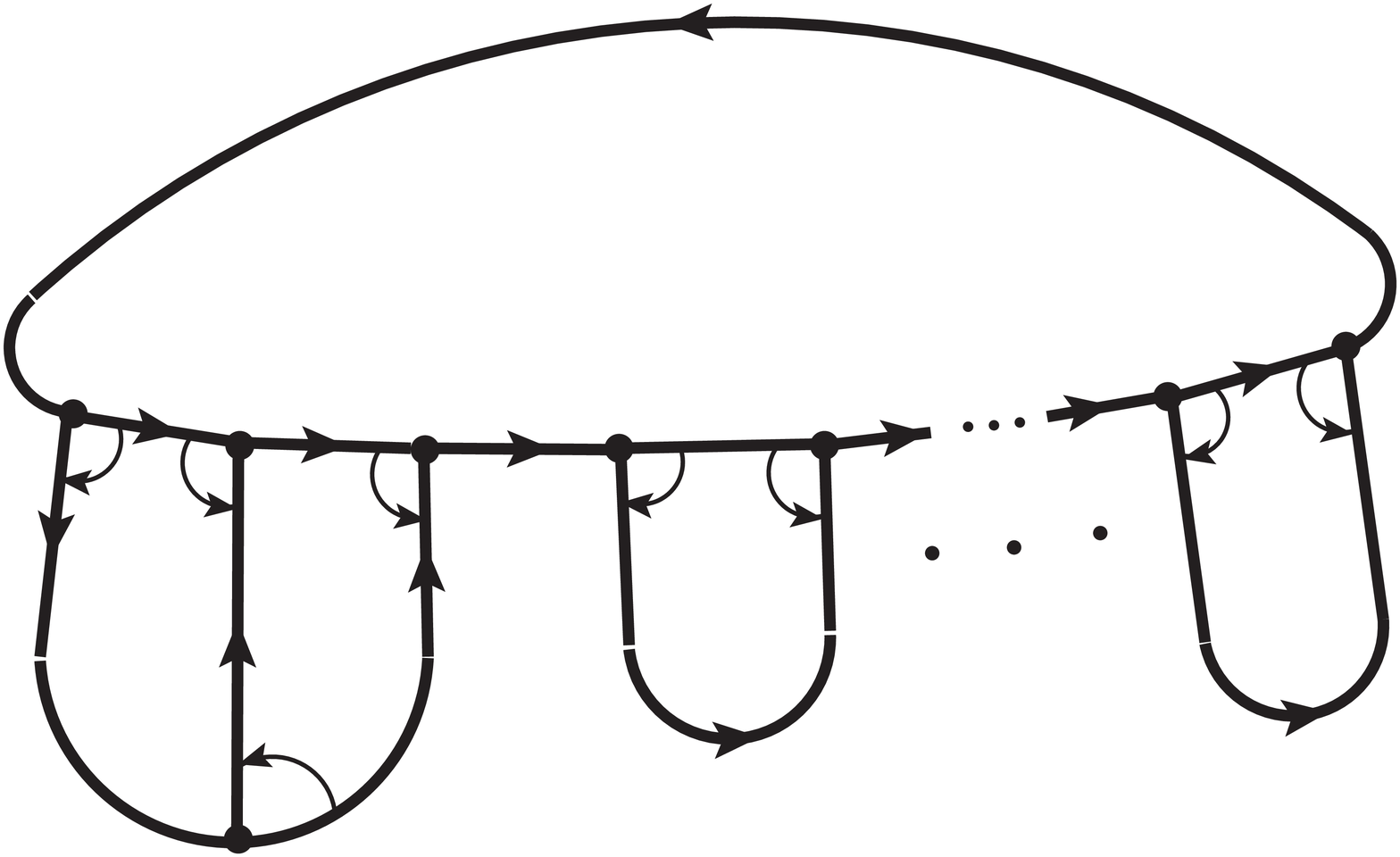}}
\>^{\pi\io}
\]
with $2n$ everywhere.
By Part 3) of Proposition \ref{piiotaprop} and Corollary \ref{negdimU}
\[
\<A|B\>_{\cH}^2=\left[\<
\parbox{2.5cm}{\psfrag{2}{$\scriptstyle{2n}$}
\includegraphics[width=2.5cm]{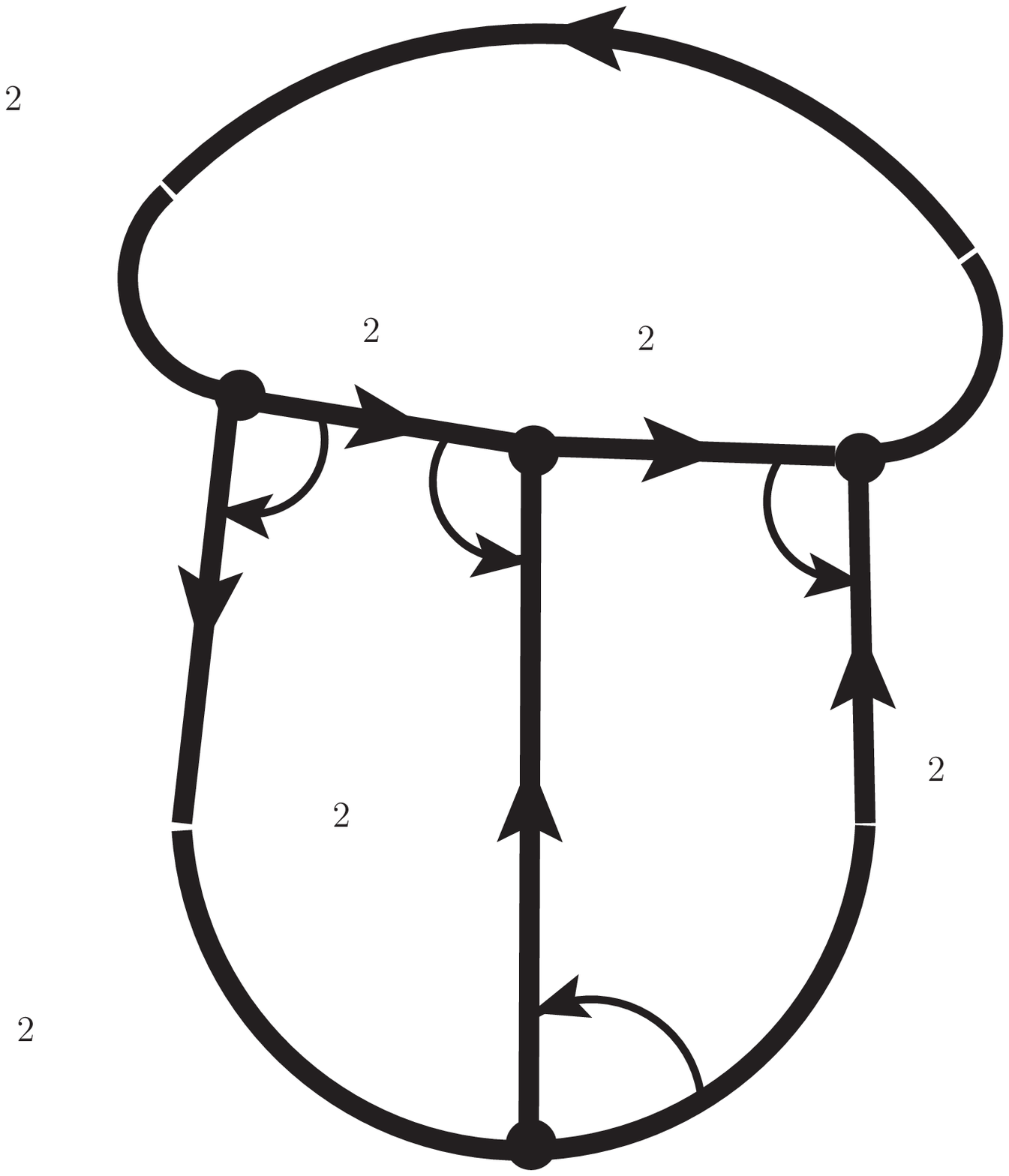}}
\>^{\pi\io}\right]^2
\]
\[
=(2n+1)^4\times\left|\<
\parbox{1.3cm}{
\includegraphics[width=1.3cm]{Fig207.eps}}
,2n\>^U\right|^2\ \ .
\]
Whereas
\[
\<B|B\>_{\cH}=(2n+1)^{\frac{s-3}{2}}\times \<
\parbox{2.5cm}{\psfrag{2}{$\scriptstyle{2n}$}
\includegraphics[width=2.5cm]{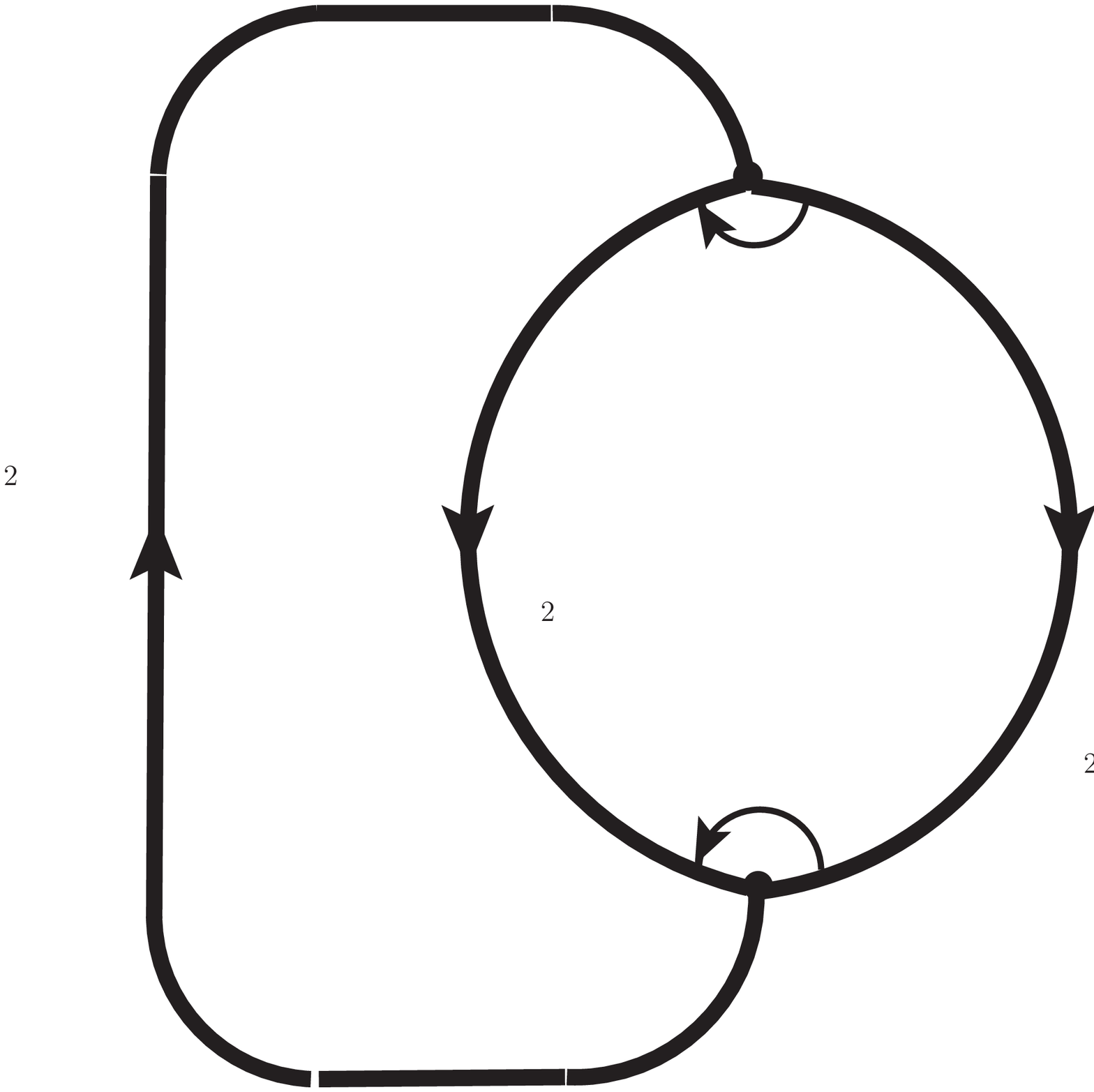}}
\>^{\pi\io}
\]
\[
=(2n+1)^{s-1}>0\ .
\]
Therefore
\[
\left|\<{\rm Drum}_s,2n\>^U\right|\ge
\frac{1}{(2n+1)^s}\times (2n+1)^4\times
\left|\<
\parbox{1cm}{
\includegraphics[width=1cm]{Fig207.eps}}
,2n\>^U\right|^2\times\frac{1}{(2n+1)^{\frac{s-1}{2}}}
\]
and $\limsup_{n\rightarrow\infty}
\left|\<{\rm Drum}_s,2n\>^U\right|^{\frac{1}{n}}\ge 1$
follows from Lemma \ref{PRlemma}.
\qed

\noindent
This completes the proof of Theorem \ref{drumthm}.

\section{Proof of the main theorem}\label{mainthmsec}

We now tackle the proof of Theorem \ref{mainthm}
using an elaboration of the Cauchy-Schwarz argument
of \S\ref{6jupsec}.
The key ingredient is the smooth orientation $\cO$.
It is enough to consider connected spin networks $(\Ga,\ga)$
without loops and which are not reduced to a trivial component.
Let $G$ be the underlying cubic graph and choose a smooth orientation
$\cO$ of $G$ by Proposition \ref{smoothprop},
as well as a gate signage $\ta$.
By Corollary \ref{negdimU}
\[
\left|\<\Ga,\ga\>^U\right|
=\left(\prod_{v\in V(G)} \frac{1}{\sqrt{{\rm dim}(v)}}\right)
\times
\left|\<G,\cO,\ta,\ga\>^{\pi\io}\right|\ \ .
\]
Now let $V_{\pi}(G)\subset V(G)$
be the set of vertices with indegree 2 and outdegree 1.
Let $V_{\io}(G)\subset V(G)$ be the complement, i.e.,
the set of vertices with indegree 1 and outdegree 2.
Since there are as many incoming half-edges as there are outgoing
half-edges
\[
2|V_{\pi}(G)|+|V_{\io}(G)|=|V_{\pi}(G)|+2|V_{\io}(G)|
\]
and thus
\[
|V_{\pi}(G)|=|V_{\io}(G)|=\frac{|V(G)|}{2}\ .
\]
By pulling the vertices of $V_{\pi}(G)$
on one side and those of $V_{\io}(G)$ on the other, we have
\[
\left|\<\Ga,\ga\>^U\right|
=\left(\prod_{v\in V(G)} \frac{1}{\sqrt{{\rm dim}(v)}}\right)
\times
\left|\<
\parbox{2.5cm}{\psfrag{i}{$V_{\io}(G)$}
\psfrag{p}{$V_{\pi}(G)$}
\includegraphics[width=2.5cm]{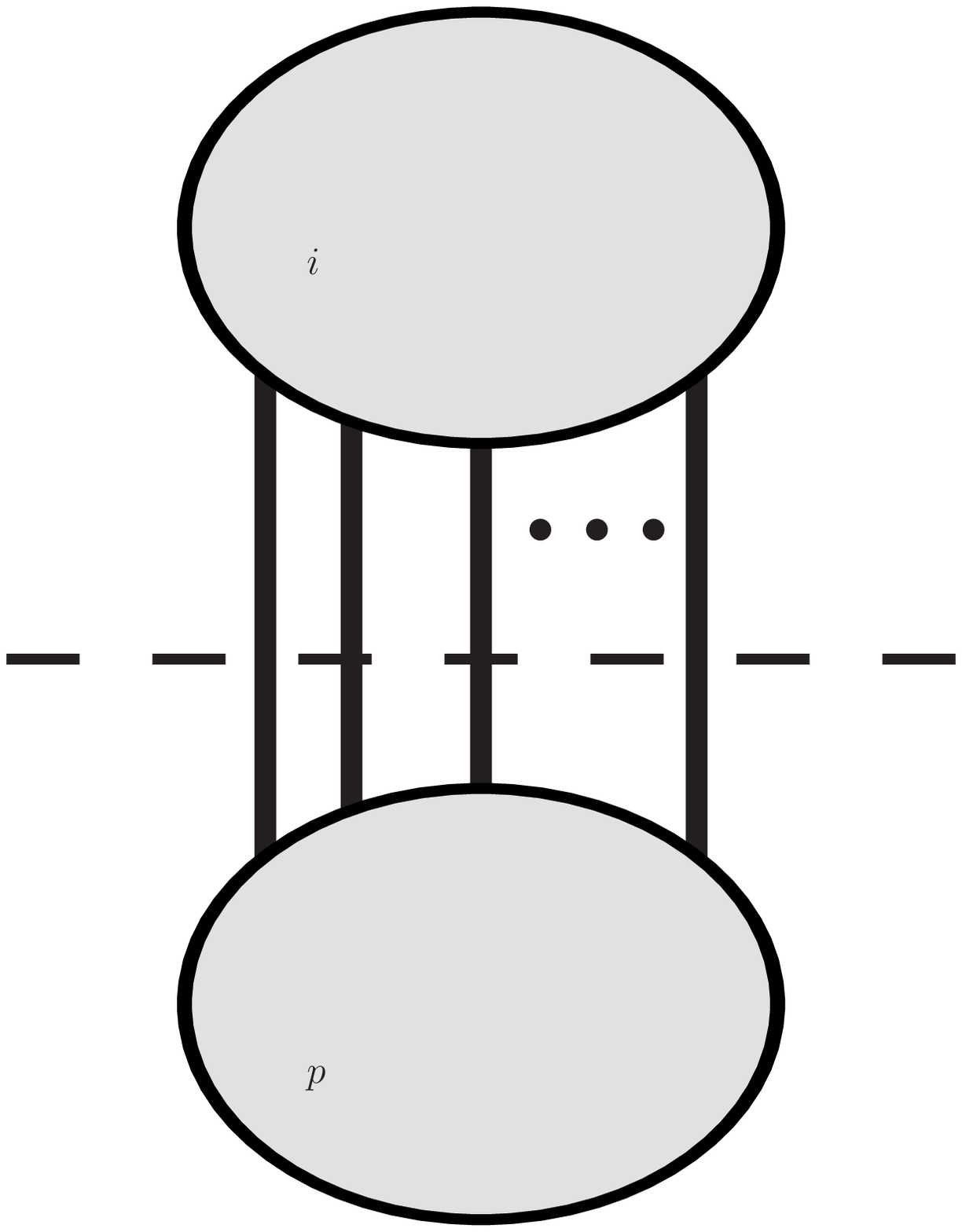}}
\>^{\pi\io}\right|
\]
where we omitted the orientations of the crossing edges in the picture,
since some edges go up and others go down.
Note that the $V_{\pi}(G)$ and $V_{\io}(G)$
induced subgraphs need not be connected.
There is a Hilbert space $\cH$ associated to the splitting indicated by the
dotted line.
One can therefore use the Cauchy-Schwarz inequality
as in \S\ref{6jupsec}, with the effect that
\[
\left|\<\Ga,\ga\>^U\right|
\le\left(\prod_{v\in V(G)} \frac{1}{\sqrt{{\rm dim}(v)}}\right)
\times
\left\{
\<
\parbox{2.5cm}{\psfrag{i}{$V_{\io}(G)$}
\psfrag{p}{$\widetilde{V_{\io}(G)}$}
\includegraphics[width=2.5cm]{Fig232.eps}}
\>^{\pi\io}\times
\<
\parbox{2.5cm}{\psfrag{i}{$\widetilde{V_{\pi}(G)}$}
\psfrag{p}{$V_{\pi}(G)$}
\includegraphics[width=2.5cm]{Fig232.eps}}
\>^{\pi\io}
\right\}^{\frac{1}{2}}
\]
where $\widetilde{V_{\io}(G)}$ subgraph is the mirror image 
of the $V_{\io}(G)$, and likewise for $\widetilde{V_{\pi}(G)}$.
This means that the decorations and gate signage are preserved, but
the edge orientations are reversed.
Now consider a graph such as
\[
\parbox{2cm}{\psfrag{p}{$V_{\pi}(G)$}
\includegraphics[width=2cm]{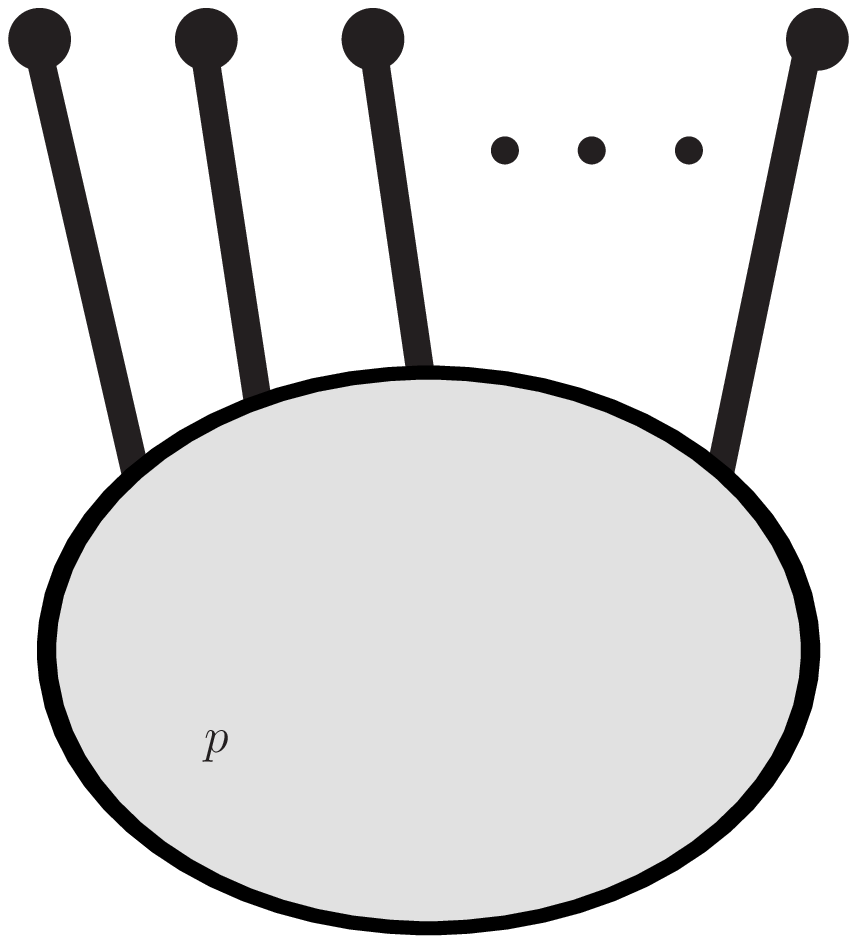}}
\]
obtained by cutting the crossing edges as in (\ref{cuttingedge})
and filling the ends by 1-valent vertices.
Such a possibly disconnected directed graph only has three types
of vertices:
\[
\parbox{1.2cm}{\includegraphics[width=1.2cm]{Fig33.eps}}
\qquad\qquad
\parbox{0.3cm}{\includegraphics[width=0.3cm]{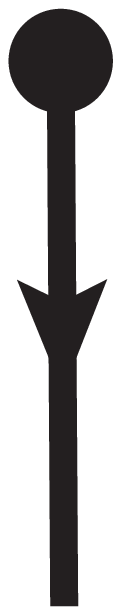}}
\qquad\qquad
\parbox{0.3cm}{\includegraphics[width=0.3cm]{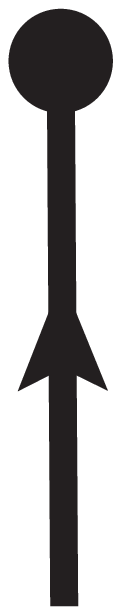}}\qquad .
\]
A (very short) moment of thought will convince the reader of the following
key observation.

\noindent{\bf Key fact:}
A connected component for such a digraph must either be a (binary)
tree coherently oriented towards the root such as
\[
\parbox{5cm}{\includegraphics[width=5cm]{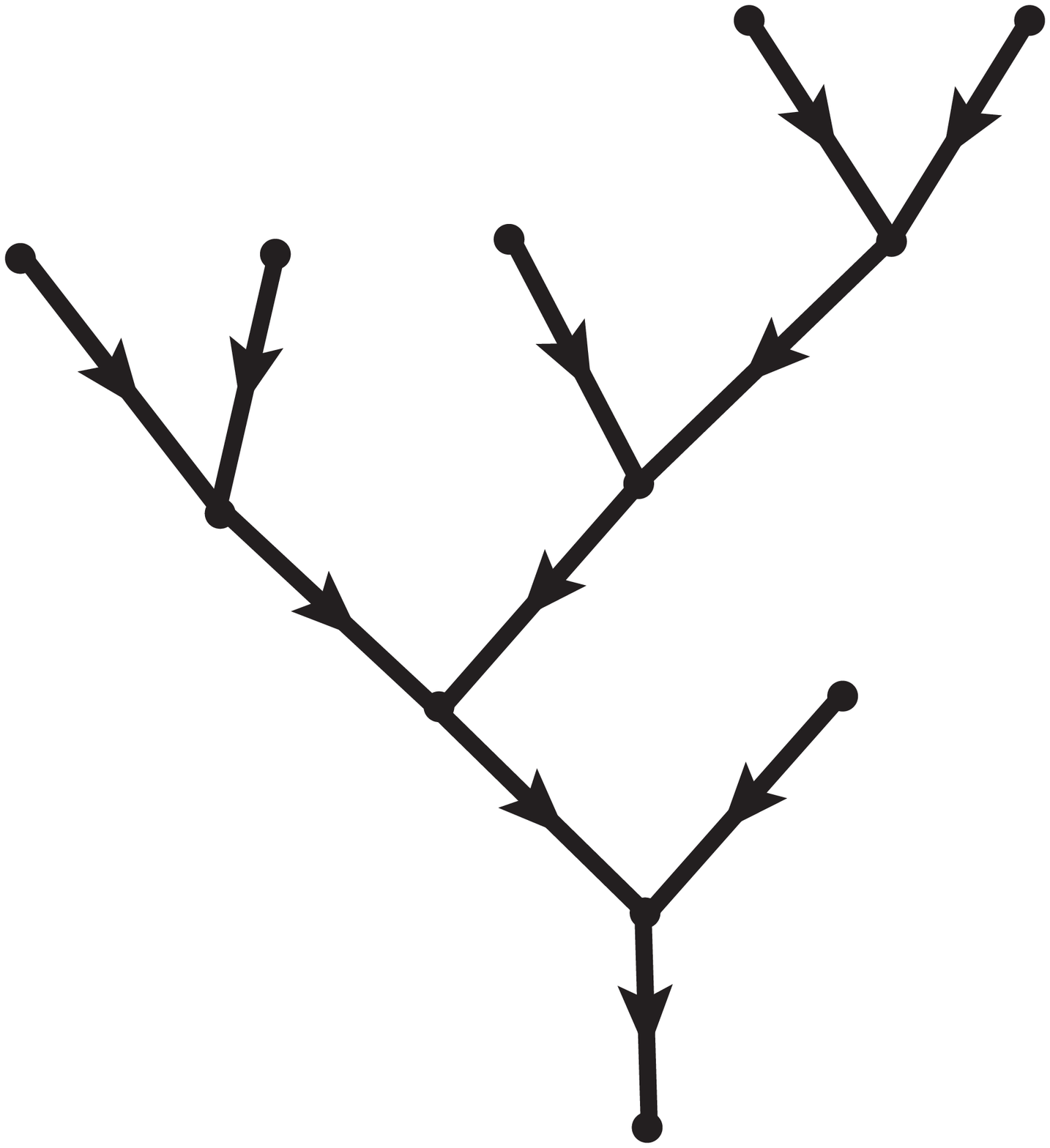}}
\qquad {\rm Type\ I}
\]
or a collection of such trees attached to a unique coherently oriented
central cycle of lenght at least 2, such as:
\[
\parbox{6cm}{\includegraphics[width=6cm]{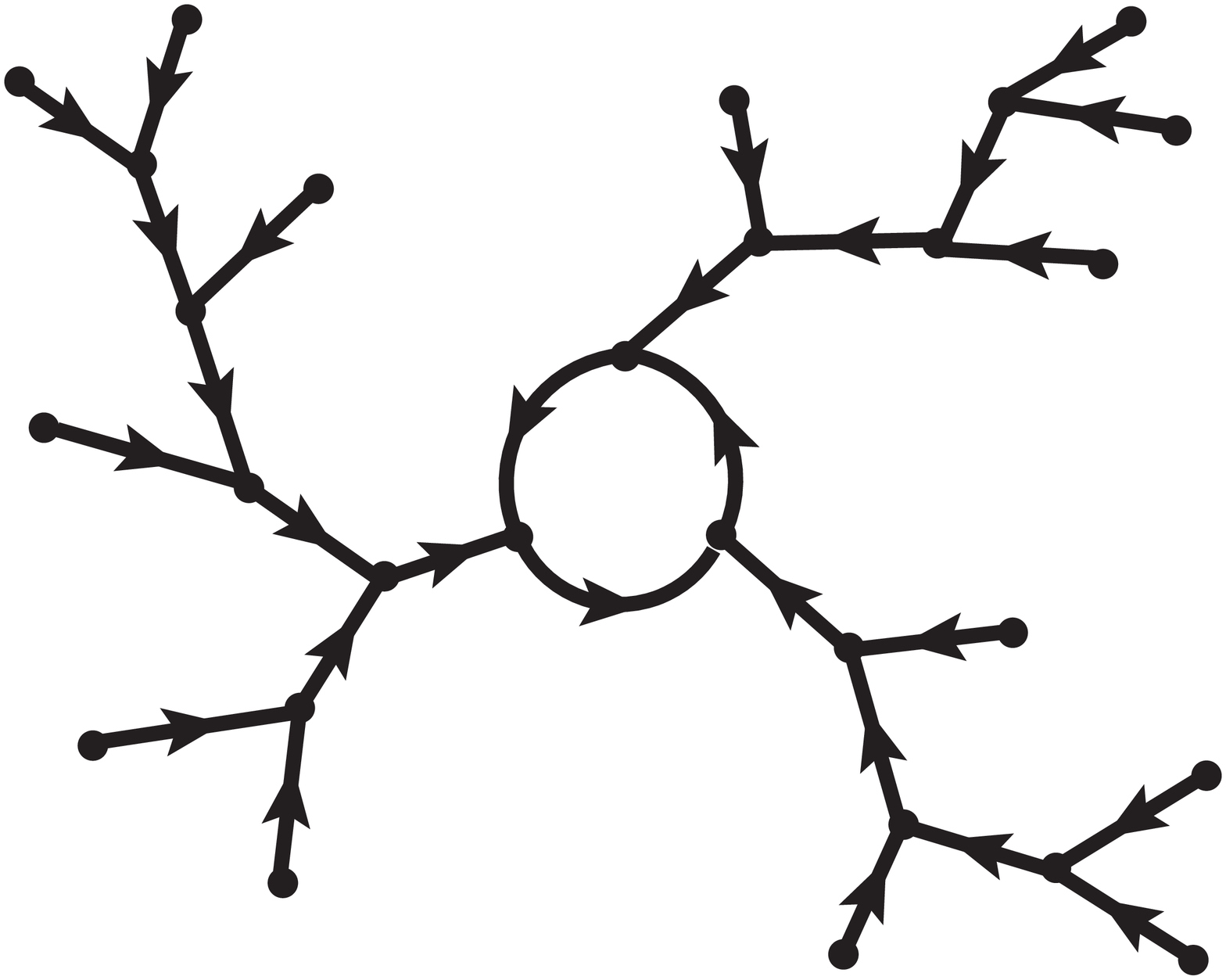}}
\qquad {\rm Type\ II}\ \ .
\]

Indeed,
such combinatorial structures are well-known to the practitioners
of the theory of species~\cite{BergeronLL}.
They correspond to the species of endofunctions where each element
has at most two preimages.
Previous experience with this type of
graphs~\cite{Abdesselamjac,AbdesselamSLC} was very useful in making this
observation.
Note that central cycles have lenght at least two since, otherwise,
the $V_{\pi}(G)$ subgraph would have a loop and therefore $G$ also which is
forbidden by the hypotheses.

Let $k_{\rm I}$ be the number of tree components of type I in $V_{\pi}(G)$.
Let $k_{\rm II}$ be the number of tree components of type II in $V_{\pi}(G)$.
For $i$, $1\le i\le k_{\rm I}$, let us denote by
$c_i$
the decoration
of the root edge of type I component number $i$ as in
\[
\parbox{3cm}{\psfrag{v}{$v_i$}
\psfrag{c}{$c_i={\rm dim}(v_i)-1$}
\includegraphics[width=3cm]{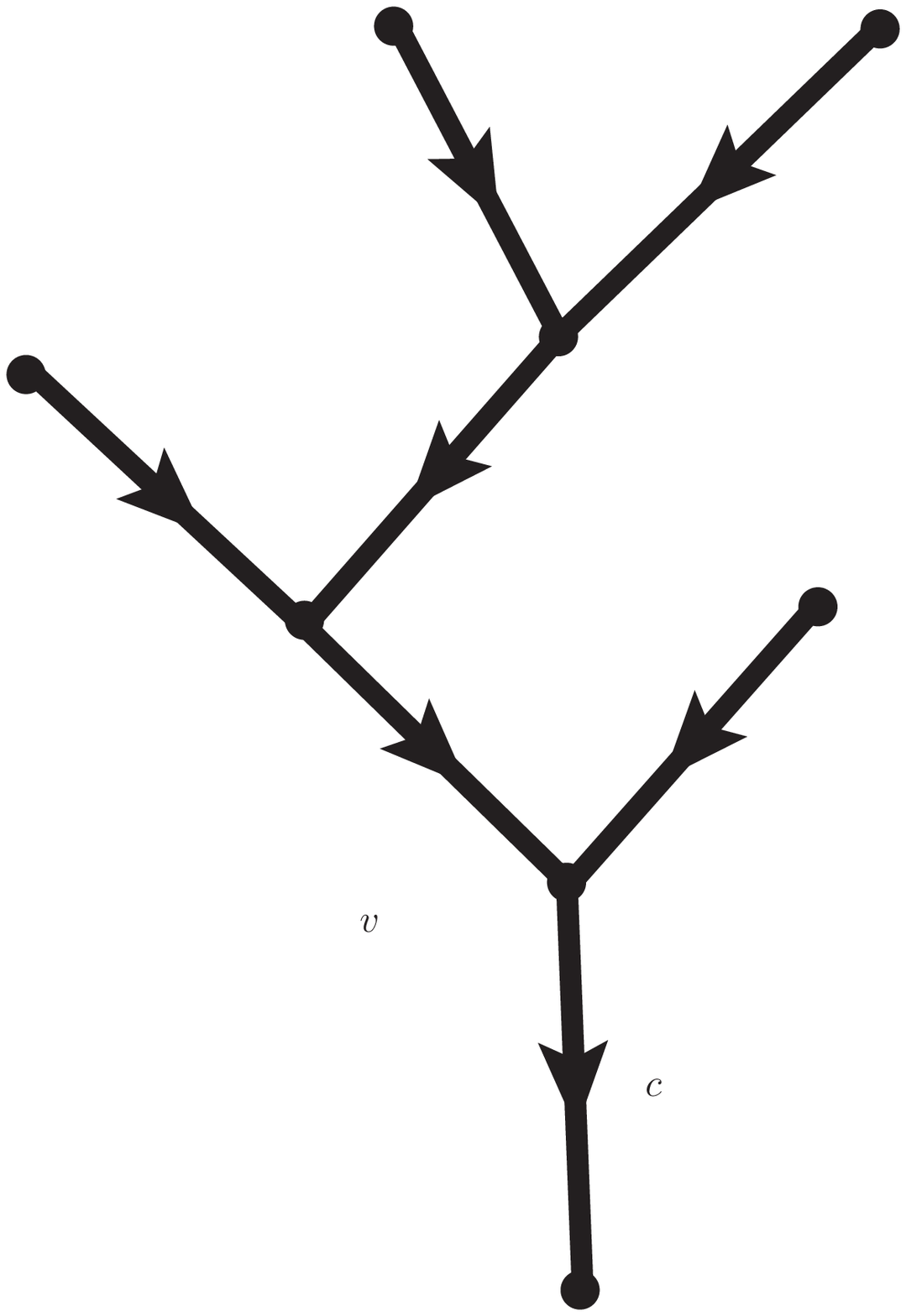}}\qquad\qquad\qquad .
\]
For $j$, $1\le j\le k_{\rm II}$, let $s_j\ge 2$
be the length of the central cycle of type II component number $j$,
and let $a_{ij}, b_{ij}$, $1\le i\le s_j$,
be the decorations around the cycle and on the edges incident to
the cycle as in
\[
\parbox{4cm}{\psfrag{1}{$a_{1j}$}
\psfrag{2}{$a_{2j}$}
\psfrag{3}{$a_{3j}$}
\psfrag{4}{$b_{1j}$}
\psfrag{5}{$b_{2j}$}
\psfrag{6}{$b_{3j}$}
\psfrag{7}{$v_{1j}$}
\psfrag{8}{$v_{2j}$}
\psfrag{9}{$v_{3j}$}
\includegraphics[width=4cm]{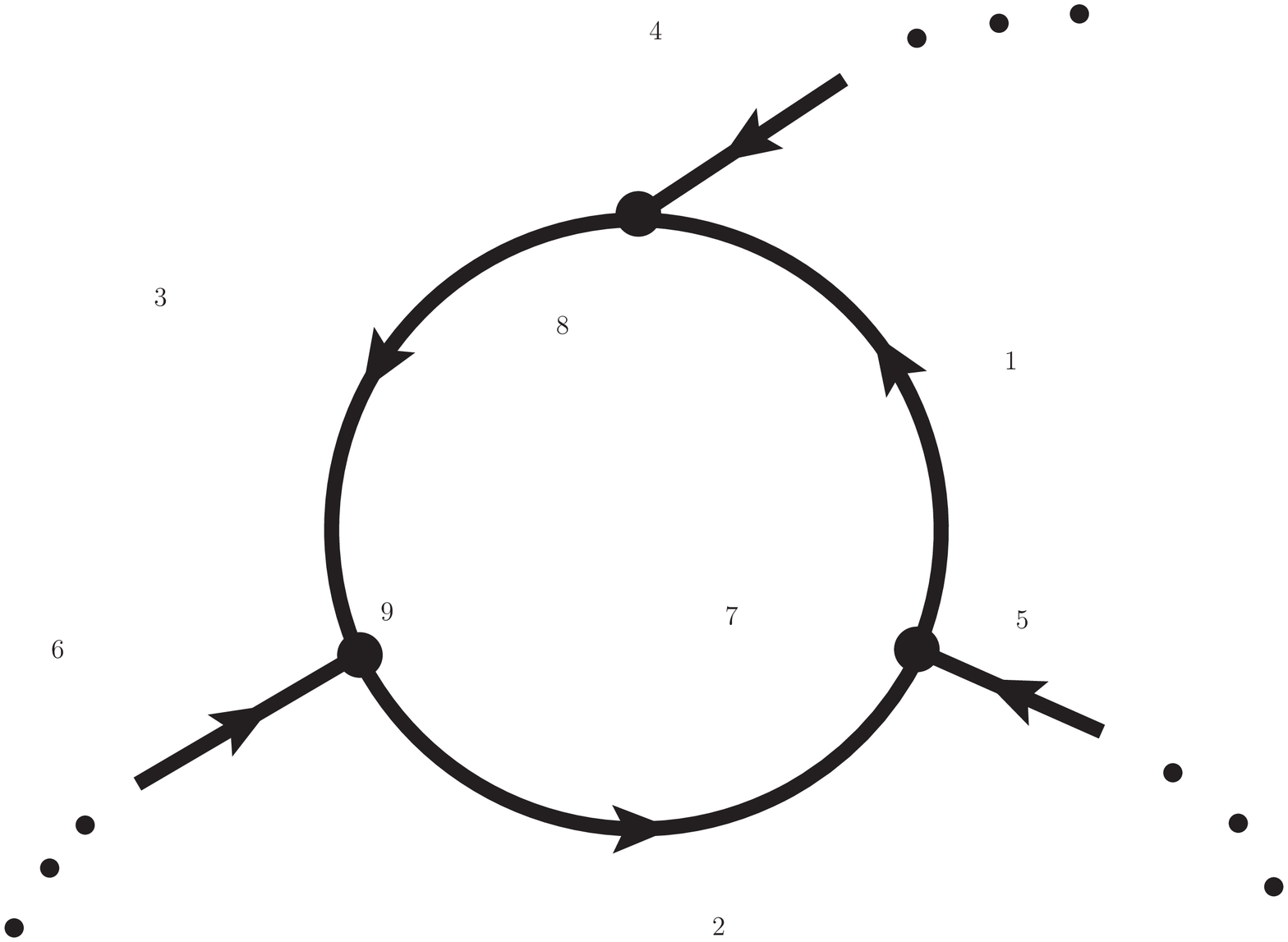}}\qquad .
\]
Also note the labelling of the vertices $v_{ij}$, so that
${\rm dim}(v_{ij})=a_{ij}+1$.
Then, by repeated application of Part 3) of Proposition \ref{piiotaprop},
we have
\[
\<
\parbox{2.5cm}{\psfrag{i}{$\widetilde{V_{\pi}(G)}$}
\psfrag{p}{$V_{\pi}(G)$}
\includegraphics[width=2.5cm]{Fig232.eps}}
\>^{\pi\io}=\prod\limits_{i=1}^{k_{\rm I}}
(c_i+1)\times
\prod\limits_{j=1}^{k_{\rm II}}
\<\ \ \ 
\parbox{5cm}{\psfrag{1}{$\scriptstyle{a_{1j}}$}
\psfrag{2}{$\scriptstyle{a_{2j}}$}\psfrag{3}{$\scriptstyle{a_{s_j,j}}$}
\psfrag{4}{$\scriptstyle{b_{1j}}$}\psfrag{5}{$\scriptstyle{b_{2j}}$}
\psfrag{6}{$\scriptstyle{b_{3j}}$}\psfrag{7}{$\scriptstyle{b_{s_j,j}}$}
\includegraphics[width=5cm]{Fig209.eps}}
\ \ \ \>^{\pi\io}
\]
with a positive sign in front. Indeed,
all $\<\cdots\>^{\pi\io}$
involved are squares of norms.
Therefore, by a coarse application of Lemma \ref{drumupperbd}
where we throw away the denominator,
\[
\<
\parbox{2.5cm}{\psfrag{i}{$\widetilde{V_{\pi}(G)}$}
\psfrag{p}{$V_{\pi}(G)$}
\includegraphics[width=2.5cm]{Fig232.eps}}
\>^{\pi\io} \le
\prod\limits_{i=1}^{k_{\rm I}}
(c_i+1)\times
\prod\limits_{j=1}^{k_{\rm II}}
\left[\min_{1\le i\le s_j}(a_{ij})+1\right]^2
\]
\[
\le \prod\limits_{i=1}^{k_{\rm I}}
(c_i+1)\times
\prod\limits_{j=1}^{k_{\rm II}}
\left[\prod\limits_{i=1}^{s_j}(a_{ij}+1)\right]
\]
because $s_j\ge 2$ for any $j$, $1\le j\le k_{\rm II}$.
As a result,
\[
\<
\parbox{2.5cm}{\psfrag{i}{$\widetilde{V_{\pi}(G)}$}
\psfrag{p}{$V_{\pi}(G)$}
\includegraphics[width=2.5cm]{Fig232.eps}}
\>^{\pi\io}\le
\prod\limits_{i=1}^{k_{\rm I}}
{\rm dim}(v_i)\times
\prod\limits_{j=1}^{k_{\rm II}}
\left[\prod\limits_{i=1}^{s_j}{\rm dim}(v_{ij})\right]
\]
\[
\le \prod\limits_{v\in V_{\pi}(G)} {\rm dim}(v)
\]
where the last bound is interpreted in the context of the original unsplit
CG network $(G,\cO,\ta,\ga)$.
Exactly the same reasoning shows
\[
\<
\parbox{2.5cm}{\psfrag{i}{$V_{\io}(G)$}
\psfrag{p}{$\widetilde{V_{\io}(G)}$}
\includegraphics[width=2.5cm]{Fig232.eps}}
\>^{\pi\io}\le
\prod\limits_{v\in V_{\io}(G)} {\rm dim}(v)
\]
and thus $\left|\<\Ga,\ga\>^U\right|\le 1$
since $V_{\pi}(G)$ and $V_{\io}(G)$ form a partition of $V(G)$.
\qed

This concludes the proof of Theorem \ref{mainthm}.
Note that an interesting question raised by the proof is that of optimizing
the choice of orientation $\cO$.
If all decorations are equal to $2n$ as in the setting of Conjecture \ref{mainconj},
and if there is only one component
of type I in $V_{\pi}(G)$ and also in $V_{\io}(G)$
one can easily see that our method provides a bound $|\<\Ga,2n\>^U|\le
(2n+1)^{1-\frac{1}{2} |V(G)|}$.
This is the case for standard 3n-j symbols obtained by comparing
two binary coupling schemes.
Not all cubic graphs can be decomposed in this
way~\cite{Yutsis}. This is related to an old conjecture by
Fran\c{c}ois Jaeger~\cite{Jaeger}.
A tantalizing question which is left for future investigations is
whether one can transpose the methods of this article to quantum spin
networks and the colored Jones polynomial of knots, with the hope of
improving known upper bounds for the volume conjecture~\cite{GaroufalidisL}.
We expect this investigation to benefit from the development
of a quantum version of CIT for binary forms by
Frank Leitenberger~\cite{Leitenberger} which could usefully complement
the more widely known techniques of~\cite{KauffmanL,MasbaumV,CarterFS}.

\bigskip
\noindent{\bf Acknowledgements:}
{\small
The drawings in this article were prepared using the JaxoDraw~\cite{BinosiT}
software.
We thank Stavros Garoufalidis, Bill Jackson, Vyacheslav Krushkal,
Gordon Royle and
Roland van der Veen for useful discussion or correspondence.
We thank Jaydeep Chipalkatti for four years of fruitful collaboration,
around CIT and its applications, which prompted us to develop the diagrammatic
formalism presented here. We thank Klaus Hoechsmann for his excellent translation
of~\cite{Gordan1} which allowed us to appreciate a sample of the combinatorial genius
of Paul Gordan.
We also thank Razvan Gurau and Vincent Rivasseau for interesting discussions
around spin networks, loop quantum gravity and noncommutative
field theory when at the conference
Combinatorial Identities and Their Applications in Statistical Mechanics,
April 7--11, 2008, at the Isaac Newton Institute for Mathematical Sciences,
Cambridge,
which we organized together with the late Pierre~Leroux.}

\end{document}